\documentclass[dvipsnames, reqno]{amsart}
\usepackage{preamble}

\renewcommand*{\emph}[1]{\textit{#1}}
\renewcommand*{\it}[1]{\textit{#1}}

\definecolor{dartmouthgreen}{rgb}{0.05, 0.5, 0.06}
\definecolor{MAGENTA}{rgb}{0.858, 0.188, 0.478}
\definecolor{DarkGreen}{rgb}{0.1,0.5,0.2}
\AtBeginDocument{%
\hypersetup{%
    colorlinks=true,
    linkcolor=RoyalPurple,
    filecolor=RedViolet,
    urlcolor=DarkGreen,
    citecolor=DarkGreen,
}}

\makeatletter
\def\@makechapterhead#1{\global\topskip 2.5pc\relax
  \begingroup
  \huge\bfseries                                                                
  \chaptername\enspace\thechapter\par                                           
  \vspace{1pc}                                                                  
     #1\par \endgroup                                                           
  \vspace{2pc}                                                                  
  }
\makeatother

\makeatletter
\DeclareRobustCommand{\ccong}{\,\mathrel{\mathpalette\@verequiv\sim}\,}
\newcommand{\@verequiv}[2]{%
  \lower.5\p@\vbox{
    \lineskiplimit\maxdimen
    \lineskip-.5\p@
    \ialign{%
      $\m@th#1\hfil##\hfil$\crcr
      #2\crcr
      \equiv\crcr
    }%
  }%
}
\makeatother

\numberwithin{equation}{section}

\renewcommand*{\theenumi}{\roman{enumi}}
\crefname{enumi}{\!}{\!}
\Crefname{enumi}{\!}{\!}
\creflabelformat{enumi}{(#2#1#3)}
\setlist[enumerate,3]{label=(\arabic*)}

\declaretheoremstyle[%
    spaceabove=1em,%
    spacebelow=1em,%
    headfont=\normalfont\bfseries,%
    notefont=\normalfont\bfseries,%
    notebraces={(}{)},%
    headpunct={.},%
    postheadspace=0.5em,%
    bodyfont=\normalfont,%
]{mplain}

\declaretheorem[style=mplain,qed=$\lrcorner$,numbered=no,name=Theorem]{theorem*}
\declaretheorem[style=mplain,qed=$\lrcorner$,numberwithin=section,name=Theorem]{theorem}
\declaretheorem[style=mplain,qed=$\lrcorner$,sibling=theorem,name=Lemma]{lemma}
\declaretheorem[style=mplain,qed=$\lrcorner$,sibling=theorem,name=Proposition]{proposition}
\declaretheorem[style=mplain,qed=$\lrcorner$,sibling=theorem,name=Corollary]{corollary}

\declaretheorem[style=mplain,qed=$\lrcorner$,sibling=theorem,name=Claim]{claim}

\declaretheorem[style=definition,qed=$\lrcorner$,sibling=theorem,name=Definition]{definition}

\declaretheorem[style=definition,qed=$\lrcorner$,sibling=theorem,name=Example]{example}
\declaretheorem[style=definition,qed=$\lrcorner$,sibling=theorem,name=Remark]{remark}

\crefname{theorem}{Thm}{Thms}
\Crefname{theorem}{Theorem}{Theorems}

\crefname{section}{Section}{Sections}
\Crefname{section}{Section}{Sections}
\crefname{subsection}{Section}{Sections}
\Crefname{subsection}{Section}{Sections}
\crefname{definition}{Def}{Defs}
\Crefname{definition}{Definition}{Definitions}
\crefname{notation}{Notation}{Notations}
\Crefname{notation}{Notation}{Notations}
\crefname{example}{Example}{Examples}
\Crefname{example}{Example}{Examples}
\crefname{remark}{Remark}{Remarks}
\Crefname{remark}{Remark}{Remarks}
\crefname{equation}{}{}
\Crefname{equation}{}{}

\crefname{corollary}{Cor}{Cors}
\Crefname{corollary}{Corollary}{Corollary}
\crefname{proposition}{Prop}{Props}
\Crefname{proposition}{Proposition}{Propositions}
\crefname{lemma}{Lem}{Lems}
\Crefname{lemma}{Lemma}{Lemmas}
\crefname{fact}{Fact}{Facts}
\Crefname{fact}{Fact}{Facts}

\crefname{tabular}{Table}{Tables}
\Crefname{tabular}{Table}{Tables}
\crefname{table}{Table}{Tables}
\Crefname{table}{Table}{Tables}
\Crefname{minipage}{Figure}{Figures}
\crefname{minipage}{Figure}{Figures}

\definecolor{RoyalPurple}{rgb}{0.47, 0.32, 0.66}

\tikzstyle{vslash} = [decoration={markings, mark=at position 0.5 with {\draw[-] (0,-1.5pt) -- (0,1.5pt);}}, postaction ={decorate}]
\tikzstyle{vslashslash} = [decoration={markings, mark=at position 0.5 with {\draw[-] (0.2pt,-1.5pt) -- (0.2pt,1.5pt);\draw[-] (-0.2pt,-1.5pt) -- (-0.2pt,1.5pt);}}, postaction ={decorate}]
\tikzstyle{vdot} = [decoration={markings, mark=at position 0.5 with {\draw [fill] (0,0) circle [radius=1pt];}}, postaction ={decorate}]
\tikzcdset{equal/.append style={double distance=1.0pt,-}}

\tikzcdset{%
    arrow style=tikz,
    diagrams={>={Straight Barb[scale=0.8]}},
    scale cd/.style={every label/.append style={scale=#1},
    cells={nodes={scale=#1}}},
    virtual/.style={row sep=small, scale cd=1},
}

\newcommand*{\sar}[1][]{%
    \ar[vslash,
        #1,
        ]
}

\newcommand*{\al}[1][]{%
    \ar[leftarrow,
        #1
        ]
}
\newcommand*{\sard}[1][]{%
    \ar[vslash, dashed,
        #1,
        ]
}

\setlist[itemize]{%
    leftmargin=1.5em,
}
\setlist[enumerate]{%
    leftmargin=2em,
}

\makeatletter
\NewDocumentCommand\Scale{O{1} m}{\scalebox{#1}{$\m@th#2$}}
\makeatother

\RenewDocumentCommand\to{O{2ex} D<>{}}{\hspace{-1.5truemm}\begin{tikzcd}[column sep=#1, ampersand replacement=\&]\!\ar[r,#2]\&\!\end{tikzcd}\hspace{-1.5truemm}}

\newcommand*{\sto}[1][]{\hspace{-1.5truemm}\begin{tikzcd}[column sep=small, ampersand replacement=\&]\!\sar[r, #1]\&\!\end{tikzcd}\hspace{-1.5truemm}}

\newcommand*{\adjointleft}{\hspace{-1.5truemm}\begin{tikzcd}[column sep=small, ampersand replacement=\&]\!%
        \al[r, shift right]\ar[r, phantom, "\mbox{\tiny $\bot$}"{marking, yshift=0.1em}]\ar[r, no head, white, shift left=2.6pt, line width=1pt]\ar[r, shift left=3.1pt]%
\&\!\end{tikzcd}\hspace{-1.5truemm}}
\NewDocumentCommand\doublecell{O{rd} m}
{
    \ar[
        #1,
        phantom, "#2"{description, inner sep=0mm}
    ]
}
\newcommand*{\pullback}[1][rd]{%
    \doublecell[#1, "\lrcorner"{description, near start, inner sep=0mm}]{}
}
\mathchardef\mhyphen="2D 

\newcommand*{\Set}{\one{S}\zero{et}}
\NewDocumentCommand{\Rel}{o o}
{%
    \dbl{R}\mathrm{el}
    \IfValueT{#2}{%
        _{#2}%
    }%
    \IfValueT{#1}{%
    (#1)%
    }%
}       
\newcommand*{\Cat}{\bi{Cat}}
\newcommand*{\CAT}{\bi{CAT}}

\newcommand*{\Fib}{\bi{Fib}}
\newcommand*{\BiFib}{\bi{BiFib}}
\newcommand*{\Reg}{\bi{Reg}}
\newcommand*{\Dbl}{\bi{Dbl}}
\newcommand*{\Eqp}{\bi{Eqp}}

\newcommand*{\Pos}{\one{P}os}
\newcommand*{\Fam}[1]{\one{F}am(#1)}

\newcommand*{\Frob}{_\bi{Frob}}

\newcommand*{\BC}{_{\bi{BC}}}

\NewDocumentCommand{\Span}{o}
{%
    \dbl{S}\mathrm{pan}
    \IfValueT{#1}{%
        (#1)%
    }%
}
\NewDocumentCommand{\Mat}{o}
{%
    \IfValueT{#1}{%
        {#1}\mhyphen
    }%
    \dbl{M}\mathrm{at}
}
\NewDocumentCommand{\Prof}{o d<>}
{%
    \IfValueT{#2}{%
        #2\mhyphen
    }%
    \dbl{C}\mathrm{at}
    \IfValueT{#1}{%
        (#1)%
    }%
}
\NewDocumentCommand{\PROF}{o d<>}
{%
    \IfValueT{#2}{%
        #2\mhyphen
    }%
    \dbl{C}\mathrm{AT}
    \IfValueT{#1}{%
        (#1)%
    }%
}
\NewDocumentCommand{\Mod}{o}
{%
    \dbl{M}\mathrm{od}
    \IfValueT{#1}{%
        (#1)%
    }%
}
    
\newcommand*{\BiSpan}{{\bi{S}\mathrm{pan}}}
\newcommand*{\BiRel}{{\bi{R}\mathrm{el}}}
\NewDocumentCommand{\PB}{o}
{%
    \one{P}\mathrm{b}
    \IfValueT{#1}{%
        (#1)%
    }%
}

\DeclareMathOperator{\ob}{ob}
\NewDocumentCommand{\Bil}{o}
{%
    \dbl{B}\mathrm{il}
    \IfValueT{#1}{%
        (#1)%
    }%
}
\NewDocumentCommand{\uni}{o}
{%
    \mf{uni}
    \IfValueT{#1}{%
        ({#1})%
    }%
}
\NewDocumentCommand{\Uni}{o}
{%
    \one{U}\mathrm{ni}
    \IfValueT{#1}{%
        ({#1})%
    }%
}
\newcommand*{\LBi}[1]{{\bi{L}}(#1)}
\newcommand*{\TBi}[1]{{\bi{T}}(#1)}
\newcommand*{\MapBi}[1]{{\bi{Map}}(#1)}
\newcommand*{\Map}[1]{\dbl{M}\mathrm{ap}(#1)}
\newcommand*{\LS}[1]{#1^>}

\newcommand*{\Cau}[1]{\overline{#1}}
\newcommand*{\Sub}[1]{\one{S}\mathrm{ub}(#1)}
\NewDocumentCommand{\Pred}{o o}
{%
    \one{P}\mathrm{red}
    \IfValueT{#2}{%
        _{#2}%
    }%
    \IfValueT{#1}{%
    (#1)%
    }%
}

\newcommand*{\Fibarr}{\zero{Fib}}
\newcommand*{\Predarr}{\zero{Prd}}

\newcommand*{\tgt}{\textsf{tgt}}
\newcommand*{\src}{\textsf{src}}

\newcommand*{\opr}{\mathsf{opRest}}

\renewcommand*{\sum}{\Sigma}

\newcommand*{\mmbox}[1]{\mbox{$#1$}}
\newcommand*{\mmboxsmall}[1]{\mbox{\small $#1$}}

\newcommand*{\id}{\mathtt{id}}

\newcommand*{\Idf}{\mathrm{Id}}
\newcommand*{\carttwo}{_{\bi{cart}}}
\newcommand*{\restc}{\mmboxsmall{\mathtt{rest}}}
\newcommand*{\cart}{\mmboxsmall{\mathtt{prn}}}
\newcommand*{\opcart}{\mmboxsmall{\mathtt{spn}}}

\newcommand*{\op}{^\textsf{op}}

\newcommand*{\lop}{^\textsf{lop}}
\newcommand*{\tiop}{^\textsf{top}}

\newcommand*{\inv}{^{-1}}
\newcommand*{\arrcat}{^{\rightarrow}}

\newcommand*{\zero}[1]{\mathrm{#1}}
\newcommand*{\one}[1]{\mathcal{#1}}
\newcommand*{\bi}[1]{{\mathbf{#1}}}
\newcommand*{\dbl}[1]{\mathbbm{#1}}

\newcommand*{\syn}[1]{\mathsf{#1}}
\newcommand*{\monunit}{\mathnormal{1}}
\newcommand*{\mf}[1]{\mathfrak{#1}}

\newcommand*{\abs}[1]{{\lvert #1\rvert}}

\newcommand*{\tpl}[1]{{\langle {#1} \rangle}}
\newcommand*{\fpl}[1]{{\,\langle\!\langle {#1} \rangle\!\rangle\,}}

\newcommand*{\smcl}{\,\fatsemi\,}
\newcommand*{\cmpr}[1]{\{\!|#1|\!\}}

\newcommand*{\ol}[1]{\overline{#1}}
\newcommand*{\ul}[1]{\underline{#1}}

\newcommand*{\wt}[1]{\widetilde{#1}}
\newcommand*{\wh}[1]{\widehat{#1}}

\NewDocumentCommand{\lc}{s t= t+}{%
    \IfBooleanT{#1}{%
    \times\land%
    }
    \IfBooleanT{#2}{%
    =
    }
    \IfBooleanT{#3}{%
    \exists
    }
}
\newcommand*{\el}{{\lc*=}}
\newcommand*{\exi}{{\lc*+}}
\newcommand*{\eltwo}{\el}
\newcommand*{\exitwo}{\exi}
\newcommand*{\eef}{{\lc*=+}}
\NewDocumentCommand{\vlc}{s t= t+}{%
    \IfBooleanT{#1}{%
    \times\land%
    }
    \IfBooleanT{#2}{%
    \ide{}
    }
    \IfBooleanT{#3}{%
    \odot
    }
}
\NewDocumentCommand{\VDbl}{t- s t= t+  t?}{%
    {\IfBooleanT{#1}{%
    \bi{F}
    }
    \bi{VDbl}_
    {\,
    \IfBooleanTF{#3}
    {\nrightarrow}{\!}
    \IfBooleanT{#4}
    {\odot}
    \IfBooleanT{#5}
    {\triangleright}
    \IfBooleanT{#2}
    {,\bi{cart}}
    }}
}
\newcommand*{\FVDbl}{\bi{FVDbl}}
\newcommand*{\VEq}{\bi{VEqp}}
\newcommand*{\FVDblP}{\FVDbl^{\odot}}

\newcommand*{\tightcat}{\text{t}}

\newcommand*{\ide}[1]{\!\nrightarrow_{#1}\!}

\newcommand*{\rest}{\textsf{rest}}

\acsetup{
    format/short=\scshape,
    format/replace=false
    }
\DeclareAcronym{DC}{short = DC, long =double category, long-plural-form = double categories}
\DeclareAcronym{FVDC}{short = FVDC, long =fibrational virtual double category, long-plural-form = fibrational virtual double categories}
\DeclareAcronym{CFVDC}{short = CFVDC, long =cartesian fibrational virtual double category, long-plural-form = cartesian fibrational virtual double categories}
\DeclareAcronym{VDC}{short = VDC, long =virtual double category, long-plural-form = virtual double categories}
\DeclareAcronym{FVDblTT}{short = FVD{\normalfont bl}TT, long =\underline{f}ibrational \underline{v}irtual \underline{d}ou\underline{bl}e \underline{t}ype \underline{t}heory}

\makeatletter
\define@key {mprset}{style}[1]{\def\TirNameStyle{#1}}
\def\l@subsection{\@tocline{1}{0pt}{2pc}{}{}}
\makeatother

\AtBeginEnvironment{mathpar}{\scriptsize}
\mprset{style={\scriptsize\scshape}}
\AtBeginEnvironment{mathparpagebreakable}{\scriptsize}
\mprset{style={\scriptsize\scshape}}

\makeatletter
\newcommand*{\colequiv}{\mathrel{\rlap{%
                    \raisebox{0.3ex}{$\m@th\cdot$}}%
                    \raisebox{-0.3ex}{$\m@th\cdot$}}%
                    \equiv}
\makeatother

\newcommand*\lon{%
        \nobreak
        \mskip6mu plus1mu
        \mathpunct{}%
        \nonscript
        \mkern-\thinmuskip
        {:}%
        \mskip2mu
        \relax
}

\tikzset{unit/.style={isosceles triangle, isosceles triangle stretches, draw, inner sep=0pt, shape border rotate=90, minimum height=0.6cm, minimum width=1.2cm, anchor = center}}
\tikzset{counit/.style={isosceles triangle, isosceles triangle stretches, draw, inner sep=0pt, shape border rotate=-90, minimum height=0.6cm, minimum width=1.2cm, anchor = center}}
\tikzset{unitw/.style={isosceles triangle, isosceles triangle stretches, draw, inner sep=0pt, shape border rotate=90, minimum width=1.6cm, minimum height=1.0cm, anchor = center}}
\tikzset{counitw/.style={isosceles triangle, isosceles triangle stretches, draw, inner sep=0pt, shape border rotate=-90, minimum width=1.6cm, minimum height=1.0cm, anchor = center}}
\tikzset{vertid/.style={rectangle, draw, inner sep=0pt, fill= black, minimum width=0.1cm}}
\tikzset{horid/.style={rectangle, draw, inner sep=0pt, fill=black, minimum height=0.1cm}}
\tikzset{unique/.style={circle, draw, inner sep=0pt, minimum size=0.5cm, fill=gray!40}}
\tikzset{first/.style={circle, inner sep=0pt, minimum size=0.5cm, top color=violet!40, bottom color = red!40, anchor=center, draw=black}}
\tikzset{second/.style={circle, inner sep=0pt, minimum size=0.5cm, top color=violet!40, bottom color = blue!40, anchor=center, draw=black}}
\tikzset{cell/.style={rectangle, draw, inner sep=0pt, minimum height=0.5cm, minimum width=0.5cm}} 

\tikzset{whf/.style={fill=white}}
\tikzset{redfill/.style={fill=red!40, draw=black}}
\tikzset{bluefill/.style={fill=blue!40, draw=black}}
\tikzset{viofill/.style={fill=violet!40, draw=black}}
\tikzset{redline/.style={draw=red}}
\tikzset{blueline/.style={draw=blue}}
\tikzset{violine/.style={draw=violet}}

\newcommand{\labeleditem}{\item\phantomsection}

\bookmarksetup{depth=3}
\usepackage[%
a4paper,margin=2cm,%
]{geometry}
\usepackage[backend=biber,style=alphabetic,sorting=nyt,giveninits=true,backref=true,isbn=false,url=true]{biblatex}
\title[Virtual double categories of relations]
{Virtual double categories of relations \\relative to fibrations}
\author{Hayato Nasu}
\thanks{
    Dalhousie University, Canada,\quad Email:hnasu@dal.ca}
\date{\today}

\begin{document}
\begin{abstract}
We study the relationship between virtual double categories and fibrations, or (hyper)doctrines, in the context of categorical logic.
We present a construction that produces a virtual double category of relations internal to a fibration with finite products.
We characterize the class of fibrations called elementary existential fibrations, a fibrational counterpart of regular logic, as those for which the induced virtual double category is a cartesian equipment.
We also show that cartesian equipments arising in this manner are exactly those that satisfy the Frobenius axiom studied by Walters and Wood in the context of cartesian bicategories.
Some logical principles, such as the comprehension axiom for predicates and the extensional equality predicate, are characterized as properties of the double categorical structure.
\end{abstract}

\maketitle
\begin{framed}
    The author's master's thesis, which contains most of the content of this paper as Chapter 2, was originally posted as this preprint (\url{https://arxiv.org/abs/2501.17869v2}).
    The current version is a standalone paper.
    The full thesis is available on their personal website (\url{https://hayatonasu.github.io/papers/logical-aspects-virtual-double-categories/}).
\end{framed}
\tableofcontents

\section{Introduction}
The aim of this paper is to propose virtual double categories as a framework for categorical semantics.
Categorical semantics provides a means to interpret logical systems and type theories in categorical structures.
The origin of this idea dates back to Lawvere's seminal work on the functorial semantics of algebraic theories \cite{lawvereFunctorialSemanticsAlgebraic1963}.
Given an algebraic theory, one can think of its models in a category with finite products by interpreting function symbols and terms as morphisms in the category and the equations as equalities of morphisms.
When one proceeds to interpret first-order logic, the interpretation of logical predicates becomes more involved,
as those predicates include quantifiers and logical connectives.
A naive way is to interpret predicates with free variables as subobjects of the product of the objects that the variables range over. 
Moreover, one must enhance the category with additional structure to interpret quantifiers, connectives, or other operators such as modalities.
For instance, to interpret regular logic, a fragment of first-order logic constituted by the equality $=$, the existential quantifier $\exists$, the conjunction $\land$, and the truth $\top$, the category in question should be a regular category.
The existential quantifier is then interpreted using the factorization system consisting of regular epimorphisms and monomorphisms.
There are various classes of categorical structures, not just 1-dimensional ones but more complex ones, that have been studied to interpret logical systems, and we aim to bring virtual double categories onto the scene.

\subsection{Fibrations and bicategories for categorical logic}
Before explaining the motivation for using virtual double categories, we review two well-known categorical approaches to interpreting first-order logic: using fibrations (or doctrines) and using bicategories.

Lawvere \cite{Law69} initiated the approach via hyperdoctrines, which can be viewed as special kinds of fibrations.
A \textit{(hyper)doctrine} is a category $\one{C}$ equipped with a contravariant pseudofunctor $\mf{P}:\one{C}^{\op}\to\Pos$ to the category of posets with additional structure.
The key idea is to remove the constraint that the interpretations of logical predicates are limited to subobjects of a fixed category, and instead allow them to be elements of indexed posets.
The category $\one{C}$ is then considered a category of contexts and terms, while predicates in a ``context'' $I$ are interpreted as elements in $\mf{P}(I)$.
The order-preserving map $\mf{P}(f):\mf{P}(J)\to\mf{P}(I)$ for a ``term'' $f:I\to J$ represents the substitution of the term, which is a fundamental concept in predicate logic.
The existential quantifier is then interpreted as the left adjoint of the substitution map, a groundbreaking idea by Lawvere in the aforementioned work.
It should be noted that a fibration over $\one{C}$ is a different formulation of the same notion, with the generalization of the target of $\mf{P}$ to the category of categories.

The study of bicategories of relations has been developed independently from the fibrational approach.
The study of \textit{cartesian bicategory} \cite{CKS84,CW87,CKWW07} is one of the most successful results in this direction.
The prototypical example is the bicategory of sets, binary relations between sets, and inclusions of relations as 0-cells, 1-cells, and 2-cells, respectively.
Binary relations are composed in this bicategory by the existential quantifier $\exists$ as follows: 
\[
\begin{aligned}
    R\colon A\sto B,\ S\colon B\sto C\quad&\mapsto\quad R\odot S\colon A\sto C,\\
    (R\odot S)(a,c)\quad&\Leftrightarrow\quad\exists b\in B.\ R(a,b)\land S(b,c).
\end{aligned}
\]
In addition, the identity relations are defined as the diagonal relations $\delta_A=\{\,(a,a')\mid a=a'\,\}$.
Another line of research in this vein was the study of \textit{allegories} by Freyd and Scedrov \cite{FS90}.
They introduced a calculus based exclusively on binary relations, but without the use of terms; this idea traces back to 19th-century work by Peirce and Schroeder and its subsequent development known as Tarski's relational calculus \cite{Tar41}.

The two approaches have their merits and drawbacks.
A merit of the fibrational approach is that it has a more transparent connection to logical systems, as indicated by its completeness theorem for first-order logic and its fragments \cite{jacobsCategoricalLogicType1999a}. 
On the other hand, it lacks binary relations as a primitive notion.
In addition, the interpretation of the equality and the existential quantifier in fibrations requires the finite products in the base category, while the bicategorical approach can handle them with their built-in composition.
Relatedly, the Beck-Chevalley condition and the Frobenius reciprocity are often imposed on fibrations to interpret logical systems, but their category-theoretic justification is not sufficiently clear; they are somewhat \textit{ad hoc} conditions primarily designed to make the interpretation sound.

A significant benefit of the bicategorical approach is the expressive power of compositionality of relations.
For instance, the exact completion of a regular category, studied broadly in the context of realizability theory \cite{Menni00} and constructive mathematics \cite{MR13}, refers to its internal equivalence relations.
The condition of an endo-relation being an equivalence relation is cumbersome to express within a regular category, but, once we construct the bicategory of the internal relations in that category, we can express it simply as a monad with symmetry.
In fact, the exact completion is described nicely through allegories in \cite{FS90,johnstoneSketchesElephantTopos2002a}.
However, not having functions as a primitive notion can be a disadvantage. 
Although one could regard internal functional relations as functions and interpret terms using them, they do not capture how terms are typically reasoned about in mathematics, particularly their operational nature.
This disadvantage is critical when interpreting logical systems on a foundation in which the principle of unique choice fails.
We need to disentangle the interpretation of functions from functional relations in orderto obtain a more flexible structure, just as the fibrational approach disentangles the interpretation of predicates from subobjects.
Thus, the two approaches are complementary to each other, and there ought to be some framework that combines the advantages of both.

\subsection{Virtual double categories of relations}
We propose that double categories are an adequate framework to achieve this goal.
The core insight is that \textit{relations are convenient tools, but functions are still indispensable}.
Double categories are categorical structures that have two kinds of 1-cells, called tight and loose arrows, and square-shaped 2-cells with two kinds of 1-cells as their boundaries, with two directions of composition.
Since both functions and relations involve composition, a double category is a natural setting to study them simultaneously.
The double category of sets, functions, and relations is an archetypal example of this idea.
Despite its intriguing properties, the double category of relations has not been studied as much as the bicategory of relations until recently.
The paper \cite{lambertDoubleCategoriesRelations2022} was the first step to filling this gap, followed by \cite{hoshinoDoubleCategoriesRelations2023}.
Dagnino and Pasquali independently developed a similar idea but with a different framework called \textit{relational doctrines} in \cite{DP23,dagnino2024cauchycompletionsruleuniquechoice}.

The reader may wonder why we go further to virtual double categories despite the fact that double categories already achieve the goal of capturing both functions and relations.
Virtual double categories are a generalization of double categories that do not assume composition of loose arrows to be defined; see \Cref{sec:virtualdoublecategories} for more details.
The reason is that the composition of relations relies on the equality and the existential quantifier, which could be unavailable in a logical system that is weaker than regular logic.
Nevertheless, even in the absence of the existential quantifier, one can still define the composition of relations when the following inclusion holds by the first line of the following equivalence:
\[
\begin{tikzcd}[virtual]
    A
    \sar[r, "R"]
    \ar[d, equal]
    \ar[drr, phantom, "\rotatebox{-90}{$\subseteq$}"]
    &
    B
    \sar[r, "S"]
    &
    C
    \ar[d, equal]
    \\
    A
    \sar[rr, "T"']
    &
    &
    C
\end{tikzcd}
\quad 
\begin{aligned}
\!
    \\
\Longleftrightarrow\quad
&
\forall a\in A.\forall b\in B.\forall c\in C.\ \left(
R(a,b)\land S(b,c)\right)\Rightarrow T(a,c)\\
(\Longleftrightarrow\quad
&
\forall a\in A.\forall c\in C.\ \left(
\exists b\in B.\ R(a,b)\land S(b,c)\right)\Rightarrow T(a,c))
\end{aligned}
\]
Moreover, we can also define inclusion with input relations comprised of a path of relations of finite length, without the existential quantifier or the equality.
Cells with multiple input relations are the basic building blocks of virtual double categories, and this observation suggests that virtual double categories can capture the structure that relations form even without the existential quantifier or equality.
Indeed, when the existential quantifier is available, the inclusion of relations can equivalently be expressed as the inclusion $R\odot S\subseteq T$, as described in the second line above, which indicates that the virtual double category of relations properly generalizes the double category of relations.

This is how virtual double categories offer a more primitive setup for binary relations than double categories.
Based on this idea, we will define a construction $\Bil$ that associates a virtual double category to a fibration with finite products both in the base and the fibers, which we call a \textit{cartesian fibration} in this paper.
This class of fibrations, or its equivalent notion of hyperdoctrines called \textit{primary doctrines}, is known to correspond to a logical fragment of first-order logic that only has conjunction $\land$ and truth $\top$.
The resulting virtual double category of ``relations'' $\Bil[\mf{p}]$ will have objects and morphisms in the base category of the fibration $\mf{p}$ as its objects and tight arrows, respectively, and its loose arrows from $I$ to $J$ will be the objects in the fiber of the product $I\times J$.
This construction recovers many known examples as summarized in \Cref{tab:primfibtovdc}.
    \begin{table}[hbtp]
        \centering
        \setlength{\aboverulesep}{0pt}
        \setlength{\belowrulesep}{0pt}
        \begin{tabular}{|c|c|}
            \toprule
            \rowcolor{gray!50}
            \textbf{Cartesian fibration} $\mf{p}$ & \textbf{Virtual double category} $\Bil[\mf{p}]$ \\
            \midrule
            \begin{tabular}{@{}c@{}}
            the codomain fibration $\one{B}^{\rightarrow}\to\one{B}$ \\
            ($\one{B}$ : a category with finite limits)
            \end{tabular} 
            & $\Span[\one{B}]$ \\
            \midrule
            \begin{tabular}{@{}c@{}}
            the subobject fibration $\Sub{\one{B}}\to\one{B}$ \\
            ($\one{B}$ : a category with finite limits)
            \end{tabular} 
            & $\Rel[\one{B}]$ \\
            \midrule
            \begin{tabular}{@{}c@{}}
            the family fibration $\Fam{L}\to\Set$ \\
            ($L$ : a meet-semilattice)
            \end{tabular} 
            & 
            \begin{tabular}{@{}c@{}}
                $\Mat[L]$ \quad
                w.r.t. \\
                the cartesian monoidal structure
            \end{tabular} 
            \\
            \bottomrule
        \end{tabular}
        \caption{Examples of the construction $\Bil$}
        \label{tab:primfibtovdc}
        \vspace{1em}
    \end{table}
    The reader may notice that the construction is reminiscent of the $\dbl{F}\mathrm{r}$ construction in \cite{shulmanFramedBicategoriesMonoidal2009}, and indeed, our construction is a generalization of its special case when the monoidal structure on the fibers is cartesian.

\subsection{Main results}
We now give an overview of the main results of this paper.

As we have seen that the existential quantifier and equality are used to define the composition and identity of relations, it is conceptually natural to speculate that the virtual double category of relations $\Bil[\mf{p}]$ has composition and identity loose arrows only when the fibration $\mf{p}$ can interpret the existential quantifier and the equality.
We prove that this is correct in the following theorem:
\begin{theorem*}[\Cref{thm:firsteefiscartdouble}]
    For a cartesian fibration $\mf{p}$ to be elementary existential, it is necessary and sufficient that the \acs{VDC} $\Bil[\mf{p}]$ is a cartesian equipment.
\end{theorem*}
\noindent Here, elementary existential fibrations are a fibrational analogue of Lawvere's elementary existentional doctrines \cite{lawvereEqualityHyperdoctrinesComprehension1970}, and one can deem them as a semantic counterpart of regular logic based on the completeness theorem for regular logic \cite{jacobsCategoricalLogicType1999a}.
Cartesian equipments are double categories that can interpret substitution of ``functions'' into ``relations'' and has a double-categorical finite-product structure, but the important point is that they have composition of loose arrows that is compatible with the finite-product structure.
Therefore, the theorem above indicates that regular logic is exactly what is needed to make the virtual double category of relations a cartesian equipment, or to get composition of relations compatible with the finite products.

The second main result is the characterization of the cartesian equipments that arise from elementary existential fibrations.
By construction, the virtual double category of relations $\Bil[\mf{p}]$ has a symmetric nature in the loose direction, namely, every object is self-dual in that direction.
In the context of cartesian bicategories, \textit{Frobenius axiom} in \cite{WW08}, which asks all objects to posess a Frobenius algebra structure in a weak sense, is known to provide such symmetry.
We prove that a double-categorical version of the Frobenius axiom is the exact condition that characterizes the image of the construction $\Bil$ on elementary existential fibrations.
\begin{theorem*}[\Cref{cor:unilateral}]
    A cartesian equipment $\dbl{D}$ is equivalent to $\Bil[\mf{p}]$ for some elementary existential fibration $\mf{p}$ if and only if it satisfies the Frobenius axiom.
\end{theorem*}
\noindent We will present this result as a 2-equivalence between the 2-category of elementary existential fibrations and the 2-category of cartesian equipments satisfying the Frobenius axiom.
We can further restrict this to a 2-equivalence between the 2-category of regular fibrations and the 2-category of cartesian equipments having Beck-Chevalley pullback squares, which is also studied in \cite{Law15} under a slightly different setting.
The former class of fibrations is also studied in papers such as \cite{Pav96}, while the Beck-Chevalley condition is undoubtedly essential in categorical logic.

More technically, an important aspect of the second main result is that it disassociates the equality and the existential quantifier from the finite-product structure in the base category using double categories, which overcomes the aforementioned limitation of fibrations.
This is possible only because the composition of loose arrows that is independent of tight arrows encodes the existential quantifier and the equality.
As a result, the two special conditions defining elementary existential fibrations, the Frobenius reciprocity and the Beck-Chevalley condition, are now understood as the conditions that make the induced double categories cartesian, which provides a more natural and conceptual understanding of them.

While these results connect fibrations and double categories, the connection between bicategories and double categories is also worth exploring.
It was proven in \cite{patterson2024transposingcartesianstructuredouble} that the loose arrows of a cartesian equipment form a cartesian bicategory in the sense of \cite{CW87,CKWW07}, and we study a partial converse of this result.
As the three structures of elementary existential fibrations, cartesian equipments, and cartesian bicategories are tied together, we will also investigate how logical principles or properties of fibrations and bicategories are reflected in the corresponding double categories.
For instance, uniqueness of identity proofs is known to be purely formulated in terms of fibrations or doctrines, and we prove that this is equivalent, in the context of cartesian equipments, to the separability condition in \cite{LWW10}. 

\subsection{Related work and future work}
This paper is a step toward understanding the capabilities of (virtual) double categories in place of the existing structures in the context of categorical logic, revealing their relationship with fibrations and bicategories.
The seed of this idea can be found in a conference talk by Par\'{e} \cite{pare2009first}, opening up the possibility to ``put logic in the realm of double categories'', as he wrote in his slides. 
Subsequently, a study based on the same motivation as ours was conducted in \cite{Law15}, but our approach is more bottom-up, starting from virtual double categories.
Our future work includes the study of the exact completion, the tripos-to-topos construction, and other logical completion procedures in the context of virtual double categories.
For instance, the existential completion \cite{Trotta20} is similar to the path construction in virtual double categories \cite{dawsonPathsDoubleCategories2006} in concept once one incorporates the idea that composition in double categories is the counterpart of the existential quantifier.
We anticipates some connection between the two.
Exact and quotient completions should also be studied, as they could be expressed elegantly as a quotient of a loose symmetric monoid in a double category, which was studied with relational doctrines in \cite{dagninoRelationalQuotientCompletion2024}.
In these ways, we hope to develop logical aspects of double categories further beyond regular logic.

\addtocontents{toc}{\protect\setcounter{tocdepth}{0}}
\subsection*{Outline}
\Cref{sec:preliminary} reviews the basic concepts of fibrations, double categories, and virtual double categories.
\Cref{sec:fromfibtoVDC} introduces the construction $\Bil[\mf{p}]$ of a virtual double category of relations from a cartesian fibration $\mf{p}$, and proves the first main result that characterizes the elementary existential fibrations as those that induce a cartesian equipment by this construction.
\Cref{sec:characterization} proves the second main result that characterizes the cartesian equipments that arise from elementary existential fibrations by the Frobenius axiom.
\Cref{sec:translating} discusses some logical principles in terms of double categories, namely, the comprehension axiom for predicates and relations (\Cref{subsec:comprehension}), function extensionality (\Cref{subsec:funcext}), the unique choice principle (\Cref{subsec:funccomp}), and extensional equality (\Cref{subsec:extensionality}). Some of the properties of double categories corresponding to these principles play a key role in understanding the relationship with other structures, such as orthogonal factorization systems (\Cref{subsec:factorization}) and cartesian bicategories (\Cref{subsec:cartesianbicategories}).

\subsection*{Notation}
We will use the fraktur font $\mf{p}$ to denote a fibration, the boldface font $\bi{K}$ to denote a bicategory, and the blackboard bold font $\dbl{D}$ to denote a double category.

For a fibration $\mf{p}\colon\one{E}\to\one{B}$, the fiber category over an object $I$ is written as $\one{E}_I$. 
We write $\alpha[f]$ for the reindexing of an object $\alpha\in\one{E}_J$ of the fiber of the fibration over $J$, along an arrow $f\colon I\to J$ in $\one{B}$.

We write $\tpl{f_0,\dots,f_{n-1}}\colon A\to B_0\times\dots\times B_{n-1}$ for the arrow induced by $f_i\colon A\to B_i$ for $i=0,\dots,n-1$.
In the case where the $f_i$'s are all projections, we adopt a more suggestive notation, and write $\tpl{i_0,\dots,i_{n-1}}\colon A_0\times\dots\times A_{m-1}\to A_{i_0}\times\dots\times A_{i_{n-1}}$ for the arrow whose $j$-th component is the $i_j$-th projection for $j=0,\dots,n-1$.
For example, we write $\tpl{0,0}\colon A\to A\times A$ for the diagonal arrow, $\tpl{0}\colon A\times B\to A$ for the first projection, and $\tpl{1}\colon A\times B\to B$ for the second projection. 
This notation facilitates calculation of the composition of arrows given by the projections.
For example, 
\[
    \left(\begin{tikzcd}
        A\times C\times C
        \ar[from=r, "\tpl{0,2,4}"]
        &
        A\times A\times C\times B\times C
        \ar[from=rr, "\tpl{\underline{0},0,\underline{2},1,\underline{2}}"]
        &&
        A\times B\times C
    \end{tikzcd}\right)
    =
    \left(\begin{tikzcd}
        A\times C\times C
        \ar[from=r, "\tpl{0,2,2}"]
        &
        A\times B\times C
    \end{tikzcd}\right).
\]
Accordingly, we have 
\[
    \alpha[\tpl{0,2,4}][\tpl{0,0,2,1,2}]\cong\alpha[\tpl{0,2,2}] \quad\text{in}\quad
    \one{E}_{A\times B\times C}
\]
for a fibration $\mf{p}\colon\one{E}\to\one{B}$ and an object $\alpha\in\one{E}_{A\times C\times C}$. 
On the other hand, we will write $!$ for the unique arrow to the terminal object, not $\tpl{}$.

\subsection*{Acknowledgements}
This paper was originally written as a part of the author's master's thesis at Kyoto University, with the title ``Logical Aspects of Virtual Double Categories''.
The author would like to express his appreciation to his former supervisor, Masahito Hasegawa for his guidance and generous support throughout the master's program, and his current supervisor, Dorette Pronk, for her helpful advice and encouragement.
The author would like to thank Nathanael Arkor, Aaron Fairbanks, Zeinab Galal, Keisuke Hoshino, Yuto Kawase, Yuki Maehara, Yutaka Maita, and Robert Morissette for their discussions, comments, and suggestions on the manuscript.
The author would also like to thank the anonymous reviewers of \cite{hoshinoDoubleCategoriesRelations2023} for their valuable comments and suggestions, which have influenced the content of this paper.
\addtocontents{toc}{\protect\setcounter{tocdepth}{1}}

\section{Preliminaries}
\label{sec:preliminary}
\subsection{Fibrations}\label{sec:fibration}
In this section, we review basic notions and results on fibrations, introducing the terminology and the notation used in this paper.
The basic theory of fibrations can be found in \cite{Her93,Her99,jacobsCategoricalLogicType1999a,johnstoneSketchesElephantTopos2002a,Pit00}.

A (Grothendieck) fibration $\mf{p}\colon\one{E}\to\one{B}$ is said to be \emph{cloven} if the prone lifts\footnote{Prone lifts are commonly called \it{cartesian} lifts in the literature.
We use the term ``prone'', adopted from \cite{Tay99,johnstoneSketchesElephantTopos2002a}, to avoid confusion with the notion of cartesian objects in a 2-category.} are chosen, and further said to be \emph{split} if the chosen prone lifts are strictly functorial: that is, $\beta[\id_{I}]=\beta$ and $\beta[g][f]=\beta[g\circ f]$ for any arrows $f\colon I\to J$ and $g\colon J\to K$ in $\one{B}$, where $\beta[g]$ is a reindexing of $\beta\in\one{E}_{K}$ along $g$ given by the chosen prone lift of $g$.

Giving a cloven fibration whose base category is $\one{B}$ is equivalent to giving an \emph{indexed category} over $\one{B}$ \cite{JP78}, that is, a pseudofunctor $\one{B}\op\to\CAT$ where $\CAT$ is the 2-category of categories, functors, and natural transformations.

In particular, an indexed category whose values are posets is sometimes called a \textit{doctrine} \cite{lawvereEqualityHyperdoctrinesComprehension1970,KR77} as a rudimentary version of a \textit{hyperdoctrine} \cite{Law69}.
In this paper, we will use the term \emph{doctrine} to refer to a fibration with each fiber being a poset, which is automatically split.

\begin{example}
    \label{ex:fibrations}
    We give some examples of fibrations.
    \begin{enumerate}
        \item The codomain functor $\one{B}^{\rightarrow}\to\one{B}$ is a fibration
        if and only if $\one{B}$ has pullbacks,
        which we call the \emph{codomain fibration} over $\one{B}$.
        \item Let $\Sub{\one{B}}$ be the category of a pair $(I,m)$ of an object $I$ in $\one{B}$ and a subobject $m$ of $I$, that is an isomorphism class of monomorphisms into $I$. 
        Then, the canonical functor $\Sub{\one{B}}\to\one{B}$ is a fibration
        if and only if $\one{B}$ has pullbacks of monomorphisms,
        which is called the \emph{subobject fibration} over $\one{B}$.
        \item For a category $\one{B}$, let $\Fam{\one{B}}$ be the small coproduct cocompletion of $\one{B}$,
        that is, the category whose objects are pairs $\left(I,(b_i)_i\right)$
        where $I$ is a set and $(b_i)_i$ is a family of objects in $\one{B}$ indexed by $I$,
        and whose arrows from $\left(I,(b_i)_i\right)$ to $\left(J,(c_j)_j\right)$ are pairs $\left(u,(f_i)_i\right)$ where
        $u\colon I\to J$ is a function and $f_i\colon b_i\to c_{u(i)}$ for each $i\in I$.
        Then the forgetful functor $\Fam{\one{B}}\to\Set$ is a fibration, called the \emph{family fibration}.
        \qedhere
    \end{enumerate}
\end{example}

We will write $\Fib$ for the 2-category of fibrations, and $\Fib_{\one{B}}$ for the 2-category of fibrations over $\one{B}$; see \cite{Her99} for the details.

\begin{proposition}[{\cite[Corollary 3.7]{Her94}, \cite[4.1]{Her99}}]
    \label{prop:cartesianfibration}
    Let $\mf{p}\colon\one{E}\to\one{B}$ be a fibration where $\one{B}$ is a category 
    with finite products.
    Then, the following conditions are equivalent:
    \begin{enumerate}
        \item $\mf{p}$ is a cartesian object in $\Fib$.
        \item $\mf{p}$ is a cartesian object in $\Fib_{\one{B}}$.
        \item $\one{E}$ has finite products and the functor $\mf{p}$ preserves them.
        \item For any object $I\in\one{B}$, the fiber $\one{E}_{I}$ has finite products,
        and for any arrow $f\colon I\to J$ in $\one{B}$, the base change functor $(-)[f]\colon\one{E}_{J}\to\one{E}_{I}$
        preserves finite products.
    \qedhere
    \end{enumerate}
\end{proposition}

\begin{definition}
    A \emph{cartesian fibration}\footnote{
    Cartesian fibrations are also used as a synonym of (Grothendieck) fibrations in the infinity-categorical setting, which is different from what we mean here.} is a fibration $\mf{p}\colon\one{E}\to\one{B}$ that satisfies the equivalent conditions in \Cref{prop:cartesianfibration}.
\end{definition}

\begin{remark}
A cartesian fibration is simply called a \it{fibration with fibered finite products whose base category has finite products} \cite{Her94,jacobsCategoricalLogicType1999a}, and the corresponding doctrine is called \emph{a primary doctrine} \cite{MR13} (or \it{prop-category} \cite{Pit00}). 
\end{remark}

\begin{example}
    \label{ex:cartesianfibration} The following are basic examples we use in this paper.
    \begin{enumerate}
    \item
    The codomain fibration $\one{B}^{\rightarrow}\to\one{B}$ for a category $\one{B}$ with finite limits is a cartesian fibration. (In fact, it is a regular fibration \Cref{ex:regularfibration}.)
    \item The subobject fibration $\Sub{\one{B}}\to\one{B}$ for a category $\one{B}$ with finite limits is a cartesian fibration. 
    \item
    The family fibration $\Fam{\one{B}}\to\Set$ is a cartesian fibration if $\one{B}$ has finite products.
    \qedhere
    \end{enumerate}
\end{example}

We now turn to the discussion of left adjoints of reindexing functors.
It is widely known that left adjoints of reindexing functors should satisfy some conditions so that they correspond to the existential quantifier and equality in logic.

We use the notation $\sum_{f}$ for the left adjoint of the reindexing functor $(-)[f]\colon\one{E}_{J}\to\one{E}_{I}$ along $f\colon I\to J$.

\begin{definition}
\label{def:BCFR}
Let $\mathfrak{p}\colon\one{E}\to\one{B}$ be a cartesian fibration.
\begin{enumerate}
    \item 
    A \emph{functorial choice of pullback squares} is a functor $\Phi\colon\one{C}\to\PB[\one{B}]$
    into the wide subcategory $\PB[\one{B}]\subseteq\one{B}^{\rightarrow}$ of $\one{B}^{\rightarrow}$
    whose arrows are pullback squares in $\one{B}$.
    For each object $c\in\one{C}$, we write $\Phi_c\colon D_c\to C_c$ for the value of $\Phi$ at $c$.

    For a cartesian fibration $\mathfrak{p}\colon\one{E}\to\one{B}$ and a functorial choice of pullback squares $\Phi\colon\one{C}\to\PB[\one{B}]$,
    we say that $\mathfrak{p}$ \emph{has $\Phi$-coproducts} if for any object $c\in\one{C}$,
    the functor $(-)[\Phi_c]\colon\one{E}_{C_c}\to\one{E}_{D_c}$ has a left adjoint $\sum_{\Phi_c}$.
    \[
        \begin{tikzcd}
            \one{E}_{D_c}
            \ar[r, bend left=20, "{\sum_{\Phi_c}}"]
            \ar[r, phantom, "\rotatebox{90}{$\vdash$}"]
            &
            \one{E}_{C_c}
            \ar[l, bend left=20, "{(-)[\Phi_c]}"]
        \end{tikzcd}
    \]
    \item We say that $\mf{p}$ satisfies \emph{the Beck-Chevalley condition (BC condition) for a pullback square with direction}
    $(g,h)\colon f\to f'$ in $\one{B}^{\rightarrow}$ as in 
    \begin{equation}
    \label{eq:typicalpullback}
    \begin{tikzcd}[column sep=small, row sep=small]
        I 
        \ar[d, "f"']
        \ar[r, "h"]
        \ar[dr, phantom, "\lrcorner", very near start]
        &
        I'
        \ar[d, "f'"]
        \\
        J
        \ar[r, "g"']
        &
        J'
    \end{tikzcd}
    \end{equation}
    with $f$ and $f'$ admitting left adjoints $\sum_f$ and $\sum_{f'}$ if the following canonical natural transformation is an isomorphism:
    \[
        \begin{tikzcd}
            \one{E}_{I}
            \ar[from=r, "h^*"']
            \ar[d, "\sum_{f'}"']
            \ar[rd, Rightarrow, shorten <= 1em, shorten >= 1em]
            &
            \one{E}_{I'}
            \ar[d, "\sum_f"]
            \\
            \one{E}_{J}
            \ar[from=r, "g^*"]
            &
            \one{E}_{J'}
        \end{tikzcd}.
    \]
    For a functorial choice of pullback squares $\Phi\colon\one{C}\to\PB[\one{B}]$ and a cartesian fibration $\mathfrak{p}\colon\one{E}\to\one{B}$ with $\Phi$-coproducts, we say that $\mathfrak{p}$ satisfies the \emph{Beck-Chevalley condition for $\Phi$} if it satisfies the Beck-Chevalley condition for any pullback square with direction $\Phi_t$ for every $t\in\one{C}$.
    \[
        \Phi_t
        =
        \begin{tikzcd}
            D_c
            \ar[d, "\Phi_c"']
            \ar[r, "D_t"]
            \ar[dr, phantom, "\lrcorner", very near start]
            &
            D_{c'}
            \ar[d, "\Phi_d"]
            \\
            C_c
            \ar[r, "C_t"']
            &
            C_{c'}
        \end{tikzcd}
    \]
    \item We say that $\mathfrak{p}$ satisfies \emph{Frobenius reciprocity} for an arrow $f\colon I\to J$ that admits a left adjoint $\sum_f$ if the following canonical natural transformation is an isomorphism:
    \[
        \begin{tikzcd}[column sep=large]
            \one{E}_{I}
            \ar[from=r, "\land"']
            \ar[d, "\sum_{f}"']
            \ar[rrd, Rightarrow, shorten <= 3em, shorten >= 3em]
            &
            \one{E}_{I}\times\one{E}_{I}
            \ar[from=r, "\id\times f^*"']
            &
            \one{E}_{I}\times\one{E}_{J}
            \ar[d, "\sum_{f}\times\id"]
            \\
            \one{E}_{J}
            \ar[from=rr, "\land"]
            &
            &
            \one{E}_{J}\times\one{E}_{J}
        \end{tikzcd}.
    \]

    For a functorial choice of pullback squares $\Phi\colon\one{C}\to\PB[\one{B}]$ and a cartesian fibration $\mathfrak{p}\colon\one{E}\to\one{B}$
    with $\Phi$-coproducts,
    we say that $\mathfrak{p}$ satisfies \emph{Frobenius reciprocity for $\Phi$} if it satisfies Frobenius reciprocity for the arrows $\Phi_c$ whenever $c\in\one{C}$.
    \qedhere
\end{enumerate}
\end{definition}

Although the Beck-Chevalley condition depends on the direction of the pullback square, we often omit reference to the direction when it is clear from the context, or when we consider the condition for both directions simultaneously.

\begin{remark}
The definition is based on \cite[Section 1.9]{jacobsCategoricalLogicType1999a}.
A standard method to choose the arrows along which the left (or right) adjoints of reindexing functors are defined is to take a subclass of arrows in the base category, as in a \it{display map category} (\cite[\S 4.3.2]{Taylor83}, \cite[Definition 10.4.1]{jacobsCategoricalLogicType1999a}).
We adopt the functorial presentation to specify the form of the pullback squares for which the Beck-Chevalley condition should hold.
\end{remark}

\begin{definition}[{\cite[Definition 2.5]{EPR21}, \cite[Definition 4.1]{EPR22}}]
    \label{def:elementaryfibration}
    Let $\mf{p}\colon\one{E}\to\one{B}$ be a cartesian fibration.
    We define a functor $\Phi_{\lc=}\colon\ob\one{B}\times\one{B}\to\PB[\one{B}]$ 
    that assigns to each pair $(I,J)$ the arrow $\tpl{0,0,1}\colon I\times J\to I\times I \times J$,
    and assigns to each arrow $f\colon (I,J)\to (I,J')$, which is simply an arrow $f\colon J\to J'$ in $\one{B}$, 
    the pullback square on the left below.
    \[
    \begin{tikzcd}
        I\times J
        \ar[d, "\tpl{0,0,1}"']
        \ar[r, "\id\times f"]
        \ar[dr, phantom, "\lrcorner", very near start]
        &
        I\times J'
        \ar[d, "\tpl{0,0,1}"]
        \\
        I\times I\times J
        \ar[r, "\id\times\id\times f"']
        &
        I\times I
        \times J
    \end{tikzcd}
    \hspace{2em}
    \begin{tikzcd}[column sep=large]
        \one{E}_{I\times J}
        \ar[d, "\sum_{\tpl{0,0,1}}"']
        \ar[from =r, "(\id\times f)^*"']
        \ar[rd, Rightarrow, shorten <= 1em, shorten >= 1em]
        &
        \one{E}_{I\times J'}
        \ar[d, "\sum_{\tpl{0,0,1}}"]
        \\
        \one{E}_{I\times I\times J}
        \ar[from = r, "(\id\times\id\times f)^*"]
        &
        \one{E}_{I\times I\times J'}
    \end{tikzcd}
    \]
    We say that $\mf{p}$ is an \emph{elementary fibration} 
    if it has $\Phi_{\lc=}$-coproducts and satisfies the Beck-Chevalley condition and Frobenius reciprocity for $\Phi_{\lc=}$.
    Here, the Beck-Chevalley condition for $\Phi_{\lc=}$ is the condition that 
    the canonical natural transformation on the right above is an isomorphism.
\end{definition}

Importantly, the object $\sum_{\tpl{0,0}}\top_I$ in $\one{E}_{I\times I}$ is seen as the equality predicate on $I$ in the fiber $\one{E}_I$, which we occationally write as $\delta_I$.

\begin{remark}
In \cite[Section 3.4]{jacobsCategoricalLogicType1999a}, a fibration with its base category having finite products  is said to have \it{(simple) equality}\footnote{
Equality defined in \cite[Definition 5.6.1]{Pit00} does not require the Beck-Chevalley condition,
but it is discussed in the following paragraphs.
} if it has $\Phi_{\lc=}$-coproducts and satisfies the Beck-Chevalley condition for $\Phi_{\lc=}$, and it is said to have \it{equality with Frobenius property} if it further satisfies Frobenius reciprocity for $\Phi_{\lc=}$.
An elementary fibration which is a doctrine is called an \emph{elementary doctrine}  (\cite{lawvereEqualityHyperdoctrinesComprehension1970}, \cite[Definition 2.1]{MREQC13}), and one with preordered fibers is called an \it{Eq-fibration} \cite[Definition 3.5.1]{jacobsCategoricalLogicType1999a}.
\end{remark}

The following two lemmas are basic properties of elementary fibrations, which we could not find in the literature.
\begin{lemma}
    \label{lem:implicationsoffrobele}
    Let $\mf{p}\colon\one{E}\to\one{B}$ be an elementary fibration.
    Then, Frobenius reciprocity for $\Phi_{\lc=}$ induces the following isomorphisms:
    \[
    \sum_{\tpl{0,0,1}}\left(\kappa_1\land\kappa_2\right)
    \overset{\iota}{\cong}
    \left(\sum_{\tpl{0,0,1}}\kappa_1\right)\land\kappa_2[\tpl{0,2}]
    \cong 
    \left(\sum_{\tpl{0,0,1}}\kappa_1\right)\land\kappa_2[\tpl{1,2}]
    \]
    for $\kappa_1,\kappa_2\in\one{E}_{I\times J}$.
    Moreover,
    the above isomorphism $\iota$ makes the following diagram commute:
    \[
    \begin{tikzcd}[row sep=small]
        \kappa_1\land\kappa_2
        \ar[r, "\eta_{\kappa_1\land\kappa_2}"]
        \ar[d, "\eta_{\kappa_1}\land\id"']
        &
        \left(\displaystyle\sum_{\tpl{0,0,1}}(\kappa_1\land\kappa_2)\right)[\tpl{0,0,1}]
        \ar[d, "{\rotatebox{90}{$\cong$}}\ {\iota[\tpl{0,0,1}]}"]
        \\
        \left(\displaystyle\sum_{\tpl{0,0,1}}\kappa_1\right)[\tpl{0,0,1}]
        \land\kappa_2
        \ar[r, "\cong"']
        &
        \left(\left(\displaystyle\sum_{\tpl{0,0,1}}\kappa_1\right)\land\kappa_2[\tpl{01}]\right)[\tpl{0,0,1}]
    \end{tikzcd},
    \]
    where $\eta$ is the unit of the adjunction $\sum_{\tpl{0,0,1}}\dashv(-)[\tpl{0,0,1}]$
    and the bottom horizontal isomorphism is given by the preservation of finite products by base change functors.
    The corresponding statement holds if we replace $\tpl{0,2}$ with $\tpl{1,2}$.
\end{lemma}
\begin{proof}
    The isomorphisms are obtained by pre-composing the base change functors $(-)[\tpl{0,2}]$ and $(-)[\tpl{1,2}]$ to Frobenius reciprocity for $\Phi_{\lc=}$, and the commutativity of the diagram is a direct consequence of the definition of the natural transformation for Frobenius reciprocity.
    \begin{align*}
    \begin{tikzcd}[ampersand replacement=\&]
        \one{E}_{I\times J}\times\one{E}_{I\times J}
        \ar[r, "\sum_{\tpl{0,0,1}}\times\id"{yshift=1ex}]
        \ar[dr,equal]
        \&
        \one{E}_{I\times I\times J}\times\one{E}_{I\times J}
        \ar[r, "{\id\times\left((-)[\tpl{0,2}]\right)}"{yshift=1ex}]
        \ar[d, "{(-)[\tpl{0,0,1}]\times\id}"]
        \ar[dr, phantom, "\rotatebox{45}{$=$}", xshift=2em]
        \&
        \one{E}_{I\times J}\times\one{E}_{I\times I\times J}
        \ar[d, "{(-)[\tpl{0,0,1}]\times\id}"]
        \ar[r, "\land"]
        \ar[dr, phantom, "\rotatebox{45}{$\cong$}", xshift=2em]
        \&
        \one{E}_{I\times I\times J}
        \ar[d, "{(-)[\tpl{0,0,1}]}"]
        \\
        \!        
        \ar[ru, Rightarrow, shorten <= 6em, "\eta\times\id", very near end]
        \&
        \one{E}_{I\times J}\times\one{E}_{I\times J}
        \ar[r, "{\id\times\left((-)[\tpl{0,2}]\right)}"'{yshift=-1ex}]
        \&
        \one{E}_{I\times J}\times\one{E}_{I\times I\times J}
        \ar[r, "{\land\circ\left(\id\times (-)[\tpl{0,0,1}]\right)}"'{yshift=-1ex}]
        \&
        \one{E}_{I\times J}
    \end{tikzcd}\\
    =
    \hspace{2em}
    \begin{tikzcd}[ampersand replacement=\&]
        \one{E}_{I\times J}\times\one{E}_{I\times J}
        \ar[r, "{\id\times\left((-)[\tpl{0,2}]\right)}"{yshift=1ex}]
        \&
        \one{E}_{I\times J}\times\one{E}_{I\times I\times J}
        \ar[r, "\sum_{\tpl{0,0,1}}\times\id"{yshift=1ex}]
        \ar[dr,equal]
        \&
        \one{E}_{I\times I\times J}\times\one{E}_{I\times I\times J}
        \ar[d, "{(-)[\tpl{0,0,1}]\times\id}"]
        \ar[r, "\land"]
        \ar[dr, phantom, "\rotatebox{45}{$\cong$}", xshift=2em]
        \&
        \one{E}_{I\times I\times J}
        \ar[d, "{(-)[\tpl{0,0,1}]}"]
        \\
        \&
        \!
        \ar[ru, Rightarrow, shorten <= 6em, "\eta\times\id", very near end]
        \&
        \one{E}_{I\times J}\times\one{E}_{I\times I\times J}
        \ar[r, "{\land\circ\left(\id\times (-)[\tpl{0,0,1}]\right)}"'{yshift=-1ex}]
        \&
        \one{E}_{I\times J}
    \end{tikzcd}\\
    =
    \hspace{2em}
    \begin{tikzcd}[ampersand replacement=\&]
        \one{E}_{I\times J}\times\one{E}_{I\times J}
        \ar[r, "{\id\times\left((-)[\tpl{0,2}]\right)}"{yshift=1ex}]
        \&
        \one{E}_{I\times J}\times\one{E}_{I\times I\times J}
        \ar[r, "{\land\circ\left(\sum_{\tpl{0,0,1}}\times\id\right)}"{yshift=1ex}]
        \ar[d,"{\land\circ\left(\id\times (-)[\tpl{0,0,1}]\right)}"']
        \&
        \one{E}_{I\times I\times J}
        \ar[d,equal]
        \\
        \&
        \one{E}_{I\times J}
        \ar[r, "\sum_{\tpl{0,0,1}}"]
        \ar[rd,equal]
        \ar[ur, "{\rotatebox{45}{$\cong$}\,(\text{FR})}"{yshift=1ex},phantom]
        \&
        \one{E}_{I\times I\times J}
        \ar[d, "{(-)[\tpl{0,0,1}]}"]
        \\
        \&
        \!
        \ar[ru, Rightarrow, shorten <= 4.5em, "\eta", very near end]
        \&
        \one{E}_{I\times J}
    \end{tikzcd}
    \end{align*}
\end{proof}

\begin{lemma}
    \label{cor:eeffrobenius}
    Let $\mf{p}\colon\one{E}\to\one{B}$ be an elementary fibration.
    Then, it satisfies the Beck-Chevalley condition for the following pullback square in $\one{B}$
    in both directions:
    \[
    \begin{tikzcd}
        I\times J
        \ar[d, "\tpl{0,0,1}"']
        \ar[r, "\tpl{0,0,1}"]
        \ar[dr, phantom, "\lrcorner", very near start]
        &
        I\times I\times J
        \ar[d, "\tpl{0,0,1,2}"]
        \\
        I\times I\times J
        \ar[r, "\tpl{0,1,1,2}"']
        &
        I\times I\times I\times J
    \end{tikzcd}.
    \]
\end{lemma}
\begin{proof}
    We only prove it in one direction
    as the other direction is similar.
    Let $\alpha\in\one{E}_{I\times I\times J}$.
    The canonical arrow that we need to show to be an isomorphism is
    at the top of the following diagram:
    \[
    \begin{tikzcd}[row sep=small]
        \sum_{\tpl{0,0,1}}\left(\alpha[\tpl{0,0,1}]\right)
        \ar[r]
        \ar[dd, "\cong"']
        &
        \left(\sum_{\tpl{0,0,1,2}}\alpha\right)[\tpl{0,1,1,2}]
        \ar[d, "\cong"', "{\sigma[\tpl{0,1,1,2}]}"]
        \\
        &
        \left(\left(\sum_{\tpl{0,0}}\top\right)[\tpl{0,1}]
        \land\alpha[\tpl{1,2,3}]
        \right)[\tpl{0,1,1,2}]
        \ar[d, "\cong"']
        \\
        \sum_{\tpl{0,0,1}}\left(\top\land\alpha[\tpl{0,0,1}]\right)
        \ar[r, "\iota"', "\cong"]
        &
        \left(\sum_{\tpl{0,0,1}}\top\right)\land\alpha[\tpl{0,0,1}][\tpl{1,2}]
    \end{tikzcd},
    \]
    where $\iota$ is the isomorphism in \Cref{lem:implicationsoffrobele},
    but for the second isomorphism.
    The arrow $\sigma$ is the canonical isomorphsim following from the
    fact that the left adjoint of the reindexing functor $(-)[\tpl{0,1,1,2}]$
    is realized by the functor $\left(\sum_{\tpl{0,0}}\top\right)[\tpl{0,1}]\land(-)[\tpl{1,2,3}]$,
    as shown in \cite{EPR22}.
    Therefore, once the commutativity of the diagram is established,
    the Beck-Chevalley condition for the pullback square in the statement follows.
    However, the commutativity of the diagram is equivalent to that of the following diagram:
    \[
    \begin{tikzcd}
        \alpha[\tpl{0,0,1}]
        \ar[r, "{\eta_{\alpha}[0,0,1]}"]
        \ar[dd, "\cong"']
        \ar[rd, "{\eta_{\alpha}'[0,0,1]}"']
        &
        \left(\sum_{\tpl{0,0,1,2}}\alpha\right)[\tpl{0,1,1,2}][\tpl{0,0,1}]
        \ar[d, "\cong"', "{\sigma[\tpl{0,1,1,2}][\tpl{0,0,1}]}"]
        \\
        &
        \left(\left(\sum_{\tpl{0,0}}\top\right)[\tpl{0,1}]
        \land\alpha[\tpl{1,2,3}]
        \right)[\tpl{0,1,1,2}][\tpl{0,0,1}]
        \ar[d, "\cong"']
        \\
        \top\land\alpha[\tpl{0,0,1}]
        \ar[r, "\eta_{\top}\land\id_{\alpha[\tpl{0,0,1}]}"]
        &
        \left(\sum_{\tpl{0,0,1}}\top\right)[\tpl{0,0,1}]\land\alpha[\tpl{0,0,1}]
    \end{tikzcd},
    \]
    where $\eta_{\alpha}'$ is the unit of the adjunction.
    The triangle is commutative by the uniqueness of the left adjoints up to isomorphism,
    and the square is commutative because 
    the unit $\eta_{\alpha}'$ is 
    isomorphic to $\eta_{\top}\land\id_{\alpha[\tpl{0,0,1}]}$.
    Therefore, the Beck-Chevalley condition for the pullback square in the statement holds.
\end{proof}

\begin{definition}
    \label{def:exisfibration}
    Let $\mf{p}\colon\one{E}\to\one{B}$ be a cartesian fibration.
    We define a functor $\Phi_{\lc+}\colon\ob\one{B}\times\one{B}\to\PB[\one{B}]$ that assigns to each pair $(I,J)$ the arrow $\tpl{1}\colon I\times J\to J$, and assigns to each arrow $f\colon (I,J)\to (I,J')$ (an arrow $f\colon J\to J'$ in $\one{B}$)
    the pullback square on the left below.
    \[
    \begin{tikzcd}
        I\times J
        \ar[d, "\tpl{1}"']
        \ar[r, "\id\times f"]
        \ar[dr, phantom, "\lrcorner", very near start]
        &
        I\times J'
        \ar[d, "\tpl{1}"]
        \\
        J
        \ar[r, "f"']
        &
        J'
    \end{tikzcd}
    \hspace{2em}
    \begin{tikzcd}[column sep=large]
        \one{E}_{I\times J}
        \ar[d, "\sum_{\tpl{1}}"']
        \ar[from =r, "(\id\times f)^*"']
        \ar[rd, Rightarrow, shorten <= 1em, shorten >= 1em]
        &
        \one{E}_{I\times J'}
        \ar[d, "\sum_{\tpl{1}}"]
        \\
        \one{E}_{J}
        \ar[from = r, "f^*"]
        &
        \one{E}_{J'}
    \end{tikzcd}
    \]
    We say that $\mf{p}$ is an \emph{existential fibration}
    if it has $\Phi_{\lc+}$-coproducts and satisfies the Beck-Chevalley condition and Frobenius reciprocity for $\Phi_{\lc+}$.
\end{definition}

\begin{remark}
The left adjoints to the reindexing functors along the product projections satisfying the Beck-Chevalley condition
are called \emph{simple coproducts} in \cite[1.9.1]{jacobsCategoricalLogicType1999a}.
We borrow the adjective ``existential'' from its doctrine counterpart 
called \emph{existential doctrine} \cite[Definition 2.11]{MR13} (\cite[Definition 3.3]{Trotta20}).
\end{remark}

\begin{definition}
    \label{def:elexifibration}
    A cartesian fibration $\mf{p}\colon\one{E}\to\one{B}$ is called an \emph{elementary existential fibration}
    if it is an elementary fibration and an existential fibration.
\end{definition}

\begin{remark}
The doctrine counterpart for this was introduced by Lawvere in \cite{lawvereEqualityHyperdoctrinesComprehension1970} under the name \emph{elementary existential doctrine}.
An elementary existential fibration with preordered fibers is called a \it{regular fibration} in \cite[Definition 4.2.1]{jacobsCategoricalLogicType1999a}. 
\end{remark}

\begin{remark}
    A well-known result (see \cite[Examples 4.3.7]{jacobsCategoricalLogicType1999a}) states that if a fibration $\mf{p}\colon\one{E}\to\one{B}$ is an elementary existential fibration, then any base change functor $(-)[f]$ has a left adjoint $\sum_f$.
    However, the Beck-Chevalley condition does not necessarily hold in general.
\end{remark}

\begin{definition}[{\cite[Definition 2.1]{Pav96}}]
    \label{def:regularfibration}
    A \emph{regular fibration} is a cartesian fibration $\mf{p}\colon\one{E}\to\one{B}$ such that the base category $\one{B}$ has finite limits and the fibration $\mf{p}$ has $\Idf_{\PB[\one{B}]}$-coproducts and satisfies the Beck-Chevalley condition and Frobenius reciprocity for $\Idf_{\PB[\one{B}]}$. 
\end{definition}

\begin{remark}
    \label{rem:regularfibrationandeed}
    Frobenius reciprocity for $\Idf_{\PB[\one{B}]}$ actually follows from the Beck-Chevalley condition for all pullback squares in $\one{B}$ and Frobenius reciprocity for $\Phi_{\lc=+}$. 
    This can be proven purely in terms of fibrations, but we will give a double-categorical proof later in \Cref{cor:FrobeniusfromBeckChevalley}.
    Therefore, a regular fibration is an elementary existential fibration that satisfies the Beck-Chevalley condition for all pullback squares in the base category.
\end{remark}

\begin{example}
    \label{ex:regularfibration}
    The codomain fibration $\one{B}^{\rightarrow}\to\one{B}$ for a category $\one{B}$ with finite limits is a regular fibration. 
    The Beck-Chevalley conditions for all pullback squares in $\one{B}$ are satisfied by virtue of \textit{the pullback pasting lemma}.
\end{example}

We have a complete characterization of when subobject fibrations and family fibrations are elementary existential fibrations and regular fibrations.

\begin{proposition}[{\cite[Examples 3.4.4,Theorem 4.4.4]{jacobsCategoricalLogicType1999a}}]
    \label{prop:subobjectfibration}
    Let $\one{B}$ be a category with finite limit.
    The subobject fibration $\Sub{\one{B}}\to\one{B}$ is always an elementary fibration.
    In addition, the following are equivalent.
        \begin{enumerate}
            \item $\one{B}$ is a regular category.
            \item The subobject fibration over $\one{B}$ is an elementary existential fibration.
            \item The subobject fibration over $\one{B}$ is a regular fibration.
            \qedhere
        \end{enumerate}
\end{proposition}

\begin{proposition}
    \label{prop:familyfibration}
    Let $\one{B}$ be a category with finite products and $\Fam{\one{B}}\to\Set$ be its family fibration.
    \begin{enumerate}
        \item The family fibration is an elementary fibration if $\one{B}$ has distributive initial objects.
        \item The following are equivalent.
        \begin{enumerate}
            \item\label{fam-1} $\one{B}$ has distributive small coproducts.
            \item The family fibration over $\one{B}$ is an elementary existential fibration.
            \item The family fibration over $\one{B}$ is a regular fibration.
            \qedhere
        \end{enumerate}
    \end{enumerate}
\end{proposition}
\begin{proof}
A proof of (i) is obtained from part of (ii), so we give a proof of (ii) here.
First, we show that (a) implies (c).
The proof for the existence of the left adjoint of reindexing functors and the Frobenius reciprocity
is given in \cite[Example 3.4.4 (iii)]{jacobsCategoricalLogicType1999a}.
The left adjoint of the reindexing functor $(-)[f]$ for a function $f\colon I\to J$ is given by
\[
\sum_f\left((\alpha_i)_{i\in I}\right) = \left(\sum_{i\in f\inv(j)}\alpha_i\right)_{j\in J},
\]
and the Beck-Chevalley condition for all pullback squares in $\Set$ follows directly from this presentation. 
Evidently, (c) implies (b). 
To show that (b) implies (a), the coproduct of a family of objects $\left(\alpha_i\right)_{i\in I}$ 
is achieved by the left adjoint of the reindexing functor for the function $!_I\colon I\to 1$.
The Frobenius reciprocity for this reindexing functor guarantees the distributivity of the coproducts.
\end{proof}

\begin{remark}
    \label{rem:terminology}
    The terms \textit{regular fibration} and \textit{elementary existential fibration} may be used in a different sense in the literature, sometimes interchangeably.
    Our terminology is based on our desire to consider the term \textit{elementary existential fibration} as a conjunction of the elementary and existential fibrations, and the term \textit{regular fibration} as a fibration with sufficiently similar properties that the subobject fibration of a regular category has. 
\end{remark}

The 2-category of the aforementioned classes of doctrines is given in \cite{MR13,MREQC13} in the style of indexed categories.
Here, we give the 2-category of the corresponding classes of fibrations.
\begin{definition}
    \label{def:2catfibration}
    Let $\mf{p}\colon\one{E}\to\one{B}$ and $\mf{p}'\colon\one{E}'\to\one{B}'$ be fibrations, and let $\mf{f}=(\mf{f}_{0},\mf{f}_1)\colon\mf{p}\to\mf{p}'$ be a morphism of fibrations.
    Here, $\mf{f}_{0}\colon\one{B}\to\one{B}'$ denotes a functor between the base categories, and $\mf{f}_1\colon\one{E}\to\one{E}'$ denotes a functor over $\mf{f}_{0}$.
    \begin{enumerate}
        \item Suppose that $\mf{p},\mf{p}'$ are cartesian fibrations.
        We say that $(\mf{f}_{0},\mf{f}_1)$ is a \emph{morphism of cartesian fibrations} if $\mf{f}_{0}$ and $\mf{f}_1$ preserve finite products.
        This is equivalent to saying that $\mf{f}_{0}$ preserves finite products and for any object $I\in\one{B}$,
        $\mf{f}_{1,I}\colon\one{E}_{I}\to\one{E}'_{\mf{f}_{0}(I)}$ preserves finite products.
        \item Suppose that we have functorial choices of pullback squares $\Phi\colon\one{C}\to\PB[\one{B}]$ and $\Phi'\colon\one{C}'\to\PB[\one{B}']$,
        that $\mf{f}_{0}$ preserves the pullback squares that arise in the image of $\Phi$,
        and that we have a functor $T\colon\one{C}\to\one{C}'$ such that $\Phi'\circ T$ is naturally isomorphic to $\PB[\mf{f}_{0}]\circ\Phi$, 
        where $\PB[\mf{f}_{0}]\circ\Phi$ makes sense by the second condition.
        Furthermore, $\mf{p}$ and $\mf{p}'$ have $\Phi$-coproducts and $\Phi'$-coproducts, respectively.
        We say that $(\mf{f}_{0},\mf{f}_1)$ \emph{sends $\Phi$-coproducts to $\Phi'$-coproducts} if for any object $c\in\one{C}$,
        the following canonical natural transformation is an isomorphism\footnote{
        We identify $KDc$ with $D'_{Tc}$ and $KCc$ with $C'_{Tc}$
        by the natural isomorphism $\PB[\mf{f}_{0}]\circ\Phi\cong\Phi'\circ T$
        since this identification does not affect the condition.
        }:
        \[
        \begin{tikzcd}
            \one{E}_{D_c}
            \ar[r, "F_{D_c}"]
            \ar[d, "\sum_{\Phi_c}"']
            &
            \one{E}'_{D'_{Tc}}
            \ar[d, "\sum_{\Phi'_{Tc}}"]
            \ar[ld, Rightarrow, shorten <= 1em, shorten >= 1em]
            \\
            \one{E}_{C_c}
            \ar[r, "F_{C_c}"']
            &
            \one{E}'_{C'_{Tc}}
        \end{tikzcd}.
        \]
        \item Supppose that $\mf{p},\mf{p}'$ are elementary fibrations.
        A morphism$(\mf{f}_{0},\mf{f}_1)$ is a \emph{morphism of elementary fibrations} if it is a morphism of cartesian fibrations
        and sends $\Phi_{\lc=}$-coproducts in $\mf{p}$ to $\Phi_{\lc=}$-coproducts in $\mf{p}'$.
        Note that when $\mf{f}_{0}$ preserves finite products, the pullback squares that arise in the image of $\Phi_{\lc=}$ are preserved by $\mf{f}_{0}$.
        \item Suppose that $\mf{p},\mf{p}'$ are existential fibrations,
        A morphism $(\mf{f}_{0},\mf{f}_1)$ is a \emph{morphism of existential fibrations} if it is a morphism of cartesian fibrations
        and sends $\Phi_{\lc+}$-coproducts in $\mf{p}$ to $\Phi_{\lc+}$-coproducts in $\mf{p}'$.
        The same note as above applies.
        \item Suppose that $\mf{p},\mf{p}'$ are elementary existential fibrations,
        we say that $(\mf{f}_{0},\mf{f}_1)$ is a \emph{morphism of elementary existential fibrations} if it is a morphism of elementary fibrations and existential fibrations.  
        \qedhere
    \end{enumerate}
\end{definition}

From now on, we will omit the indices $0$ and $1$ in the notation of morphisms of fibrations. 

\begin{definition}
    \label{def:2catfibration2}
    We define the 2-category $\Fib\carttwo$ (resp. $\Fib_\eltwo$, $\Fib_\exitwo$, $\Fib_{\eef}$) of cartesian fibrations (resp. elementary fibrations, existential fibrations, elementary existential fibrations) as follows:
    \begin{enumerate}
        \item The objects are cartesian fibrations (resp. elementary fibrations, existential fibrations, elementary existential fibrations).
        \item The morphisms are morphisms of cartesian fibrations (resp. elementary fibrations, existential fibrations, elementary existential fibrations).
        \item The 2-cells are natural transformations between morphisms of fibrations.
        \qedhere
    \end{enumerate}
\end{definition}
Note that the way we define the 2-category of cartesian fibrations gives exactly the same 2-category as the 2-category of cartesian objects (c.f. \Cref{def:cartesianobj}) in $\Fib$ because of \Cref{prop:cartesianfibration}.

Let $\BiFib$ be the 2-category of bifibrations, fibered functors preserving the left adjoints of reindexing functors, and arbitrary fibered natural transformations.

\begin{lemma}
    \label{lem:eeftobifib}
    An elementary existential fibration is a bifibration.
    This gives rise to a fully faithful 
    2-functor $\Fib_{\eef}\to\BiFib$ that sends an elementary existential fibration to its associated bifibration.
\end{lemma}
\begin{proof}
    The first statement is a classical result as explained in \Cref{rem:regularfibrationandeed}.
    Since the left adjoint of $(-)[f]\colon\one{E}_{J}\to\one{E}_{I}$ is achieved by 
    \[
    \begin{tikzcd}[column sep=large]
        \one{E}_{I}
        \ar[r, "{\left\langle\delta_I,(-)[\tpl{0}]\right\rangle}"]
        &
        \one{E}_{I\times I}\times\one{E}_{I\times J}
        \ar[r, "{(-)[\id\times f]\times\Idf}"]
        &
        \one{E}_{I\times J}\times\one{E}_{I\times J}
        \ar[r, "\land"]
        &
        \one{E}_{I\times J}
        \ar[r,"\tpl{1}"]
        &
        \one{E}_J
    \end{tikzcd},
    \]
    all of which are preserved by a morphism of elementary existential fibrations, the 2-functor is well-defined.
    This also implies that the condition for a fibered functor to be a 1-cell in $\BiFib$ and in $\Fib_{\eef}$ are the same, that is, the preservation of the left adjoints of reindexing functors for all arrows in the base category,
    and the 2-functor is fully faithful.
\end{proof}

\subsection{Double categories}
\label{subsec:doublevirt}
We will use basic concepts of double categories in this paper.
We refer the reader to \cite{hoshinoDoubleCategoriesRelations2023} for the basic definitions; see also \cite{Grandis20} for a comprehensive account.

\begin{remark}[Terminology]
    \label{rem:terminologydoublecategory}
    We mostly follow the terminology of \cite{hoshinoDoubleCategoriesRelations2023} for double categories with the following exceptions.
    \begin{enumerate}
        \item For a (pseudo) double category, we use the term \emph{tight arrows} for the arrows that compose strictly associatively and unitally, and \emph{loose arrows} for the arrows that compose associatively and unitally up to coherent isomorphisms.
        Correspondingly, the opposites of a double category are called the \emph{tight opposite} and the \emph{loose opposite}, denoted by $\dbl{D}\tiop$ and $\dbl{D}\lop$, respectively.
        The bicategory of objects, loose arrows, and \emph{globular cells}, meaning cells whose source and target tight arrows are identities, is called the \emph{loose bicategory} of the double category $\dbl{D}$, and denoted by $\LBi{\dbl{D}}$.
        The 2-category of objects, tight arrows, and cells whose source and target loose arrows are identities is called the \emph{tight 2-category} of the double category $\dbl{D}$, and denoted by $\TBi{\dbl{D}}$.
        The loose identity arrow on an object $I$ is denoted by $\delta_I$.
        \item An \emph{equipment} is a double category $\dbl{D}$ such that the functor $\src\times\tgt\colon\dbl{D}_1\to\dbl{D}_0\times\dbl{D}_0$ is a fibration, or equivalently, an opfibration.
        We use the term \emph{prone cells} and \emph{supine cells} for the cells that are prone and supine lifts of a pair of tight arrows with respect to that functor.
        The top loose arrow (resp. bottom loose arrow) of a prone (resp. supine) cell is called \emph{the restriction} (resp. \emph{the oprestriction}) of the bottom (resp. top) loose arrow along the tight arrows, and denoted by $\alpha[f,g]$ (resp. $\opr(\beta;f,g)$).
    The restriction of an identity loose arrow along a pair of tight arrows with the one on the left (resp. right) being an identity gives the way of turning a tight arrow into a loose arrow, and the resulting loose arrow is called the \emph{companion} (resp. \emph{conjoint}) of the tight arrow, denoted by $f_*$ (resp. $f^*$).
    A \emph{cartesian equipment} is a cartesian object in the 2-category of equipments, pseudo functors, and tight transformations; see \cite{aleiferiCartesianDoubleCategories2018} for details.
    \qedhere
    \end{enumerate}
\end{remark}

\begin{remark}[Diagrammatic notation]
	\label{rem:notationtriangle}
	We often use the convention that identity arrows are contracted to vertices of objects.
	They are also drawn as arrows with double lines like $\to<equal>$ and $\sto[equal]$.
	By an alignment of arrows, we mean the composite of them.
	The following exemplifies these conventions:
	\[
		\begin{tikzcd}[virtual,column sep=small]
		&
		I
		\ar[ld, equal]
		\ar[rd, "f"]
		\\
		I
		\sar[rr,"\alpha"']
		\ar[rr, phantom, "\varphi", shift left=2ex]
		&&
		J
		\end{tikzcd}
		\quad
		\coloneqq
		\quad
		\begin{tikzcd}[virtual,column sep=small]
			I
			\sar[r, "\delta_I"]
			\ar[d, "\id_I"']
			\ar[rd, "\varphi", phantom]
			&
			I
			\ar[d, "f"]
			\\
			I
			\sar[r, "\alpha"']
			&
			J
		\end{tikzcd},
		\qquad
		\begin{tikzcd}[virtual,column sep=small]
			I
			\ar[dr, "g"']
			\sar[r, "\alpha"]
			&
			J
			\sar[r,equal]
			\ar[d, phantom, "\varphi"]
			&
			J
			\ar[dl, "h"]
			\\
			&
			K
			&
			\!
		\end{tikzcd}
		\quad
		\coloneqq
		\quad
		\begin{tikzcd}[virtual]
			I
			\sar[r, "{\alpha\odot\delta_J(\cong\alpha)}"]
			\ar[d, "g"']
			\ar[rd, "\varphi"{description}, phantom]
			&
			J
			\ar[d, "h"]
			\\
			K
			\sar[r, "\delta_K"']
			&
			K
		\end{tikzcd}.
	\]
\end{remark}	
\begin{example}
	\label{ex:doublecat}
	We give some basic examples of double categories.
	\begin{enumerate}
		\item 	The double category $\Rel[\Set]$ 
		of relations between sets is defined as follows.
		Its tight category $\Rel[\Set]_0$ is the category $\Set$ of sets and functions.
		The loose arrows are binary relations between sets, \text{i.e.,} a loose arrow $\alpha\colon A\sto B$ is a subset of $A\times B$.
		There is at most one cell of the form 
        \[
			\begin{tikzcd}[virtual]
				A
				\ar[d, "f"']
				\sar[r,"\alpha"]
				\ar[rd, phantom,"\varphi"]
				&
                B
				\ar[d, "g"]
				\\
				C
				\sar[r,"\beta"']
				&
				D
			\end{tikzcd}
        \]
        and it exists if and only if for any $a\in A$ and $b\in B$ such that $(a,b) \in \alpha$, we have $(f(a),g(b))\in \beta$.
		The composite $\alpha\odot\beta$ of relations $\alpha\colon A\sto B$ and $\beta\colon B\sto C$ is defined as the following relation.
		For $a\in A$ and $c\in C$, we have $(a,c)\in \alpha\odot\beta$ if and only if there exists $b\in B$ such that $(a,b)\in \alpha$ and $(b,c)\in \beta$.
		The identity loose arrow $\delta_A$ is defined by the diagonal $\{\,(a,a)\mid a\in A\,\}$.
		This construction is generalized to relations in a regular category in \cite{lambertDoubleCategoriesRelations2022}.
		\item For a category $\one{C}$ with pullbacks, one can form the double category of spans in $\one{C}$, whose tight category is $\one{C}$, whose loose arrows are spans in $\one{C}$, that is, pairs of arrows with the same source, and whose cells are arrows from the vertex of the top loose arrow to the vertex of the bottom loose arrow.
		\[
			\begin{tikzcd}[virtual]
				I 
				\ar[d, "f"']
				\sar[r,"\alpha"]
				\ar[rd, phantom,"\varphi"]
				&
				J
				\ar[d, "g"]
				\\
				K
				\sar[r,"\beta"']
				&
				L
			\end{tikzcd}	
			\quad
			=
			\quad
			\begin{tikzcd}[virtual]
				&
				\abs{\alpha}
				\ar[ld, "\ell_\alpha"']
				\ar[rd, "r_\alpha"]
				\ar[d,"\varphi"]
				\\
				I
				\ar[d, "f"']
				\ar[r,phantom,"\circlearrowright"]
				&
				\abs{\beta}
				\ar[ld, "\ell_\beta"]
				\ar[rd, "r_\beta"']
				\ar[r,phantom,"\circlearrowright"]
				&
				J
				\ar[d, "g"]
				\\
				K
				&
				&
				L
			\end{tikzcd}
		\]
		The identity loose arrow $\delta_I$ is the span $(\id_I,\id_I)$,
		and the composition of spans is defined by pullback in $\one{C}$.
        \item The examples (i) and (ii) are generalized to the double category $\Rel[\one{C}][(\zero{E},\zero{M})]$ of $\zero{M}$-relations for a stable orthogonal factorization system $(\zero{E},\zero{M})$ on a category $\one{C}$ with finite limits in \cite[Definition 3.2.1]{hoshinoDoubleCategoriesRelations2023}.
        Here, $\zero{M}$-relations are spans $(\ell\colon R\to I, r\colon R\to J)$ such that $\tpl{\ell,r}\colon R\to I\times J$ is in $\zero{M}$.
        \qedhere
	\end{enumerate}
\end{example}

In the joint work \cite{hoshinoDoubleCategoriesRelations2023}, the author and Hoshino made the following small observation, which turns out to be a convenient and powerful tool for calcuations in equipments.

\begin{lemma}[Sandwich Lemma, {\cite[Lemma 2.1.8]{hoshinoDoubleCategoriesRelations2023}}]
	\label{lem:Sandwich}
	Let $\dbl{D}$ be an equipment.
	Given a sequence of loosewise composable cells
	\begin{equation}
		\label{dgm:cartopcartcart}
	\begin{tikzcd}
		\cdot 
		\sar[r] 
		\ar[d] 
		\doublecell{\cart}
		&
		\cdot 
		\sar[r] 
		\ar[d] 
		\doublecell{\opcart}
		& 
		\cdot 
		\sar[r] 
		\ar[d] 
		\doublecell{\cart}
		& 
		\cdot 
		\ar[d] 
		\\
		\cdot 
		\sar[r] 
		& 
		\cdot 
		\sar[r] 
		& 
		\cdot 
		\sar[r] 
		& 
		\cdot 
	\end{tikzcd}
	\end{equation}
	with the supine cell sandwiched between two prone cells,
	the composition of these cells is prone.
	The same thing holds when swapping the roles of `prone' and `supine'.
\end{lemma}

The notion of \textit{Beck-Chevalley pullback squares} in a double category plays a fundamental role in this paper, as well as in \cite{hoshinoDoubleCategoriesRelations2023}.
It is a double categorical counterpart of the \text{Beck-Chevalley condition} in a bicategory \cite[\S 2.4]{WW08}.

\begin{definition}[{\cite[Definition 3.1.1]{hoshinoDoubleCategoriesRelations2023}}]
	\label{defn:BeckChevalley}
	Let $\dbl{D}$ be a double category.
	A \emph{diamond cell} in $\dbl{D}$ is a quadruple of tight arrows together with a tight cell $\alpha$ of the form on the left below.
	A diamond cell is called an \emph{identity diamond cell} if the tight cell is an identity cell.
	We say that a diamond cell satisfies the \emph{Beck-Chevalley condition} if there exists a loose arrow $\alpha\colon B\sto C$ and $\alpha$ factors as an opcartesian cell followed by a cartesian cell as shown in the right below:
	\begin{equation}
		\begin{tikzcd}[column sep=small]
				&
				I
				\ar[ld, "g"']
				\ar[rd, "f"]
				\doublecell[dd]{\alpha}
					&
			\\
			J
			\ar[rd, "h"']
				&
					&
					K
					\ar[ld, "k"]
			\\
				&
				L
					&
		\end{tikzcd}
		\hspace{1ex}
		=
		\hspace{1ex}
		\begin{tikzcd}[column sep=small]
				&
				I
				\ar[ld, "g"']
				\ar[rd, "f"]
					&
			\\
			J
			\sar[rr, "\alpha"{near end, description, inner sep=0.1mm}]
			\ar[rd, "h"']
			\doublecell[rr, shift left=2.2ex]{\opcart}
			\doublecell[rr, shift right=2.2ex]{\cart}
				&
					&
					K
					\ar[ld, "k"]
			\\
				&
				L
					&
		\end{tikzcd}
	\end{equation}
	Although this condition is defined for a diamond cell, we often abuse the terminology and say that a cell $\alpha$ satisfies the Beck-Chevalley condition when the quadruple of vertical arrows is evident from the context.
\end{definition}

\begin{definition}[{\cite[Definition 3.1.2]{hoshinoDoubleCategoriesRelations2023}}]
	A \emph{Beck-Chevalley pullback square} in $\dbl{D}$ is a pullback square in $\dbl{D}_0$
	as presented on the left below for which the two identity diamond cells placed in both directions
	as in the diagrams in the middle and right below
	satisfy the Beck-Chevalley condition.
	\[
	\begin{tikzcd}
		P 
		\ar[r, "s"]
		\ar[d, "t"']
		\pullback[rd]
		&
		I
		\ar[d, "f"]
		\\
		J
		\ar[r, "g"']
		&
		K
	\end{tikzcd}	
	\hspace{2ex}
	,
	\hspace{2ex}
	\begin{tikzcd}[column sep=small,virtual]
		&
		P
		\ar[ld, "s"']
		\ar[rd, "t"]
		\doublecell[dd]{=}
			&
	\\
	I
	\ar[rd, "f"']
		&
			&
			J
			\ar[ld, "g"]
	\\
		&
		K
			&
	\end{tikzcd}
	\hspace{2ex}
	,
	\hspace{2ex}
	\begin{tikzcd}[column sep=small,virtual]
		&
		P
		\ar[ld, "t"']
		\ar[rd, "s"]
		\doublecell[dd]{=}
			&
	\\
	J
	\ar[rd, "g"']
		&
			&
			I
			\ar[ld, "f"]
	\\
		&
		K
			&
	\end{tikzcd}
	\]
	We say a double category $\dbl{D}$ \emph{has Beck-Chevalley pullbacks} if the tight category $\dbl{D}_0$ has pullbacks and its pullback squares are all Beck-Chevalley pullback squares.
\end{definition}

\begin{example}
    \label{ex:BeckChevalley}\ 
    \begin{enumerate}
    \item For a category with finite limits $\one{C}$, the double category of spans $\Span[\one{C}]$ has Beck-Chevalley pullbacks \cite[Example 5.2.2]{aleiferiCartesianDoubleCategories2018}.
    \item For a regular category $\one{C}$, the double category of relations $\Rel[\one{C}]$ has Beck-Chevalley pullbacks.
    \item More generally, the double category of $\zero{M}$-relations $\Rel[\one{C}][(\zero{E},\zero{M})]$ for a stable orthogonal factorization system $(\zero{E},\zero{M})$ on a category $\one{C}$ with finite limits has Beck-Chevalley pullbacks \cite[Proposition 3.2.5]{hoshinoDoubleCategoriesRelations2023}.
    \qedhere
    \end{enumerate}
\end{example}

\begin{lemma}
    \label{lem:BeckChevalleyclosure}
    Let $\dbl{D}$ be a cartesian equipment.
    \begin{enumerate}
        \item The pullback of an identity arrow along any arrow gives a Beck-Chevalley pullback square in $\dbl{D}$.
        \item Beck-Chavalley pullback squares in $\dbl{D}$ are closed under finite products.
		\item Beck-Chevalley pullback squares in $\dbl{D}$ are closed under pasting.
	\end{enumerate}
    \[
    \text{(i)}\ 
    \begin{tikzcd}
        I
        \ar[r, "f"]
        \ar[d, "\id_I"']
        \ar[rd, phantom, "\text{(BC)}"]
        &
        K
        \ar[d, "\id_K"]
        \\
        I
        \ar[r, "f"']
        &
        K
    \end{tikzcd} 
    \hspace{2ex}
    \text{(ii)}\ 
    \begin{tikzcd}
        I
        \ar[r, "f"]
        \ar[d, "g"']
        \ar[rd, phantom, "\text{(BC)}"]
        &
        K
        \ar[d, "h"]
        \\
        J
        \ar[r, "k"']
        &
        L
    \end{tikzcd}
    ,
    \begin{tikzcd}
        I'
        \ar[r, "f'"]
        \ar[d, "g'"']
        \ar[rd, phantom, "\text{(BC)}"]
        &
        K'
        \ar[d, "h'"]
        \\
        J'
        \ar[r, "k'"']
        &
        L'
    \end{tikzcd}
    \Rightarrow
    \begin{tikzcd}
        I\times I'
        \ar[r, "f\times f'"]
        \ar[d, "g\times g'"']
        \ar[rd, phantom, "\text{(BC)}"]
        &
        K\times K'
        \ar[d, "h\times h'"]
        \\
        J\times J'
        \ar[r, "k\times k'"']
        &
        L\times L'
    \end{tikzcd}
    \]
	\[
	\text{(iii)}\
	\begin{tikzcd}
		I
		\ar[r, "f"]
		\ar[d, "g"']
		\ar[rd, phantom, "\text{(BC)}"]
		&
		K
		\ar[d, "n"]
		\\
		J
		\ar[r, "k"']
		&
		L
	\end{tikzcd}
	,
	\begin{tikzcd}
		K 
		\ar[r, "h"]
		\ar[d, "n"']
		\ar[rd, phantom, "\text{(BC)}"]
		&
		M
		\ar[d, "l"]
		\\
		L
		\ar[r, "m"']
		&
		N
	\end{tikzcd}
	\Rightarrow
	\begin{tikzcd}
		I
		\ar[r, "f"]
		\ar[d, "g"']
		\ar[rrd, phantom, "\text{(BC)}"]
		&
		K
		\ar[r, "h"]
		&
		M
		\ar[d, "l"]
		\\
		J
		\ar[r, "k"']
		&
		L
		\ar[r, "m"']
		&
		N
	\end{tikzcd}
	\]
\end{lemma}
\begin{proof}
	\ 
    \begin{enumerate}
        \item The two identity diamond cells for this pullback square 
        are given by the companion and the conjoint of the arrow $f$.
        \item Since every double functor preserves prone and supine cells,
        the Beck-Chevalley condition on diamond cells
        is preserved under the product functor $\times\colon\dbl{D}\times\dbl{D}\to\dbl{D}$.
		\item In the following diagram, the big top triangle is a supine cell
		and the big bottom triangle is a prone cell
		because of the sandwich lemma \Cref{lem:Sandwich}.
		\[
			\begin{tikzcd}[column sep=small, row sep=small]
				&
				&
				I
				\ar[ld,equal]
				\ar[rd, "f"]
				\ar[d, phantom, "\opcart"]
				\\
				&
				I
				\sar[rr]
				\ar[ld, "g"']
				\ar[rd, "f"{description}]
				\ar[d, phantom, "\opcart"]
				&
				\!
				\ar[d, phantom, "\cart"]
				&
				K
				\ar[ld, equal]
				\ar[rd, "h"]
				\ar[d, phantom, "\opcart"]
				\\
				J
				\sar[rr]
				\ar[rd, "h"']
				&
				\!
				\ar[d, phantom, "\cart"]
				&
				K
				\sar[rr]
				\ar[ld, "n"{description}]
				\ar[rd, "h"{description}]
				\ar[d, phantom, "\opcart"]
				&
				\!
				\ar[d, phantom, "\cart"]
				&
				M
				\ar[ld, equal]
				\\
				&
				L
				\ar[rd, "m"']
				\sar[rr]
				&
				\!
				\ar[d, phantom, "\cart"]
				&
				M
				\ar[ld, "l"]
				&
				\\
				&
				&
				N
			\end{tikzcd}
		\]
	\end{enumerate}
\end{proof}

\begin{lemma}
    \label{lem:BeckChevalleycell}
	Let $\dbl{D}$ be a cartesian equipment.
    Suppose we have
	a pullback square
    and a loose arrow as follows.
	\[
		\begin{tikzcd}[virtual]
			&
		I
			\ar[ld, "s"']
			\ar[rd, "t"]
			\ar[dd, phantom, very near start, "\rotatebox{-45}{$\lrcorner$}"description]
			&
		\\
		J
			\ar[rd, "f"']
			&
			&
		K
			\ar[ld, "g"]
			\sar[r, "\alpha"]
			&
		M
		\\
			&
		L
			&
		\end{tikzcd}
	\]
	Then the canonical cell $\sigma$ on the right below is an isomorphism if the pullback square is a Beck-Chevalley pullback square.
	\[
		\begin{tikzcd}[row sep=1.5em]
		I
			\ar[d, "t"']
			\sar[r, "t_*\alpha"]
			\doublecell[rd]{\cart}
				&
		M
			\ar[d, equal]
		\\
		K 
			\ar[d, "g"']
			\sar[r, "\alpha"]
			\doublecell[rd]{\opcart}
				&
		M
			\ar[d, equal]
		\\
		L
			\sar[r, "g^*\alpha"']
				&
		M
		\end{tikzcd}
		=
		\hspace{2ex}
		=
		\begin{tikzcd}[row sep=1.5em]
			I
				\ar[d, "s"']
				\sar[r, "t_*\alpha"]
					&
			M
				\ar[d, equal]
			\\
			J
				\ar[d, "f"']
				\sar[r, bend left = 25, "s^*t_*\alpha"]
                \sar[r, bend right = 25, "f_*g^*\alpha"']
				\ar[r, phantom, "\cart"{yshift=-5.5ex}]
                \ar[r, phantom, "\opcart"{yshift=4.5ex}]
                \ar[r, phantom, "\sigma"]
					&
			M
				\ar[d, equal]
			\\
			L
				\sar[r, "g^*\alpha"']
					&
			M
			\end{tikzcd}
	\]
\end{lemma}
\begin{proof}
	Applying the sandwich lemma \Cref{lem:Sandwich} to the diagram below,
	we obtain the desired result.
	\[
		\begin{tikzcd}[virtual]
			&
		I
			\ar[ld, "s"']
			\ar[rd, "t"]
			\sar[rr, "t_*\alpha"]
			&
			&
		M 
			\ar[d,equal]
		\\
		J
			\ar[rd, "f"']
			\sar[rr]
			\ar[rr, phantom, "\opcart"{yshift=2ex}]
			\ar[rr, phantom, below, "\cart"{yshift=-2ex}]
			&
			&
		K
			\ar[ld, "g"]
			\sar[r, "\alpha"]
			\ar[r, phantom, "\cart"{yshift=2.5ex, xshift=-1.5ex}]
			\ar[r, phantom, "\opcart"{yshift=-2.5ex, xshift=-1.5ex}]
			&
		M
			\ar[d,equal]
		\\
			&
		L
			\sar[rr, "g^*\alpha"']
			&
			&
		M
		\end{tikzcd}
	\]
\end{proof}

\begin{lemma}[{\cite[Proposition 3.1.10]{hoshinoDoubleCategoriesRelations2023}}]
    \label{lem:Frobenius}
    Let $\dbl{D}$ be a cartesian equipment.
    Suppose we have the following data in a double category $\dbl{D}$.
    \[
        \begin{tikzcd}[virtual]
            I
            \sar[r, "\alpha"]
            \ar[d, "f"']
                &
                K
            \\
            J
            \sar[r, "\beta"']
                &
                K
        \end{tikzcd}
    \]
    Then the canonical cell $\sigma$ on the right below is an isomorphism
    if the pullback of $\tpl{0,0}$ and $f\times\id_K$, which is always given by the span $f$ and $\tpl{\id,f}$,
    is a Beck-Chevalley pullback square.
    \[
        \begin{tikzcd}[row sep=1.5em]
            I
            \sar[r, "\alpha\land f_*\beta"]
            \ar[d, "\mmbox{\langle\id, f\rangle}"']
                &
                K
                \ar[d, "\tpl{0,0}"]
            \\
            I\times J
            \sar[r, "\alpha\times \beta"]
            \ar[d, "f\times\id"']
            \ar[r, phantom, "\opcart"{yshift=-3ex}]
            \ar[r, phantom, "\cart"{yshift=3ex}]
                &
                K\times K
                \ar[d, equal]
            \\
            J\times J
            \sar[r, "f^*\alpha\times \beta"']
                &
                K\times K
        \end{tikzcd}
        \hspace{1ex}
        =
        \hspace{1ex}
        \begin{tikzcd}[row sep=1.5em]
            I
                \sar[r, "\alpha\land f_*\beta"]
                \ar[d, "f"']
                &
            K
                \ar[d, equal]
            \\
            J
                \ar[d, "\tpl{0,0}"']
                \sar[r, bend left = 20, "f^*(\alpha\land f_*\beta)"]
                \sar[r, bend right = 20, "f^*\alpha\land\beta"']
                \ar[r, phantom, "\cart"{yshift=-5ex}]
                \ar[r, phantom, "\opcart"{yshift=5ex}]
                \ar[r, phantom, "\sigma"]
                &
            K
                \ar[d, "\tpl{0,0}"]
            \\
            J\times J
                \sar[r, "f^*\alpha\times \beta"']
                &
            K\times K
        \end{tikzcd}
    \]
\end{lemma}

\subsection{Virtual double categories}
\label{sec:virtualdoublecategories}
We turn to the basic concepts of (fibrational) virtual double categories. 
Virtual double categories were first introduced by Burroni in \cite{Bur71} under the name of \textit{multicat\'egories} as an example of $T$-categories.
Since then, this concept has turned up in several papers under different names, such as \textit{$\mathbf{fc}$-multicategories} in \cite{leinsterHigherOperadsHigher2004}, or \textit{lax double categories} in \cite{dawsonPathsDoubleCategories2006}.
The more standard name ``virtual double categories'' nowadays was introduced by Cruttwell and Shulman (\cite{cruttwellUnifiedFrameworkGeneralized2010}).

\begin{definition}[{\cite[Definition 2.1]{cruttwellUnifiedFrameworkGeneralized2010}}]
    A \emph{\acl{VDC}} (\ac{VDC}) $\dbl{X}$ is a structure consisting of the following data and axioms:
    \begin{enumerate}
        \item A category $\dbl{X}_{\tightcat}$. Its objects are simply called \emph{objects}, and its arrows are called \emph{tight arrows}, which are depicted vertically in this paper.
        \item A class of \emph{loose arrows} $\dbl{X}(I,J)_0$ for each pair of objects $I,J \in \dbl{X}_{\tightcat}$. These arrows are depicted horizontally with slashes as $\alpha\colon I \sto J$.
        \item A class of \emph{(virtual) cells} 
        \begin{equation}
            \label{eq:cell1}
            \begin{tikzcd}[virtual]
                I_0
                \ar[d, "s"']
                \sar[r, "\alpha_1"]
                \ar[rrrd, phantom, "\mu"]
                & I_1
                \sar[r]
                & \cdots
                \sar[r, "\alpha_n"]
                & I_n
                \ar[d, "t"] \\
                J_0
                \sar[rrr, "\beta"']
                & & & J_1
            \end{tikzcd}
        \end{equation}
        for each dataset consisting of $n\geq 0$, objects $I_0, \ldots, I_n, J_0, J_1 \in \dbl{X}_{\tightcat}$, tight arrows $s\colon I_0 \to J_0$ and $t\colon I_n \to J_1$, and loose arrows $\alpha_1, \ldots, \alpha_n, \beta$ with source and target as shown in \Cref{eq:cell1}.
        To specify the number $n$ of loose arrows, we call the cell $\mu$ an \emph{$n$-ary cell}.
        We will write the finite sequence of loose arrows as $\ol\alpha = \alpha_1,\dots,\alpha_n$.
        When $s$ and $t$ are identities, $\mu$ is said to be \emph{globular} and let $\mu\colon \ol\alpha \Rightarrow \beta$ denote the cell.
        The class of globular cells $\ol\alpha\Rightarrow\beta$ is denoted by $\dbl{X}(\ol{I})(\ol\alpha,\beta)$ in which $\ol{I} = I_0,\dots,I_n$.
        \item A composition operation on cells that assigns to each family of cells
        \[
            \begin{tikzcd}[column sep=4em,virtual]
                I_{1,0}
                \ar[d, "s_0"']
                \sard[r, "\ol\alpha_1"]
                \ar[dr, phantom, "\mu_1"]
                & I_{1,m_1}
                \ar[d, "s_1"']
                \sard[r, "\ol\alpha_{2}"]
                \ar[dr, phantom, "\mu_2"]
                & I_{2,m_2}
                \ar[d, "s_2"']
                \sard[r]
                & \cdots
                \sard[r, "\ol\alpha_{n}"]
                \ar[dr, phantom, "\mu_n"]
                & I_{n,m_n}
                \ar[d, "s_n"] \\
                J_{0}
                \ar[d, "t_0"']
                \sar[r, "\beta_1"'] 
                \ar[drrrr, phantom, "\nu"]
                & J_{1}
                \sar[r, "\beta_2"']
                & J_{2}
                \sar[r]
                & \cdots
                \sar[r, "\beta_n"' ]  
                & J_{n}
                \ar[d, "t_1"] \\
                K_{0}
                \sar[rrrr, "\gamma"']
                & & & & K_{1}
            \end{tikzcd}
        \]
        a cell 
        \[
            \begin{tikzcd}[column sep=4em,virtual]
                I_{1,0}
                \ar[d, "s_0"']
                \sard[r, "\ol\alpha_1"]
                \ar[ddrrrr, phantom, "{\nu\{\mu_1,\dots,\mu_n\}}"]
                & I_{1,m_1}
                \sard[r, "\ol\alpha_{2}"]
                & I_{2,m_2}
                \sard[r]
                & \cdots
                \sard[r, "\ol\alpha_{n}"]
                & I_{n,m_n}
                \ar[d, "s_n"] \\
                J_{0}
                \ar[d, "t_0"']
                &&&& J_{n}
                \ar[d, "t_1"] \\
                K_{0}
                \sar[rrrr, "\gamma"']
                & & & & K_{1}
            \end{tikzcd},
        \]
        where the dashed arrows represent finite sequences of loose arrows.
        We will write the finite sequence of cells as $\ol\mu = \mu_1,\dots,\mu_n$.
        \item An identity cell for each loose arrow $\alpha\colon I \sto J$
        \[
            \begin{tikzcd}[virtual]
                I
                \ar[d, "\id_I"']
                \sar[r, "\alpha"]
                \ar[dr, phantom, "\id_\alpha"]
                & J
                \ar[d, "\id_J"] \\
                I
                \sar[r, "\alpha"']
                & J
            \end{tikzcd}.
        \]
        (Henceforth, we will just write $=$ for the identity tight arrows.)
        \item The composition of cells is associative and unital with respect to the identity cells.
        \qedhere
    \end{enumerate}
\end{definition}

\begin{remark}
    \label{remark:vdcat}
    As already mentioned, we will often use dashed horizontal arrows to represent sequences of loose arrows.
    Correspondingly, we will use the expression on the left below to represent a sequence of identity cells as on the right below:
    \[
        \begin{tikzcd}[virtual]
            I_0
            \ar[d, equal]
            \sar[rr, "\ol\alpha", dashed]
            \ar[rrd, phantom, "\rotatebox{90}{$=$}"]
            &&
            I_n
            \ar[d, equal]
            \\
            I_0
            \sar[rr, "\ol\alpha"', dashed]
            &&
            I_n
        \end{tikzcd}
        \quad \coloneqq \quad
        \begin{tikzcd}[virtual]
            I_0
            \ar[d, "\id_{I_0}"']
            \sar[r, "\alpha_1"]
            \ar[rd, phantom, "\id_{\alpha_1}"]
            & I_1
            \sar[r]
            \ar[d, "\id_{I_1}"]
            & \cdots
            \sar[r, "\alpha_n"]
            \ar[rd, phantom, "\id_{\alpha_n}"]
            & I_n
            \ar[d, "\id_{I_n}"] \\
            I_0
            \sar[r, "\alpha_1"']
            & I_1
            \sar[r]
            & \cdots
            \sar[r, "\alpha_n"']
            & I_n
        \end{tikzcd}
    \]
    
    We also note that a cell whose top sequence of loose arrows is 
    the empty sequence is depicted as a triangle:
    \[
        \begin{tikzcd}[virtual, column sep=small]
            &
            I 
            \ar[dl, "s"']
            \ar[dr, "t"]
            &
            \\
            J_0
            \sar[rr, "\beta"']
            \ar[rr, phantom, "\mu", shift left=1em]
            &
            & J_1
        \end{tikzcd}.
    \]
\end{remark}

\begin{example}
\label{example:dblcat}
A double category can be seen as a virtual double category in the following way.
A cell \cref{eq:cell1} is defined as a virtual cell 
\[
    \begin{tikzcd}[virtual, column sep=large]
        I_0
        \ar[d, "s"']
        \sar[r, "\alpha_1\odot\dots\odot\alpha_n"]
        \ar[rd, phantom, "\mu"]
        & I_n
        \ar[d, "t"] \\
        J_0
        \sar[r, "\beta"']
        & J_1
    \end{tikzcd}
\]
where $\odot$ is the composition of loose arrows in the double category. 
The composition of cells is given by first composing cells in the loose direction on each row and then composing the resulting cells in the tight direction.
\end{example}

\begin{definition}[{\cite[Definition 3.1]{cruttwellUnifiedFrameworkGeneralized2010}}]
    \label{def:vdfunc}
    A \emph{virtual double functor} $F\colon \dbl{X} \to \dbl{Y}$ between virtual double categories $\dbl{X}$ and $\dbl{Y}$ 
    consists of the following data and conditions:
    \begin{itemize}
        \item A functor $F_{\tightcat}\colon \dbl{X}_{\tightcat} \to \dbl{Y}_{\tightcat}$.
        \item A family of functions $F_{I,J}\colon \dbl{X}(I,J)_0 \to \dbl{Y}(F_{\tightcat}(I), F_{\tightcat}(J))_0$ for each pair of objects $I,J$ of $\dbl{X}$.
        \item A family of functions sending each cell $\mu$ of $\dbl{X}$ on the left below to a cell $F(\mu)$ of $\dbl{Y}$ on the right below:
        \begin{equation}
            \label{eq:vdfunc}
            \begin{tikzcd}[virtual, column sep=small]
                I_0
                \ar[d, "s_0"']
                \sar[r, "\alpha_1"]
                \ar[rrrd, phantom, "\mu"]
                & I_1
                \sar[r]
                & \cdots
                \sar[r, "\alpha_n"]
                & I_n
                \ar[d, "s_1"] \\
                J_0
                \sar[rrr, "\beta"']
                & & & J_1
            \end{tikzcd}
            \mapsto 
            \begin{tikzcd}[virtual, column sep=small]
                F_{\tightcat}(I_0)
                \ar[d, "F_{\tightcat}(s_0)"']
                \sar[r, "F_{I_0,I_1}(\alpha_1)"]
                \ar[rrrd, phantom, "F(\mu)"]
                & F_{\tightcat}(I_1)
                \sar[r]
                & \cdots
                \sar[r, "F_{I_{n-1},I_n}(\alpha_n)"]
                & F_{\tightcat}(I_n)
                \ar[d, "F_{\tightcat}(s_1)"] \\
                F_{\tightcat}(J_0)
                \sar[rrr, "F_{J_0,J_1}(\beta)"']
                & & & F_{\tightcat}(J_1)
            \end{tikzcd}
            .
        \end{equation}
        \item The identity cells are preserved.
        \item Composition of cells is preserved.
    \end{itemize}
    As usual, we will often omit the subscripts of the functor and functions $F_{\tightcat}$ and $F_1$.

    A \emph{tight transformation} $\theta\colon F \to G$ between virtual double functors $F, G\colon \dbl{X} \to \dbl{Y}$
    consists of the following data and conditions:
    \begin{itemize}
        \item A natural transformation $\theta_0\colon F_{\tightcat} \to G_{\tightcat}$.
        \item A cell $\theta_{1,\alpha}$ for each loose arrow $\alpha\colon I \sto J$ of $\dbl{X}$:
        \[
            \begin{tikzcd}[virtual]
                FI 
                \ar[d, "\theta_{0,I}"']
                \sar[r, "F\alpha"]
                \ar[dr, phantom, "\theta_{1,\alpha}"]
                & FJ
                \ar[d, "\theta_{0,J}"] \\
                GI
                \sar[r, "G\alpha"']
                & GJ
            \end{tikzcd}
        \]
        \item The naturality condition for cells:
        \[
            \begin{tikzcd}[virtual]
                FI_0
                \ar[d, "Fs_0"']
                \sard[r, "F\ol\alpha"]
                \ar[dr, phantom, "F\mu"]
                & FI_n
                \ar[d, "Fs_n"] \\
                FJ_0
                \ar[d, "\theta_{J_0}"']
                \sar[r, "F\beta"']
                \ar[dr, phantom, "\theta_\beta",yshift=-0.3em]
                & FJ_1
                \ar[d, "\theta_{J_1}"] \\
                GJ_0
                \sar[r, "G\beta"']
                & GJ_1
            \end{tikzcd}
            \quad = \quad
            \begin{tikzcd}[virtual]
                FI_0
                \ar[d, "\theta_{I_0}"']
                \sard[r, "F\ol\alpha"]
                \ar[dr, phantom, "\theta_{\ol\alpha}"]
                & FI_n
                \ar[d, "\theta_{I_n}"] \\
                GI_0
                \ar[d, "Gs_0"']
                \sard[r, "G\ol\alpha"']
                \ar[dr, phantom, "G\mu",yshift=-0.3em]
                & GI_n
                \ar[d, "Gs_n"] \\
                GJ_0
                \sar[r, "G\beta"']
                & GJ_1
            \end{tikzcd}.
        \]
    \end{itemize}
    We denote by $\VDbl$ the 2-category of virtual double categories, virtual double functors, and tight transformations.
\end{definition}

Next, we review the notion of cartesian fibrational virtual double categories, which generalizes cartesian equipments.

\begin{definition}[{\cite[Definition 7.1]{cruttwellUnifiedFrameworkGeneralized2010}}]
    \label{def:fibrational}
    Let $\dbl{X}$ be a virtual double category. 
    A \emph{restriction} of a loose arrow $\alpha\colon I \sto J$ 
    along a pair of tight arrows $s\colon I' \to I$ and $t\colon J' \to J$ is the loose arrow $\alpha[s, t]\colon I' \sto J'$ 
    equipped with a cell 
    \[
        \begin{tikzcd}[virtual]
            I'
            \ar[d, "s"']
            \sar[r, "{\alpha[s, t]}"]
            \ar[dr, phantom, "\restc"]
            & J'
            \ar[d, "t"] \\
            I
            \sar[r, "\alpha"']
            & J
        \end{tikzcd}
    \]
    with the following universal property: any cell $\mu$ of the form on the left below factors uniquely through the cell 
    $\restc$ as on the right below.
    \begin{equation}
        \label{eq:rest}
        \begin{tikzcd}[virtual]
            K 
            \ar[d, "u"']
            \sard[r, "\ol\beta"]
            \ar[ddr, phantom, "\mu"]
            & 
            L 
            \ar[d, "v"] \\
            I'
            \ar[d, "s"']
            &
            J'
            \ar[d, "t"] \\
            I
            \sar[r, "\alpha"']
            & J
        \end{tikzcd}
        \quad = \quad
        \begin{tikzcd}[virtual]
            K 
            \ar[d, "u"']
            \sard[r, "\ol\beta"]
            \ar[dr, phantom, "\widehat{\mu}"]
            &
            L
            \ar[d, "v"] \\
            I'
            \ar[d, "s"']
            \sar[r, "{\alpha[s, t]}"']
            \ar[dr, phantom, "{\restc}", yshift=-0.3em]
            & J'
            \ar[d, "t"] \\
            I
            \sar[r, "\alpha"']
            & J
        \end{tikzcd}
    \end{equation}
    In this case, we call the cell $\restc$ a \emph{restricting cell}.
    If the restrictions exist for all triples $(\alpha, s, t)$, then we say that $\dbl{X}$ is a \emph{\ac{FVDC}}\footnote{
    The term ``fibrational'' is not standard in the literature.
    If we follow the terminology of \cite{aleiferiCartesianDoubleCategories2018}, we should call this a \textit{fibrant virtual double category}, but we prefer to use the term because it has nothing to do with any model structure, at least \textit{a priori}.
    }

    A \emph{fibrational virtual double functor} $F\colon \dbl{X} \to \dbl{Y}$ between fibrational virtual double categories $\dbl{X}$ and $\dbl{Y}$ is a virtual double functor that preserves restrictions.
    We denote by $\FVDbl$ the 2-category of fibrational virtual double categories, fibrational virtual double functors, and tight transformations.
\end{definition}

\begin{example}
    \label{example:vdcat}
    A double category is an equipment if and only if it is fibrational as a virtual double category.
\end{example}

\begin{lemma}[{\cite[Lemma 2.5]{nasu2024internallogicvirtualdouble}}]
    \label{lemma:fibvdblequiv}
    A virtual double functor $F\colon \dbl{X} \to \dbl{Y}$ is an equivalence in ${\VDbl}$ if and only if
    \begin{enumerate}
        \labeleditem  the functor $F_{\tightcat}\colon \dbl{X}_{\tightcat} \to \dbl{Y}_{\tightcat}$ for $F$ is an equivalence of categories, \label{lem:fibvdblequiv1}
        \labeleditem \label{lem:fibvdblequiv2} for any loose arrow $\alpha\colon I \sto J$ in $\dbl{Y}$, there exists
        a loose arrow $\beta\colon I' \sto J'$ in $\dbl{X}$ and an isomorphism cell $\mu$ as below:
        \[
            \begin{tikzcd}[virtual]
                FI'
                \ar[d, "\rotatebox{90}{$\cong$}"']
                \sar[r, "F\beta"]
                \ar[dr, phantom, "\mu\ \rotatebox{90}{$\cong$}"]    
                & FJ'
                \ar[d, "\rotatebox{90}{$\cong$}"] \\
                I
                \sar[r, "\alpha"']
                & J
            \end{tikzcd},
            \quad\text{and}
        \]
        \labeleditem \label{lem:fibvdblequiv3} for any quadruple $(s,t,\ol\alpha,\beta)$,
        the function $F$ on the cells \cref{eq:vdfunc} is a bijection.
    \end{enumerate}
    A fibrational virtual double functor $F\colon \dbl{X} \to \dbl{Y}$ is an equivalence in $\FVDbl$ if and only if
    \Cref{lem:fibvdblequiv1},\Cref{lem:fibvdblequiv2}, and the special case of \Cref{lem:fibvdblequiv3} where $s$ and $t$ are identities are satisfied.
\end{lemma}

Next, we explicitly describe the notion of \ac{CFVDC}, although it is already implicitly defined because we have the 2-category of fibrational virtual double categories $\FVDbl$, which has strict finite products. 

\begin{proposition}[{\cite[Proposition 2.9]{nasu2024internallogicvirtualdouble}}]
    \label{prop:FibVDblCart}
    An \ac{FVDC} $\dbl{X}$ is cartesian if and only if the following conditions are satisfied:
    \begin{enumerate}
        \item $\dbl{X}_{\tightcat}$ has finite products.
        \item $\dbl{X}$ locally has finite products, that is, for each $I,J\in\dbl{X}_{\tightcat}$.
        \begin{enumerate}
        \item For any loose arrows
        $\alpha,\beta\colon I\sto J$ in $\dbl{X}$,
        there exists a loose arrow $\alpha\land\beta\colon I\sto J$ and globular cells $\pi_0\colon\alpha\land\beta\Rightarrow\alpha$, $\pi_1\colon\alpha\land\beta\Rightarrow\beta$ such that for any finite sequence of loose arrows $\ol\gamma$ where $\gamma_i\colon I_{i-1}\sto I_i$ for $1\leq i\leq n$ where $I_0=I$ and $I_n=J$, the function 
        \[ 
            \dbl{X}(\ol{I})(\ol\gamma,\alpha\land\beta) \to
            \dbl{X}(\ol{I})(\ol\gamma,\alpha)\times\dbl{X}(\ol{I})(\ol\gamma,\beta)\quad;
            \quad \mu\mapsto (\pi_0\circ\mu,\pi_1\circ\mu)
        \]
        is a bijection.
        \item There exists a loose arrow $\top\colon I\sto J$ such that $\dbl{X}(\ol{I})(\ol\gamma,\top)_0$ is
        a singleton for any finite sequence of loose arrows $\ol\gamma$;
        \end{enumerate}
        \item The local finite products are preserved by restrictions.
    \end{enumerate} 
    A morphism between cartesian \acp{FVDC} is a cartesian morphism if and only if the underlying tight functor preserves finite products
    and the morphism preserves local finite products.
\end{proposition}

\begin{remark}
    \Cref{prop:FibVDblCart} is a virtual counterpart of \cite[Proposition 4.12]{aleiferiCartesianDoubleCategories2018}.
    The third condition in \cref{prop:FibVDblCart} is necessary for \acp{FVDC} but not for equipments as in \cite{aleiferiCartesianDoubleCategories2018} since the latter has oprestrictions of loose arrows.
\end{remark}

\begin{example}
        \label{example:rel}
        Many examples of \acp{CFVDC} can be found in formal category theory \cite{cruttwellUnifiedFrameworkGeneralized2010,arkorFormalTheoryRelative2024}, but here we give examples with more logical flavor.
        \begin{enumerate}
            \item In \Cref{ex:doublecat}, we defined a double category $\Rel[\Set]$ of sets, functions, and relations, and mentioned that it is generalizable to $\Rel[\one{B}]$ for a regular category $\one{B}$.
            We can drop the regularity condition and define a virtual double category $\Rel[\one{B}]$ of objects, arrows, and internal relations, meaning subobjects of the product of two objects in a category with finite limits $\one{B}$.
            This is possible because without the regularity condition, we can interpret Horn sentences internally in $\one{B}$ when it has finite limits.
            More concretely, it is a virtual double category with at most one cell for each frame where a cell exists precisely when the corresponding Horn sentence is valid. 
            For instance, nullary and binary cells are respectively defined as follows:
            \begin{align*}
                \begin{tikzcd}[virtual, column sep=small, ampersand replacement=\&]
                    \&
                    I
                    \ar[ld, "s"']
                    \ar[rd, "t"]
                    \&
                    \\
                    J
                    \sar[rr, "\alpha"']
                    \ar[rr, phantom, "{\rotatebox{-90}{$\leq$}}"{yshift=1em}]
                    \&
                    \& K
                \end{tikzcd}
                &
                \Longleftrightarrow
                \qquad
                \begin{tikzcd}[column sep=small, ampersand replacement=\&]
                    I
                    \ar[d, "\tpl{0,0}"']
                    \ar[r, dashed, "\exists"]
                    \ar[rd, phantom, "\circlearrowleft"]
                    \& \alpha 
                    \ar[d, hook] \\
                    I\times I
                    \ar[r, "s\times t"']
                    \& J\times K
                \end{tikzcd},\\
                \begin{tikzcd}[virtual, column sep=small, ampersand replacement=\&]
                    I
                    \ar[d, "s"']
                    \sar[r, "\alpha"]
                    \ar[rrd, phantom, "{\rotatebox{-90}{$\leq$}}"]
                    \& J
                    \sar[r, "\beta"]
                    \& K
                    \ar[d, "t"] 
                    \\
                    L
                    \sar[rr, "\gamma"']
                    \& \& M
                \end{tikzcd}
                &
                \Longleftrightarrow
                \qquad
                \begin{tikzcd}[ampersand replacement=\&]
                    (\alpha\times K) \cap (I\times\beta)
                    \ar[d, hook]
                    \ar[r, dashed, "\exists"]
                    \ar[rd, phantom, "\circlearrowleft"]
                    \& \gamma
                    \ar[d, hook] \\
                    I\times J\times K
                    \ar[r, "{(s\times t)\circ\tpl{0,2}}"']
                    \& L\times M
                \end{tikzcd}
            \end{align*}
            A restriction of a relation along a pair of functions is given by the pullback
            of the relation along the product of the functions.
            This double category is
            a \ac{CFVDC}.
            \item For a monoidal category $\one{V}$, we can define a fibrational virtual double category
            $\Mat[\one{V}]$ as follows.
            Its tight category is $\Set$, and the loose arrows $I\sto J$ are
            matrices $\left(A_{i,j}\right)_{i\in I,j\in J}$ of objects in $\one{V}$.
            A cell of the form on the left below, for instance, is
            a family of morphisms in $\one{V}$ on the right below:
            \[
                \begin{tikzcd}[virtual]
                    I
                    \ar[d, "s"']
                    \sar[r, "\left(A_{i,j}\right)_{i,j}"]
                    \ar[rrd, phantom, "\mu"]
                    & J
                    \sar[r, "\left(B_{j,k}\right)_{j,k}"]
                    & K
                    \ar[d, "t"] 
                    \\
                    L
                    \sar[rr, "\left(C_{l,m}\right)_{l,m}"']
                    & & M
                \end{tikzcd}
                \quad\vline\!\vline\quad
                \left(\mu_{i,j,k}\colon A_{i,j}\otimes B_{j,k}\to C_{s(i),t(k)}\right)
                _{i,j,k}
            \]
            Defining general cells and composition of cells involves the monoidal structure of $\one{V}$.
            A restriction of a matrix along a pair of functions $s\colon I'\to I$ and $t\colon J'\to J$ is given by 
            the matrix $\left(A_{s(i),t(j)}\right)_{i\in I',j\in J'}$.
            It is a \ac{CFVDC} if $\one{V}$ is cartesian monoidal.
            \qedhere
        \end{enumerate}
\end{example}

It is fundamental that pseudo double categories are seen as virtual double categories with structure defined by a universal property, called \emph{composites}.
We now review the basic results based on \cite{dawsonPathsDoubleCategories2006,cruttwellUnifiedFrameworkGeneralized2010}.

\begin{definition}[{\cite[Definition 2.7]{dawsonPathsDoubleCategories2006},\cite[Definition 5.2]{cruttwellUnifiedFrameworkGeneralized2010}}]
    \label{def:composite}
    A \emph{composite} of a given sequence of loose arrows
    $\ol\alpha=\left(I_0\sto["\alpha_1"] I_{1}\sto\cdots\sto["\alpha_m"] I_m\right)$ in a virtual double category 
    is a loose arrow $\odot\ol\alpha$ from $I_0$ to $I_m$ equipped with a cell
    \[
        \begin{tikzcd}[virtual, column sep=small]
            I_0
            \sar[r, "\alpha_1"]
            \ar[d, equal]
            \ar[rrrd, phantom, "\varkappa_{\ol\alpha}"]
            &
            I_1
            \sar[r]
            &
            \cdots
            \sar[r, "\alpha_m"]
            &
            I_m
            \ar[d, equal]
            \\
            I_0
            \sar[rrr, "\odot\ol\alpha"']
            &&&
            I_m
        \end{tikzcd}
    \]
    with the following universal property:
    given any cell $\nu$ on the left below where $\ol\beta,\ol\beta'$ 
    are arbitrary sequences of loose arrows,
    it uniquely factors through the sequence of the identity cells with $\mu_{\ol\alpha}$ as on the right below.
    \begin{equation}
        \label{eq:composite}
        \begin{tikzcd}[virtual, column sep=small]
            J_0
            \sard[r, "\ol\beta"]
            \ar[d, "f"']
            \ar[phantom,rrrrd, "\nu"]
            &
            I_0
            \sard[rr, "\ol\alpha"]
            &&
            I_m
            \sard[r, "\ol\beta'"]
            &
            J'_{n'}
            \ar[d, "f'"]
            \\
            K
            \sar[rrrr, "\gamma"']
            &&&&
            K'
        \end{tikzcd}
=
        \begin{tikzcd}[virtual, column sep=small]
            J_0
            \sard[r, "\ol\beta"]
            \ar[d, equal]
            \ar[dr, phantom, description, "{\rotatebox{90}{=}}"]
            &
            I_0
            \ar[d, equal]
            \ar[drr,phantom, "\varkappa_{\ol\alpha}"]
            \sar[rr, "\ol\alpha",dashed]
            &
            &
            I_m 
            \sard[r, "\ol\beta'"]
            \ar[d, equal]
            \ar[dr, phantom, description, "{\rotatebox{90}{=}}"]
            &
            J'_{n'}
            \ar[d, equal]
            \\
            J_0
            \sard[r, "\ol\beta"']
            \ar[d, "f"']
            \ar[phantom,rrrrd, "\wt{\nu}",description, yshift=-1ex]
            &
            I_0
            \sar[rr, "\odot\ol\alpha"']
            &&
            I_m
            \sard[r, "\ol\beta'"']
            &
            J'_{n'}
            \ar[d, "f'"]
            \\
            K
            \sar[rrrr, "\gamma"']
            &&&&
            K'
        \end{tikzcd}
\end{equation}
We call the cell $\varkappa_{\ol\alpha}$ the \emph{composing cell} of $\ol\alpha$.
In particular, a composite of the empty sequence of loose arrows on $I$
is called a \emph{unit} on $I$ and denoted by $\delta_I$.

A virtual double functor is said to \emph{preserve the composite $\odot\ol\alpha$} 
if it sends the composing cell of $\ol\alpha$ to
the cell that exhibits the image of $\odot\ol\alpha$ 
as the composite of the images of $\ol\alpha$.
It is said to \emph{preserve composition} if it preserves all composites.
\end{definition}

For the purpose of this paper, we will separately discuss 
composability of sequences of loose arrows of positive length
and those of length zero.
\begin{definition}
    \label{def:composable}\ 
    \begin{enumerate}
    \item A \ac{VDC} is said to be \emph{unital} if it has composites of sequences of length zero, or equivalently, if it has units on all objects.
    In particular, a \emph{virtual equipment} is an \ac{FVDC} that is also unital.
    \item A \ac{VDC}is called \emph{positively representable} if it has composites of any sequence of loose arrows of positive length. 
    \item A \ac{VDC} is called \emph{representable} if it has composites of any sequence of loose arrows of any finite length, including zero.
    \qedhere
    \end{enumerate}
\end{definition}

The following proposition is a fundamental result in this context.
The proof can be found in \cite{Her00}, which is mentioned in \cite{dawsonPathsDoubleCategories2006}.
\begin{proposition}[{\cite[Proposition 2.8]{dawsonPathsDoubleCategories2006}}]
    A representable virtual double category is presented as a double category seen as a virtual double category in the way described in \Cref{example:dblcat}.
    More precisely, the 2-category of representable virtual double categories is biequivalent to the 2-category of pseudo double categories.
\end{proposition}

For this reason, we will treat representable virtual double categories as pseudo double categories in the following, and write $\Dbl$ for the 2-category of representable virtual double categories, virtual double functors preserving composites (up to isomorphisms), and tight transformations.

\begin{definition}
    \label{def:2catofvdbl}
    We define the following 2-categories of \acp{FVDC} with various structures.
    \begin{enumerate}
        \item $\VEq$ is the 2-category of virtual equipments, fibrational virtual double functors preserving units, and tight transformations.
        \item $\FVDblP$ is the 2-category of positively representable virtual double categories with restrictions, fibrational virtual double functors preserving composites of positive-length sequences of loose arrows, and tight transformations.
        \item $\Eqp$ is the 2-category of equipments (\text{i.e.,} representable virtual double categories with restrictions), fibrational virtual double functors preserving composites, and tight transformations.
        \qedhere
    \end{enumerate}
\end{definition}

\begin{remark}
    \label{rem:equipmentsasfVDC}
    It was shown in \cite[Theorem 7.24]{cruttwellUnifiedFrameworkGeneralized2010} 
    that a unital virtual double functor between virtual equipments
    automatically preserves restrictions.
    Therefore, the 2-category $\VEq$ can be seen as a full sub-2-category of the 2-category of unital virtual double categories, and the 2-category $\Eqp$ as a full sub-2-category of $\Dbl$.
\end{remark}

\begin{remark}
    \label{rem:composabilityinFVDC}
    When a \ac{VDC} is fibrational, the universal property of a composite of a sequence of loose arrows can be reduced to a simpler condition.
    In an \ac{FVDC} $\dbl{X}$, a globular cell $\varkappa$ is a composing cell if and only if it has the universal property that for any cells $\nu$ in \Cref{eq:composite} with $f$ and $f'$ being identities, there is a unique cell $\wt{\nu}$ that makes the equation hold.
\end{remark}

The 2-categories $\VEq$, $\FVDblP$ have strict finite products, so we can discuss cartesianness in these 2-categories in the sense of \Cref{def:cartesianobj}, which is unravelled in the following proposition.

\begin{proposition}[{\cite[Proposition A.2]{nasu2024internallogicvirtualdouble}}]
    \label{prop:cartesianunital}
    A unital \ac{FVDC} $\dbl{X}$ is cartesian in $\VEq$ if and only if:
    \begin{enumerate}
        \item $\dbl{X}$ is a cartesian \ac{FVDC},
        \item $\delta_1\cong\top_{1,1}$ in $\dbl{X}(1,1)$ canonically, and
        \item for any $I,J\in\dbl{X}$, $\delta_{I,J}\cong \delta_I\times \delta_J$ canonically in $\dbl{X}(I\times J,I\times J)$.
        \qedhere
    \end{enumerate}
\end{proposition}

\begin{proposition}[{\cite[Proposition A.3]{nasu2024internallogicvirtualdouble}}]
    \label{prop:cartesiancompose}
    Let $\FVDblP$ be the locally full sub-2-category of $\FVDbl$ spanned by the \acp{FVDC} with composites 
    of sequences of loose arrows of positive length and functors preserving those composites.
    A \ac{VDC} $\dbl{X}$ in $\FVDblP$ is cartesian in this 2-category if and only if
    \begin{enumerate}
        \item $\dbl{X}$ is a cartesian \ac{FVDC},
        \item $\top_{1,1}\odot\top_{1,1}\cong\top_{1,1}$ canonically in $\dbl{X}(1,1)$, and
        \item for any paths of positive length 
        \[
        \begin{tikzcd}
            I_0\sar["\alpha_1",r] & I_1\sar[r] & \dots\sar[r,"\alpha_n"] & I_n
        \end{tikzcd}  
        \quad
        \text{and}
        \quad
        \begin{tikzcd}
            J_0\sar["\beta_1",r] & J_1\sar[r] & \dots\sar[r,"\beta_n"] & J_n
        \end{tikzcd}
        \]
        in $\dbl{X}$,
        we have
        \[
        (\alpha_1\odot\dots\odot\alpha_n)\times(\beta_1\odot\dots\odot\beta_n)\cong(\alpha_1\times\beta_1)\odot\dots\odot(\alpha_n\times\beta_n)
        \]
        canonically in $ \dbl{X}(I_0\times J_0,I_n\times J_n)$.
    \qedhere
    \end{enumerate}
\end{proposition}

\begin{example}
    \label{example:composablity} \ 
    \begin{enumerate}
        \item $\Rel[\one{B}]$ is unital, as its unit on $I$ is given as the equality relation on $I$,
        namely, $\tpl{0,0}\colon I\to I\times I$.
        It is positively representable if and only if $\one{B}$ is regular.
        This follows from our main result \Cref{thm:firsteefiscartdouble}.
        \item $\Mat[\one{V}]$ is unital if $\one{V}$ has an initial object $0$ that is 
        preserved by the tensor product.
        In this case, the unit on $I$ is given by $\left(\delta_{i,j}\right)_{i,j\in I}$
        where $\delta_{i,j}$ is the monoidal unit of $\one{V}$ when $i=j$ and the initial object otherwise. 
        It is positively representable if and only if
        $\one{V}$ has small coproducts that are preserved by the tensor product.
        In this case, the composing cell of a sequence of matrices is given as follows:
        \[
        \begin{tikzcd}
            I
            \sar[r, "\left(A_{i,j}\right)_{i,j}"]
            & J
            \sar[r, "\left(B_{j,k}\right)_{j,k}"]
            & K
        \end{tikzcd};
        \qquad 
        \left(A_{i,j}\right)_{i,j}\odot\left(B_{j,k}\right)_{j,k}=
        \left(\sum_{j\in J}A_{i,j}\otimes B_{j,k}\right)_{i\in I,k\in K}.
        \]
        See \cite[Example 5.3]{cruttwellUnifiedFrameworkGeneralized2010}.
    \end{enumerate}
    See \Cref{ex:regularfibrationtoequip} for another explanation for representability of these virtual double categories. 
\end{example}

\section{Construction of a virtual double category from a cartesian fibration}
\label{sec:fromfibtoVDC}
In this section, we will construct a virtual double category out of a cartesian fibration, and show that left adjoints of base change functors along tight arrows in the fibration satisfying the Beck-Chevalley condition and Frobenius reciprocity induce composition of loose arrows in the resulting virtual double category that is compatible with the cartesian structure.
When we see objects in thefibers of a cartesian fibration as predicates, the loose arrows in the resulting virtual double category are the binary relations described by these predicates.

\begin{definition}
    \label{def:primfibtovdc}
    Let $\mf{\mathfrak{p}}\colon\one{E}\to \one{B}$ be a cartesian fibration.
    Then the following data form a \acs{VDC} $\Bil[\mf{p}]$\footnote{
        $\Bil$ stands for \textit{bilateral}, suggesting that this construction turns a fibration, which has unary predicates in its fibers, into a virtual double category, which has binary relations as its loose arrows.
    }:
    \begin{enumerate}
        \item The tight part of $\Bil[\mf{p}]$ is the category $\one{B}$.
        \item The loose arrows from $I$ to $J$ in $\Bil[\mf{p}]$ are objects $\alpha$ in $\one{E}$ over $I\times J$.
        \item The cells of the form
        \begin{equation}
            \label{eq:generalcell}
        \begin{tikzcd}
            I_0
            \ar[phantom,rrrd, "\xi"]
            \sar[r, "\alpha_1"]
            \ar[d, "s_0"']
            &
            I_1
            \sar[r]
            &
            \cdots
            \sar[r, "\alpha_{n}"]
            &
            I_n
            \ar[d, "s_1"]
            \\
            J_0
            \sar[rrr, "\beta"']
            &&
            &
            J_n
        \end{tikzcd}    
        \end{equation}
        in $\Bil[\mf{p}]$
        are arrows $\xi\colon\alpha_1[\tpl{0,1}]\land\dots\land\alpha_n[\tpl{n-1,n}]\to\beta[(s_0\times s_1)\circ\tpl{0,n}]$ 
        in $\one{E}_{I_0\times\dots\times I_n}$,
        where $\tpl{i,j}$ denotes the pairing of the $i$-th and $j$-th projections $I_0\times\dots\times I_n\to I_i\times I_j$. 
        This is equivalent to the data of an arrow 
        $\xi\colon\bigwedge_{1\le i\le n}\alpha_i[\tpl{i-1,i}]\to\beta[s_0\times s_1]$ 
        over the projection $I_0\times\dots\times I_n\to I_0\times I_n$.
        \item The composite of the following cells
        \[
            \begin{tikzcd}[column sep=4em,virtual]
                I_{1,0}
                \ar[d, "s_0"']
                \sard[r, "\ol\alpha_1"]
                \ar[dr, phantom, "\xi_1"]
                & I_{1,m_1}
                \ar[d, "s_1"']
                \sard[r]
                & \cdots
                \sard[r, "\ol\alpha_{n}"]
                \ar[dr, phantom, "\xi_n"]
                & I_{n,m_n}
                \ar[d, "s_n"] \\
                J_{0}
                \ar[d, "t_0"']
                \sar[r, "\beta_1"'] 
                \ar[drrr, phantom, "\zeta"]
                & J_{1}
                \sar[r]
                & \cdots
                \sar[r, "\beta_n"']  
                & J_{n}
                \ar[d, "t_1"] \\
                K_{0}
                \sar[rrr, "\gamma"']
                & & & K_{1}
            \end{tikzcd}
        \]
        in $\Bil[\mf{p}]$
        with the cells $\xi_i$ and $\zeta$ given by the arrows
        \[
        \begin{aligned} 
            \xi_1&\colon\bigwedge_{1\le i\le m_1}\alpha_1^i[\tpl{i-1,i}]\to\beta_1[(s_0\times s_1)\circ\tpl{0,m_1}] && \text{in}\quad\one{E}_{\prod_{0\le i\le m_1}I_{1,i}},\\
            &\vdots\\
            \xi_n&\colon\bigwedge_{1\le i\le m_n}\alpha_n^i[\tpl{i-1,i}]\to\beta_n[(s_{n-1}\times s_n)\circ\tpl{0,m_n}] && \text{in}\quad\one{E}_{\prod_{0\le i\le m_n}I_{n,i}}, \text{ and}\\
            \zeta&\colon\bigwedge_{1\le j\le n}\beta_j[\tpl{j-1,j}]\to \gamma[(t_0\times t_1)\circ\tpl{0,n}] && \text{in}\quad\one{E}_{\prod_{0\le j\le n}J_j},\\
            &(\text{where}\quad I_{j,0}:=I_{j-1,m_{j-1}}\quad\text{for}\quad 1< j\le n).
        \end{aligned}
        \]
        is the cell 
        \[
        \begin{aligned}
            &\bigwedge_{\substack{1\le j\le n\\1\le i\le m_j}}\alpha_j^i[\scriptstyle\tpl{\sum_{1\le k<j}m_k+i-1,\sum_{1\le k<j}m_k+i}]
            \\
            \to[3ex]<"\cong">&
            \bigwedge_{1\le j\le n}\left(\bigwedge
            _{1\le i\le m_j}\alpha_j^i[\tpl{i-1,i}]
            \right)
            [\scriptstyle\tpl{\sum_{1\le k<j}m_k,\dots,\sum_{1\le k\le j}m_k}]
            \\
            \to[33ex]<"{\bigwedge_{1\le j\le n}\xi_j[\tpl{\sum_{1\le k<j}m_k,\sum_{1\le k\le j}m_k}]}">&
            \bigwedge_{1\le j\le n}\beta_j[(s_{j-1}\times s_j)\circ\tpl{0,m_j}][\tpl{\scriptstyle\sum_{1\le k<j}m_k,\dots,\sum_{1\le k\le j}m_k}]\\
            \to[3ex]<"\cong">&
            \bigwedge_{1\le j\le n}\beta_j[(s_{j-1}\times s_j)\circ\tpl{\scriptstyle\sum_{1\le k<j}m_k,\sum_{1\le k\le j}m_k}]
            \\
            \to[3ex]<"\cong">&
            \left(\bigwedge_{1\le j\le n}\beta_j[\tpl{j-1,j}]\right)[(s_0\times\dots\times s_n)\circ\tpl{0,m_1,\dots,\scriptstyle\sum_{1\le j\le n}m_j}]
            \\
            \to[28ex]<"{\zeta[(s_0\times\dots\times s_n)\circ\tpl{0,\dots,\scriptstyle\sum_{1\le j\le n}m_j}]}">&
            \gamma[(t_0\times t_1)\circ\tpl{0,n}][(s_0\times\dots\times s_n)\circ\tpl{0,\dots,\scriptstyle\sum_{1\le j\le n}m_j}]\\
            \to[3ex]<"\cong">&
            \gamma[\left((t_0\circ s_0)\times (t_1\circ s_1)\right)\circ\tpl{0,\scriptstyle\sum_{1\le j\le n}m_j}]
            \text{ in }\one{E}_{\prod_{\substack{0\le j\le n\\0\le i\le m_j}}I_{j,i}}.
        \end{aligned}
        \]
        \item The identity cell for a loose arrow $\alpha\colon I\sto J$ is 
        the canonical isomorphism $\alpha\to \alpha[\id_{I\times J}]=\alpha[(\id_I\times\id_J)]$ in $\one{E}_{I\times J}$.
    \end{enumerate}
\end{definition}

It is straightforward to check that the data in \Cref{def:primfibtovdc} form a virtual double category. 
One way to see this is to use the internal language of fibrations, as in \Cref{rem:intlang}.
Loose arrows of $\Bil[\mf{p}]$ correspond to propositions $\syn{\alpha}(\syn{x},\syn{y})$ in the context $\syn{x}:\syn{I},\syn{y}:\syn{J}$, and its cells correspond to proofs $\syn{\xi}$ as follows:
\[
\syn{x}_0:\syn{I}_0,\dots,\syn{x}_n:\syn{I}_n\mid
\syn{\alpha}_1(\syn{x}_0,\syn{x}_1),\dots,\syn{\alpha}_n(\syn{x}_{n-1},\syn{x}_n)\vdash^{\syn{\xi}}
\syn{\beta}(\syn{f}_0(\syn{x}_0),\syn{f}_1(\syn{x}_n)).
\]
Composing cells in $\Bil[\mf{p}]$ means constructing a new proof from given proofs.
Suppose we have the following proofs:
\[
\begin{aligned} 
    \syn{x}_{1,0}:\syn{I}_{1,0},\dots,\syn{x}_{1,m_1}:\syn{I}_{1,m_1}
    &\mid\,
    \syn{\alpha}_1^1(\syn{x}_{1,0},\syn{x}_{1,1}),\dots,\syn{\alpha}_1^{m_1}(\syn{x}_{1,m_1-1},\syn{x}_{1,m_1})\\
    &\vdash^{\syn{\xi}_1} \syn{\beta}_1(\syn{s}_0(\syn{x}_{1,0}),\syn{s}_1(\syn{x}_{1,m_1})),\\
    &\vdots\\
    \syn{x}_{n-1,m_{n-1}}:\syn{I}_{n-1,m_{n-1}},\dots,\syn{x}_{n,m_n}:\syn{I}_{n,m_n}
    &\mid\,
    \syn{\alpha}_n^1(\syn{x}_{n-1,m_{n-1}},\syn{x}_{n,1}),\dots,\syn{\alpha}_n^{m_n}(x_{n,m_n-1},\syn{x}_{n,m_n})\\
    &\vdash^{\syn{\xi}_n} \syn{\beta}_n(\syn{s}_{n-1}(\syn{x}_{n-1,m_{n-1}}),\syn{s}_n(\syn{x}_{n,m_n})),\\
    \syn{y}_0:\syn{J}_0,\dots,\syn{y}_n:\syn{J}_n
    &\mid\,
    \syn{\beta}_1(\syn{y}_0,\syn{y}_1),\dots,\syn{\beta}_n(\syn{y}_{n-1},\syn{y}_n)\\
    &\vdash^{\syn{\zeta}} \syn{\gamma}(\syn{t}_0(\syn{y}_0),\syn{t}_1(\syn{y}_n)).
\end{aligned}
\]
Then, we can construct the proof 
\[
\begin{aligned}
    \syn{x}_{1,0}:\syn{I}_{1,0},\syn{x}_{1,m_1}:\syn{I}_{1,m_1},\dots,\syn{x}_{n,m_n}:\syn{I}_{n,m_n} 
    &\mid\,
    \syn{\beta}_1(\syn{s}_0(\syn{x}_{1,0}),\syn{s}_1(\syn{x}_{1,m_1})),\dots,\syn{\beta}_n(\syn{s}_{n-1}(\syn{x}_{n-1,m_{n-1}}),\syn{s}_n(\syn{x}_{n,m_n}))\\ 
    &\vdash^{\syn{\zeta}[\ol{\syn{s}_{\ul{j}}(\syn{x}_{\ul{j},m_j})}]} \syn{\gamma}(\syn{t}_0(\syn{s}_0(\syn{x}_{1,0})),\syn{t}_1(\syn{s}_n(\syn{x}_{n,m_n}))), 
\end{aligned}
\]
by substituting $\syn{s}_j(\syn{x}_{j,m_j})$'s for $\syn{y}_j$'s in the proof of $\syn{\zeta}$.
Subsequently, we can combine the proofs $\syn{\xi}_1,\dots,\syn{\xi}_n$ with $\syn{\zeta}[\ol{\syn{s}_{\ul{j}}(\syn{x}_{\ul{j},m_j})}]$ to obtain a proof of 
\[
\begin{aligned}
    \syn{x}_{1,0}:\syn{I}_{1,0},\syn{x}_{1,1}:\syn{I}_{1,1},\dots,\syn{x}_{n,m_n}:\syn{I}_{n,m_n} 
    &\mid\,
    \syn{\alpha}_1^1(\syn{x}_{1,0},\syn{x}_{1,1}),\dots,\syn{\alpha}_1^{m_1}(\syn{x}_{1,m_1-1},\syn{x}_{1,m_1}),
    \dots, \syn{\alpha}_n^{m_n}(\syn{x}_{n,m_n-1},\syn{x}_{n,m_n})\\
    &\vdash \syn{\gamma}(\syn{t}_0(\syn{s}_0(\syn{x}_{1,0})),\syn{t}_1(\syn{s}_n(\syn{x}_{n,m_n}))).
\end{aligned}
\]
This corresponds to the composite of the corresponding cells in $\Bil[\mf{p}]$.
The proof of the associativity of the composition in $\Bil[\mf{p}]$
is almost the same as the proof of \cite[Proposition 4.6]{nasu2024internallogicvirtualdouble},
but without the bilaterality of the propositions and the proofs.

\begin{remark}
    This $\Bil$ construction is a generalization of Shulman's $\dbl{F}\mathrm{r}$ construction of framed bicategories from cartesian fibrations with additional structure \cite[Theorem 14.4]{shulmanFramedBicategoriesMonoidal2009}.
    The construction assumes these structures on the fibration so that the resulting entity is a framed bicategory, which in our terminology is an equipment; we will revisit this in \Cref{rem:comparisonwithshulman}.
    On the other hand, the paper deals with more general fibrations than cartesian fibrations; the construction $\dbl{F}\mathrm{r}$ takes a \textit{monoidal fibration} with the monoidal structure on the base category cartesian as input.
    We could follow the same path and defined a virtual double category from a monoidal fibration with cartesian base, but we do not proceed in this direction in this paper.
\end{remark}

\begin{proposition}
    For a cartesian fibration $\mf{p}\colon\one{E}\to\one{B}$,
    the \ac{VDC} $\Bil[\mf{p}]$ is a \ac{CFVDC}.
\end{proposition}
\begin{proof}
    The restriction of a loose arrow $\alpha\colon I\sto J$,
    which is an object in $\one{E}$ over $I\times J$,
    along a pair of tight arrows $f\colon I'\to I$ 
    and $g\colon J'\to J$ is 
    given by $\alpha[f\times g]\in \one{E}_{I'\times J'}$,
    and the restricting cell $\rest$ is the identity arrow on $\alpha[f\times g]$. 
    One can check this by seeing that
    a cell on the left-hand side of the equation
    corresponds to a morphism $\xi$ in $\one{E}_{\prod_{0\le i\le n}K_i}$
    on the right diagram below:
    \[
    \begin{tikzcd}[column sep=1.5em, row sep=1em]
        K_0
        \sar[rr, dashed, "\ol\beta"]
        \ar[d, "h"']
        \ar[ddrr, phantom, "\xi"]
        &
        &
        K_n
        \ar[d, "k"]
        \\
        I'
        \ar[d, "f"']
        &
        &
        J'
        \ar[d, "g"]
        \\
        I
        \sar[rr, "\alpha"']
        &
        &
        J
    \end{tikzcd}
    =
    \begin{tikzcd}[column sep=1.5em, row sep=1em]
        K_0
        \ar[d, "h"']
        \ar[drr, phantom, "\wt\xi"{yshift=1ex}]
        \sar[rr, dashed, "\ol\beta"]
        &
        &
        K_n
        \ar[d, "k"]
        \\
        I'
        \ar[d, "f"']
        \sar[rr, "{\alpha[f\times g]}"]
        \ar[rrd, phantom, "\restc"]
        &
        &
        J'
        \ar[d, "g"]
        \\
        I
        \sar[rr, "\alpha"']
        &
        &
        J
    \end{tikzcd},
    \hspace{1em}
    \begin{tikzcd}[column sep=1em, row sep=1.5em]
        \displaystyle\bigwedge_{1\le i\le n}\beta_i[\tpl{i-1,i}]
        \ar[r, "\wt\xi"]
        \ar[rd,"\xi"', near start]
        &
        \alpha[f\times g][(h\times k)\circ\tpl{0,n}]
        \ar[d, "\cong"]
        \\
        &
        \alpha[(f\times g)\circ(h\times k)\circ\tpl{0,n}]
    \end{tikzcd}
    \hspace{-2em}
    \text{in }\one{E}_{\prod_{0\le i\le n}K_i}.
    \]
    Here, post-composing the canonical isomorphism on the rightmost diagram
    represents post-composing the cell $\restc$ on the leftmost diagram. 
    This shows that $\Bil[\mf{p}]$ has restrictions.

    To show that $\Bil[\mf{p}]$ is cartesian,
    let us recall \Cref{prop:FibVDblCart}, which provides an explicit description 
    of the cartesian structure on an \ac{FVDC}.
    The vertical part of $\Bil[\mf{p}]$ has finite products by definition,
    and for each pair of objects $I$ and $J$ in $\one{B}$,
    the finite products in $\one{E}_{I\times J}$ give 
    the local finite products in $\Bil[\mf{p}](I,J)$.
    Finally, these are preserved by restriction along tight arrows
    because it is given by base change in $\one{E}$,
    which preserves finite products in a cartesian fibration.
\end{proof}

\begin{proposition} 
    \label{prop:primfibtovdc}
    The assignment of a \ac{CFVDC} $\Bil[\mf{p}]$ to a cartesian fibration $\mf{p}\colon\one{E}\to\one{B}$ 
    extends to a 2-functor
    $\Bil\colon\Fib\carttwo\to{\VDbl-}\carttwo$.
\end{proposition}
\begin{proof} 
    A morphism of cartesian fibrations $\mf{f}$
    from $\mf{p}$ to $\mf{p}'$
    induces a morphism of \acp{CFVDC} $\Bil[\mf{f}]\colon\Bil[\mf{p}]\to\Bil[\mf{p}']$.
    This is because the morphism $\mf{f}$ preserves the structure of the fibration
    including the cartesian structure, 
    and hence the structure of the \ac{CFVDC} $\Bil[\mf{p}]$.
    The assignment of a \ac{CFVDC} $\Bil[\mf{p}]$ to a cartesian fibration $\mf{p}\colon\one{E}\to\one{B}$ 
    is functorial accordingly.

    A 2-cell $\theta\colon\mf{f}\to\mf{g}\colon\mf{p}\to\mf{p}'$ is given by the pair $\theta=(\theta_0,\theta_1)$ of natural transformations as follows.
\[\begin{tikzcd}[ampersand replacement=\&]
	{\one{E}} \&\& {\one{E}'} \\
	\\
	{\one{B}} \&\& {\one{B}'}
	\arrow[""{name=0, anchor=center, inner sep=0}, "{\mf{f}_1}", curve={height=-12pt}, from=1-1, to=1-3]
	\arrow[""{name=1, anchor=center, inner sep=0}, "{\mf{g}_1}"', curve={height=12pt}, from=1-1, to=1-3]
	\arrow["{\mf{p}}"', from=1-1, to=3-1]
	\arrow["{\mf{p}}", from=1-3, to=3-3]
	\arrow[""{name=2, anchor=center, inner sep=0}, "{\mf{f}_0}", curve={height=-12pt}, from=3-1, to=3-3]
	\arrow[""{name=3, anchor=center, inner sep=0}, "{\mf{g}_0}"', curve={height=12pt}, from=3-1, to=3-3]
	\arrow["{\theta_1}", between={0.2}{0.8}, Rightarrow, nfold, from=0, to=1]
	\arrow["{\theta_0}", between={0.2}{0.8}, Rightarrow, nfold, from=2, to=3]
\end{tikzcd}\]
    To define a tight transformation $\Bil[\theta]\colon\Bil[\mf{f}]\Rightarrow\Bil[\mf{g}]$, we use the components of $\theta_0$ as the components of $\Bil[\theta]$ at an object $I$ in $\one{B}$, and the component $\Bil[\theta]_\alpha$ at a loose arrow $\alpha\colon I\sto J$ in $\Bil[\mf{p}]$ is given by the unique arrow $\wt\theta_\alpha$ that makes the following diagram commute:
    \[
\begin{tikzcd}[ampersand replacement=\&]
	{\mf{f}_0(I)} \& {\mf{f}_0(J)} \\
	{\mf{g}_0(I)} \& {\mf{g}_0(J)}
	\arrow[""{name=0, anchor=center, inner sep=0}, "{\Bil[\mf{f}](\alpha)}"{inner sep=.8ex}, "\shortmid"{marking}, from=1-1, to=1-2]
	\arrow["{\theta_{0,I}}"', from=1-1, to=2-1]
	\arrow["{\theta_{0,J}}", from=1-2, to=2-2]
	\arrow[""{name=1, anchor=center, inner sep=0}, "{\Bil[\mf{g}](\alpha)}"'{inner sep=.8ex}, "\shortmid"{marking}, from=2-1, to=2-2]
	\arrow["{\Bil[\theta]_\alpha}"{description}, draw=none, from=0, to=1]
\end{tikzcd},
    \hspace{4em}
\begin{tikzcd}[ampersand replacement=\&]
	{\mf{f}_1(\alpha)} \&\& {\one{E}'} \\
	{(\mf{g}_1(\alpha))[\theta_I\times\theta_J]} \& {\mf{g}_1(\alpha)} \\
	{\mf{f}_0({I})\times\mf{f}_0({J})} \& {\mf{g}_0({I})\times\mf{g}_0({J})} \& {\one{B}'}
	\arrow["\begin{array}{c} \exists!\wt{\theta}_\alpha \end{array}"', dashed, from=1-1, to=2-1]
	\arrow["{\theta_{1,\alpha}}", from=1-1, to=2-2]
	\arrow["{\mf{p}'}", from=1-3, to=3-3]
	\arrow["{\text{prone}}"', from=2-1, to=2-2]
	\arrow[between={0.2}{0.8}, maps to, from=2-1, to=3-1]
	\arrow[between={0.2}{0.8}, maps to, from=2-2, to=3-2]
	\arrow["{\theta_I\times\theta_J}"', from=3-1, to=3-2]
\end{tikzcd}
    \] 
    The naturality conditions of $\Bil[\theta]$ follow from
    the naturality of $\theta$.
\end{proof}

\begin{example}
    \label{ex:primfibtovdc}
    Some examples of \acp{CFVDC} we have seen so far can be obtained through the construction $\Bil$.
    \begin{enumerate}
    \item The subobject fibration over a category with finite limits $\one{C}$ gives rise to $\Rel[\one{C}]$, the virtual double category of relations internal to $\one{C}$. When $\one{C}$ is not regular, this virtual double category is not a double category. This generalizes double categories studied in \cite{lambertDoubleCategoriesRelations2022}.
    \item The codomain fibration over a category with finite limits $\one{C}$ gives rise to the double category of spans $\Span[\one{C}]$.
    \item The family fibration $\Fam{\one{B}}\to\Set$ for a cartesian category $\one{B}$ gives rise to the virtual double category of $\one{B}$-valued matrices $\Mat[\one{B}]$. This can be properly virtual when $\one{B}$ does not have distributive coproducts. 
    \item For an orthogonal factorization system $(\zero{E},\zero{M})$ on a category $\one{C}$ with finite limits, the fibration of $\zero{M}$-subobjects gives rise to $\Rel[\one{C}][(\zero{E},\zero{M})]$, the virtual double category of $\zero{M}$-relations internal to $\one{C}$. This can be properly virtual when the factorization system is not stable. This generalizes double categories studied in \cite{hoshinoDoubleCategoriesRelations2023}.
    \qedhere
\end{enumerate}
\end{example}

Now we have the construction of a \ac{CFVDC} from a cartesian fibration.
The next step is to show how the properties of the fibration 
are reflected in the resulting \ac{CFVDC}.
The primary interest is in coproducts in fibrations,
since they are the key ingredient to interpret regular logic in fibrations.

\begin{lemma} 
    \label{lem:unitalvdc}
    For an elementary fibration $\mf{p}$,
    the \ac{VDC} $\Bil[\mf{p}]$ is unital.
\end{lemma}
\begin{proof}
    We will prove that a unit on an object $I$ in $\one{B}$ is given by the object $\delta_I\coloneqq\sum_{\tpl{0,0}}\top_I$ in $\one{E}_{I\times I}$, where $\top_I$ is the terminal object in $\one{E}_I$.
    Here, the unit cell $\eta_I$ is the component of the unit $\eta$ of the adjunction $\sum_{\tpl{0,0}}\colon\one{E}_{I}\adjointleft\one{E}_{I\times I}\lon(-)[\tpl{0,0}]$ at the object $\top_I$.
    The universal property of the unit cell $\eta_I$ that we want to show is stated as follows: Any cell $\nu$ on the left below uniquely factors through the unit cell $\eta_I$ as on the right below:
    \[
        \begin{tikzcd}[virtual, column sep=small]
            J_m
            \sar[r, "\alpha_m"{yshift=1ex}]
            \ar[d, equal]
            \ar[rrrrrrd, phantom, "\nu",description]
            &
            J_{m-1}
            \sar[rr,dashed,"{\alpha_{m-1},\dots,\alpha_1}"{yshift=1ex}]
            &&
            I
            \sar[rr, "{\beta_1,\dots,\beta_{n-1}}"{yshift=1ex}, dashed]
            &
            &
            K_{n-1}
            \sar[r, "\beta_{n}"{yshift=1ex}]
            &
            K_n
            \ar[d, equal]
            \\
            J_m
            \sar[rrrrrr, "\gamma"']
            &&&&\!&&
            K_n
        \end{tikzcd}
        =
        \begin{tikzcd}[virtual, column sep=small]
            J_m
            \sar[r, "\alpha_m"]
            \ar[d, equal]
            \ar[dr, phantom, description, "{\rotatebox{90}{=}}"]
            &
            J_{m-1}
            \sar[rr,dashed,"{\alpha_{m-1},\dots,\alpha_1}"]
            \ar[d,equal]
            \ar[dr, phantom, description, "{\rotatebox{90}{=}}", xshift=1ex]
            &&
            I
            \ar[ld, equal]
            \ar[rd, equal]
            \ar[d,phantom, "\eta_I"]
            \sar[rr, "{\beta_1,\dots,\beta_{n-1}}", dashed]
            &
            &
            K_{n-1}
            \ar[d, equal]
            \ar[dl, phantom, description, "{\rotatebox{90}{=}}", xshift=-1ex]
            \ar[dr, phantom, description, "{\rotatebox{90}{=}}"]
            \sar[r, "\beta_{n}"]
            &
            K_n
            \ar[d, equal]
            \\
            J_m
            \sar[r, "\alpha_m"']
            \ar[d, equal]
            \ar[phantom,rrrrrrd, "\wt{\nu}",description, yshift=-1ex]
            &
            J_{m-1}
            \sar[r,dashed]
            &
            I
            \sar[rr, "\delta_I"']
            &\!&
            I
            \sar[r, dashed]
            &
            K_{n-1}
            \sar[r, "\beta_{n}"']
            &
            K_n
            \ar[d, equal]
            \\
            J_m
            \sar[rrrrrr, "\gamma"']
            &&&&\!&&
            K_n
        \end{tikzcd}.
    \]
    This amounts to saying that for any arrow $\nu$ as below,
    there is a unique arrow $\wt{\nu}$ 
    for which $\wt{\nu}[\Delta]$ makes the following diagram commute:
    \[
    \begin{tikzcd}[row sep=small]
        \kappa\coloneqq\displaystyle\bigwedge_{1\le i\le m}\alpha_i[{\scriptstyle\tpl{m-i,m-i+1}}]
        \land  
        \bigwedge_{1\le j\le n}\beta_j[{\scriptstyle\tpl{m+j-1,m+j}}]
        \ar[rd,"\nu"]
        \ar[d,"\cong"']
        \ar[ddd, bend right=30, shift right = 12em, "\tau"']
        \\
        \displaystyle\bigwedge_{1\le i\le m}\alpha_i[{\scriptstyle\tpl{m-i,m-i+1}}]
        \land  
        \top_I[{\scriptstyle\tpl{m}}]
        \land  
        \bigwedge_{1\le j\le n}\beta_j[{\scriptstyle\tpl{m+j-1,m+j}}]
        \ar[d, "{\id\land\eta_{I,\top_I}[\tpl{m}]\land\id}"']
        \ar[rd, phantom, "\circlearrowright"]
        &
        \gamma[{\scriptstyle\tpl{0,m+n}}]
        \ar[d,"\cong"]
        \\
        \displaystyle\bigwedge_{1\le i\le m}\alpha_i[{\scriptstyle\tpl{m-i,m-i+1}}]
        \land  
        \delta_I[{\scriptstyle\tpl{0,0}}][{\scriptstyle\tpl{m}}]
        \land
        \bigwedge_{1\le j\le n}\beta_j[{\scriptstyle\tpl{m+j-1,m+j}}]
        \ar[d,"\cong"']
        &
        \gamma[{\scriptstyle\tpl{0,m+n+1}}][\Delta]
        \\
        \left(
        \displaystyle\bigwedge_{1\le i\le m}\alpha_i[{\scriptstyle\tpl{m-i,m-i+1}}]
        \land  
        \delta_I[{\scriptstyle\tpl{m,m+1}}]
        \land
        \bigwedge_{1\le j\le n}\beta_j[{\scriptstyle\tpl{m+j,m+j+1}}]
        \right)
        [\Delta]
        \ar[ru, dashed, "{\wt{\nu}[\Delta]}"']
    \end{tikzcd}
    \]
    Here, $J_0\coloneqq I$, $K_0\coloneqq I$, and $\Delta\coloneqq\scriptstyle\tpl{0,\dots,m-1,m,m,m+1,\dots,m+n}$.
    We now show that the composite $\tau$ of the arrows in the left column of the above diagram coincides with the component of the unit $\eta'$ of the adjunction $\sum_{\Delta}\dashv\,(-)[\Delta]$ at the object $\kappa$.
    Let us consider the Beck-Chevalley condition for the following diagram:
    \begin{equation}
    \label{eq:beckchevalleyinst}
    \begin{tikzcd}[row sep=2.5em]
        \one{E}_{\prod_{1\le i\le m}J_i\times I\times\prod_{1\le j\le n}K_j}
        \ar[r, "\scriptstyle\sum_{\Delta}", bend left=10]
        \ar[r, phantom, "\rotatebox{90}{$\vdash$}"]
        &
        \one{E}_{\prod_{1\le i\le m}J_i\times I\times I\times\prod_{1\le j\le n}K_j}
        \ar[l, "{(-)[\Delta]}", bend left=10]
        \\
        \one{E}_{I}
        \ar[u, "{(-)[\tpl{m}]}"]
        \ar[r, "\scriptstyle\sum_{\tpl{0,0}}", bend left=10]
        \ar[r, phantom, "\rotatebox{90}{$\vdash$}"]
        &
        \one{E}_{I\times I}
        \ar[l, "{(-)[\tpl{0,0}]}", bend left=10]
        \ar[u, "{(-)[\tpl{m,m+1}]}"']
    \end{tikzcd}.
    \end{equation}
    Looking at the components of the units at the terminal objects,
    we have the following commutative triangle:
    \[
    \begin{tikzcd}[column sep=1em]
        \top_{I}[{\scriptstyle\tpl{m}}]
        \ar[dr, "{\eta'_{\top_{I}[{\scriptstyle\tpl{m}}]}}"']
        \ar[r,"{(\eta_{I,\top_{I}})[{\scriptstyle\tpl{m}}]}"{yshift=1ex}]
        &
        \left(\delta_I[\tpl{0,0}]\right)[\tpl{m}]
        \ar[d, "\cong"]
        \\
        &
        \left(\sum_\Delta(\top_{I}[{\scriptstyle\tpl{m}}])\right)[\Delta]    
    \end{tikzcd}
    \text{ in }
    \one{E}_{\prod_{1\le i\le m}J_i\times I\times\prod_{1\le j\le n}K_j}.
    \]
    In addition, it follows from \Cref{lem:implicationsoffrobele} that the following commutes
    in $\one{E}_{\prod_{1\le i\le m}J_i\times I\times\prod_{1\le j\le n}K_j}$:
    \[
    \begin{tikzcd}[column sep=1em]
        \kappa
        \ar[r, "\eta'_{\kappa}"]
        \ar[d, "\cong"']
        &
        \left(\sum_\Delta\kappa\right)[\Delta]
        \ar[d, "\cong"]
        \\
        \top_{I}[{\scriptstyle\tpl{m}}]\land\kappa
        \ar[r, "{\eta'_{\top_{I}[{\scriptstyle\tpl{m}}]}}"]
        \ar[d, "{\eta'_{\top_{I}[{\scriptstyle\tpl{m}}]}\land\id_\kappa}"']
        &
        \left(\sum_\Delta\left(\top_{I}[{\scriptstyle\tpl{m}}]\land\kappa\right)\right)[\Delta]
        \ar[d, "\cong"]
        \\
        \left(\sum_\Delta\left(\top_{I}[{\scriptstyle\tpl{m}}]\right)\right)[\Delta]\land\kappa
        \ar[r, "\cong"]
        &
        \left(\sum_\Delta\left(\top_{I}[{\scriptstyle\tpl{m}}]\right)\land\kappa[{\scriptstyle\id\times\tpl{0}\times\id}]\right)[\Delta] 
    \end{tikzcd}
    \]
    With these diagrams, we can show that the composite $\tau$ is equal 
    to the arrow $\eta'_{\kappa}$ up to the isomorphisms that are given by the structures of $\mf{p}$
    as an elementary fibration.
    Note that the subexpression in the codomain of $\tau$ 
    to which the base change $(-)[\Delta]$ is applied
    can be identified with the object $\left(\sum_\Delta\kappa\right)[\Delta]$
    via the isomorphisms in the argument above;
    this can be summarized as the following sequence of isomorphisms:
    \begin{align*}
        \sum_\Delta\kappa
        &\cong
        \sum_\Delta\left(\top_I[{\scriptstyle\tpl{m}}]\land\kappa\right)
        &
        \text{by the preservation of finite products by base change}
        \\
        &\cong
        \sum_\Delta\left(\top_I[{\scriptstyle\tpl{m}}]\right)\land\kappa[{\scriptstyle\id\times\tpl{0}\times\id}] 
        &\text{by Frobenius reciprocity, specifically \Cref{lem:implicationsoffrobele}}
        \\
        &\cong
        \delta_I[{\scriptstyle\tpl{m,m+1}}]\land\kappa[{\scriptstyle\id\times\tpl{0}\times\id}]
        &\text{by the Beck-Chevalley condition \Cref{eq:beckchevalleyinst}}
        \\
        &\cong
        \displaystyle\bigwedge_{1\le i\le m}\alpha_i[{\scriptstyle\tpl{m-i,m-i+1}}]
        \land  
        \delta_I[{\scriptstyle\tpl{m,m+1}}] \land
        \bigwedge_{1\le j\le n}\beta_j[{\scriptstyle\tpl{m+j,m+j+1}}].
        \hspace{-13em}
    \end{align*}
\end{proof}

\begin{proposition}
    \label{prop:bil-elem-fib-is-cart-veq}
    For an elementary fibration $\mf{p}\colon\one{E}\to\one{B}$,
    the \ac{VDC} $\Bil[\mf{p}]$ is a cartesian virtual equipment.
\end{proposition}
\begin{proof}
    We have shown that $\Bil[\mf{p}]$ is a unital \ac{FVDC} in the previous lemma.
    Considering \Cref{prop:cartesianunital}, it remains to show that units are compatible with the cartesian structure
    in the sense of $(ii)$ and $(iii)$ there.
    The second condition is easily satisfied because the diagonal $\tpl{00}$ 
    on the terminal object $1$ in $\one{B}$ is an isomorphism.
    The third condition follows from the Beck-Chevalley condition
    and the Frobenius reciprocity.
    From the construction of the units in $\Bil[\mf{p}]$
    and the translation between the fiberwise and the total finite products \Cref{prop:cartesianfibration},
    the canonical arrow in $(iii)$ of \Cref{prop:cartesianunital} 
    is given by the arrow
    \[
        \delta_{I\times J}
        \to[5ex]<"{\tpl{\mu,\nu}}">
        \delta_I[\tpl{0_I,2_I}]\land\delta_J[\tpl{1_J,3_J}]
        \quad \text{in }\one{E}_{I\times J\times I\times J},
    \]
    where $\mu$ and $\nu$ are the arrows that correspond respectively 
    to the following arrows in $\one{E}_{I\times J}$:
    \[
    \begin{aligned}    
        \top_{I\times J} \ \cong \ \top_I[\tpl{0_I}]
        \to[5ex]<"{\eta_I[\tpl{0_I}]}">
        \delta_I[\tpl{0_I,0_I}][\tpl{0_I}]
        \ \cong \
        \delta_I[\tpl{0_I,2_I}][\tpl{0_I,1_J,0_I,1_J}],
        \\
        \top_{I\times J} \ \cong \ \top_J[\tpl{1_J}]
        \to[5ex]<"{\eta_J[\tpl{0_J}]}">
        \delta_J[\tpl{0_J,0_J}][\tpl{1_J}]
        \ \cong \
        \delta_J[\tpl{1_J,3_J}][\tpl{0_I,1_J,0_I,1_J}].
    \end{aligned}
    \]
    However, it is not hard to see that
    this coincides with the composite of the arrows as follows:
    \[
    \begin{aligned}
        \delta_{I\times J}
        &\to<"\cong">
        \sum_{\scriptstyle\tpl{0_I,1_J,0_I,2_J}}\sum_{\tpl{0_I,1_J,1_J}}\top_{I\times J}
        &
        \text{\small by the uniqueness of the left adjoint}
        \\
        &\to<"\cong">
        \sum_{\tpl{0_I,1_J,0_I,2_J}}\left(\top_{I\times J\times J}\land
        \sum_{\tpl{0_I,1_J,1_J}}\left(\top_{I\times I\times J}[{\scriptstyle\tpl{0_I,0_I,1_J}}]\right)\right)
        \\
        &\to<"\cong">
        \sum_{\tpl{0_I,1_J,0_I,2_J}}\left(\top_{I\times J\times J}\land
        \left(\sum_{\scriptstyle\tpl{0_I,2_J,1_I,2_J}}\top_{I\times I\times J}\right)[{\scriptstyle\tpl{0_I,1_J,0_I,2_J}}]\right)
        \hspace{-2em}
        &
        \text{\small by the Beck-Chevalley condition}
        \\
        &\to<"\cong">
        \sum_{\tpl{0_I,1_J,0_I,2_J}}\top_{I\times J\times J}
        \land \sum_{\tpl{0_I,2_J,1_I,2_J}}\top_{I\times I\times J}
        &
        \text{\small by the Frobenius reciprocity}
        \\
        &\to<"\cong">
        \delta_I[\tpl{0_I,2_I}]\land\delta_J[\tpl{1_J,3_J}]
        &
        \text{\small by the Beck-Chevalley condition.}
    \end{aligned}
    \]
    This completes the proof.
\end{proof}

\begin{corollary}
    The assignment of a cartesian unital \ac{FVDC} $\Bil[\mf{p}]$
    to an elementary fibration $\mf{p}\colon\one{E}\to\one{B}$
    extends to a 2-functor
    $\Bil\colon\Fib_{\lc*=}\to\VEq\carttwo$.
\end{corollary}
\begin{proof}
    A morphism between elementary fibrations $\mf{f}\colon\mf{p}\to\mf{p}'$
    preserves the left adjoints of base change functors,
    which define the units and the unit cells of the \ac{FVDC} $\Bil[\mf{p}]$.
    Hence, the assignment of an \ac{FVDC} $\Bil[\mf{p}]$ to an elementary fibration $\mf{p}\colon\one{E}\to\one{B}$
    extends to a functor.
    Since $\Fib_{\lc*=}$ and $\VEq\carttwo$ are both locally fully embedded in $\Fib\carttwo$ and $\FVDbl\carttwo$, respectively, the assignment extends to a 2-functor.
\end{proof}

\begin{lemma}
    For an existential fibration $\mf{p}\colon\one{E}\to\one{B}$,
    the \ac{VDC} $\Bil[\mf{p}]$ is positively representable.
\end{lemma}
\begin{proof}
    The proof is analogous to the proof of \Cref{lem:unitalvdc}.
    It is enough to show that a composite of two arrows 
    $\alpha_0\colon J_0\to I_0$ and $\beta_0\colon I_0\to K_0$
    is given by the object $\sum_{\tpl{0,2}}\left(\alpha_0[\tpl{0,1}]\land\beta_0[\tpl{1,2}]\right)$
    in $\one{E}_{J_0\times K_0}$.
    The cell $\varkappa_{\alpha_0;\beta_0}\colon\alpha_0;\beta_0\Rightarrow\alpha_0\odot\beta_0$ 
    in the definition of composition of paths is, in this case, 
    the unit component of the adjunction
    $\sum_{\tpl{0,2}}\dashv\,(-)[\tpl{0,2}]$ at the object $\alpha_0[\tpl{0,1}]\land\beta_0[\tpl{1,2}]$
    in $\one{E}_{J_0\times I_0\times K_0}$.

    Let $p$ be the prodcut projection $\prod_{0\le i\le m}J_i\times I_0\times\prod_{0\le j\le n}K_j\to \prod_{0\le i\le m}J_i\times\prod_{0\le j\le n}K_j$. 
    The universal property of the cell $\varkappa_{\alpha_0;\beta_0}$ is
    that for any cell $\nu\colon \alpha_m;\dots;\alpha_0;\beta_0;\dots;\beta_n\Rightarrow\gamma$,
    there is a unique cell $\wt{\nu}\colon\alpha_m;\dots;\alpha_0\odot\beta_0;\dots;\beta_n\Rightarrow\gamma$ 
    for which $\wt{\nu}[p]$ makes the composite of the following arrows
    equal to $\nu$:
    \[
        \begin{tikzcd}[row sep=small]
            \theta\coloneqq\displaystyle\bigwedge_{1\le i\le m}\alpha_i[{\scriptstyle\tpl{m-i,m-i+1}}]
            \land  
            \alpha_0[{\scriptstyle\tpl{m,m+1}}]
            \land
            \beta_0[{\scriptstyle\tpl{m+1,m+2}}]
            \land
            \bigwedge_{1\le j\le n}\beta_j[{\scriptstyle\tpl{m+j+1,m+j+2}}]
            \ar[d,"\cong"']
            \\
            \displaystyle\bigwedge_{1\le i\le m}\alpha_i[{\scriptstyle\tpl{m-i,m-i+1}}]
            \land  
            \left(\alpha_0[{\scriptstyle\tpl{0,1}}]\land\beta_0[{\scriptstyle\tpl{1,2}}]\right)[{\scriptstyle\tpl{m,m+1,m+2}}]
            \land  
            \bigwedge_{1\le j\le n}\beta_j[{\scriptstyle\tpl{m+j+1,m+j+2}}]
            \ar[d, "{\id\land\varkappa_{\alpha_0;\beta_0}[\tpl{m,m+1,m+2}]\land\id}"']
            \\
            \displaystyle\bigwedge_{1\le i\le m}\alpha_i[{\scriptstyle\tpl{m-i,m-i+1}}]
            \land  
            \left(\sum_{\scriptstyle\tpl{0,2}}\left(\alpha_0[{\scriptstyle\tpl{0,1}}]\land\beta_0[{\scriptstyle\tpl{1,2}}]\right)\right)[{\scriptstyle\tpl{0,2}}][{\scriptstyle\tpl{m,m+1,m+2}}]
            \land
            \bigwedge_{1\le j\le n}\beta_j[{\scriptstyle\tpl{m+j+1,m+j+2}}]
            \ar[d,"\cong"']
            \\
            \!
        \end{tikzcd}
        \]
        \[
        \begin{tikzcd}[row sep=small]
            \left(
            \displaystyle\bigwedge_{1\le i\le m}\alpha_i[{\scriptstyle\tpl{m-i,m-i+1}}]
            \land  
            \left(\sum_{\scriptstyle\tpl{0,2}}\left(\alpha_0[{\scriptstyle\tpl{0,1}}]\land\beta_0[{\scriptstyle\tpl{1,2}}]\right)\right)[{\scriptstyle\tpl{m,m+1}}]
            \land
            \bigwedge_{1\le j\le n}\beta_j[{\scriptstyle\tpl{m+j,m+j+1}}]
            \right)
            [p]
            \ar[d, dashed, "{\wt{\nu}[p]}"']
            \\
            \gamma[{\scriptstyle\tpl{0,m+n+1}}][p]
            \ar[d, "\cong"']
            \\
            \gamma[{\scriptstyle\tpl{0,m+n+2}}]
        \end{tikzcd}.
        \]

        The following shows that the domain of the $\wt{\nu}$
        is isomorphic to the image of $\theta$
        under $\sum_{p}$:
        \begin{align*}
            \sum_p\theta
            &\cong
            \sum_p\left(
            \displaystyle\bigwedge_{1\le i\le m}\alpha_i[{\scriptstyle\tpl{m-i,m-i+1}}] 
            \land
            \bigwedge_{1\le j\le n}\beta_j[{\scriptstyle\tpl{m+j+1,m+j+2}}]
            \land
            \left(
            \alpha_0[{\scriptstyle\tpl{0,1}}]\land\beta_0[{\scriptstyle\tpl{1,2}}] 
            \right)[{\scriptstyle\tpl{m,m+1,m+2}}]
            \right)
            \\
            &\cong
            \displaystyle\bigwedge_{1\le i\le m}\alpha_i[{\scriptstyle\tpl{m-i,m-i+1}}] 
            \land
            \bigwedge_{1\le j\le n}\beta_j[{\scriptstyle\tpl{m+j,m+j+1}}]
            \land
            \sum_p\left(\left(
                \alpha_0[{\scriptstyle\tpl{0,1}}]\land\beta_0[{\scriptstyle\tpl{1,2}}] 
                \right)[{\scriptstyle\tpl{m,m+1,m+2}}]\right)
            \\
            &\cong
            \displaystyle\bigwedge_{1\le i\le m}\alpha_i[{\scriptstyle\tpl{m-i,m-i+1}}] 
            \land
            \bigwedge_{1\le j\le n}\beta_j[{\scriptstyle\tpl{m+j,m+j+1}}]
            \land
            \left(
            \sum_{\tpl{0,2}}
                \alpha_0[{\scriptstyle\tpl{0,1}}]\land\beta_0[{\scriptstyle\tpl{1,2}}] 
                \right)[{\scriptstyle\tpl{m,m+1}}]
        \end{align*}
        From the second to the third line, we used the Frobenius reciprocity,
        and from the third to the fourth line, we used the Beck-Chevalley condition. 
        What remains is to show the image of this isomorphism under the base change $(-)[p]$
        precomposed with the unit component of the adjunction $\sum_p\dashv\,(-)[p]$ at $\theta$ 
        is indeed the same as the upper part of the diagram above,
        which is a straightforward calculation as in the proof of \Cref{lem:unitalvdc}. 
\end{proof}

\begin{proposition}
    For an existential fibration $\mf{p}\colon\one{E}\to\one{B}$,
    the \ac{VDC} $\Bil[\mf{p}]$ is a cartesian positively representable \ac{FVDC}.
\end{proposition}    

\begin{proof}
    The proof is analogous to the proof of \Cref{prop:bil-elem-fib-is-cart-veq}
    besides that we consult \Cref{prop:cartesiancompose} instead.
    Again, we only give the canonical isomorphisms, 
    but not the explicit calculations for verifying that 
    the canonical isomorphisms indeed arise as the universal properties
    of composition of paths.

    For (ii), we have 
    \[
        \top_{1}\odot\top_{1}
        =
        \sum_{\tpl{0,2}}\left(\top_{1}[\tpl{0,1}]\land\top_{1}[\tpl{1,2}]\right) 
        \cong
        \top_{1}
    \]
    since $\tpl{i,j}$ is the identity on $1$ for all $i$ and $j$.
    For (iii),
    if we have $\alpha_i\colon I_{i-1}\sto I_i$ and $\beta_i\colon J_{i-1}\sto J_i$ for $i=1,2$,
    then we obtain the following isomorphisms.
    Here, we do not write the projections by numbers but by the names of the objects
    for the sake of readability.
    \begin{align*}
        &\left(\alpha_1\times\beta_1\right)\odot\left(\alpha_2\times\beta_2\right)
        \\
        &=
        \sum_{\tpl{I_0,J_0,I_2,J_2}}\left(
        (\alpha_1\times\beta_1)[{\scriptstyle\tpl{I_0,J_0,I_1,J_1}}]
        \land
        (\alpha_2\times\beta_2)[{\scriptstyle\tpl{I_1,J_1,I_2,J_2}}]
        \right)
        \\
        &\cong
        \sum_{\tpl{I_0,J_0,I_2,J_2}}\left(
        \alpha_1[{\scriptstyle\tpl{I_0,I_1}}]\land\alpha_2[{\scriptstyle\tpl{I_1,I_2}}]
        \land
        \beta_1[{\scriptstyle\tpl{J_0,J_1}}]\land\beta_2[{\scriptstyle\tpl{J_1,J_2}}]
        \right)
        \\
        &\cong
        \sum_{\tpl{I_0,J_0,I_2,J_2}}\sum_{\tpl{I_0,J_0,J_1,I_2,J_2}}\left(
        (\alpha_1[{\scriptstyle\tpl{I_0,I_1}}]\land\alpha_2[{\scriptstyle\tpl{I_1,I_2}}])\land
        \left(
        \beta_1[{\scriptstyle\tpl{J_0,J_1}}]\land\beta_2[{\scriptstyle\tpl{J_1,J_2}}]
        \right)[{\scriptstyle\tpl{I_0,J_0,J_1,I_2,J_2}}]
        \right)
        \\
        &\cong
        \sum_{\tpl{I_0,J_0,I_2,J_2}}\left(
        \left(
        \sum_{\tpl{I_0,J_0,J_1,I_2,J_2}}
        (\alpha_1[{\scriptstyle\tpl{I_0,I_1}}]\land\alpha_2[{\scriptstyle\tpl{I_1,I_2}}])
        [{\scriptstyle   \tpl{I_0,J_0,I_1,I_2,J_2}}]
        \right)
        \land
        \beta_1[{\scriptstyle\tpl{J_0,J_1}}]\land\beta_2[{\scriptstyle\tpl{J_1,J_2}}]
        \right)
        \\
        &\cong
        \sum_{\tpl{I_0,J_0,I_2,J_2}}\left(
            \left(
            \sum_{\tpl{I_0,J_0,I_2,J_2}}
            (\alpha_1[{\scriptstyle\tpl{I_0,I_1}}]\land\alpha_2[{\scriptstyle\tpl{I_1,I_2}}])
            \right)
            [{\scriptstyle\tpl{I_0,J_0,I_2,J_2}}]
            \land
            \beta_1[{\scriptstyle\tpl{J_0,J_1}}]\land\beta_2[{\scriptstyle\tpl{J_1,J_2}}]
            \right)
        \\
        &\cong
        \left(
            \sum_{\tpl{I_0,J_0,I_2,J_2}}
            (\alpha_1[{\scriptstyle\tpl{I_0,I_1}}]\land\alpha_2[{\scriptstyle\tpl{I_1,I_2}}])
        \right)
        \land
        \left(
            \sum_{\tpl{I_0,J_0,I_2,J_2}}
            \beta_1[{\scriptstyle\tpl{J_0,J_1}}]\land\beta_2[{\scriptstyle\tpl{J_1,J_2}}]
        \right)\\
        &\cong
        \left(
        \sum_{\tpl{I_0,I_2}}
        (\alpha_1[{\scriptstyle\tpl{I_0,I_1}}]\land\alpha_2[{\scriptstyle\tpl{I_1,I_2}}])
        \right)
        [{\scriptstyle\tpl{I_0,I_2}}]
        \land
        \left(
        \sum_{\tpl{J_0,J_2}}
        \beta_1[{\scriptstyle\tpl{J_0,J_1}}]\land\beta_2[{\scriptstyle\tpl{J_1,J_2}}]
        \right)
        [{\scriptstyle\tpl{J_0,J_2}}]\\
        &\cong
        (\alpha_1\odot\alpha_2)\times(\beta_1\odot\beta_2).
        \end{align*}
\end{proof}

\begin{corollary}
    The assignment of a cartesian positively representable \ac{FVDC} $\Bil[\mf{p}]$
    to an existential fibration $\mf{p}\colon\one{E}\to\one{B}$
    extends to a 2-functor
    $\Bil\colon\Fib_{\lc*+}\to\FVDblP\carttwo$.
\end{corollary}
\begin{proof}
    The preservation of the left adjoints of base change functors along 
    product projections leads to the preservation of composition of paths of positive length.
    Hence, the assignment of an \ac{FVDC} $\Bil[\mf{p}]$ to an existential fibration $\mf{p}\colon\one{E}\to\one{B}$ 
    extends to a functor, 
    and further to a 2-functor by the same argument as in the proof of the previous corollary. 
\end{proof}

\begin{corollary}
    The assignment of a cartesian representable \ac{FVDC} $\Bil[\mf{p}]$
    to an elementary existential fibration $\mf{p}\colon\one{E}\to\one{B}$
    extends to a 2-functor $\Bil\colon\Fib_{\lc*=+}\to\Eqp\carttwo$.
\end{corollary}

\begin{remark}
    \label{rem:comparisonwithshulman}
    The paper \cite{shulmanFramedBicategoriesMonoidal2009} gives two sufficient conditions 
    for a monoidal fibration to induce a framed bicategory.
    Since our definition concentrates on cartesian fibrations,
    let us compare the two conditions with our results in this case,
    although there seem no obstacles to generalizing our results to monoidal fibrations
    and the corresponding comparison.

    In both conditions, the given cartesian fibration is required to be a bifibration,
    which holds true for elementary existential fibrations (cf. discussion in \Cref{rem:regularfibrationandeed}).
    The first condition in \textit{loc. cit.} is 
    that it is \textit{strongly BC}, which means that the Beck-Chevalley condition holds for all pullback squares.
    By \Cref{cor:FrobeniusfromBeckChevalley}, this is a stronger assumption than that the fibration is elementary existential. 
    
    The second condition is that the fibration is \textit{internally closed},
    which means that the fibration is fiberwise cartesian closed
    and base changes preserve the internal homs,
    and is \textit{weakly BC}, 
    which means that the Beck-Chevalley condition holds for pullback squares with one of the legs
    being a projection of a binary product.
    Examining the proofs in \cite[\S 17]{shulmanFramedBicategoriesMonoidal2009},
    we see that the internal closure of the fibration is used to show
    that it satisfies what we call the Frobenius reciprocity for arbitrary arrows (cf. (16.7) and (16.8) 
    in \textit{loc. cit.}),
    and the Beck-Chevalley conditions they use are the ones
    with the concerning left adjoints are of the forms $\sum_{\tpl{0,0,1}}$ or $\sum_{\tpl{0}}$.
\end{remark}

The main theorem for this chapter is the following.
\begin{theorem}
    \label{thm:firsteefiscartdouble}
    For a cartesian fibration $\mf{p}$ to be elementary existential, it is necessary and sufficient that the \ac{VDC} $\Bil[\mf{p}]$ is a cartesian representable \ac{FVDC}, \text{i.e.,} a cartesian equipment.
\end{theorem}

To state the theorem more precisely in 2-categorical terms,
we need to introduce the following definition.

\begin{definition}
    A 2-functor $\bi{f}\colon\bi{K}\to\bi{L}$ is a local inclusion if
    for each pair of 0-cells $k$ and $k'$ in $\bi{K}$,
    the (1-)functor 
    \[
    \bi{f}_{k,k'}\colon\bi{K}(k,k')\to\bi{L}(\bi{f}(k),\bi{f}(k'))
    \]
    that $\bi{f}$ induces is fully faithful and injective on objects.
\end{definition}

\begin{lemma}
    \label{lem:localinclusion}
    For a commutative square of 2-functors
    \[
    \begin{tikzcd}
        \bi{K}'
        \ar[r, "\bi{f}'"]
        \ar[d, "\bi{i}"']
        &
        \bi{L}'
        \ar[d, "\bi{j}"]
        \\
        \bi{K}
        \ar[r, "\bi{f}"']
        &
        \bi{L}
    \end{tikzcd}
    \]
    where $\bi{i}$ and $\bi{j}$ are local inclusions and isofibrations,
    then it is a pullback if and only if
    \begin{enumerate}
        \item a 0-cell $k$ in $\bi{K}$ is (essentially) in the image of $\bi{i}$ 
        precicely when $\bi{f}(k)$ is (essentially) in the image of $\bi{j}$, and
        \item a 1-cell $t\colon \bi{i}(k'_0)\to\bi{i}(k'_1)$ in $\bi{K}'$,
        in which $k'_0$ and $k'_1$ are 0-cells in $\bi{K}'$,
        is in the image of $\bi{i}$
        precicely when $\bi{f}'(t)\colon \bi{j}(\bi{f}'(k'_0))\to\bi{j}(\bi{f}'(k'_1))$ is
        in the image of $\bi{j}$.
        \qedhere
    \end{enumerate}
\end{lemma}
\begin{proof}
    Since the 2-category of 2-categories allows a concrete description of strict pullbacks, 
    what we need to show is that 
    the induced 2-functor $\langle\bi{i},\bi{f}'\rangle\colon\bi{K}'\to\bi{K}\times_{\bi{L}}\bi{L}'$ is an equivalence
    if and only if the two conditions hold.
    Here, we use the Whitehead theorem for 2-categories \cite[\S 1.11]{kellyBasicConceptsEnriched2005}:
    a 2-functor is an equivalence (in the 2-category of 2-categories)
    if and only if it is essentially surjective on 0-cells and 
    locally isomorphic (not an equivalence).

    The essential surjectivity on 0-cells for $\langle\bi{i},\bi{f}'\rangle$
    states that for each pair of 0-cells $k$ and $l'$ in $\bi{K}'$ and $\bi{L}'$
    such that $\bi{f}'(k)=\bi{j}(l')$,
    there is a 0-cell $k'$ in $\bi{K}$ such that $\langle\bi{i},\bi{f}'\rangle(k')\cong(k,l')$.
    This is equivalent to the first condition in the lemma since $\bi{i}$ and $\bi{j}$ are isofibrations.
    For the second condition, note that we have 
    \[
    \begin{tikzcd} 
        \bi{K}'(k'_0,k'_1)
        \ar[r, "{\langle\bi{i},\bi{f}'\rangle}"]
        \ar[dr,"\bi{i}"']
        \ar[dr, phantom, "\circlearrowright", very near end, shift left=1em]
        &
        \bi{K}\left(\bi{i}(k'_0),\bi{i}(k'_1)\right)\times_{\bi{L}\left(\bi{j}(\bi{f}'(k'_0)),\,\bi{j}(\bi{f}'(k'_1))\right)}\bi{L}'\left(\bi{f}'(k'_0),\bi{f}'(k'_1)\right) 
        \ar[d, "\cong"]
        \\
        &
        \left\{\,
        t\colon\bi{i}(k'_0)\rightarrow\bi{i}(k'_1)
        \mid
        \bi{f}'(t)\text{ is in the image of }\bi{j}\,
        \right\}
        \underset{\text{full}}{\subseteq}
        \bi{K}\left(\bi{i}(k'_0),\bi{i}(k'_1)\right)
    \end{tikzcd},
    \]
    since $\bi{j}$ induces the inclusion between the hom-categories of $\bi{L}'$ and $\bi{L}$. 
    As $\bi{i}$ is also a local inclusion, 
    the $\bi{i}$ in the diagram above is isomorphic if and only if it is bijective
    on objects (1-cells in the 2-categories), which is (ii) in the lemma.
\end{proof}

\begin{theorem}[{Restatement of \Cref{thm:firsteefiscartdouble} functorially}]
    \label{thm:eefiscartdouble}
    We have the following pullback square of 2-functors:
    \[
    \begin{tikzcd}
        \Fib_{\lc*=+}
        \ar[r]
        \ar[d]
        \ar[dr, phantom, "\lrcorner", very near start]
        &
        \Eqp\carttwo
        \ar[d]
        \\
        \Fib\carttwo
        \ar[r, "\Bil"']
        &
        \FVDbl\carttwo.
    \end{tikzcd}
    \]
\end{theorem}

\begin{proof}
    We have shown that the assignment of $\Bil[\mf{p}]$ in $\Eqp\carttwo$ to an elementary existential fibration $\mf{p}$
    extends to a 2-functor $\Bil\colon\Fib_{\lc*=+}\to\Eqp\carttwo$.
    Since $\Fib_{\lc*=+}\to\Fib\carttwo$ and $\Eqp\carttwo\to\VDbl-\carttwo$ are local inclusions and isofibrations,
    we can apply \Cref{lem:localinclusion} to show that the square is a pullback.

    Therefore, what we need to show is
    \begin{enumerate}
        \item for a cartesian fibration $\mf{p}\colon\one{E}\to\one{B}$,
            if $\Bil[\mf{p}]$ is a cartesian equipment\footnote{
            Since the 2-functor $\Eqp\carttwo\to\Eqp\carttwo$ is essentially surjective on 0-cells,
            $\Bil[\mf{p}]$ being a cartesian equipment is equivalent to saying that it is in the image of
            the 2-functor $\Eqp\carttwo\to\FVDbl\carttwo$.
            The only reason we phrase it in this way is to make the statement more readable.
            },
            then $\mf{p}$ is an elementary existential fibration.
        \item for elementary existential fibrations $\mf{p}$ and $\mf{p}'$
        and $\mf{f}\colon\mf{p}\to\mf{p}'$ in $\Fib_{\lc*}$,
        if $\Bil[\mf{f}]$ is a morphism of cartesian equipments,
        then $\mf{f}$ is a morphism of elementary existential fibrations.
    \end{enumerate}
    The second statement is easy.
    If $\Bil[\mf{f}]$ is a morphism of cartesian equipments,
    it preserves any spine cells,
    which means every coproducts including $\Phi_=$- and $\Phi_\exists$-coproducts.
    We will show the first statement.

    Since an equipment $\Bil[\mf{p}]$ gives a bifibration 
    $\Bil[\mf{p}]_1\to\one{B}\times\one{B}$,
    we can recover $\mf{p}$ as a base change of this along the functor 
    $\one{B}\to\one{B}\times\one{B}$ that sends an object $I$ to $(I,1)$,
    and hence, we have that $\mf{p}$ is a bifibration.
    Let us explain this in explicit terms. 
    Recall that a cell $\xi$ in $\Bil[\mf{p}]$ as in the left diagram in \Cref{eq:cellinBil}
    corresponds to an arrow $\xi\colon\alpha\to\beta[f]$ in $\one{E}_{I}$.
    Since $\Bil[\mf{p}]$ is a cartesian equipment,
    we know that there is a oprestriction $f^*\alpha$ of $\alpha$ along $f$ and $\id_1$.
    Using this, we find that the arrow $\xi$ corresponds to the arrow $\wt\xi$ 
    in the right diagram in \Cref{eq:cellinBil}.
    \begin{equation}
        \label{eq:cellinBil}
    \begin{tikzcd}[virtual]
        I
        \sar[r, "\alpha"]
        \ar[d, "f"']
        \ar[dr, phantom, "\xi"]
        &
        1
        \ar[d, equal]
        \\
        J
        \sar[r, "\beta"']
        &
        1
    \end{tikzcd}
    \hspace{2em}
    \vline\,\vline 
    \hspace{2em}
    \begin{tikzcd}[virtual]
        J
        \sar[r, "f^*\alpha"]
        \ar[d, equal]
        \ar[dr, phantom, "\wt\xi"]
        &
        1
        \ar[d, equal]
        \\
        J
        \sar[r, "\beta"']
        &
        1
    \end{tikzcd}
    \end{equation}
    This shows that the oprestriction functor $f^*(-)\colon\one{E}_{I}\to\one{E}_{J}$ gives a left adjoint 
    to the base change functor $(-)[f]\colon\one{E}_{J}\to\one{E}_{I}$.
    Therefore, it remains to show that $\mf{p}$ satisfies the Beck-Chevalley condition
    and the Frobenius reciprocity for the base change functors along 
    arrows in $\Phi_{\exists}$ and $\Phi_{=}$.

    For the Beck-Chevalley condition, note that the canonical arrow which we need to show is an isomorphism
    is the one in \Cref{lem:BeckChevalleycell} for the cartesian equipment $\Bil[\mf{p}]$,
    where we take 
    the pullback square
    in $\Phi_{=}$ and $\Phi_{\exists}$
    and $M=1$.
    Since those pullback squares are the products of the pullback squares in which one side is the identity arrow,
    we can apply \Cref{lem:BeckChevalleyclosure} to show that these are Beck-Chevalley pullback squares,
    and then apply \Cref{lem:BeckChevalleycell} to show that the canonical arrow is an isomorphism.

    For the Frobenius reciprocity for $\Phi_{\exists}$, 
    apply \Cref{lem:Frobenius} to the case where $f$ is the
    projection $\tpl{1}\colon I\times J\to J$
    and $K$ is the terminal object $1$.
    Here, the assumption that the pullback of $\tpl{0,0}\colon J\to J\times J$ and $\tpl{1}\times\id_J=\tpl{1,2}\colon I\times J\times J\to J\times J$ is a Beck-Chevalley pullback square
    follows again from \Cref{lem:BeckChevalleyclosure}:
    \[
        \begin{tikzcd}[virtual]
            I\times J
            \ar[r, "\tpl{0,1,1}"]
            \ar[d, "\tpl{1}"']
            \ar[dr, phantom, "\lrcorner", very near start]
            &
            I\times J\times J
            \ar[d, "\tpl{1,2}"]
            \\
            J
            \ar[r, "\tpl{0,0}"']
            &
            J\times J
        \end{tikzcd}
        \hspace{2em}
        \cong 
        \hspace{2em}
        \begin{tikzcd}[virtual]
            I
            \ar[r, "\id_I"]
            \ar[d, "\tpl{}"']
            \ar[dr, phantom, "\lrcorner", very near start]
            &
            I
            \ar[d, "\tpl{}"]
            \\
            1
            \ar[r, "\id_1"']
            &
            1
        \end{tikzcd}
        \hspace{2em}
        \times 
        \hspace{2em}
        \begin{tikzcd}[virtual]
            J
            \ar[r, "\tpl{0,0}"]
            \ar[d, "\id_J"']
            \ar[dr, phantom, "\lrcorner", very near start]
            &
            J\times J
            \ar[d, "\id_{J\times J}"]
            \\
            J
            \ar[r, "\tpl{0,0}"']
            &
            J\times J
        \end{tikzcd}
    \]

    The Frobenius reciprocity for $\Phi_{=}$ is not as straightforward as the other cases. 
    Recall that we need to show that the canonical arrow from 
    $\sum_{\tpl{0,0,1}}(\alpha\land\beta[{\scriptstyle\tpl{0,0,1}}])$ to $(\sum_{\tpl{0,0,1}}\alpha)\land\beta$ is an isomorphism. 
    The strategy we take is as follows.
    {
    \renewcommand{\theenumi}{\textbf{\Roman{enumi}}}
    \begin{enumerate}
        \item We reduce the problem to the case where $\alpha$ is the terminal
        and simultaneously present another equivalent form of the canonical arrow up to isomorphism. 
        \item We show that the canonical arrow has a retraction.
        \item We show that the retraction in turn has a retraction,
        which implies that the original arrow is an isomorphism.
    \end{enumerate}
    }
    In the following, 
    when we index the isomorphisms with the symbol $\cong$ without defining them explicitly,
    we mean that the isomorphisms are the canonical ones
    that are induced by the iterated base changes and the Beck-Chevalley condition.

    \begin{flushleft}
        \textbf{(I)} Reduction to a special case.
    \end{flushleft}
    We have already shown that $\mf{p}$ is a bifibration,
    and since the composites including the units in an \ac{FVDC} are unique up to isomorphism,
    we can identify the unit on $I$ with the object $\sum_{\tpl{0,0}}\top_I$ in $\one{E}_{I\times I}$
    and the composite of $\alpha\colon I\sto J$ and $\beta\colon J\sto K$ with the object 
    $\sum_{\tpl{0,2}}(\alpha[{\scriptstyle\tpl{0,1}}]\land\beta[{\scriptstyle\tpl{1,2}}])$ in $\one{E}_{I\times K}$,
    with the composing cells defined by the units of the adjunctions.
    
    We postpone the proof of the following claim to the end of this proof.
    \begin{claim}
        \label{claim:Frobenius}
        The natural transformation $\sum_{\tpl{0,0,1}}\Rightarrow(\sum_{\tpl{0,0,1}}\top_{I\times K})\land(-)[{\scriptstyle\tpl{1,2}}]$ 
        defined by the universal property of the local binary product $\land$ from the natural transformations
        \[
        \begin{tikzcd}[virtual]
            &
            \one{E}_{I\times K}
            \ar[rd, "{\sum_{\tpl{0,0,1}}}"]
            &
            \\
            \one{E}_{I\times K}
            \ar[rd, "{(-)[{\scriptstyle\tpl{1,2}}]}"']
            \ar[ru, equal]
            \ar[r, phantom, "{\rotatebox{-45}{$\cong$}}", near start]
            &
            \!
            \ar[r, phantom, "{\rotatebox{-45}{$\Rightarrow$}}", near end]
            &
            \one{E}_{I\times I\times K}
            \\
            &
            \one{E}_{I\times I\times K}
            \ar[uu, "{(-)[{\scriptstyle\tpl{0,0,1}}]}"{description}]
            \ar[ur, equal]
        \end{tikzcd}
        \hspace{2em},
        \hspace{2em}
        \begin{tikzcd}[virtual]
            &
            \one{E}_{I\times K}
            \ar[rd, "{\sum_{\tpl{0,0,1}}}"]
            &
            \\
            \one{E}_{I\times K}
            \ar[rd, "!"']
            \ar[ru, equal]
            \ar[r, phantom, "{\rotatebox{-45}{$\Rightarrow$}}", near start]
            &
            \!
            \ar[r, phantom, "{\rotatebox{-45}{$\cong$}}", near end]
            &
            \one{E}_{I\times I\times K}
            \\
            &
            \bf{1}
            \ar[uu, "{\top_{I\times K}}"{description}]
            \ar[ur, "{\sum_{\tpl{0,0.1}}\top_{I\times K}}"']
        \end{tikzcd}
        \]
        is an isomorphism.
        Here, the 2-cells in the diagram are the unit and counit of the adjunctions.
    \end{claim}
    Using this claim, we have the following commutative diagram
    for a pair of objects $\alpha\in\one{E}_{I\times K}$ and $\beta\in\one{E}_{I\times I\times K}$
    with the vertical arrows being isomorphisms.
    The whole diagram is natural in $\alpha$ and $\beta$.
    \[
    \begin{tikzcd}
        \sum_{\tpl{0,0,1}}(\alpha\land\beta[{\scriptstyle\tpl{0,0,1}}])
        \ar[d, "{\text{\Cref{claim:Frobenius}}}\quad\cong"']
        \ar[r, "(i)"]
        \ar[rd, phantom, "\circlearrowright"]
        &
        (\sum_{\tpl{0,0,1}}\alpha)\land\beta 
        \ar[d, "\cong\quad\text{\Cref{claim:Frobenius}}"]
        \\
        (\sum_{\tpl{0,0,1}}\top_{I\times K})\land 
        \left(\alpha\land\beta[{\scriptstyle\tpl{0,0,1}}]\right)[{\scriptstyle\tpl{1,2}}]
        \ar[r]
        \ar[d, "\cong"']
        \ar[rd, phantom, "\circlearrowright"]
        &
        \left(\sum_{\tpl{0,0,1}}\top_{I\times K}\land\alpha[{\scriptstyle\tpl{1,2}}]\right)\land\beta
        \ar[d, "\cong"]
        \\
        (\sum_{\tpl{0,0,1}}\top_{I\times K})\land  
        \alpha[{\scriptstyle\tpl{1,2}}]\land\beta[{\scriptstyle\tpl{1,1,2}}]
        \ar[r, "(ii)"']
        &
        (\sum_{\tpl{0,0,1}}\top_{I\times K})\land\alpha[{\scriptstyle\tpl{1,2}}]\land\beta
    \end{tikzcd}
    \]
    The Frobenius reciprocity for $\Phi_{=}$ states that 
    the canonical arrow $(i)$ is an isomorphism,
    so it suffices to show that the canonical arrow $(ii)$ is an isomorphism.
    By looking at how $(i)$ are defined and the naturality on the unique arrow $!\colon\beta\to\top_{I\times K\times K}$ in $\one{E}_{I\times I\times K}$, 
    we can see that $(ii)$ is compatible with the projections to $\sum_{\tpl{0,0,1}}\top_{I\times K}$
    and $\alpha[{\scriptstyle\tpl{1,2}}]$.
    Therefore, the proof reduces to showing that $(ii)$ is an isomorphism in the case 
    where $\alpha$ is the terminal object $\top_{I\times K}$.
    In this case, the above diagram that defines $(ii)$ 
    becomes the following diagram:
    \begin{equation}
        \label{eq:reduction}
    \begin{tikzcd}[column sep=2.5em]
        \sum_{\tpl{0,0,1}}(\beta[{\scriptstyle\tpl{0,0,1}}])
        \ar[d, "{\text{\Cref{claim:Frobenius}}}\quad\cong"',"\iota"]
        \ar[r, "\fpl{\sum_{\tpl{0,0,1}}!,\ \varepsilon_{\beta}}"]
        &
        (\sum_{\tpl{0,0,1}}\top_{I\times K})\land\beta
        \\
        (\sum_{\tpl{0,0,1}}\top_{I\times K})\land
        \beta[{\scriptstyle\tpl{1,1,2}}]
        \ar[ru, "\zeta"']
    \end{tikzcd}
    \end{equation}
    By definition, the counterpart of $(i)$ is the pairing 
    of the image of the unique arrow $!\colon\beta\to\top_{I\times K\times K}$ by $\sum_{\tpl{0,0,1}}$ and
    the counit component $\varepsilon_{\beta}$ of the adjunction $\sum_{\tpl{0,0,1}}\dashv(-)[{\scriptstyle\tpl{0,0,1}}]$
    at $\beta$.
    Here $\fpl{-,-}$ is the pairing of arrows in the fiber categories.

    \begin{flushleft}
        \textbf{(II)} Construction of a retraction of $\zeta$.
    \end{flushleft}
    Let $\nu$ be the arrow that corresponds to the following cell $\nu$
    in the virtual double category $\Bil[\mf{p}]$,
    where $\delta_I\coloneqq\sum_{\tpl{0,0}}\top_I$ is the unit on $I$. 
    \begin{equation}
        \label{eq:nu}
    \begin{tikzcd}
        \delta_I[{\scriptstyle\tpl{0,1}}]
        \land\beta[{\scriptstyle\tpl{1,2,3}}]
        \ar[d, "\nu"']
        \\
        \beta[{\scriptstyle\tpl{0,2,3}}]
    \end{tikzcd}
    \hspace{-2em}
    \text{in}\ \one{E}_{I\times I\times I\times K}
    \hspace{1em}
    \vline 
    \hspace{1em}
    \begin{tikzcd}[column sep=2.5em,virtual]
        \!
        \ar[rd, phantom, "\eta_I", very near end]
        &
        I
        \sar[r, "\beta"]
        \ar[ld, equal]
        \ar[d, equal]
        \ar[dr, phantom, "\rotatebox{90}{$=$}"]
        &
        I\times K
        \ar[d, equal]
        \\
        I
        \sar[r, "\delta_I"']
        \ar[d, equal]
        \ar[rrd, "\nu", phantom]
        &
        I
        \sar[r, "\beta"']
        &
        I\times K
        \ar[d, equal]
        \\
        I
        \sar[rr, "\beta"']
        &
        &
        I\times K
    \end{tikzcd}
    =
    \begin{tikzcd}[column sep=2.5em]
        I
        \sar[r, "\beta"]
        \ar[d, equal]
        \ar[rd, phantom, "\rotatebox{90}{$=$}"]
        &
        I\times K
        \ar[d, equal]
        \\
        I
        \sar[r, "\beta"']
        &
        I\times K
    \end{tikzcd}
    \end{equation}
    Then, let $\xi\colon(\sum_{\tpl{0,0,1}}\top_{I\times K})\land\beta\to\beta[{\scriptstyle\tpl{1,1,2}}]$
    be the following composite
    \[
    \begin{tikzcd}
        (\sum_{\tpl{0,0,1}}\top_{I\times K})\land\beta
        \ar[r, "\cong"]
        &
        \left(
        \delta_I[{\scriptstyle\tpl{0,1}}]
        \land\beta[{\scriptstyle\tpl{1,2,3}}]
        \right)
        [{\scriptstyle\tpl{1,0,1,2}}]
        \ar[r, "{\nu[{\scriptstyle\tpl{1,0,1,2}}]}"{yshift=1ex}]
        &
        \beta[{\scriptstyle\tpl{0,2,3}}][{\scriptstyle\tpl{1,0,1,2}}]
        \ar[r, "\cong"]
        &
        \beta[{\scriptstyle\tpl{1,1,2}}]
    \end{tikzcd}.
    \]
    The first isomorphism uses the Beck-Chevalley condition for $\Phi_{=}$.
    Our claim is that $\fpl{\pi_0,\xi}$ is the inverse of $\zeta$,
    where $\pi_0$ is the 0-th projection.

    We show that $\fpl{\pi_0,\xi}$ is a retraction of $\zeta$.
    By \Cref{eq:reduction}, this is equivalent to showing that
    $\fpl{\pi_0,\xi}\circ\fpl{\sum_{\tpl{0,0,1}}!,\ \varepsilon_{\beta}}=\iota$.
    On the 0-th projection, the equation follows directly
    from the construction of $\iota$ in \Cref{claim:Frobenius}.
    The 1-st projection of the equation is the commutativity of the following diagram:  
    \[
    \begin{tikzcd}[column sep=4em,virtual]
        \sum_{\tpl{0,0,1}}(\beta[{\scriptstyle\tpl{0,0,1}}])
        \ar[dr, "\varepsilon"']
        \ar[r, "\fpl{\sum_{\tpl{0,0,1}}!,\ \varepsilon_{\beta}}"]
        \ar[rd, phantom, "\circlearrowright"{yshift=1ex, xshift=3ex}]
        &
        (\sum_{\tpl{0,0,1}}\top_{I\times K})\land\beta
        \ar[d, "\xi"]
        \\
        &
        \beta[{\scriptstyle\tpl{1,1,2}}]
    \end{tikzcd}.
    \]
    The adjunction $\sum_{\tpl{0,0,1}}\dashv(-)[{\scriptstyle\tpl{0,0,1}}]$ translates this to the commutativity of the diagram 
    \[
    \begin{tikzcd}[column sep=4em,virtual]
        \beta[{\scriptstyle\tpl{0,0,1}}]
        \ar[dr, "\cong"']
        \ar[drr, bend right=20, equal]
        \ar[r]
        \ar[rr,bend left=10, "\fpl{\eta\circ !,\ \id_{\beta[{\scriptstyle\tpl{0,0,1}}]}}"] 
        &
        \left((\sum_{\tpl{0,0,1}}\top_{I\times K})\land\beta\right)
        [{\scriptstyle\tpl{0,0,1}}]
        \ar[d, "{\xi[{\scriptstyle\tpl{0,0,1}}]}"]
        \ar[r, "\cong"]
        &
        \delta_I[{\scriptstyle\tpl{0,0}}]
        \land\beta[{\scriptstyle\tpl{0,0,1}}]
        \ar[d, "{\nu[\tpl{0,0,0,1}]}"]       
        \\
        &
        \beta[{\scriptstyle\tpl{1,1,2}}][{\scriptstyle\tpl{0,0,1}}]
        \ar[r, "\cong"']
        &
        \beta[{\scriptstyle\tpl{0,0,1}}]
    \end{tikzcd},
    \]
    which holds true by the definition of $\nu$ in \Cref{eq:nu}.

    \begin{flushleft}
        \textbf{(III)} Construction of a retraction of $\fpl{\pi_0,\xi}$.
    \end{flushleft}
    Considering the pairing with the 0-th projection
    and how the arrow $\xi$ is defined,
    it is enough to show that there exists an arrow $\lambda$
    that makes the following triangle commute:
    \begin{equation}
        \label{eq:lambda}
    \begin{tikzcd}
        (\delta_I[{\scriptstyle\tpl{0,1}}]\land\beta[{\scriptstyle\tpl{1,2,3}}])
        [{\scriptstyle\tpl{1,0,1,2}}]
        \ar[d, "\cong"']
        \ar[dr, phantom, "\circlearrowright"]
        \ar[r, "{\fpl{\pi_0,\nu}[{\scriptstyle\tpl{1,0,1,2}}]}"{yshift=1ex}] 
        &
        \left(
        \delta_I[{\scriptstyle\tpl{0,1}}]\land\beta[{\scriptstyle\tpl{0,2,3}}]
        \right)
        [{\scriptstyle\tpl{1,0,1,2}}]
        \ar[d, "\cong"]
        \\
        \delta_I[{\scriptstyle\tpl{1,0}}]\land\beta
        \ar[rd, "\pi_1"']
        \ar[r,"{\fpl{\pi_0,\nu[{\scriptstyle\tpl{1,0,1,2}}]}}"]
        &
        \delta_I[{\scriptstyle\tpl{1,0}}]\land\beta[{\scriptstyle\tpl{1,1,2}}]
        \ar[d, "\lambda"]
        \\
        &
        \beta
    \end{tikzcd}
    \text{in}\ \one{E}_{I\times I\times K}.
    \end{equation}
    Before the construction, we observe that we have the following cell $\mu$
    and the corresponding arrow $\mu\colon\delta_I[\tpl{0,1}]\land\delta_I[\tpl{1,2}]\to\delta_I[\tpl{0,2}]$ in $\one{E}_{I\times I\times I}$
    satisfying the following commutative diagram:
    \[
        \begin{tikzcd}[virtual,column sep=2em]
        &
        I
        \ar[ld, equal]
        \ar[d, equal]
        \ar[rd,equal]
        \\
        I
        \sar[r, "\delta_I"']
        \ar[d, equal]
        \ar[ru, phantom, "\eta_I"{yshift=-1ex, xshift=1ex}]
        \ar[rrd, phantom, "\mu"]
        &
        I
        \sar[r, "\delta_I"']
        &
        I
        \ar[lu, phantom, "\eta_I"{yshift=-1ex, xshift=-1ex}]
        \ar[d, equal]
        \\
        I
        \sar[rr, "\delta_I"']
        &
        &
        I
        \end{tikzcd}
        =
        \begin{tikzcd}[virtual,column sep=1em]
            &
            I
            \ar[ld, equal]
            \ar[rd, equal]
            \ar[d, "\eta_I", phantom]
            \\
            I
            \sar[rr, "\delta_I"']
            &
            \!
            &
            I
        \end{tikzcd}
        ,
        \hspace{0.5em}
        \begin{tikzcd}
            \top_{I}
            \ar[d, "\fpl{\eta_I,\eta_I}"']
            \ar[rd, "\eta_I"]
            \ar[rdd, phantom, "\circlearrowright"]
            \\
            \delta_I[{\scriptstyle\tpl{0,0}}]\land\delta_I[{\scriptstyle\tpl{0,0}}]
            \ar[d, "\cong"']
            &
            \delta_I[{\scriptstyle\tpl{0,0}}]
            \ar[d, "\cong"]
            \\
            \left(\delta_I[\tpl{0,1}]\land\delta_I[\tpl{1,2}]\right)[{\scriptstyle\tpl{0,0,0}}]
            \ar[r, "{\mu[{\scriptstyle\tpl{0,0,0}}]}"']
            &
            \delta_I[{\scriptstyle\tpl{0,2}}][{\scriptstyle\tpl{0,0,0}}]
        \end{tikzcd}
    \]
    Let $\sigma\colon \delta_I\to \delta_I[\tpl{1,0}]$ be the canonical isomorphism,
    which follows from the observation that $\delta_I[\tpl{1,0}]$ exhibits the same
    universal property as $\delta_I$.
    From the natural isomorphism between $(-)[\tpl{0,0,0}]$ and $(-)[\tpl{1,0,1}][\tpl{0,0}]$
    and from the above diagram, 
    we obtain the following commutative diagrams in $\one{E}_{I\times I}$: 
    \begin{equation}
        \label{eq:lambda2}
        \begin{tikzcd}[column sep=1em,row sep=2em]
        \delta_I
        \ar[drr, phantom, "\circlearrowright"]
        \ar[d, "\fpl{\sigma,\id_{\delta_I}}"{description}]
        \ar[r, "!"]
        &
        \top_I
        \ar[r, "\cong"]
        &
        \top_I[\tpl{1}]
        \ar[d, "{\eta_I[\tpl{1}]}"]
        \\
        \delta_I[\tpl{1,0}]\land\delta_I
        \ar[r, "{\mu[{\tpl{1,0,1}}]}"'{yshift=-1ex}]
        &
        \delta_I[\tpl{1,1}]
        \ar[r, "\cong"'{yshift=-1ex}]
        &
        \delta_I[\tpl{0,0}][\tpl{1}]
        \end{tikzcd}
        \begin{tikzcd}[column sep=tiny]\!\ar[r,maps to,"{(-)[\tpl{10}]}"{yshift=1ex}]&\!\end{tikzcd}
        \begin{tikzcd}[column sep=1em,row sep=2em]
            \delta_I[\tpl{1,0}]
            \ar[drr, phantom, "\circlearrowright"]
            \ar[d, "\fpl{\sigma\inv,\id_{\delta_I}}"{description}]
            \ar[r, "!"']
            &
            \top_I
            \ar[r, "\cong"]
            &
            \top_I[\tpl{0}]
            \ar[d, "{\eta_I[\tpl{0}]}"]
            \\
            \delta_I\land\delta_I[\tpl{1,0}]
            \ar[r, "{\mu[{\tpl{0,1,0}}]}"'{yshift=-1ex}]
            &
            \delta_I[\tpl{0,0}]
            \ar[r, "\cong"'{yshift=-1ex}]
            &
            \delta_I[\tpl{0,0}][\tpl{0}]
        \end{tikzcd}
    \end{equation}
    Note that the vertical arrow $!$ is also the counit component at $\top_I$.
    Now, we get back to the construction of $\lambda$.
    By the associativity of loose composition, or equivalently by the universal property,
    we have the following commutative diagram:
    \[
    \begin{tikzcd}
        \delta_I[{\scriptstyle\tpl{0,1}}]\land\delta_I[{\scriptstyle\tpl{1,2}}]\land\beta[{\scriptstyle\tpl{2,3,4}}] 
        \ar[r, "{\mu[\scriptstyle\tpl{0,1,2}]\land\id_{\beta[{\scriptstyle\tpl{2,3,4}}]}}"{yshift=1ex}]
        \ar[rd, phantom, "\circlearrowright"]
        \ar[d, "{\id_{\delta_I[{\scriptstyle\tpl{0,1}}]}\land\nu[{\scriptstyle\tpl{1,2,3,4}}]}"']
        &
        \delta_I[{\scriptstyle\tpl{0,2}}]\land\beta[{\scriptstyle\tpl{2,3,4}}]
        \ar[d, "{\nu[{\scriptstyle\tpl{0,2,3,4}}]}"]
        \\
        \delta_I[{\scriptstyle\tpl{0,1}}]\land\beta[{\scriptstyle\tpl{1,3,4}}]
        \ar[r, "{\nu[{\scriptstyle\tpl{0,1,3,4}}]}"']
        &
        \beta[{\scriptstyle\tpl{0,3,4}}]
    \end{tikzcd}
    \quad\text{in}\quad \one{E}_{I\times I\times I\times I\times K}.
    \] 
    Sending this diagram by $(-)[{\scriptstyle\tpl{0,1,0,1,2}}]$,
    we obtain the following commutative diagram: 
    \begin{equation}
        \label{eq:lambda3}
    \begin{tikzcd}[column sep=4em]
        \delta_I[{\scriptstyle\tpl{0,1}}]\land\delta_I[{\scriptstyle\tpl{1,0}}]\land\beta 
        \ar[r, "{\mu[\scriptstyle\tpl{0,1,0}]\land\id_{\beta}}"{yshift=1ex}]
        \ar[rd, phantom, "\circlearrowright"]
        \ar[d, "{\id_{\delta_I[{\scriptstyle\tpl{0,1}}]}\land\nu[{\scriptstyle\tpl{1,0,1,2}}]}"']
        &
        \delta_I[{\scriptstyle\tpl{0,0}}]\land\beta
        \ar[d, "{\nu[{\scriptstyle\tpl{0,0,1,2}}]}"]
        \\
        \delta_I[{\scriptstyle\tpl{0,1}}]\land\beta[{\scriptstyle\tpl{1,1,2}}]
        \ar[r, "{\nu[{\scriptstyle\tpl{0,1,1,2}}]}"']
        &
        \beta
    \end{tikzcd}
    \quad\text{in}\quad \one{E}_{I\times I\times K}.
    \end{equation}
    These diagrams can be combined into the desired diagram \Cref{eq:lambda}:
    \[
    \begin{tikzcd}[column sep=6em, row sep=2.5em]
        \delta_I[{\scriptstyle\tpl{1,0}}]\land\beta 
        \ar[r, "\pi_1"]
        \ar[d, "\fpl{\sigma\inv\circ\pi_0,\,\pi_0,\,\pi_1}"{description}]
        \ar[ddd, bend right=50, "\fpl{\pi_0,\,\nu[{\scriptstyle\tpl{1,0,1,2}}]}"', shift right=4em, near start]
        \ar[rd, phantom, "\text{\Cref{eq:lambda2}}"]
        &
        \beta 
        \ar[d,"{\eta_I[\tpl{0}]\land\id_{\beta}}"{description}]
        \ar[dd, bend left=90, equal, shift left=2.5em]
        \ar[dd, phantom, shift left=4em, "\text{\Cref{eq:nu}}"{yshift=2ex}]
        \\
        \delta_I[{\scriptstyle\tpl{0,1}}]\land\delta_I[{\scriptstyle\tpl{1,0}}]\land\beta
        \ar[d, "{\id_{\delta_I[{\scriptstyle\tpl{0,1}}]}\land\nu[{\scriptstyle\tpl{1,0,1,2}}]}"{description}]
        \ar[r, "{\mu[\scriptstyle\tpl{0,1,0}]\land\id_{\beta}}"{yshift=0.5ex}]
        \ar[d, phantom, "="{xshift=-6em}]
        \ar[rd, phantom, "\text{\Cref{eq:lambda3}}"]
        &
        \delta_I[{\scriptstyle\tpl{0,0}}]\land\beta
        \ar[d, "{\nu[{\scriptstyle\tpl{0,0,1,2}}]}"{description}]
        \\
        \delta_I[{\scriptstyle\tpl{0,1}}]\land\beta[{\scriptstyle\tpl{1,1,2}}]
        \ar[r, "{\nu[{\scriptstyle\tpl{0,1,1,2}}]}"']
        \ar[d, "\sigma\land\id_{\beta[{\scriptstyle\tpl{1,1,2}}]}"{description}]
        \ar[rd, phantom, "\rotatebox{135}{$\coloneqq$}"]
        &
        \beta
        \\
        \delta_I[{\scriptstyle\tpl{1,0}}]\land\beta[{\scriptstyle\tpl{1,1,2}}]
        \ar[ru, bend right=20, "\lambda"']
        &
        \!
    \end{tikzcd}
    \]
    Consequently, we have shown that $\fpl{\pi_0,\xi}$ has a retraction.
\end{proof}

\begin{proof}[{Proof of \Cref{claim:Frobenius}}]
    In a cartesian equipment, we have the following two  canonical isomorphisms:
    \[
    \begin{tikzcd}[column sep=6em,virtual]
        I\times I
        \sar[r, "{\delta_{I\times I}[\id\smcl\tpl{0,0}]}"]
        \ar[d, equal]
        \ar[rd, phantom, "\rotatebox{90}{$\cong$}"]
        &
        I
        \ar[d, equal]
        \\
        I\times I
        \sar[r, "{\delta_{I}[\tpl{0}\smcl\id]\land\delta_{I}[\tpl{1}\smcl\id]}"']
        &
        I
    \end{tikzcd},
    \hspace{2em}
    \begin{tikzcd}[virtual]
        I
        \sar[r, "{\delta_{I}}"]
        \ar[d, equal]
        \ar[rrd, phantom, "\rotatebox{90}{$\cong$}"]
        &
        I
        \sar[r, "{\gamma}"]
        &
        I\times K
        \ar[d, equal]
        \\
        I
        \sar[rr, "{\gamma}"']
        &&
        I\times K
    \end{tikzcd}
    \quad (\gamma\ \text{is arbitrary.})
    \]
    The first isomorphism derives from the cartesian condition for virtual double categories with units, 
    and the second isomorphism is the unitality of the horizontal composition. 
    These lead to the following isomorphisms in $\one{E}_{I\times I\times K}$:
    \[
    \begin{aligned} 
        \sum_{\tpl{0,0,1}}\alpha
        &\cong
        \sum_{\tpl{0,2,3}}\left(\delta_{I\times I}[\tpl{\tpl{0,2},\tpl{1,1}}]\land\alpha[\tpl{1,3}]\right)
        &
        &\text{by the presentation of $\sum_{\tpl{0,0,1}}$ in \Cref{eq:cellinBil}}\\
        &\cong
        \sum_{\tpl{0,2,3}}\left(\delta_{I}[\tpl{0,1}]\land\delta_{I}[\tpl{1,2}]\land\alpha[\tpl{1,3}]\right) 
        &
        &\text{by the first isomorphism above}\\
        &\cong
        \delta_{I}[\tpl{0,1}]\land\alpha[\tpl{1,2}]
        &
        &\text{by the second isomorphism above}\\
        &\cong
        \sum_{\tpl{0,0,1}}\top_{I\times K}\land\alpha[\tpl{1,2}]
        &
        &\text{by the Beck-Chevalley condition for $\Phi_{=}$}.
    \end{aligned}
    \]
    By tracing the isomorphisms, 
    one can see that this isomorphism is the desired one.
\end{proof}

See \Cref{ex:regularfibrationtoequip} for examples of the use of \Cref{thm:firsteefiscartdouble}.

\begin{remark}
    \label{remark:handinhand}
    What we have proved is that the fibration $\mf{p}$ is elementary existential when $\Bil[\mf{p}]$ has units and binary loose composition that are compatible with the cartesian structure. 
    It is worth noting that both the units and the binary loose composition are crucial for each of the properties of the fibration.
    Neither the existence of units nor the existence of binary loose composition alone would make the proof here work.
\end{remark}

\section{Elementary existential fibrations and Frobenius cartesian equipments}
\label{sec:characterization}
\subsection{Regular fibrations and cartesian equipments with Beck-Chevalley pullbacks}
\label{subsec:BeckChevalleyRegular}
A stricter class of fibrations than elementary existential fibrations is that of regular fibrations.
It will turn out to be the relationship between the Beck-Chevalley conditions for fibrations and the Beck-Chevalley pullbacks in cartesian equipments that characterizes the regular fibrations by a condition on the corresponding cartesian equipments.

We start by recalling some consequences of Beck-Chevalley pullbacks in cartesian equipments.
Applying \Cref{lem:BeckChevalleycell,lem:Frobenius}, we have the following corollaries.

\begin{proposition}
    \label{cor:BeckChevalleydeduction}
    Let $\mf{p}\colon\one{E}\to\one{B}$ be an elementary existential fibration. If a pullback square in $\one{B}$ is a Beck-Chevalley pullback square in $\Bil[\mf{p}]$, then $\mf{p}$ satisfies the Beck-Chevalley condition for the pullback square in both directions.
\end{proposition}
\begin{proof}
    In the statement of \Cref{lem:BeckChevalleycell},
    take $M$ to be the terminal object in the double category $\Bil[\mf{p}]$.
    Then the canonical cell $\sigma$ in \Cref{lem:BeckChevalleycell} 
    reduces to the component of the 
    canonical natural transformation in the definition of the Beck-Chevalley condition
    at $\alpha$.
\end{proof}

\begin{proposition}
    \label{cor:FrobeniusfromBeckChevalley}
    Let $\mf{p}\colon\one{E}\to\one{B}$ be an elementary existential fibration.
    If the pullback of $f\times\id_J\colon I\times J\to J\times J$ along the diagonal $\tpl{0,0}\colon J\to J\times J$ is a Beck-Chevalley pullback square in $\Bil[\mf{p}]$ for $f$ in $\one{B}$, then $\mf{p}$ satisfies the Frobenius reciprocity for $f$.

    In particular, an elementary existential fibration $\mf{p}$ with the Beck-Chevalley condition for all pullback squares in $\one{B}$ satisfies the Frobenius reciprocity for all arrows in $\one{B}$.
\end{proposition}
\begin{proof}
    In the statement of \Cref{lem:Frobenius}, take $K$ to be the terminal object in the double category $\Bil[\mf{p}]$.
\end{proof}

\begin{proposition}
    \label{lemma:BeckChevalleyfibration}
    Let $\mf{p}\colon\one{E}\to\one{B}$ be an elementary existential fibration.
    A pullback square in $\one{B}$ on the left below is a Beck-Chevalley pullback in $\Bil[\mf{p}]$ if $\mf{p}$ satisfies the Beck-Chevalley condition for the square as right below, when we see it as an arrow from $\tpl{f,g}$ to $\tpl{0,0}$ in $\one{B}^{\rightarrow}$.
    \[
        \begin{tikzcd}
            I
            \ar[r, "f"]
            \ar[d, "g"']
            \pullback
            &
            J
            \ar[d, "h"]
            \\
            K
            \ar[r, "k"']
            &
            L
        \end{tikzcd}
        \hspace{4em}
        \begin{tikzcd}
            I
            \ar[d, "\tpl{f,g}"']
            \ar[r]
            \pullback
            &
            L
            \ar[d, "\tpl{0,0}"]
            \\
            J\times K
            \ar[r, "h\times k"']
            &
            L\times L
        \end{tikzcd}
    \]
\end{proposition}
\begin{proof}
    The original square is a Beck-Chevalley pullback in $\Bil[\mf{p}]$
    when 
    the following canonical arrow is an isomorphism:
    \[
        \begin{tikzcd}
            \sum_{f\times g}(\delta_I)
            \ar[r,"\cong"]
            &
            \sum_{\tpl{f,g}}\top_I
            \ar[r]
            &
            \delta_L[k\times h]
        \end{tikzcd}.
    \]
    Since $\delta_L=\sum_{\tpl{0,0}}\top_L$,
    the above isomorphism is equivalent to the component at $\top_L$
    of the canonical transformation
    for the Beck-Chevalley condition for the second square in the lemma.
\end{proof}

\begin{remark}
    \Cref{cor:BeckChevalleydeduction} states that the Beck-Chevalley condition on a pullback square in $\one{B}$ for a fibration $\mf{p}$ implies the Beck-Chevalley pullback condition on the same pullback square for the corresponding cartesian equipment $\Bil[\mf{p}]$.
    \Cref{lemma:BeckChevalleyfibration} may seem to state the converse of this statement, but the pullback squares in question are not the same.
    However, this is enough to state the equivalence of the Beck-Chevalley condition and the Beck-Chevalley pullback condition globally, as we will see in \Cref{cor:BeckChevalleyall}.

    We could not find exactly how the class of pullback squares in $\one{B}$ for which the Beck-Chevalley condition holds in $\mf{p}$ corresponds to the class of the Beck-Chevalley pullbacks in $\Bil[\mf{p}]$ in general.
    Another delicate point is that the latter is closed under taking products but the former is not in general.
    Classes of the pullback squares in $\one{B}$ with respect to these two conditions should be investigated further.
    Not much is known about the subtleties of the Beck-Chevalley conditions but \cite{Seely83,Law15} give detailed discussions on the related topics.
\end{remark}

\begin{corollary}
    \label{cor:BeckChevalleyall}
    Let $\mf{p}\colon\one{E}\to\one{B}$ be an elementary existential fibration
    with $\one{B}$ having finite limits.
    Then, the following are equivalent:
    \begin{enumerate}
        \item $\mf{p}$ is a regular fibration,
        that is, the Beck-Chevalley condition holds for all pullback squares in $\one{B}$
        (see \Cref{cor:FrobeniusfromBeckChevalley}).
        \item $\Bil[\mf{p}]$ has the Beck-Chevalley pullbacks.
    \qedhere
    \end{enumerate}
\end{corollary}
\begin{proof}
    The implication (i)$\Rightarrow$(ii) follows from \Cref{lemma:BeckChevalleyfibration},
    and the implication (ii)$\Rightarrow$(i) follows from \Cref{cor:BeckChevalleydeduction}.
\end{proof}

\begin{example}
    \label{ex:regularfibrationtoequip}
    Combining the classical results on the examples of fibrations in \Cref{ex:regularfibration}
    and the results in this section,
    we have the following characterizations of the \acp{CFVDC},
    which are mentioned in \Cref{example:composablity}.
    \begin{enumerate}
        \item For a category $\one{B}$ with finite limits, 
        the \ac{CFVDC} $\Span[\one{B}]$ is a cartesian equipment with Beck-Chevalley pullbacks 
        since $\one{B}^{\rightarrow}\to\one{B}$ is a regular fibration.
        \item For a category $\one{B}$ with finite limits,
        the \ac{CFVDC} $\Rel[\one{B}]$ is a cartesian virtual equipment. It is a cartesian equipment 
        if and only if $\one{B}$ is a regular category.
        In this case, all pullbacks in $\one{B}$ are Beck-Chevalley pullbacks in $\Rel[\one{B}]$.
        \item For a cartesian monoidal category $\one{B}$,
        the \ac{CFVDC} $\Mat[\one{B}]$ is a cartesian equipment if and only if $\one{B}$ has distributive small coproducts.
        In this case, all pullbacks in $\one{B}$ are Beck-Chevalley pullbacks in $\Mat[\one{B}]$.
        \qedhere
    \end{enumerate}
\end{example}

\subsection{The Frobenius axiom and recovering fibrations from cartesian equipments}
\label{subsec:Frobenius}

We turn our attention to determining when a cartesian equipment is of the form $\Bil[\mf{p}]$ for an elementary existential fibration $\mf{p}$.
To this end, we introduce the Frobenius axiom for cartesian equipments, and show that it is a characteristic property of cartesian equipments of the form $\Bil[\mf{p}]$.

\begin{equation}
    \label{eq:Frobenius}
    \begin{tikzcd}[column sep=small, row sep=small]
        &
        I
        \ar[dl, "\tpl{0,0}"']
        \ar[dr, "\tpl{0,0}"]
        \\
        I\times I
        \ar[dr, "\tpl{0,0,1}"']
        &&
        I\times I
        \ar[dl, "\tpl{0,1,1}"]
        \\
        &
        I\times I\times I
    \end{tikzcd}
\end{equation}
Using the notation $\Delta=\tpl{0,0}$,
the above square being a Beck-Chevalley pullback means the canonical isomorphism
$\Delta^*\Delta\cong(\Delta\times\delta_I)_*(\delta_I\times\Delta)^*$ 
and $\Delta^*\Delta\cong(\delta_I\times\Delta)_*(\Delta\times\delta_I)^*$
hold.
In the context of cartesian bicategories,
this condition is known as the Frobenius axiom in \cite{WW08}\footnote{
This is called \textit{discreteness} in \cite{CW87,lambertDoubleCategoriesRelations2022}.
}.

\begin{definition}
    \label{def:Frobenius}
    Let $\dbl{D}$ be a cartesian equipment.
    An object $I$ in $\dbl{D}$ is said to be \emph{Frobenius} 
    if the pullback square \Cref{eq:Frobenius} is a Beck-Chevalley pullback.
    A cartesian equipment is said to be \emph{Frobenius} if every object in it is Frobenius.
\end{definition}

\begin{proposition}
    \label{prop:FrobeniusforBil}
    Let $\mf{p}\colon\one{E}\to\one{B}$ be an elementary existential fibration.
    Then, $\Bil[\mf{p}]$ is Frobenius.
\end{proposition}
\begin{proof}
    To see that $\Bil[\mf{p}]$ is Frobenius,
    by \Cref{lemma:BeckChevalleyfibration},
    it suffices to show that $\mf{p}$ satisfies the Beck-Chevalley condition
    for the following square:\footnote{we omit the product symbols $\times$ and the commas ``,'' in $\tpl{-}$ for the sake of readability.}
    \[
        \begin{tikzcd}
            I
            \ar[r, "\tpl{000}"]
            \ar[d, "\tpl{0000}"']
            \pullback
            &
            III
            \ar[d, "\tpl{012012}"]
            \\
            IIII 
            \ar[r, "\tpl{011223}"']
            &
            IIIIII
        \end{tikzcd}
        =
        \begin{tikzcd}[row sep=small]
            I
            \ar[r]
            \ar[d]
            \pullback
            &
            II
            \ar[d]
            \ar[r]
            \pullback
            &
            III
            \ar[d, "\tpl{0,1,2,0}"]
            \\
            II
            \ar[r]
            \ar[d]
            \pullback
            &
            III
            \ar[d]
            \ar[r]
            \pullback
            &
            IIII
            \ar[d, "\tpl{0,1,2,3,1}"]
            \\
            III
            \ar[d]
            \ar[r]
            \pullback
            &
            IIII
            \ar[d]
            \ar[r]
            \pullback
            &
            IIIII
            \ar[d, "\tpl{0,1,2,3,4,2}"]
            \\
            IIII
            \ar[r, "\tpl{0,1,1,2,3}"']
            &
            IIIII
            \ar[r, "\tpl{0,1,2,3,3,4}"']
            &
            IIIIII
        \end{tikzcd}.
    \]
    However, this can be decomposed into the six pullback squares as above,
    where the right bottom square belongs to $\Phi_{=}$,
    and the other five squares are of the form presented in \Cref{cor:eeffrobenius}.
    Since the Beck-Chevalley conditions are closed under pasting of pullback squares,
    the square above is a Beck-Chevalley pullback.
\end{proof}

\begin{remark}
    For any pushout square 
    \[
        \begin{tikzcd}
            N_1
            \ar[from=r, "f"']
            \ar[from=d, "g"]
            \pullback
            &
            N_2
            \ar[from=d, "h"']
            \\
            N_3
            \ar[from=r, "k"]
            &
            N_4
        \end{tikzcd} 
    \]
    in the category of finite sets
    with $N_1+N_4=N_2+N_3$ as natural numbers,
    the pullback square
    \[
        \begin{tikzcd}
            I^{N_1}
            \ar[r,"I^f"]
            \ar[d,"I^g"']
            \pullback
            &
            I^{N_2}
            \ar[d,"I^h"]
            \\
            I^{N_3}
            \ar[r,"I^k"']
            &
            I^{N_4}
        \end{tikzcd}
    \]
    can be decomposed as in the proof of \Cref{prop:FrobeniusforBil}.
\end{remark}

What makes the Frobenius axiom interesting is that it leads to the self-dual structure on the objects in the cartesian equipment. 
This was observed in \cite{CW87,WW08} in the context of cartesian bicategories, and in \cite{hoshinoDoubleCategoriesRelations2023} in the context of equipments\footnote{
In \cite{hoshinoDoubleCategoriesRelations2023}, the authors assume discreteness, which is a stronger condition than the Frobenius axiom,
to obtain the self-dual structure,
but the Frobenius axiom suffices for most of the results in the paper
except for Lemma 3.1.16.
}.
To discuss further results on the self-duality in Frobenius cartesian equipments,
we give a brief summary of the results in the paper\footnote{
We do not define what the self-dual structure is in this paper,
because the author is not confident at the moment
that a definition we have is of the correct generality.
A tentative definition is given \cite[Definition 3.1.11]{hoshinoDoubleCategoriesRelations2023}.
}.

\begin{proposition}
    \label{prop:Frobeniusselfdual}
    Suppose that $\dbl{D}$ is a Frobenius cartesian equipment.
    For any object $I$ in $\dbl{D}$,
    let $\iota_I\colon 1\sto I\times I$ and $\epsilon_I\colon I\times I\sto 1$
    be defined by the following oprestrictions:
    \[
        \begin{tikzcd}[column sep=small]
            &
            I
            \ar[dl, "!"']
            \ar[dr, "\tpl{0,0}"]
            \\
            1
            \sar["\iota_I"',rr]
            \ar[rr, "\opcart"{yshift=1em}, phantom]
            &&
            I\times I
        \end{tikzcd}
        \hspace{4em}
        \begin{tikzcd}[column sep=small]
            &
            I
            \ar[dl, "\tpl{0,0}"']
            \ar[dr, "!"]
            \\
            I\times I
            \ar[rr, "\opcart"{yshift=1em}, phantom]
            \sar["\epsilon_I"',rr]
            &&
            1
        \end{tikzcd}
    \]
    Then, the following hold:
    \begin{enumerate}
        \item They come equipped with isomorphsms $\zeta_I\colon\delta_I\Rightarrow(\iota_I\times\delta_I)(\delta_I\times\epsilon_I)$
        and $\theta_I\colon\delta_I\Rightarrow(\delta_I\times\iota_I)(\epsilon_I\times\delta_I)$.
        \item These data extend to the functors $\iota,\epsilon\colon\dbl{D}_0\to\dbl{D}_1$
        such that $\src\circ\iota=\tgt\circ\epsilon=1\circ !$
        and $\tgt\circ\iota=\src\circ\epsilon=\times\circ\Delta$.
        \item Objects in $\dbl{D}$ are self-dual in the sense of \cite[Definition 4.11]{Stay16},
        thus, $\LBi{\dbl{D}}$ is a compact closed bicategory.
        \item The dagger structure induced from the above self-dual structure
        extends to the whole of $\dbl{D}$
        as a indentity-on-tight-parts double functor $(-)^{\dagger}\colon\dbl{D}\lop\to\dbl{D}$.
        \item For a tight arrow $f\colon I\to J$ in $\dbl{D}$,
        $(f_*)^{\dagger}\cong f^*$ canonically,
        in particular, $(\delta_I)^{\dagger}\cong\delta_I$.
        \item Let $f\colon I_0\to I_1$, $g\colon J_0\to J_1$ be tight arrows in $\dbl{D}$,
        and $\left(\alpha_i\colon I_i\sto J_i, \beta_i\colon I_i\times J_i\sto 1, \gamma_i\colon J_i\sto I_i\right)$
        be triples of loose arrows in $\dbl{D}$ for $i=0,1$
        where $\alpha_i,\beta_i,\gamma_i$ correspond to each other using $\iota,\epsilon,\zeta,\theta$.
        Then, we also have the bijective correspondence between the cells below:
        \[
		\begin{tikzcd}
			I_0
			\sar[r, "\alpha_0"]
			\ar[d, "f"']
			\doublecell[rd]{\alpha}
				&
				J_0
				\ar[d, "g"]
			\\
			I_1
			\sar[r, "\alpha_1"']
				&
                J_1
		\end{tikzcd}
		\vline
		\,
		\vline
		\begin{tikzcd}
			I_0\times J_0
			\sar[r, "\beta_0"]
			\ar[d, "f\times g"']
			\doublecell[rd]{\beta}
				&
				1
				\ar[d, equal]
			\\
			I_1\times J_1
			\sar[r, "\beta_1"']
				&
				1
		\end{tikzcd}
		\vline
		\,
		\vline
		\begin{tikzcd}
			J_0
			\sar[r, "\gamma_0"]
			\ar[d, "g"']
			\doublecell[rd]{\gamma}
				&
				I_0
				\ar[d, "f"]
			\\
            J_1
			\sar[r, "\gamma_1"']
				&
                I_1
		\end{tikzcd}.
	\]
    This correspondence is functorial with respect to tightwise composition.
    In other words, when we define the category $(\LS{\dbl{D}})_1$ 
    by the pullback 
    \[
        \begin{tikzcd}[column sep=huge]
            (\LS{\dbl{D}})_1
            \ar[r]
            \ar[d, "\tpl{\src',\tgt'}"']
            \ar[rd, phantom, "\lrcorner", very near start]
            &
            \dbl{D}_1
            \ar[d, "\tpl{\src,\tgt}"]
            \\
            \dbl{D}_0\times\dbl{D}_0
            \ar[r, "{(I,J)\mapsto(I\times J,1)}"']
            &
            \dbl{D}_0\times\dbl{D}_0
        \end{tikzcd},
    \]
    the above correspondence gives
    the fibered equivalence $(\LS{\dbl{D}})_1\simeq\dbl{D}_1$
    over $\dbl{D}_0\times\dbl{D}_0$.
    \qedhere
    \end{enumerate}
\end{proposition}    
\begin{proof}
    Most of these are presented in \cite[Section 3.1]{hoshinoDoubleCategoriesRelations2023} except for (iii).
    The definition of the self-dual structure in \textit{loc.cit.} lacks
    the last condition in \cite[Definition 4.11]{Stay16}
    \footnote{
        This was pointed out by Zeinab Galal in a private communication with Keisuke Hoshino and the author.
        },
    which is the following axiom called ``swallowtail equation'':
    \[
    \begin{tikzcd}[column sep=large, row sep=small]
        &
        I\times I
        \sar[rd, "\iota_I\times\delta_I\times\delta_I", near start]
        \sar[rrd, bend left, "\delta_{I\times I}"]
        \ar[rr, phantom, "\Downarrow\ \theta_I\times\delta_I"]
        &&
        \!
        \\
        1
        \sar[ru, "\iota_I"]
        \sar[rd, "\iota_I"']
        \ar[rr, phantom, "\rotatebox{90}{$\cong$}"]
        &&
        I\times I\times I\times I
        \sar[r, "\delta_I\times\epsilon_I\times\delta_I"]
        &
        I\times I
        \\
        &
        I\times I
        \sar[ru, "\delta_I\times\delta_I\times\iota_I"', near start]
        \sar[rru, bend right, "\delta_{I\times I}"']
        \ar[rr, phantom, "\Downarrow\ \delta_I\times\zeta_I\inv"']
        &&
        \!
    \end{tikzcd}
    =
    \begin{tikzcd}
        1
        \sar[r, "\iota_I"]
        &
        I\times I
    \end{tikzcd}.
    \]
    This axiom is satisfied because all the 2-cells in the diagram
    are induced from the supine cells in the cartesian equipment.
    More precisely, since we have the supine cell whose bottom face is $\iota_I$,
    it suffices to show the equation composed with the supine cell,
    which is the following:
    (again, we omit the product symbols $\times$ and the commas ``,'' for the sake of readability)
    \[
    \begin{tikzcd}[column sep=small, row sep=small]
        &&&
        I
        \ar[lldd, equal]
        \ar[rd, "\tpl{00}"]
        \\
        &&&&
        II
        \ar[lldd, equal]
        \ar[rrdd, equal]
        \\
        &
        I
        \ar[dl, "!"']
        \ar[dr, "\tpl{00}"]
        \\
        1
        \sar[rr, "\iota_I"']
        \ar[rr, phantom, "\opcart"{yshift=0.75em}]
        &&
        II
        \sar[rr, "\iota_I\delta_I\delta_I"']
        \ar[rrrr, phantom, "\theta_I\delta_I"{yshift=2em}]
        &&
        I I I I
        \sar[rr, "\delta_I\epsilon_I\delta_I"']
        &&
        I I
    \end{tikzcd}
    =
    \begin{tikzcd}[column sep=small, row sep=small]
        &&&
        I
        \ar[lldd, equal]
        \ar[rd, "\tpl{00}"]
        \\
        &&&&
        I I
        \ar[lldd, equal]
        \ar[rrdd, equal]
        \\
        &
        I
        \ar[dl, "!"']
        \ar[dr, "\tpl{00}"]
        \\
        1
        \sar[rr, "\iota_I"']
        \ar[rr, phantom, "\opcart"{yshift=0.75em}]
        \sar[rrd, "\iota_I"']
        &&
        I I
        \sar[rr, "\delta_I\delta_I\iota_I"]
        \ar[rrrr, phantom, "\delta_I\zeta_I"{yshift=2em}]
        \ar[d, phantom, "\rotatebox{90}{$\cong$}"]
        &&
        I I I I
        \sar[rr, "\delta_I\epsilon_I\delta_I"']
        &&
        I I
        \\
        &&
        I I I I
        \sar[rru, "\iota_I\delta_I\delta_I"']
    \end{tikzcd}.
    \]
    Note that the square in the diagram above is a Beck-Chevalley pullback.
    The left-hand side of the equation can be computed as follows:
    \begin{align*}
    \text{(LHS)}
    &=
    \begin{tikzcd}[row sep=small, ampersand replacement=\&]
        \&\&\&
        I
        \ar[rd, "\tpl{00}"]
        \ar[ld, "\tpl{00}"']
        \\
        \&\&
        I I
        \ar[ld,"\tpl{1}"']
        \ar[rd, "\tpl{011}"]
        \&\&
        I I
        \ar[ld, "\tpl{001}"']
        \ar[rd, "\tpl{001}"]
        \\
        \&
        I
        \ar[dl, "!"']
        \ar[dr, "\tpl{00}"]
        \&\&
        I I I
        \ar[ld, "\tpl{12}"']
        \ar[rd, "\tpl{0012}"]
        \&\&
        I I I
        \ar[ld, "\tpl{0112}"']
        \ar[rd, "\tpl{02}"]
        \\
        1
        \sar[rr, "\iota_I"']
        \ar[rr, phantom, "\opcart"{yshift=0.75em}]
        \&\&
        I I
        \sar[rr, "\iota_I\delta_I\delta_I"']
        \ar[rr, phantom, "\opcart"{yshift=0.75em}]
        \&\&
        I I I I
        \sar[rr, "\delta_I\epsilon_I\delta_I"']
        \ar[rr, phantom, "\opcart"{yshift=0.75em}]
        \&\&
        I I
    \end{tikzcd}
    \\
    &=
    \begin{tikzcd}[row sep=small, ampersand replacement=\&]
        \&\&\&
        I
        \ar[rrdd, "\tpl{000}"]
        \ar[ld, "\tpl{00}"']
        \\
        \&\&
        I I
        \ar[lldd, "!"']
        \ar[rrdd, "\tpl{0011}"]
        \&\&
        \\
        \&
        \&
        \&
        \&
        \&
        I I I
        \ar[ld, "\tpl{0112}"']
        \ar[rd, "\tpl{02}"]
        \\
        1
        \sar[rr, "\iota_I"']
        \ar[rrrr, phantom, "\opcart"{yshift=2em}]
        \&\&
        I I
        \sar[rr, "\iota_I\delta_I\delta_I"']
        \&\&
        I I I I
        \sar[rr, "\delta_I\epsilon_I\delta_I"']
        \ar[rr, phantom, "\opcart"{yshift=0.75em}]
        \&\&
        I I
    \end{tikzcd}
    ,
        \\
    &=
    \begin{tikzcd}[row sep=small, ampersand replacement=\&]
        \&\&\&
        I
        \ar[rd, "\tpl{00}"]
        \ar[ld, "\tpl{00}"']
        \\
        \&\&
        I I
        \ar[ld,"\tpl{0}"']
        \ar[rd, "\tpl{001}"]
        \&\&
        I I
        \ar[ld, "\tpl{011}"']
        \ar[rd, "\tpl{011}"]
        \\
        \&
        I
        \ar[dl, "!"']
        \ar[dr, "\tpl{00}"]
        \&\&
        I I I
        \ar[ld, "\tpl{01}"']
        \ar[rd, "\tpl{0122}"]
        \&\&
        I I I
        \ar[ld, "\tpl{0112}"']
        \ar[rd, "\tpl{02}"]
        \\
        1
        \sar[rr, "\iota_I"']
        \sar[rrd, "\iota_I"']
        \ar[rrrr, phantom, "\rotatebox{90}{$\cong$}"{yshift=-1em}]
        \ar[rr, phantom, "\opcart"{yshift=0.75em}]
        \&\&
        I I
        \sar[rr, "\delta_I\delta_I\iota_I"']
        \ar[rr, phantom, "\opcart"{yshift=0.75em}]
        \&\&
        I I I I
        \sar[rr, "\delta_I\epsilon_I\delta_I"']
        \ar[rr, phantom, "\opcart"{yshift=0.75em}]
        \&\&
        I I
        \\
        \&
        \&
        I I
        \sar[rru, "\iota_I\delta_I\delta_I"']
    \end{tikzcd}
    = \text{(RHS)},
    \end{align*}
    where all pullback squares are Beck-Chevalley pullbacks
    by the trivial reasons \Cref{lem:BeckChevalleyclosure} and the Frobenius axiom
    so that we can use the sandwich lemma \Cref{lem:Sandwich} iteratively 
    to deduce that all triangles pointing upwards are supine cells.
\end{proof}

The paper \cite{hoshinoDoubleCategoriesRelations2023} does not explicitly mention
how the loose composition relates to the compact closed structure,
although the connection to the dagger structure is discussed in (iv) of \Cref{prop:Frobeniusselfdual}.

\begin{proposition}
    \label{prop:Frobeniusdaggercomposition}
    Suppose that $\dbl{D}$ is a Frobenius cartesian equipment,
    and let $\iota_I$ and $\epsilon_I$ be as in \Cref{prop:Frobeniusselfdual}.
    For loose arrows $\beta\colon I\times J\sto 1$ and $\beta'\colon J\times K\sto 1$,
    let $\alpha\colon I\sto J$ and $\alpha'\colon J\sto K$ be the corresponding loose arrows
    induced by $\iota,\epsilon$.
    Then, $\alpha\circ\alpha'$ is given by the following composite:
    \begin{equation}
    \label{eq:Frobeniusdaggercomposition}
    \begin{tikzcd}[column sep=huge]
        I
        \sar[r, "\delta_I\times\iota_K"]
        &
        I\times K\times K
        \sar[r, "(\delta_I\times\iota_J\times\delta_K)\times\delta_K"]
        &
        I\times J\times J\times K \times K
        \sar[r, "(\beta\times\beta')\times\delta_K"]
        &
        K
    \end{tikzcd}.
    \end{equation}
    This assignment is functorial with respect to the tightwise composition.
\end{proposition}
\begin{proof}
    The first statement follows from the general properties of
    compact closed bicategories.
    The second statement is a straightforward calculation,
    as in the proof of \cite[Proposition 3.1.15]{hoshinoDoubleCategoriesRelations2023}. 
\end{proof}

\begin{corollary}
    \label{cor:Frobeniusdaggercomposition}
    Suppose the same setting as in \Cref{prop:Frobeniusdaggercomposition}.
    When we define the category $(\LS{\dbl{D}})_1$ as in \Cref{prop:Frobeniusselfdual},
    define the composition functor $\odot'\colon(\LS{\dbl{D}})_1\times_{\dbl{D}_0}(\LS{\dbl{D}})_1\to(\LS{\dbl{D}})_1$ 
    by
    \[
        \left(
        \beta\colon I\times J\sto 1,
        \beta'\colon J\times K\sto 1
        \right)
        \mapsto
        \begin{tikzcd}[column sep=large]
            I\times K
            \sar[r, "\delta_I\times\iota_J\times\delta_K"]
            &
            I\times J\times J\times K
            \sar[r, "\beta\times\beta'"]
            &
            1
        \end{tikzcd}.
    \]
    Then, this data together with $\tpl{\src',\tgt'}$ gives rise to
    a double category $\LS{\dbl{D}}$:
    \[
    \begin{tikzcd}
        (\LS{\dbl{D}})_1\times_{\dbl{D}_0}(\LS{\dbl{D}})_1
        \ar[r, "\odot'"] 
        &
        (\LS{\dbl{D}})_1
        \ar[r, shift left =2 , "\src'"]
        \ar[r, shift right =2 , "\tgt'"']
        &
        \dbl{D}_0
        \ar[l, "\epsilon"{description}]
    \end{tikzcd}\text{.}
    \]
    The fibered equivalence in (vi) of \Cref{prop:Frobeniusselfdual} 
    is lifted to an equivalence $\LS{\dbl{D}}\simeq\dbl{D}$
    as double categories, \text{i.e.,} an equivalence in the 2-category $\Dbl$.
\end{corollary}

Therefore, with the Frobenius axiom, 
the compact closed structure behaves well with respect
to both the tight and loose compositions.

\begin{definition}
    \label{def:unilateral}
    Let $\dbl{D}$ be a cartesian equipment.
    We define a \emph{unilateral} fibration 
    to be a fibration $\uni[\dbl{D}]\colon\Uni[\dbl{D}]\to\dbl{D}_0$ 
    defined by the pullback
    \[
        \begin{tikzcd}
            \Uni[\dbl{D}]
            \ar[r]
            \ar[d, "{\uni[\dbl{D}]}"']
            \ar[rd, phantom, "\lrcorner", very near start]
            &
            \dbl{D}_1
            \ar[d, "\tpl{\src,\tgt}"]
            \\
            \dbl{D}_0
            \ar[r, "{I\mapsto(I,1)}"']
            &
            \dbl{D}_0\times\dbl{D}_0
        \end{tikzcd}.
    \]
\end{definition}
Although we have defined the unilateral fibration for arbitrary cartesian equipments, 
the resulting fibration loses the information of the original equipment
when the Frobenius axiom is not assumed.
For instance, from the equipment $\PROF$ of profunctors,
we obtain the fibration of the presheaves over the category of categories, 
which no longer remembers all profunctors,
in particular, copreshaves.
Note that, with the Frobenius axiom, the fibration $\uni[\dbl{D}]$ 
is equivalent to the pullback of $(\LS{\dbl{D}})_1$ over $\dbl{D}_0\times\dbl{D}_0$
along the diagonal functor by construction.
The extraordinary symmetric nature that the Frobenius axiom gives to the objects 
enables the equipment to be reconstructed from only one side of the loose arrows and cells.

\begin{proposition}
    \label{lemma:unilateral}
    Let $\dbl{D}$ be a Frobenius cartesian equipment.
    Then, the unilateral fibration $\uni[\dbl{D}]$ is a cartesian fibration,
    and its bilateral virtual double category $\Bil[\uni[\dbl{D}]]$
    is equivalent to $\LS{\dbl{D}}$, 
    hence to $\dbl{D}$.
    In particular, $\uni[\dbl{D}]$ is an elementary existential fibration. 
\end{proposition}

\begin{proof}
    The fibration $\tpl{\src,\tgt}$ is a cartesian fibration 
    since $\dbl{D}$ is a cartesian equipment.
    Since base change along finite-product-preserving functors
    preserves cartesian fibrations,
    the unilateral fibration $\uni[\dbl{D}]$ is a cartesian fibration.
    For the second statement, $\Bil[\uni[\dbl{D}]]$ has 
    $\dbl{D}_0$ as its tight part, 
    and a loose arrow $\alpha\colon I\sto J$ in $\Bil[\uni[\dbl{D}]]$
    is given by a loose arrow $\alpha\colon I\times J\sto 1$ in $\dbl{D}$
    by construction,
    which is a loose arrow of $I\sto J$ in $\LS{\dbl{D}}$.
    In the same way, the unary cells in $\Bil[\uni[\dbl{D}]]$ are in bijection with 
    those in $\LS{\dbl{D}}$.
    However, we do not know that $\Bil[\uni[\dbl{D}]]$ is a double category,
    requiring to check the correspondence of the general $n$-ary globular cells. 
    See \Cref{lemma:fibvdblequiv} for the condition for the equivalence of fibrational virtual double categories. 
    We will only show this in the case of $n=2$,
    and the general case is similar.
    \[
        \begin{tikzcd}[virtual]
            I
            \sar[r, "\alpha"]
            \ar[d, equal]
            \doublecell[rrd]{\tau}
            &
            J
            \sar[r, "\beta"]
            &
            K
            \ar[d, equal]
            \\
            I
            \sar[rr, "\gamma"']
            &&
            I
        \end{tikzcd}
    \]
    In $\Bil[\uni[\dbl{D}]]$, the above cell is given by a cell
    \[
        \begin{tikzcd}[virtual, column sep=8em]
            I\times J \times K
            \sar[r, "{\wh{\alpha}[\tpl{0,1}\fatsemi\id]}\land{\wh{\beta}[\tpl{1,2}\fatsemi\id]}"]
            \ar[d, "\tpl{0,2}"']
            \ar[dr, phantom, "\tau"]
            &
            1
            \ar[d, equal]
            \\
            I\times K
            \sar[r, "\wh{\gamma}"']
            &
            1
        \end{tikzcd}
    \]
    where $\wh{\alpha},\wh{\beta},\wh{\gamma}$ are the corresponding loose arrows in $\dbl{D}$
    by the compact closed structure.
    Here, the loose arrow at the top is isomorphic to the restriction of $\wh{\alpha}\times\wh{\beta}$
    along the tight arrow $\tpl{0,1,1,2}$.
    In other words, the cell above is equivalently given by the cell
    \[
        \begin{tikzcd}[virtual, row sep=tiny]
            &
            I\times J\times K
            \sar[r, "{\tpl{0,1,1,2}_*}"]
            &
            I\times J\times J\times K
            \sar[rd, "{\wh{\alpha}\times\wh{\beta}}"]
            \\
            I\times K
            \sar[rrr, "\wh{\gamma}"']
            \ar[rrr, phantom, "\tau"{yshift=1em}]
            \sar[ru, "{\tpl{0,2}^*}"]
            &&&
            1
        \end{tikzcd},
    \]
    but the composite on the top row is exactly the composite of $\wh{\alpha}$ 
    and $\wh{\beta}$ in $\LS{\dbl{D}}$.
    Therefore, the virtual cells in $\Bil[\uni[\dbl{D}]]$ are in bijection with those in $\LS{\dbl{D}}$.

    The last statement follows from \Cref{thm:firsteefiscartdouble}.
\end{proof} 

The equivalence $\Bil[\uni[\dbl{D}]]\simeq\dbl{D}$ in \Cref{lemma:unilateral}
is natural in $\dbl{D}$ in the following sense.
\begin{lemma}
    \label{lem:frobeqp-to-eef}
    The assignment of $\uni[\dbl{D}]$ to $\dbl{D}$ gives rise to a 2-functor 
    $\uni\colon\Eqp\Frob\to\Fib_{\lc*=+}$.
\end{lemma}
\begin{proof}
    We have the 2-functor $\Eqp\Frob\to\BiFib$
    that sends an equipment $\dbl{D}$ to the associated bifibration $\tpl{\src,\tgt}\colon\dbl{D}_1\to\dbl{D}_0\times\dbl{D}_0$,
    and since the base change preserves bifibrations,
    we have the 2-functor $\uni\colon\Eqp\to\BiFib$\footnote{
    This should be better explained in terms of 2-fibrations \cite{Her99}.
    }.
    On $\Eqp\Frob$, this factors through the 2-functor $\Fib_{\lc*=+}\to\BiFib$ by \Cref{lemma:unilateral}
    on the level of 0-cells,
    and also as a 2-functor because $\Fib_{\lc*=+}\to\BiFib$ is fully faithful by
    \Cref{lem:eeftobifib}. 
\end{proof}
\begin{proposition}
    \label{prop:uniisnaturalequivalence}
    The equivalence $\Bil[\uni[\dbl{D}]]\simeq\dbl{D}$ 
    is pseudo-natural in $\dbl{D}$.
\end{proposition}
\begin{proof}[Sketch of proof]
    Cartesian double functors preserve loose composition,
    finite-product structure, and (op)restrictions,
    and in particular, $\iota_I$ and $\epsilon_I$
    as in \Cref{prop:Frobeniusselfdual} are preserved.
    This ensures that the equivalence constructed above is 
    pseudo-natural because the data of the equivalence is determined 
    by those structures.
\end{proof}

On the other hand, we have the canonical isomorphism $\uni[\Bil[\mf{p}]]\cong\mf{p}$ for any elementary existential fibration $\mf{p}$
by the construction, and it is 2-natural in $\mf{p}$.
Combining the arguments above, we have the following corollary.

\begin{corollary}
    \label{cor:unilateral}
    The 2-functor $\Bil\colon\Fib_{\lc*=+}\to\Eqp\carttwo$ is locally an equivalence,
    and the essential image of $\Bil$ up to equivalence is $\Eqp\Frob$.
    The inverse 2-functor is given by $\uni\colon\Eqp\Frob\to\Fib_{\lc*=+}$.
\end{corollary}

It is worth mentioning that the Frobenius axiom is an instance of the Beck-Chevalley pullback condition, so the biequivalence restricts to give the following corollary. 
Let $\Eqp\BC$ be the full sub-2-category of $\Eqp\carttwo$ consisting of cartesian equipments with Beck-Chevalley pullbacks, and let $\Reg\Fib$ be the full sub-2-category of $\Fib_{\lc*=+}$ consisting of regular fibrations.

\begin{corollary}
    \label{cor:equivregbc}
    The 2-functor $\Bil\colon\Reg\Fib\to\Eqp\BC$, obtained by restricting the 2-functor $\Bil$ to the full sub-2-category $\Reg\Fib$, is a biequivalence,
    with the 2-functor $\uni$ as its inverse.
\end{corollary}
\begin{proof}
    \Cref{cor:BeckChevalleyall} ensures that the biequivalence
    restricts to the full sub-2-categories of regular fibrations and 
    cartesian equipments with Beck-Chevalley pullbacks.
\end{proof}

\begin{remark}
    A similar result to these corollaries is presented in \cite[\S 4.2.2]{Law15}.
    What is called \textit{regular fibrations} in \cite{Law15} lies between \textit{elementary existential fibrations} and \textit{regular fibrations} in our terminology, and they are assumed to satisfy the Beck-Chevalley condition for \textit{product-absolute} pullbacks only.
    Thus, this is another restriction of the equivalence in \Cref{cor:unilateral}.
\end{remark}

\section{Properties of fibrations and double categories and relationships to other structures}\label{sec:translating}
In this section, we discuss some properties of both fibrations and cartesian equipments, and we can translate between them under the correspondence in the previous section.
These properties reflect what kind of logical systems are available when we view fibrations and cartesian equipments as semantic fields.
Based on the correspondence, we will show that the properties that have been studied in the context of fibrations are translated to natural properties of cartesian equipments.
We will also give comparisons with other categorical structures, such as factorization systems (\Cref{subsec:factorization}) and bicategories (\Cref{subsec:cartesianbicategories}).
The relationship among the structures we have discussed in this paper is summarized in \Cref{fig:struc}.

\pgfdeclarelayer{background layer}
\pgfdeclarelayer{middle layer}
\pgfdeclarelayer{foreground layer}
\pgfsetlayers{background layer,middle layer,main,foreground layer}
\begin{figure}[htbp]
\centering
\begin{tikzpicture}[auto]
    \tikzstyle{every node}=[font=\tiny]
    \tikzstyle{rel}=[fill=white, align=center, rounded corners]
    \tikzstyle{span}=[fill=white, align=center, rounded corners]
    \tikzstyle{ele}=[dashed, line width=1pt, fill=white, align=left, rounded corners]
    \tikzstyle{exi}=[dashed, line width=1pt, fill=white, align=left, rounded corners]
    \tikzstyle{trans}=[->,shorten >= 3pt, shorten <= 3pt]
    \tikzstyle{char}=[draw,double, fill=red!10, align=center, rounded corners=3pt,auto=false]
    \tikzstyle{ind}=[draw, fill=red!10, align=center, rounded corners=3pt,auto=false]
    \tikzstyle{ope}=[draw, fill=white, align=center, rounded corners=5pt,auto=false] 

    \begin{pgfonlayer}{background layer}
    \fill[gray!40] (-1,0.9) rectangle (16,-1.9);
    \end{pgfonlayer}
    \begin{pgfonlayer}{foreground layer}
    \node[anchor=west] (B0) at (-1,0.5) {\small\hyperref[subsec:factorization]{\textbf{Factorization systems}}\ };
    \end{pgfonlayer}
    \node[draw, rel] (F0) at (1,-0.25) {$(\zero{RegEpi},\zero{Mono})$\\ on a regular category $\one{B}$};
    \node[draw, span] (F1) at (3,-1.25) {$(\zero{Iso},\zero{All})$\\ on a lex category $\one{C}$};
    \node[draw, fill=white, align=left] (F2) at (6,-0.75) {Stable orthogonal\\ factorization system}; 
    \node[draw, fill=white, align=left] (F3) at (11,-0.75) {Orthogonal\\ factorization system};
    \draw[-] (F0) -- (F2);
    \draw[-] (node cs:name=F1, anchor=east) -- (node cs:name=F2, anchor=base west);
    \draw[-] (F2) -- (F3);

    \begin{pgfonlayer}{background layer}
        \fill[gray!20] (-1,-2.1) rectangle (16,-5.9);
    \end{pgfonlayer}
    \begin{pgfonlayer}{foreground layer}
    \node[anchor=west, fill=gray!20] (B1) at (-1,-2.5) {\small\hyperref[sec:fibration]{\textbf{Fibrations}}\ };
    \end{pgfonlayer}
    \node[draw, rel] (G0) at (1,-3.5) {Subobject fibration\\ $\Pred[\one{B}]\rightarrow\one{B}$};
    \node[draw, span] (G1) at (3,-5) {Codomain fibration\\ $\one{C}^{\rightarrow}\rightarrow \one{C}$}; 
    \node[draw, fill=white, align=left] (G2) at (6,-4.25) {Regular\\ fibration};
    \node[draw, fill=white, align=left] (G3) at (11,-2.5) {Bifibration};
    \node[draw, fill=white, align=left] (G4) at (9,-4.25) {Elementary\\ existential\\ fibration};
    \node[draw, fill=white, align=left] (G5) at (12,-3.25) {Elementary \\fibration}; 
    \node[draw, fill=white, align=left] (G6) at (11,-5.25) {Existential\\ fibration};
    \node[draw, fill=white, align=left] (G7) at (14,-4.25) {Cartesian\\ fibration};
    \draw[-] (G0) -- (G2);
    \draw[-] (node cs:name=G1, anchor=east) -- (node cs:name=G2, anchor=base west);
    \draw[-] (G2) -- (G4);
    \draw[-, rounded corners=10pt] (G4) -- (9.3,-2.5) -- (G3);
    \draw[-] (G4) -- (G5);
    \draw[-] (G4) -- (G6);
    \draw[-] (G5) -- (G7);
    \draw[-] (G6) -- (G7);

    \begin{pgfonlayer}{middle layer}
    \draw[trans,|-|] (node cs:name=F0) -- (node cs:name=G0);
    \draw[trans,|-|] (F1) -- (G1);
    \draw[trans] (F2) -- node [char] {\cite{HJ03}\\(\cref{lem:stablefibration}\\\cref{thm:factorizationsystem})} (G2);
    \draw[trans] (F3) -- node [char] {\cite{HJ03} (\cref{def:predfibration},\cref{thm:factorizationsystem})} (G3);
    \end{pgfonlayer}

    \begin{pgfonlayer}{background layer}
        \fill[gray!40] (-1,-6.1) rectangle (16,-10.4);
    \end{pgfonlayer}
    \begin{pgfonlayer}{foreground layer}
        \node[anchor=west, fill=gray!40] (B2) at (-1,-6.5) {\small\hyperref[def:fibrational]{\textbf{Fibrational virtual double categories}}}; 
    \end{pgfonlayer}
    \node[draw, rel] (H0) at (1,-7.75) {Double category of\\ relations $\Rel[\one{B}]$};
    \node[draw, span] (H1) at (3,-9.25) {Double category of\\ spans $\Span(\one{C})$};
    \node[draw, fill=white, align=left] (H2) at (6,-8) {Cartesian equipment\\ with BC pullbacks};
    \node[draw, fill=white, align=left] (H4) at (9,-8) {Cartesian\\ equipment};
    \node[draw, fill=white, align=left] (H5) at (12,-7) {Cartesian \\ virtual equipment}; 
    \node[draw, fill=white, align=left] (H6) at (11,-9.25) {Cartesian positively-representable\\ fibrational \acs{VDC}}; 
    \node[draw, fill=white, align=left] (H7) at (14,-8) {Cartesian \\ fibrational \acs{VDC}}; 
    \node[draw, fill=white, align=left] (H8) at (7, -9) {Unit-pure map-complete\\ cartesian equipment};
    \draw[-] (H0) -- (H2);
    \draw[-] (node cs:name=H1, angle=60) -- (node cs:name=H2, anchor=base west);
    \draw[-] (H2) -- (H4);
    \draw[-] (H4) -- (H5);
    \draw[-] (H4) -- (H6);
    \draw[-] (H5) -- (H7);
    \draw[-] (H6) -- (H7);
    \draw[-, rounded corners=2pt] (H8) -- (8, -8) -- (H4);
    \draw[-] (H0) -- (H8);
    \draw[-] (node cs:name=H1, anchor=east) -- (H8);

    \begin{pgfonlayer}{middle layer}
    \draw[trans,|-|] (node cs:name=G0) -- (node cs:name=H0);
    \draw[trans,|-|] (G1) -- (H1);
    \draw[trans,<->] (G2) -- node [char, near start] {\cref{cor:BeckChevalleyall}\\\cref{cor:equivregbc}} (H2);
    \draw[trans,->] (node cs:name=F2, angle=340) -- node [char, near start,right] {$\Rel[-][-]$\\\cite{hoshinoDoubleCategoriesRelations2023}} (node cs:name=H2, angle=20);
    \draw[trans] (G4) -- node [char] {\cref{thm:firsteefiscartdouble}\\\cref{cor:unilateral}} (H4);
    \draw[trans] (G5) -- node [ind,near end] {\cref{prop:cartesianunital}} (H5);
    \draw[trans] (node cs:name=G6, angle=220) -- node [ind, near end,right] {\cref{prop:cartesiancompose}} (node cs:name=H6, angle=140);
    \draw[trans] (G7) -- node [ind] {\cref{def:primfibtovdc}\\\cref{prop:primfibtovdc}} (H7); 
    \end{pgfonlayer}

    \begin{pgfonlayer}{background layer}
        \fill[gray!20] (-1,-10.6) rectangle (11.5,-14.9);
    \end{pgfonlayer}
    \begin{pgfonlayer}{foreground layer}
    \node[anchor=west,fill=gray!20] (B3) at (-1,-11) {\small\hyperref[subsec:cartesianbicategories]{\textbf{Bicategories}}\ };
    \end{pgfonlayer}
    \node[draw, rel] (I0) at (1,-12.5) {Bicategory of\\ relations $\BiRel(\one{B})$};
    \node[draw, span] (I1) at (3,-11.25) {Bicategory of\\ spans $\BiSpan(\one{C})$};
    \node[draw, fill=white, align=left] (I2) at (7,-12.5) {Map-discrete\\ cartesian \\bicategory}; 
    \node[draw, fill=white, align=left] (I3) at (9.5,-12.5) {Cartesian \\ bicategory \\ \cite{CW87,CKWW07}};
    \node[draw, fill=white, align=left] (I4) at (4,-12.5) {Frobenius \\locally posetal\\ cartesian bicategory};
    \node[draw, fill=white, align=left, rel] (I5) at (1,-14.25) {Unitary tabular \\ allegory};
    \node[draw, fill=white, align=left] (I6) at (4,-14.25) {Unitary pretabular \\ allegory};
    \node[draw, fill=white, align=left] (I7) at (8, -14.25) {Allegory \\ \cite{FS90,johnstoneSketchesElephantTopos2002a}};

    \draw[-] (I0) -- (I4);
    \draw[-, rounded corners=2pt] (I1) -- (5.75, -11.25) -- (5.75,-12.5) -- (I2);
    \draw[-] (I2) -- (I4);
    \draw[-] (I2) -- (I3);
    \draw[-] (I5) -- (I6);
    \draw[-] (I6) -- (I7);
    \draw[<->] (I0) -- node [char] {\cite{FS90}} (I5);
    \draw[<->] (I4) -- node [char] {\cite[\S 2.2]{Law15}} (I6);

    \begin{pgfonlayer}{middle layer}
    \draw[trans,|-|] (H0) -- (I0);
    \draw[trans,|-|] (H1) -- (I1);
    \draw[trans,<->] (H8) -- node [char] {\cref{thm:eefiscartdouble}} (I2);
    \draw[trans] (node cs:name=H4, angle=240) -- node [ind] {\cite{patterson2024transposingcartesianstructuredouble}\\\cref{prop:CartdoubleCartbicat}} (node cs:name=I3, angle=140); 
    \end{pgfonlayer}

    \draw [->,line width=1.5pt] (14.5,-1.5) to[edge node={node[pos=0.5,ope]{\small$\Pred[-][-]$}}, bend left=20](14.5,-3);
    \draw [->,line width=1.5pt] (15.2,-5.25) to[edge node={node[pos=0.5,ope]{\small$\Bil$}}, bend left=20](15.2,-6.75); 
    \draw [->,line width=1.5pt] (10.5,-10.25) to[edge node={node[pos=0.5,ope]{\small$\LBi{-}$}}, bend left=20](10.5,-11.75); 

    \begin{pgfonlayer}{foreground layer}
    \draw [->] (12.5, 0.25-13.5) -- (13, 0.25-13.5);
    \draw [<->] (12.5, 0.75-13.5) -- (13, 0.75-13.5);
    \draw [|-|] (12.5, 1.25-13.5) -- (13, 1.25-13.5);
    \draw [-] (12.5, 1.75-13.5) -- (13, 1.75-13.5);
    \node [anchor=west] at (13, 0.25-13.5) {: Construction};
    \node [anchor=west] at (13, 0.75-13.5) {: Equivalence};
    \node [anchor=west] at (13, 1.25-13.5) {:\begin{tabular}{l}Corresponding\\ \ examples\end{tabular}}; 
    \node [anchor=west] at (13, 1.75-13.5) {: Implication};

    \node at (13,2.25-13.5) {\small \textbf{Legend}};

    \node (MG) at (14,0.5) {\small\textbf{General Classes}};
    \node (MS) at (8,0.5) {\small\textbf{Specific Classes}};
    \draw [<->,line width=1pt] (MS) -- (MG);
\end{pgfonlayer}
    \begin{pgfonlayer}{background layer}
        \draw[line width=1pt, fill=white] (12,2.5-13.5) rectangle (16,-0.2-13.5);
    \end{pgfonlayer}
\end{tikzpicture}
\caption{Relationships between fibrational virtual double categories and other structures that have been studied}
\label{fig:struc}
\end{figure}

\subsection{The predicate comprehension and tabulators}
\label{subsec:comprehension}

In set theory, the comprehension axiom states that for any set $I$ and
any predicate $\alpha(x)$, there exists a subset $\{\alpha\}$ of $I$ such that
\[
    \forall x : I.\  \alpha(x) \iff x \in \{\alpha\}
\]
The following definition is a categorical reformulation of this axiom.
\begin{definition}
    \label{def:predicatecomprehension}
    Let $\mf{p}\colon\one{E}\to\one{B}$ be a fibration with fiberwise terminal objects.
    For an object $\alpha$ in $\one{E}_I$, a \textbf{predicate comprehension} of $\alpha$ is
    a terminal object in the comma category $\top_{-}\downarrow\alpha$,
    where $\top_{-}$ is the functor that sends an object $I$ to the terminal object in $\one{E}_I$. 
    We say that $\mf{p}$ has \textbf{predicate comprehension} if every object in $\one{E}$ has a predicate comprehension, 
    or equivalently, if the functor $\top_{-}$ has a right adjoint.
\end{definition}

The statement that $\alpha$ has a predicate comprehension is unwound as follows.
Let $\tau_\alpha\colon\top_{\{\alpha\}}\to\alpha$ be a terminal object in $\top_{-}\downarrow\alpha$,
and let $c_\alpha\colon\{\alpha\}\to I$ be its image under $\mf{p}$.
Then, for any object $J$ in $\one{B}$ with an arrow $\varphi\colon\top_{J}\to\alpha$, 
there exists a unique arrow $u$ such that $\tau_\alpha\circ \top_{u} = \varphi$.
Using the internal language of fibrations, the data of $\tau_\alpha$ and $c_\alpha$ can be encoded as a term $\syn{x}:\cmpr{\syn{\alpha}}\vdash\syn{c}(\syn{x}):\syn{I}$ and a proof of $\syn{x}:\cmpr{\syn{\alpha}}\mid\top\vdash\syn{\alpha}(\syn{c}(\syn{x}))$, respectively.
The statement above is then translated into the following.
For any tuple of a type $\syn{J}$, a term $\syn{y}:\syn{J}\vdash\syn{v}(\syn{y}):\syn{I}$,
and a proof of $\syn{y}:\syn{J}\mid\top\vdash\syn{\alpha}(\syn{v}(\syn{y}))$,
there exists a unique term $\syn{y}:\syn{J}\vdash\syn{u}(\syn{y}):\cmpr{\syn{\alpha}}$ such that
$\syn{v}(\syn{y}) = \syn{c}(\syn{u}(\syn{y}))$ and 
the proof of $\syn{\alpha}(\syn{c}(\syn{x}))$ for $\syn{x}$ replaced by $\syn{u}(\syn{y})$
is the same as the given proof.
If one takes $\syn{J}$ to be the terminal type $\syn{1}$ and goes down to proof irrelevance, 
the above statement looks quite similar to the comprehension axiom in set theory. 

The predicate comprehension is called the \textit{subset (type)}
in \cite[Definition 4.6.1]{jacobsCategoricalLogicType1999a},
and just \textit{comprehension} in \cite[\S 4]{MREQC13}.
It should be noted that this notion is a special case of what 
are called \textit{comprehension structures} in several contexts,
such as in \cite{PN20},
where a comprehension structure with section is defined as a section of the fibration $\mf{p}$, not necessarily terminal, with a right adjoint.

\begin{definition}
    \label{def:fullcomprehension}
    Let $\mf{p}\colon\one{E}\to\one{B}$ be a fibration 
    with fiberwise terminal objects and
    predicate comprehension.
    We say that $\mf{p}$ has \textbf{full predicate comprehension} if
    the functor $c_{-}\colon\one{E}\to\one{B}\arrcat$ is fully faithful.
    Since this is shown to be a fibered functor,
    we can equivalently say that the functor $c_{-}\colon\one{E}_I\to
    \one{B}/I$ is fully faithful for every object $I$ in $\one{B}$.
    We write $\Predarr(\mf{p})$ for the class of arrows in $\one{B}$ 
    that arise, up to isomorphism,
    as $c_{\alpha}\colon\{\alpha\}\to I$ for some $\alpha\in\one{E}_I$.
\end{definition}
This condition is required to ensure that the behaviors of predicates are completely determined by their comprehension.
Note that, when $\mf{p}$ is fiberwise preordered, the faithfulness of $\{-\}$ is automatically satisfied\footnote{
    This seems to be why only the term \textit{full} is used in the definition, not \textit{fully faithful}.
}
since the counit components are necessarily epimorphisms.

In the view that a double category is a generalization of an elementary existential fibration,
a double categorical analog of predicate comprehension is the notion of tabulators.
\begin{definition}[Tabulators \cite{grandisLimitsDoubleCategories1999}]
    \label{def:tabulators}
    A \emph{(1-dimensional) tabulator} of a loose arrow $\alpha\colon I\sto J$ is an object $\{\alpha\}$ equipped
    with a pair of tight arrows $\ell_{\alpha}\colon \{\alpha\}\to I$ and $r_{\alpha}\colon \{\alpha\}\to J$
    and a cell
    \[
        \begin{tikzcd}[virtual, column sep=small]
            &
            \{\alpha\}
            \ar[dl, "\ell_{\alpha}"']
            \ar[dr, "r_{\alpha}"]
            &
            \\
            I
            \sar[rr, "\alpha"']
            \ar[rr, phantom, "\tau_\alpha", yshift=2ex]
            &&
            J
        \end{tikzcd}
    \]
    such that,
    for any cell $\nu$ on the left below, there exists a unique tight arrow $t_{\nu}\colon X\to \{\alpha\}$ 
    that makes the following two cells equal.
    \[
        \begin{tikzcd}[virtual, column sep=small]
            &
            X
            \ar[dl, "h"']
            \ar[dr, "k"]
            &
            \\
            I
            \sar[rr, "\alpha"']
            \ar[rr, phantom, "\nu", yshift=2ex]
            &&
            J
        \end{tikzcd}
        =
        \begin{tikzcd}[virtual, column sep=small]
            &
            X
            \ar[ddl, "h"', bend right=30]
            \ar[ddl, phantom, "\circlearrowleft"]
            \ar[ddr, "k", bend left=30]
            \ar[ddr, phantom, "\circlearrowleft"]
            \ar[d, "t_{\nu}"]
            \\
            &
            \{\alpha\}
            \ar[dl, "\ell_{\alpha}"']
            \ar[dr, "r_{\alpha}"]
            &
            \\
            I
            \sar[rr, "\alpha"']
            \ar[rr, phantom, "\tau_\alpha", yshift=2ex]
            &&
            J
        \end{tikzcd}
    \]
    Henceforth, we call the cell $\nu$ the \emph{tabulating cell} of $\alpha$. 
    In other words, a tabulator is a terminal object in the comma category 
    $\delta_{-}\downarrow\alpha$,
    where $\delta_{-}\colon\dbl{D}_0\to\dbl{D}_1$ is the functor that sends an object $I$ 
    to the loose arrow $\delta_I$,
    and $\alpha$ is seen as an object in $\dbl{D}_1$.
    A tabulator is called \emph{effective}\footnote{
    In \cite{hoshinoDoubleCategoriesRelations2023}, the authors use the term \textit{strong tabulator} 
    to mean the same thing according to the definition in \cite{aleiferiCartesianDoubleCategories2018}.
    We adopt the term \textit{effective} because it is a fixed point of 
    the adjunction between the category of spans between $I$ and $J$ 
    and the category of loose arrows from $I$ to $J$,
    as an effective epimorphism from $I$ is a fixed point of the adjunction between the category of
    parallel pairs into $I$ and the category of arrows from $I$.
    }
    if the tabulating cell is supine.

    In a double category $\dbl{D}$ with a terminal object $1$ in the tight category, 
    tight arrow $f\colon I\to J$ is called a \emph{fibration} 
    if there exists a loose arrow $\alpha\colon I\sto 1$ and a tabulating cell as follows.
    \[
        \begin{tikzcd}[virtual, column sep=small]
            &
            I
            \ar[dl, "f"']
            \ar[dr, "!"]
            &
            \\
            J
            \sar[rr, "\alpha"']
            \ar[rr, phantom, "\tau_f", yshift=2ex]
            &&
            1
        \end{tikzcd}
    \]
    We write $\Fibarr(\dbl{D})$ for the class of fibrations in $\dbl{D}$.
\end{definition}

\begin{proposition}
    \label{prop:tabulatorspredicatecomp}
    Let $\mf{p}\colon\one{E}\to\one{B}$ be an elementary existential fibration.
    Then, the following are equivalent:
    \begin{enumerate}
        \item $\mf{p}$ has (full) predicate comprehension,
        \item $\Bil[\mf{p}]$ has (effective) tabulators, and
        \item $\Bil[\mf{p}]$ has left-sided (effective) tabulators.
    \end{enumerate}
    Here, the additional conditions of fullness and effectiveness
    are satisfied simultaneously.
    Furthermore, $\Predarr(\mf{p})$ coincides with $\Fibarr(\Bil[\mf{p}])$.
\end{proposition}
\begin{proof}
    In the double category $\Bil[\mf{p}]$,
    the comma category $\delta_{-}\downarrow\alpha$ for 
    a loose arrow $\alpha\colon I\sto J$ 
    is equivalent to the comma category $\top_{-}\downarrow\alpha$ 
    where $\alpha$ is seen as an object in $\one{E}_{I\times J}\subseteq\one{E}$.
    Since the first condition is equivalent to the existence of a terminal object
    in $\top_{-}\downarrow\alpha$ for every $I$ and $\alpha\in\one{E}_I$,
    and the third condition is the statement for the cases where $J$ is terminal, 
    we obtain the equivalence of the three conditions.
    The last statement is immediate when one observes how
    the tabulator and the predicate comprehension are 
    related by the above equivalence.

    In an equipment $\Bil[\mf{p}]$, the operation of taking tabulators
    gives the right adjoints of the functors that sends each span to
    its oprestriction: 
    \[
        \begin{tikzcd}[column sep=large, row sep=large]
            \Span(\one{B})(I,J)
            \ar[r, shift left=2, "{(f,g)}\mapsto f^*g_*"]
            \ar[r, phantom, "\rotatebox{90}{$\vdash$}"]
            \ar[d,"\cong"']
            &
            \Bil[\mf{p}](I,J)
            \ar[l, shift left=2, "\tpl{\ell_{\alpha},r_{\alpha}}\mapsfrom\,\alpha"]
            \ar[d,"\cong"]
            \\
            \one{B}/I\times J
            \ar[r, shift left=2, "\tpl{f,g}\mapsto \sum_{\tpl{f,g}}\top"]
            \ar[r, phantom, "\rotatebox{90}{$\vdash$}"]
            &
            \one{E}_{I\times J}
            \ar[l, shift left=2, "c_{\alpha}\mapsfrom\,\alpha"]
        \end{tikzcd}.
    \]
    The effectiveness of tabulators is equivalent to the counit components of the adjunction 
    being isomorphisms, 
    while the fullness of predicate comprehension is equivalent to 
    the right adjoint $c_{-}$ being fully faithful for every $I$ and $J$.
    Therefore, these conditions are satisfied simultaneously.
\end{proof}

\subsection{Double categories of relations relative to factorization systems}
\label{subsec:factorization}
In the paper \cite{HJ03}, the authors show the characterization of factorization systems as special bifibrations using the notion of predicate comprehension (\Cref{thm:factorizationsystem}).
Meanwhile, \cite{hoshinoDoubleCategoriesRelations2023} provides a characterization of cartesian equipments of relations relative to stable factorization systems using strong tabulators (\Cref{cor:DCR}).
As we have established the correspondence between predicate comprehension and tabulators, we can connect these two results, and recover the latter from the former.

\begin{definition}[{\cite[Definition 2.12]{HJ03}}]
    \label{def:strongproductsalongpredicatecomprehensionprojection}
    Let $\mf{p}\colon\one{E}\to\one{B}$ be a bifibration with fiberwise terminal objects and full predicate comprehension. 
    We say that $\mf{p}$ has \textbf{strong products along predicate projections} if for any object $I$ in $\one{B}$,
    any object $\alpha$ in $\one{E}_I$, and any object $\beta$ in $\one{E}_{\{\alpha\}}$, 
    the supine lift of $c_\alpha\colon\{\alpha\}\to I$ to $\beta$ 
    induces an isomorphism $\{\beta\}\to\{\sum_{c_\alpha}\beta\}$.
    \[
    \begin{tikzcd}
        &
        \beta
            \ar[r, dashed, "\text{supine lift}"]           
        &
        \sum_{c_\alpha}\beta
            \ar[dd, mapsto, shorten >= 0.5em, shorten <= 0.5em, bend left=20] 
        \\
        \{\beta\} 
            \ar[r] 
            \ar[rd, "c_{\beta}"']
            &
        \{\sum_{c_\alpha}\beta\} 
            \ar[rd, "c_{\sum_{c_\alpha}\beta}"]
        \\
        &
        \{\alpha\}
            \ar[r, "c_\alpha"']
            \ar[from=uu, mapsto, shorten >= 0.5em, shorten <= 0.5em, bend left=50, crossing over]
            &
        I
    \end{tikzcd}
    \]
\end{definition}

\begin{lemma}
    \label{lem:predicational}
    Let $\mf{p}\colon\one{E}\to\one{B}$ be a bifibration with fiberwise terminal objects and full predicate comprehension, and suppose that $\one{B}$ has finite limits.
    Then, the following are equivalent:
    \begin{enumerate}
        \item $\mf{p}$ has strong products along subset projections, and 
        \item $\Predarr(\mf{p})$ is closed under composition.
        \qedhere
    \end{enumerate}
\end{lemma}
\begin{proof}
    The implication (i)$\Rightarrow$(ii) is immediate from the definition of strong products along subset projections.
    Suppose we have a triple $\alpha,\gamma\in\one{E}_I$ and $\beta\in\one{E}_{\{\alpha\}}$
    such that the diagram on the left below commutes and the top arrow is an isomorphism. 
    \[
    \one{B}\arrcat
    \ni\quad
    \begin{tikzcd}
        \{\beta\}
            \ar[r, "\cong"]
            \ar[d, "c_{\beta}"']
        &
        \{\gamma\}
            \ar[d, "c_{\gamma}"]
        \\
        \{\alpha\}
            \ar[r, "c_{\alpha}"']
        &
        I
    \end{tikzcd}
    \quad
    \overset{c_{-}}{\longmapsfrom}
    \quad
    \begin{tikzcd}
        \beta
            \ar[r, dashed, "\exists!\ \mu"]
        &
        \gamma
    \end{tikzcd}
    \quad\in\one{E}
    \]
    By the fullness of predicate comprehension, 
    this is the image of an arrow $\mu\colon\beta\to\gamma$ in $\one{E}$
    under the functor $c_{-}$.
    In the codomain fibration $\one{B}\arrcat\to\one{B}$,
    the arrow $c_{\mu}$ is supine since the top arrow is an isomorphism.
    Since a fully faithful fibered functor reflects supine arrows,
    the arrow $\mu$ is supine.
    Thus, (ii)$\Rightarrow$(i) holds.
\end{proof}

\begin{definition}
	\label{def:predfibration}
    Let $\one{B}$ be a category and $(\zero{E},\zero{M})$ be an orthogonal factorization system on $\one{B}$.
    Suppose that $\one{B}$ admits pullbacks of arrows in $\zero{M}$ along arbitrary arrows.
    The fibration $\Pred[\one{B}][\zero{M}]\to\one{B}$ is defined as follows.
    \begin{itemize}
    \item The total category $\Pred[\one{B}][\zero{M}]$ is the full subcategory of $\one{B}^{\rightarrow}$ spanned by the arrows in $\zero{M}$.
    \item The functor $\Pred[\one{B}][\zero{M}]\to\one{B}$ is the functor that sends an arrow to its codomain.
    \end{itemize}
	By abuse of notation, we denote the fibration $\Pred[\one{B}][\zero{M}]\to\one{B}$ by $\Pred[\one{B}][\zero{M}]$.
\end{definition}

\begin{lemma}[{\cite[Lemma 2.8]{HJ03}}]
    \label{lem:stablefibration}
    Let $\one{B}$ be a category and $(\zero{E},\zero{M})$ be an orthogonal factorization system on $\one{B}$.
    If $\one{B}$ admits finite limits, then the fibration $\Pred[\one{B}][\zero{M}]\to\one{B}$ is a cartesian fibration.
    If is is also stable, meaning that the class $\zero{E}$ is closed under pullbacks, then the fibration $\Pred[\one{B}][\zero{M}]\to\one{B}$ is a regular fibration.
\end{lemma}

\begin{theorem}[{\cite[\S 3]{HJ03}}]
    \label{thm:factorizationsystem}
    Let $\mf{p}\colon\one{E}\to\one{B}$ be a bifibration with
    fiberwise terminal objects,
    and suppose that $\one{B}$ has all pullbacks.
    Then, the following are equivalent:
    \begin{enumerate}
        \item $\mf{p}$ has full predicate comprehension
        and strong products along subset projections.
        \item $\mf{p}$ is equivalent to $\Pred[\one{B}][\zero{M}]$
        for some factorization system
        $(\zero{E},\zero{M})$ on $\one{B}$.
    \end{enumerate}
    The factorization system $(\zero{E},\zero{M})$ is 
    uniquely determined by $\mf{p}$.
    The bifibration is regular if and only if the factorization system is stable.
\end{theorem}

\begin{corollary}[{\cite[Theorem 3.3.20]{hoshinoDoubleCategoriesRelations2023}}]
    \label{cor:DCR}
    The following are equivalent for a double category $\dbl{D}$:
    \begin{enumerate}
        \item $\dbl{D}$ is equivalent to $\Rel[\one{B}][(\zero{E},\zero{M})]$ for some category 
        with finite limits $\one{B}$ and a stable factorization system $(\zero{E},\zero{M})$ on $\one{B}$.
        \item $\Fibarr(\dbl{D})$ is closed under composition,
        and $\dbl{D}$ is a cartesian equipment with Beck-Chevalley pullbacks and effective tabulators.
        \item $\Fibarr(\dbl{D})$ is closed under composition,
        and $\dbl{D}$ is a cartesian equipment with Beck-Chevalley pullbacks and left-sided effective tabulators.
        \qedhere
    \end{enumerate}
\end{corollary}
\begin{proof}[An alternative proof of \Cref{cor:DCR} using \Cref{thm:factorizationsystem}]
    By \Cref{cor:equivregbc},
    the three conditions are subsumed by the condition that $\dbl{D}$ is
    of the form $\Bil[\mf{p}]$ for some regular fibration $\mf{p}$.
    Then, the equivalence follows from \Cref{prop:tabulatorspredicatecomp,lem:predicational,thm:factorizationsystem}.
\end{proof}

\subsection{Function extensionality and unit-pureness}
\label{subsec:funcext}
    
Function extensionality is a principle that states that two functions are equal if they are equal at every point.
\[
\forall f, g : I \to J. \  \left( \forall x : I. \  f(x) = g(x) \right) \Longrightarrow f = g.
\]
In the context of doctrines or fibrations, this principle is formulated
by interpreting the equality of functions as the equality in the base category
and the equality of elements as the predicate expressed in the fiber category.
An elementary preordered fibration with this property,
meaning that the inequality $\top_I\le \delta_J[\tpl{f,g}]$ implies $f = g$
for every $f,g\colon I\to J$ in the base category,
is said to have \textit{very strong equality}
in \cite[Section 3.4]{jacobsCategoricalLogicType1999a}.
This property is equivalent to the property that 
the diagonal arrow $\tpl{0,0}\colon I\to I\times I$ 
is a predicate comprehension of $\delta_I$ for every object $I$ in the base category
(see \cite[Proposition 2.12]{MPR17}).
From this observation, an elementary doctrine with this property is said
to have \textit{comprehensive diagonals} in \cite{MREQC13}, 
and the combination of this property with the property of full predicate comprehension
is called \textit{m-variational} in \cite[Definition 2.16]{MPR17}.
In the paper \cite{DP23} introducing relational doctrines,
the authors use the term \textit{extensional} for the corresponding property
to the very strong equality.
We adopt the term \emph{comprehensive diagonals} in this paper.

Note that in the context of allegories, extensionality is a trivial notion since functions are defined as maps there and pointwise equality leads to the equality of the maps themselves
(cf. \cite[Proposition A 3.2.3]{johnstoneSketchesElephantTopos2002a}).

In the context of double categories, the corresponding property should be
\textit{unit-pureness}.
\begin{definition}[{\cite[Definition 4.3.7]{aleiferiCartesianDoubleCategories2018}}]
    A double category (or a unital virtual double category) is called
    \textbf{unit-pure} if a cell of the form
    \[
        \begin{tikzcd}[virtual]
            I
            \sar[r, "\delta_I"]
            \ar[d, "f"']
            \ar[dr, "\mu",phantom]
            &
            I
            \ar[d, "g"]
            \\
            J
            \sar[r, "\delta_J"']
            &
            J
        \end{tikzcd}
    \]
    is necessarily the identity cell $\delta_f$ with $f = g$.
    In other words, the functor $\delta_{-}\colon\dbl{D}_0\to\dbl{D}_1$ is fully faithful, or the tight 2-category $\TBi{\dbl{D}}$ is locally discrete.
\end{definition}
Similarly to the case of fibrations,
the unit-pureness of a double category is equivalent to the property 
that the span $(\id_I,\id_I)$ exhibits $I$ as a tabulator of $\delta_I$ for every object $I$. 

The following is an easy observation.
\begin{proposition}
    \label{prop:unitpureness}
    Let $\mf{p}\colon\one{E}\to\one{B}$ be an elementary existential fibration.
    Then, the following are equivalent:
    \begin{enumerate}
        \item $\Bil[\mf{p}]$ is unit-pure.
        \item For each parallel pair $f,g\colon I\to J$ in $\one{B}$ and each arrow $\mu\colon \top_I\to\delta_J[\tpl{f,g}]$ in $\one{E}_I$, it holds that $f = g$ and $\mu$ is the composite of the following:
        \[
            \begin{tikzcd}
            \top_I
                \ar[r, "\cong"]
            &
            \top_J[\tpl{f}]
                \ar[r, "{\eta_J[f]}"]
            &
            \delta_J[\tpl{0,0}][f]
                \ar[r, "\cong"]
            &
            \delta_J[\tpl{f,f}]
            \end{tikzcd}
        \]
    \end{enumerate}
    In particular, when $\mf{p}$ is fiberwise preordered, the unit-pureness of $\Bil[\mf{p}]$ is equivalent to very strong equality of $\mf{p}$.
\end{proposition}

\begin{example}
    \label{ex:unit-pure}
    \begin{enumerate}
    \item For the double category $\Rel[\one{B}][(\zero{E},\zero{M})]$ for a stable factorization system $(\zero{E},\zero{M})$ on $\one{B}$, the unit-pureness is equivalent to the property that the left class $\zero{E}$ is included in the class of epimorphisms in $\one{B}$ \cite[Theorem 4.1.6]{hoshinoDoubleCategoriesRelations2023}.
    \item For a category with finite products $\one{C}$, the loosewise chaotic double category $\dbl{C}\mathrm{h}(\one{C})$, the double category with the source, target, and identity functors in $\bi{Cat}$ being the identity on $\one{C}$, is not unit-pure unless $\one{C}$ is the trivial category.
    Correspondingly, the identity functor $\one{C}\to\one{C}$, seen as a fibration, does not have comprehensive diagonals unless $\one{C}$ is trivial.
    \qedhere
    \end{enumerate}
\end{example}

\subsection{The unique choice principle and map-completeness}
\label{subsec:funccomp}

The unique choice principle, or functional comprehension, states that if a predicate $\alpha(x,y)$ is total and single-valued in the sense that
\[
    \begin{aligned}
        \forall x : I. \exists y : J. \ & \alpha(x,y) \\
        \forall x : I. \forall y, y' : J. \ & 
        \left( \alpha(x,y) \land \alpha(x,y') \right) \Longrightarrow y = y', 
    \end{aligned}
\]
then there exists a function $f\colon I\to J$ such that
\[
    \forall x : I. \forall y : J. \left( \alpha(x,y) \Longleftrightarrow f(x) = y \right). 
\]
The totalness and single-valuedness of $\alpha$ are equivalently 
expressed as the following when $\beta(y,x):\equiv\alpha(x,y)$:
\[
    \begin{aligned}
    \forall x , x' : I.  &
    \left( x = x' \Longrightarrow \exists y:J. \alpha(x,y) \land\beta(y,x') \right) \\
    \forall x : I. \forall y, y' : J. \ &
    \left( \beta(y,x) \land \alpha(x,y') \right) \Longrightarrow y = y'.
    \end{aligned},
\]
which serve as the unit and the counit of the adjunction $\alpha\dashv\beta$. 
The triangle identities are only meaningful when we respect the proofs of these implications,
but they should be satisfied proof-theoretically, as explained in \cite{Pav95}.
This shows the importance of left adjoints (or maps) in a bicategory
as an appropriate categorical counterpart of functional relations.
If we follow this view, the unique choice principle is a statement that a left adjoint relation
always arises from a function.

There have been several studies on this notion in various contexts, such as \cite{Pav95,Pav96} in fibrations, \cite{MPR17} in hyperdoctrines, and \cite{dagnino2024cauchycompletionsruleuniquechoice} in relational doctrines.
In this paper, we follow the terminology in \cite{Pav96} and call an elementary existential fibration with this property \emph{function comprehensive}. 
A double categorical counterpart of this notion is the map-completeness of a double category, which is defined as follows.

\begin{definition}[{\cite[Definition 19]{Par21},\cite[\S 4.2]{hoshinoDoubleCategoriesRelations2023}}]
	A double category $\dbl{D}$ is \emph{map-complete}\footnote{
    This condition is called \textit{Cauchy} in \cite{Par21,hoshinoDoubleCategoriesRelations2023}, due to the fact that a category is Cauchy-complete if and only if every left adjoint profunctor into it is representable.
    } if any adjunction $\alpha\colon I\adjointleft J\lon \beta$ in the bicategory $\LBi{\dbl{D}}$ is representable, namely, is of the form $f_*\colon I\adjointleft J\lon f^*$ for some tight arrow $f\colon I\to J$.
\end{definition}

Maps in a bicategory do not form a category in general, but they do form a bicategory.

\begin{definition}[{\cite[Definition 1.5]{CW87}}]
	For a bicategory $\bi{B}$,
	$\MapBi{B}$ is the locally full sub-bicategory of $\bi{B}$ spanned by the left adjoint 1-cells.
	We call it the \emph{map bicategory} of $\bi{B}$\footnote{
	This terminology comes from the convention to call a left adjoint 1-cell a \textit{map} \cite{lawvereMetricSpacesGeneralized1973a}.}.
	By abuse of notation, we write $\MapBi{\dbl{D}}$ for $\MapBi{\LBi{\dbl{D}}}$ where $\dbl{D}$ is a double category.
\end{definition}

\begin{remark}
    \label{rem:map-complete}
    It is worth rephrasing the map-completeness of a double category $\dbl{D}$ in terms of the map bicategory $\MapBi{\LBi{\dbl{D}}}$.
        For an equipment $\dbl{D}$, there is a functor $(-)_*\colon\TBi{\dbl{D}}(I,J)\to\MapBi{\LBi{\dbl{D}}}(I,J)$ for every pair of objects $I$ and $J$ in $\dbl{D}$, which sends a tight arrow to its companion.
    This functor is always fully faithful.  
    $\dbl{D}$ is map-complete if and only if the functor $(-)_*$ is essentially surjective, and hence an equivalence.
\end{remark}

\begin{remark}
    \label{rem:unique-choice}
    What is called unique choice in \cite[Definition 4.9.1]{jacobsCategoricalLogicType1999a} is different from the unique choice principle we are discussing here, since the condition therein involves coproducts along product projections when a predicate is single-valued, but not necessarily total.
    \cite[Proposition 5.3]{MPR17} seems to be proven by assuming that these two conditions are equivalent, but it is not clear whether this is the case.
\end{remark}

The Frobenius condition on cartesian equipments has some implications on the behavior of maps in the bicategory of loose arrows, which was discussed in \cite{WW08} and similarly in \cite{Pav96} in the context of regular fibrations.
We only state the results we will use later.
\begin{proposition}[{\cite[Theorems 3.11, 3.12]{WW08}}]
    \label{prop:frobenius-map}
    Let $\dbl{D}$ be a Frobenius cartesian equipment.
    \begin{enumerate}
        \item If a loose arrow $\alpha\colon I\sto J$ has a right adjoint, the dual $\alpha^\dagger$ is its right adjoint.
        \item The map bicategory $\MapBi{\LBi{\dbl{D}}}$ is locally groupoidal, meaning that every 2-cell is invertible.
        \qedhere
    \end{enumerate}
\end{proposition}
\begin{proof}
    We know $\dbl{D}$ is of the form $\Bil[\mf{p}]$ for some elementary existential fibration $\mf{p}$ by \Cref{cor:unilateral}, and the statements are given in \cite[Proposition 4.2]{Pav96} and a discussion before it.
    Note that the proofs there assume the fibration to be regular, but nevertheless, the proofs are valid for any elementary existential fibration, similarly to the arguments in \cite{Pav95}.
    
    Alternatively, we can use the fact that $\LBi{\dbl{D}}$ is a cartesian bicategory by \Cref{prop:CartdoubleCartbicat} to reduce the statements to the results in \cite{WW08}.
\end{proof} 

\subsection{Extensional equality and separability}
\label{subsec:extensionality}
As we do not assume preorderedness on the fibers of fibrations to achieve proof-irrelevance, equality in elementary fibrations is not extensional \textit{a priori}.
Although fiberwise preorderedness implies existential equality, we do not have the converse.
Therefore, it is worth discussing this condition independently as a property of fibrations.

The reflexivity of equality     
\[
    \forall x:I.\quad (\top \Longrightarrow x=x)
\]
corresponds to the unit component of the adjunction $\sum_{\tpl{0,0}}\dashv(-)[\tpl{0,0}]$ at the terminal in the fiber, as the equality $x=y$ corresponds to $\sum_{\tpl{0,0}}\top$. 
This unit component always has a retraction, but this is not the inverse in general.
One form of extensional equality is to ask this to be the inverse.
This is logically equivalent to the principle known as \textit{uniqueness of identity proofs (UIP)}, which states that a proof of equality $x=y$ is unique if it exists.
In categorical terms, this states that the equality object is subterminal, that is, it admits at most one arrow from each object into it.

Now, we formulate these equivalent conditions as a property of fibrations.
Recall that we write $\delta_I$ for the equality predicate in the fiber $\one{E}_{I\times I}$, which is defined as $\sum_{\tpl{0,0}}\top_I$.

\begin{definition}
    \label{def:extensional_equality}
    Let $\mf{p}\colon\one{E}\to\one{B}$ be an elementary fibration.
    An object $I$ in $\one{B}$ \emph{admits extensional equality} in $\mf{p}$ if the object $\delta_I[\tpl{0,0}]$ is a terminal object in $\one{E}_{I}$.
    An elementary fibration is said to be \emph{extensional} if every object in the base category has extensional equality.
\end{definition}

\begin{proposition}
    \label{prop:extensional_equality}
    For an elementary fibration $\mf{p}\colon\one{E}\to\one{B}$ and an object $I$ in $\one{B}$, the following are equivalent:
    \begin{enumerate}
        \item $I$ admits extensional equality in $\mf{p}$.
        \item $\delta_I[\tpl{0,0}]$ is a subterminal object in $\one{E}_I$.
        \item $\delta_I$ is a subterminal object in $\one{E}_{I\times I}$.
        \item $\sum_{\tpl{0,0}}\colon\one{E}_I\to\one{E}_{I\times I}$ preserves subterminal objects.
        \qedhere
    \end{enumerate}
\end{proposition}
\begin{proof}
    The equivalence (i)$\iff$(ii) follows from the fact that there are always mutual morphisms between $\top_I$ and $(\sum_{\tpl{0,0}}\top_I)[\tpl{0,0}]$ in $\one{E}_I$.
    The implications (iv)$\Rightarrow$(iii)$\Rightarrow$(ii) is immediate.
    For (i)$implies$(iii), we have the following commutative diagram where the vertical morphism is an isomorphism by Frobenius reciprocity:
\[\begin{tikzcd}[ampersand replacement=\&]
	\& {\sum_{\tpl{0,0}}(\top_I\land\delta_I[\tpl{0,0}])} \&\&\& \\
	\&\&\& {\sum_{\tpl{0,0}}\top_I} \& {\delta_I} \\
	{\delta_I\land\delta_I} \& {\sum_{\tpl{0,0}}\top_I\land\delta_I}
	\arrow["{\sum_{\tpl{0,0}}\pi_0}", from=1-2, to=2-4]
	\arrow[""{name=0, anchor=center, inner sep=0}, "\cong"', from=1-2, to=3-2]
	\arrow["{=}"{description}, draw=none, from=2-4, to=2-5]
	\arrow["{=}"{description}, draw=none, from=3-1, to=3-2]
	\arrow["{\pi_0}"', from=3-2, to=2-4]
	\arrow["\circlearrowright"{description}, draw=none, from=0, to=2-4]
\end{tikzcd}\]
    However, the assumption that $(\sum_{\tpl{0,0}}\top_I)[\tpl{0,0}]$ is terminal implies that the morphism $\sum_{\tpl{0,0}}\top_I$ is also an isomorphism, and hence the morphism $\pi_0$ is an isomorphism as well.
    This means $\sum_{\tpl{0,0}}\top_I$ is a subterminal object.
    For (iii)$\Rightarrow$(iv), we apply \Cref{lem:implicationsoffrobele} to deduce that for any object $\alpha$ in $\one{E}_I$, we have the canonical isomorphism 
    \[
    \sum_{\tpl{0,0}}\alpha\cong\sum_{\tpl{0,0}}(\top_I\land\alpha)\cong(\sum_{\tpl{0,0}}\top_I)\land\alpha[\tpl{0,0}]=\delta_I\land\alpha[\tpl{0,0}].
    \]
    Thus, when $\alpha$ is a subterminal object, so is this object.
\end{proof}

We introduce a double-categorical expression of extensional equality.
\begin{definition}
	\label{def:CoverInclusion}
	Let $\dbl{D}$ be a double category.
	We say a tight arrow $f\colon I\to J$ is an \emph{inclusion} if
	the loose identity cell on $f$ is prone.
	\[
		f\text{ is an inclusion }
		\Longleftrightarrow
		\begin{tikzcd}[virtual]
			I
			\ar[d, "f"']
			\sar[r, "\delta_I"]
			\doublecell{\cart}
				&
				I
				\ar[d, "f"]
			\\
			J
			\sar[r, "\delta_J"']
				&
				J
		\end{tikzcd}
	\]
\end{definition}

\begin{definition}[{\cite[Definition 1.2]{LWW10}}]
    Let $\dbl{D}$ be a cartesian equipment.
    An object $I$ in $\dbl{D}$ is called \emph{separable} if $\tpl{0,0}\colon I\to I\times I$ is an inclusion.
    A cartesian equipment is called \emph{separable} if every object in it is separable.
\end{definition}

\begin{remark}
    In \cite[Definition 3.10]{LWW10} and \cite[Definition 3.1.7]{hoshinoDoubleCategoriesRelations2023}, an object that is Frobenius and separable is called \textit{discrete}, modifying the terminology in \cite{CW87,lambertDoubleCategoriesRelations2022}. However, this terminology is still unfortunate, as Frobenius and separable objects in the double category of categories and profunctors are not exactly discrete categories, but rather preordered groupoids. In this paper, we do not introduce any terminology for such objects, and simply call them Frobenius and separable.
\end{remark}

In \cite{LWW10}, separability of cartesian bicategories is studied, which we reformulate in terms of cartesian equipments.
\begin{proposition}[{\cite[Propostion 3.4]{LWW10}}]
    \label{prop:separable}
    For a cartesian equipment $\dbl{D}$ and an object $J$ in $\dbl{D}$, the following are equivalent:
    \begin{enumerate}
        \item $J$ is separable.
        \item The identity loose arrow on $J$ is a subterminal object in $\LBi{\dbl{D}}(J,J)$.
        \item Every map $\alpha\colon I\sto J$ in $\dbl{D}$ is a subterminal object in $\LBi{\dbl{D}}(I,J)$.
    \end{enumerate}
    Therefore, for a cartesian equipment $\dbl{D}$, the following are equivalent:
    \begin{enumerate}
        \item $\dbl{D}$ is separable.
        \item The identity loose arrow on every object in $\dbl{D}$ is a subterminal object in the loose hom-category.
        \item Every map in $\dbl{D}$ is a subterminal object in the loose hom-category.
    \end{enumerate}
    In particular, when $\dbl{D}$ is separable, the map bicategory $\MapBi{\LBi{\dbl{D}}}$ is locally preordered.
\end{proposition}
\begin{proof}
    We prove the first statement, and the rest follows from it.
    The identity loose arrow on $J\times J$ is isomorphic to $\delta_J\times\delta_J$, and its restriction along the diagonals is $\delta_J\land\delta_J$, with the diagonal globular cell $\delta_J\Rightarrow\delta_J\land\delta_J$ induced by the universal property of the restriction. 
    Therefore, $J$ is separable if and only if this diagonal cell is an isomorphism, which means $\delta_J$ is a subterminal object. This proves the equivalence of (i) and (ii).

    (ii) is an instance of (iii), so it suffices to prove that (ii) implies (iii). 
    For a map $\alpha$, the precomposition with $\alpha$ gives a right adjoint $\LBi{\dbl{D}}(J,J)\to\LBi{\dbl{D}}(I,J)$, with its left adjoint given by the precomposition of the right adjoint of $\alpha$. 
    Thus, a subterminal object is preserved by the functor, and $\alpha$ as the image of the identity is also subterminal.
\end{proof}

\begin{proposition}
    \label{prop:separable-extensional}
    For a Frobenius cartesian equipment $\dbl{D}$, the following are equivalent:
    \begin{enumerate}
        \item $\dbl{D}$ is separable.
        \item $\dbl{D}$ is equivalent to $\Bil[\mf{p}]$ for an extensional elementary existential fibration $\mf{p}$.
        \qedhere
    \end{enumerate}
\end{proposition}
\begin{proof}
    By \Cref{cor:unilateral}, a Frobenius cartesian equipment can be expressed as $\Bil[\mf{p}]$ for an elementary existential fibration $\mf{p}$ over $\dbl{D}_0$.
    Note that the loose identity on $I$ corresponds to the equality predicate $\delta_I$ in $\one{E}_{I\times I}$.
    Thus, (iii) in \Cref{prop:extensional_equality} and (ii) in \Cref{prop:separable} are equivalent for each object.
\end{proof}

\begin{example}
    \label{ex:separable}\
    \begin{enumerate}
        \item A regular fibration is extensional because the Beck-Chevalley condition for the pullback of $\tpl{0,0}$ along itself implies the extensional equality. It also follows from the fact that regular fibrations correspond to cartesian equipments with Beck-Chevalley pullbacks (\Cref{cor:BeckChevalleyall}). 
        \item A doctrine is extensional as every object in fibers is subterminal.
        \item In $\PROF$, an object is separable if and only if it is thin, \text{i.e.,} a preordered set.
        Thus, $\PROF$ is not separable as a whole.
        When we restrict the double category to groupoids, it becomes a Frobenius cartesian equipment.
        This corresponds to a non-extensional elementary existential fibration.
        \qedhere
    \end{enumerate}
\end{example}

\subsection{Cartesian bicategories and map-completion}
\label{subsec:cartesianbicategories}
Since tight arrows are completely recoverable from the data of the loose bicategories in a map-complete cartesian equipment, it is natural to seek an equivalent formulation of map-complete cartesian equipments in terms of bicategories.
Roughly speaking, one can say that the framework of cartesian bicategories \cite{CKS84,CW87,CKWW07,LWW10} is this bicategorical formulation, although it was found earlier than the notion of cartesian equipments.
In this section, we discuss the relationship between map-complete cartesian equipments and cartesian bicategories, revisiting results from the literature.

\begin{remark}[Allegories: an alternative bicategorical approach]
    \label{rem:allegory}
    Although we focus on cartesian bicategories in this paper, one could also consider the framework of allegories \cite{FS90,johnstoneSketchesElephantTopos2002a}\footnote{
According to \cite{Bunge2017IntroductionAP}, the original content of \cite{FS90} had been already presented in the 1970s, in the unpublished paper ``On Canonizing Category Theory'', or ``On Functorializing Model Theory'' by Freyd in 1974.
} as a bicategorical approach to capture the notion of relations.
    Unlike cartesian bicategories, allegories are defined as locally posetal bicategories with an involution structure subject to some axioms.
    A similar approach is taken in the theory of ordered categories with involution, as in \cite{CGR84,Lam99}.
    More recently, allegories of relations with respect to  factorization systems were studied in \cite{HNSTY22}.
    
    Bicategories of relations in regular categories are axiomatized in terms of both allegories and cartesian bicategories, so they are related in this sense.
    A broader connection between them was studied in \cite{Law15}, where unitary pre-tabular allegories were proved to be equivalent to bicategories of relations in the sense of \cite{CW87} (\cite[Proposition 2.2.33, Theorem 2.2.34]{Law15}), or in more explicit terms, locally posetal Frobenius cartesian bicategories.
\end{remark}

First, we recall the definition of cartesian bicategories, and see how we can obtain a cartesian bicategory from a cartesian equipment.

\begin{definition}[{\cite{CKWW07}}]
	\label{def:cartbicat}
	A bicategory $\bi{B}$ is a \emph{cartesian bicategory} if it satisfies the following:
	\begin{enumerate}
		\item Each hom-category $\bi{B}(I,J)$ has finite products for any pair of objects $I,J$ in $\bi{B}$.
		\item $\MapBi{\bi{B}}$ has finite products in the sense of bilimits.
		More explicitly, 
		\begin{enumerate}
			\item there is an object $1$ such that for every object $K$ in $\bi{B}$, $\MapBi{\bi{B}}(K,1)$ is equivalent to the terminal category, and 
			\item for every pair of objects $I,J$ in $\bi{B}$,
			there is an object $I\times J$ in $\bi{B}$ such that
			for every object $K$ in $\bi{B}$,
			$\MapBi{\bi{B}}(K,I\times J)$ is equivalent to the product category $\MapBi{\bi{B}}(K,I)\times\MapBi{\bi{B}}(K,J)$. 
		\end{enumerate}
		\item The lax functors $1\colon\bi{1}\to\bi{B}$ and $\times\colon\bi{B}\times\bi{B}\to\bi{B}$ induced by 
		the terminal object and the product object in the method described in \cite{CKWW07}
		are pseudofunctors.
        \qedhere
	\end{enumerate}
\end{definition}

Most examples of cartesian bicategories arise from cartesian equipments by taking the loose part of the double category.
In fact, for a cartesian equipment $\dbl{D}$, the loose bicategory $\LBi{\dbl{D}}$ is a cartesian bicategory, proven in \cite{patterson2024transposingcartesianstructuredouble}.

\begin{theorem}[{\cite[Corollary 5.8]{patterson2024transposingcartesianstructuredouble}}]
	\label{prop:CartdoubleCartbicat}
	The loose bicategory $\LBi{\dbl{D}}$ of a cartesian equipment
	$\dbl{D}$ is a cartesian bicategory.
\end{theorem}

We have also obtained a proof of this fact independently, which is based on string diagrammatic reasoning and is included in \Cref{sec:diagrammatic-proof}.

\begin{remark}
        Patterson's proof is based on Trimble's reformulation of cartesian bicategories \cite{nlab:cartesian_bicategory} and is more abstract than our approach.
    Another proof of this proposition was given in the locally posetal case in \cite{lambertDoubleCategoriesRelations2022}, but it contains an error in the proof.
    In the last line of the proof, the author applies the universal property of the binary product to loose arrows, which is not valid in general.
	Indeed, the bicategory of relations on sets has disjoint unions as the binary product, not the same as a category of functions.
\end{remark}

For the converse direction, namely, to construct a cartesian equipment from a cartesian bicategory, it is natural to take functional relations as the tight arrows, as in \cite{JW00}. As we showed in the previous section, functional relations are categorically characterized as left adjoints in the bicategory of relations, and for this reason, in the context of cartesian bicategories, left adjoints, or maps, are used as ``functions'' defined in the bicategory. 
Accordingly, one would try to construct a double category from a cartesian bicategory by taking the maps as tight arrows and the 1-cells as loose arrows. 

One problem of this construction is that the composition of tight arrows in a double category is strictly associative and unital, and hence, when we take maps in a bicategory to construct a double category, one needs to quotient, which brings about a coherence issue.
One way to avoid this is to consider doubly-pseudo double categories, or \textit{double bicategories} introduced in \cite{Ver11}.
In Example 1.5.19 of this paper, the author constructs a double-bicategorical equivalent of an equipment from a bicategory, and characterizes cartesian bicategories in the sense of \cite{CKWW07} as those that induce a cartesian object in the bicategory-enriched category of (double-bicategorical) equipments and homomorphisms \cite[{(p.152)}]{Ver11}.

On the other hand, we can avoid this coherence issue by restricting to a special case where maps in a bicategory essentially form a category.
\begin{definition}
    \label{def:mapdiscrete}
    A bicategory is called \emph{map-discrete} if the locally full sub-bicategory $\MapBi{\bi{B}}$ of $\bi{B}$ is locally equivalent to discrete categories, namely, for every pair of objects $I$ and $J$ in $\bi{B}$, $\MapBi{\bi{B}}(I,J)$ is equivalent to a discrete category.
    A double category $\dbl{D}$ is called \emph{map-discrete} if its loose bicategory $\LBi{\dbl{D}}$ is map-discrete.
\end{definition}

There are two important classes of cartesian equipments that are map-discrete, which we will present in the following.

\begin{proposition}
    \label{lem:mapdiscrete}
    A unit-pure and map-complete equipment $\dbl{D}$ is map-discrete.
\end{proposition}
\begin{proof}
    From \Cref{rem:map-complete}, we see that if a double category $\dbl{D}$ is map-complete, then the functor from $\TBi{\dbl{D}}(I,J)$ to $\MapBi{\LBi{\dbl{D}}}(I,J)$ that sends a tight arrow $f$ to its companion $f_*$ is an equivalence for every pair of objects $I$ and $J$ in $\dbl{D}$.
    Since $\dbl{D}$ is also unit-pure, $\TBi{\dbl{D}}(I,J)$ is discrete for each $I$ and $J$.
    Therefore, $\MapBi{\LBi{\dbl{D}}}(I,J)$ is also discrete for each $I$ and $J$, which means that $\LBi{\dbl{D}}$ is map-discrete.
\end{proof}

\begin{lemma}[{\cite[Proposition 3.14]{LWW10}}]
    \label{lem:separableFrobeniusmapdiscrete}
    A separable Frobenius cartesian equipment is map-discrete.
\end{lemma}
\begin{proof}
    By \Cref{prop:separable}, for a separable cartesian equipment $\dbl{D}$, the map bicategory $\MapBi{\LBi{\dbl{D}}}$ is locally preordered.
    By \Cref{prop:frobenius-map}, it is also locally groupoidal, and hence, locally equivalent to discrete categories.
\end{proof}

\begin{definition}
    \label{def:mapdoublecategory}
    Let $\bi{B}$ be a map-discrete bicategory.
    We define the \textbf{double category of maps} $\Map{\bi{B}}$ as follows:
    \begin{itemize}
        \item Its objects are the same as the 0-cells of $\bi{B}$.
        \item Its tight arrows are the isomorphism classes of maps in $\bi{B}$.
        The composition of tight arrows is the composition as in $\bi{B}$.
        \item  Its loose arrows are the same as the 1-cells of $\bi{B}$.
        The composition of loose arrows is the composition as in $\bi{B}$.
        \item Its cells of the form on the left below are the 2-cells in $\bi{B}$.
        \[
            \begin{tikzcd}
                I
                    \sar[r, "\alpha"]
                    \ar[d, "{[f]}"']
                    \ar[dr, "\mu",phantom]
                &
                J
                    \ar[d, "{[g]}"]
                \\
                K
                    \sar[r, "\beta"']
                &
                L
            \end{tikzcd}
            \hspace{2ex}
            \vline\!\vline
            \hspace{2ex}
            \begin{tikzcd}[row sep=tiny]
                &
                J
                \sar[dr, "g"]
                \\
                I
                \sar[ur, "\alpha"]
                \sar[dr, "f"']
                \ar[rr, "\mu", phantom]
                &&
                K
                \\
                &
                L
                \sar[ur, "\beta"']
            \end{tikzcd}
        \]
        Note that the choice of representatives $f$ and $g$ is arbitrary.
        \item The composition of cells is defined by the composition in $\bi{B}$ as follows:
        \[
            \begin{tikzcd}
                I
                    \sar[r, "\alpha"]
                    \ar[d, "{[f]}"']
                    \ar[dr, "\mu",phantom]
                &
                J
                    \ar[d, "{[g]}"]
                    \sar[r, "\beta"]
                    \ar[dr, "\nu",phantom]
                &
                K
                    \ar[d, "{[h]}"]
                \\
                I'
                    \sar[r, "\alpha'"']
                &
                J'
                    \sar[r, "\beta'"']
                &
                K'
            \end{tikzcd}
            \hspace{1ex}
            \vline\!\vline
            \hspace{1ex}
            \begin{tikzcd}[row sep=tiny, column sep=small]
                &&
                K
                \sar[dr, "h"]
                \\
                &
                J
                \sar[dr, "g"]
                \ar[rr, "\nu\Downarrow", phantom]
                \sar[ur, "\beta"]
                &&
                K'
                \\
                I
                \sar[ur, "\alpha"]
                \ar[rr, "\mu\Downarrow", phantom]
                \sar[dr, "f"']
                &&
                J'
                \sar[ur, "\beta'"']
                \\
                &
                I'
                \sar[ur, "\alpha'"']
            \end{tikzcd}
            ,
            \hspace{2ex}
            \begin{tikzcd}
                I
                    \sar[r, "\alpha"]
                    \ar[d, "{[f]}"']
                    \ar[dr, "\mu",phantom]
                &
                J
                    \ar[d, "{[g]}"]
                \\
                I'
                    \sar[r, "\alpha'"']
                    \ar[d, "{[f']}"']
                    \ar[dr, "\mu'",phantom]
                &
                J'
                    \ar[d, "{[g']}"]
                \\
                I''
                    \sar[r, "\alpha''"']
                &
                J''
            \end{tikzcd}
            \hspace{1ex}
            \vline\!\vline
            \hspace{1ex}
            \begin{tikzcd}[row sep=tiny, column sep=small]
                &
                J
                \sar[dr, "g"]
                \\
                I
                \sar[dr, "f"']
                \sar[ur, "\alpha"]
                \ar[rr, "\mu\Downarrow", phantom]
                &&
                J'
                \sar[dr, "g'"]
                \\
                &
                I'
                \sar[ur, "\alpha'"']
                \sar[dr, "f'"']
                \ar[rr, "\mu'\Downarrow", phantom]
                &&
                J''
                \\
                &&
                I''
                \sar[ur, "\alpha''"']
            \end{tikzcd}.
        \]
        These compositions are well-defined regardless of the choice of representatives by the map-discreteness of $\bi{B}$.
        In particular, the horizontal composition of cells requires the representative of the tight arrow in the middle by composing the unique isomorphism, but the result is independent of the choice of the representative.
        \qedhere
    \end{itemize}
\end{definition}

\begin{lemma}
    \label{lem:mapdoublecategory}
    The map double category $\Map{\bi{B}}$ for a map-discrete bicategory $\bi{B}$ is a unit-pure map-complete equipment,
    and the loose bicategory $\LBi{\Map{\bi{B}}}$ is isomorphic to $\bi{B}$.
\end{lemma}
\begin{proof}
    Since the cells with the loose arrows at the top and the bottom being identities correspond to
    the 2-cells between the maps on the left and the right,
    the two maps are equal up to unique isomorphism and hence the cell is the identity.
    This shows the unit-pureness of $\Map{\bi{B}}$.
    The companion and the conjoint for a map $f$,
    which has a right adjoint $f^*$,
    are given by the following cells:
    \begin{align*}
    &\text{(companion)}
    &
    \begin{tikzcd}[row sep=small, column sep=small, ampersand replacement=\&]
        I
        \sar[rr, "f_*=f"]
        \ar[rd, "f"']
        \ar[rr, phantom, "\cart"{yshift=-1.5ex}]
        \&\&
        J
        \ar[ld, equal]
        \\
        \&
        J
        \&
    \end{tikzcd}
    &
    \hspace{0.5ex}
    \vline\!\vline
    \hspace{0.5ex}
    \begin{tikzcd}[row sep=tiny, column sep=small, ampersand replacement=\&]
        \&
        J
        \sar[rd, equal]
        \\
        I
        \sar[rr, bend right=15, "f"']
        \sar[ru, "f"]
        \ar[rr, phantom, "\rotatebox{90}{$=$}"{yshift=1ex}]
        \&\&
        J
    \end{tikzcd}
    \hspace{0.5ex}
    ,
    \hspace{0.5ex}
    &\begin{tikzcd}[row sep=small, column sep=small, ampersand replacement=\&]
        \&
        I
        \ar[ld, equal]
        \ar[rd, "f"]
        \&
        \\
        I
        \sar[rr, "{f_*=f}"']
        \ar[rr, phantom, "\opcart"{yshift=1.5ex}]
        \&\&
        J
    \end{tikzcd}
    &
    \hspace{0.5ex}
    \vline\!\vline
    \hspace{0.5ex}
    \begin{tikzcd}[row sep=tiny, column sep=small, ampersand replacement=\&]
        I
        \sar[rr, bend left=15, "f"]
        \sar[rd, equal]
        \ar[rr, phantom, "\rotatebox{90}{$=$}"{yshift=-1ex}]
        \&\&
        J
        \\
        \&
        I
        \sar[ru, "f"']
        \&
    \end{tikzcd}
    \\
    &\text{(conjoint)}
    &
    \begin{tikzcd}[row sep=small, column sep=small, ampersand replacement=\&]
        J
        \sar[rr, "f^*"]
        \ar[rd,equal]
        \ar[rr, phantom, "\cart"{yshift=-1.5ex}]
        \&\&
        I
        \ar[ld, "f"]
        \\
        \&
        J
        \&
    \end{tikzcd}
    &
    \hspace{0.5ex}
    \vline\!\vline
    \hspace{0.5ex}
    \begin{tikzcd}[row sep=tiny, column sep=small, ampersand replacement=\&]
        \&
        I
        \sar[rd, "f"]
        \\
        J
        \sar[rr, bend right=15, equal]
        \sar[ru, "f^*"]
        \ar[rr, phantom, "\varepsilon_{f}"{yshift=1ex}]
        \&\&
        J
    \end{tikzcd}
    \hspace{0.5ex}
    ,
    \hspace{0.5ex}
    &\begin{tikzcd}[row sep=small, column sep=small, ampersand replacement=\&]
        \&
        I
        \ar[ld, "f"']
        \ar[rd, equal]
        \&
        \\
        J
        \sar[rr, "{f^*}"']
        \ar[rr, phantom, "\opcart"{yshift=1.5ex}]
        \&\&
        I
    \end{tikzcd}
    &
    \hspace{0.5ex}
    \vline\!\vline
    \hspace{0.5ex}
    \begin{tikzcd}[row sep=tiny, column sep=small, ampersand replacement=\&]
        I
        \sar[rr, bend left=15, equal]
        \sar[rd, "f"']
        \ar[rr, phantom, "\eta_{f}"{yshift=-1ex}]
        \&\&
        I
        \\
        \&
        J
        \sar[ru, "f^*"']
        \&
    \end{tikzcd}
    \end{align*}
    Thus, $\Map{\bi{B}}$ is an equipment.
    The map-completeness also follows immediately from the above construction.
    The isomorphism of the loose bicategory of $\Map{\bi{B}}$ and $\bi{B}$ is clear from the construction. 
\end{proof}

\begin{proposition}
    \label{prop:mapdiscreteandmap-complete}\ 
    \begin{enumerate}
        \item For a map-discrete cartesian bicategory $\bi{B}$,
        the double category of maps $\Map{\bi{B}}$ is a unit-pure map-complete cartesian equipment,
        with $\LBi{\Map{\bi{B}}}\cong\bi{B}$.
        \item For a unit-pure map-complete cartesian equipment $\dbl{D}$, $\LBi{\dbl{D}}$ is a map-discrete cartesian bicategory, and $\Map{\LBi{\dbl{D}}}\cong\dbl{D}$.
    \end{enumerate}
    Therefore, the categories of map-discrete cartesian bicategories and unit-pure map-complete cartesian equipments are equivalent.
\end{proposition}
\begin{proof}
    \begin{enumerate}
        \item By \Cref{lem:mapdoublecategory}, $\Map{\bi{B}}$ is a unit-pure map-complete equipment.
        The cartesianness of $\Map{\bi{B}}$ is proven in the same way as in \cite[Example 1.5.19]{Ver11}. 
        Recall the characterization of cartesian equipments given in \cite[Corollary 4.3.3]{aleiferiCartesianDoubleCategories2018}
        The category $\Map{\bi{B}}_0$ is biequivalent to $\MapBi{\bi{B}}$ as a bicategory, and since the latter has finite biproducts, so does the former.
        However, it is locally discrete by the assumption, hence it has strict finite products.
        The finite products in the loose hom-categories follow from the assumption, and the pseudo-ness condition follows from the corresponding condition in the definition of cartesian bicategories.
        \item The first part of the statement follows from \Cref{lem:mapdiscrete,prop:CartdoubleCartbicat}.
        We take the isomorphism classes of maps in $\LBi{\dbl{D}}$ as the tight arrows of $\Map{\LBi{\dbl{D}}}$, but these are in one-to-one correspondence with the tight arrows of $\dbl{D}$ by the unit-pureness and map-completeness of $\dbl{D}$, and hence we have an isomorphism of double categories $\Map{\LBi{\dbl{D}}}\cong\dbl{D}$.
        \qedhere
    \end{enumerate}
\end{proof}

We now turn our attention to the universal property that the cartesian equipment of maps $\Map{\bi{B}}$ satisfies.
Let $\ol{\dbl{D}}$ denote $\Map{\LBi{\dbl{D}}}$ from now on.

\begin{lemma}
    \label{lem:map-completion}
    Let $\dbl{D}$ be a map-discrete cartesian equipment.
    Then, it admits a double functor $C\colon\dbl{D}\to\ol{\dbl{D}}$ that is identity on objects and loose arrows, and sends each tight arrow $f$ to $f_*$.
    The codomain $\ol{\dbl{D}}$ is a unit-pure map-complete cartesian equipment.
\end{lemma}
\begin{proof}
    Since a cell $\tau$ in $\dbl{D}$ of the form below is in one-to-one correspondence with a globular cell
	$\widetilde{\tau}\colon \alpha g_*\Rightarrow f_*\beta$ in $\LBi{\dbl{D}}$,
	so $C$ sends such a cell to the cell on the right below.
	Note that a pseudo double functor preserves companions.
	\[
	\begin{tikzcd}
		I
			\ar[d,"f"']
			\sar[r,"\alpha"]
			\ar[dr, phantom, "\tau"]
		&
			J
			\ar[d,"g"]
		\\
			K
			\sar[r,"\beta"']	
		&
			L
	\end{tikzcd}
	\hspace{2ex},
	\hspace{2ex}
	\begin{tikzcd}[column sep = large, row sep = small]
		A 
			\ar[d,equal]
			\sar[r,"\alpha"]
			\doublecell[rd]{\rotatebox{90}{$=$}}
		&
		B 
			\ar[rd,"Cg",bend left=15]
			\ar[d,equal]
		\\
		I
			\sar[r,"\alpha"']
		    \ar[d,equal]
			\ar[drr, phantom, "\widetilde{\tau}"]
		&
		J
			\sar[r,"g_*"']
			\ar[r, phantom, shift left=1.5ex,xshift=-1ex, "\opcart"]
		&
		L
			\ar[d,equal]
		\\
		I
			\sar[r,"f_*"]
			\ar[rd, "Cf"', bend right=15]
			\ar[r, phantom, shift right=1.5ex, xshift=1ex,"\cart"]
			&
		K
			\sar[r,"\beta"]
			\ar[d,equal]
			\doublecell[rd]{\rotatebox{90}{$=$}}
		&
		L
			\ar[d,equal]
		\\
		&
		K
			\sar[r,"\beta"']
		&
		L
	\end{tikzcd}
	\]
	Thus defined $C$ is easily shown to be a double functor.

    The last statement is a direct consequence of \Cref{lem:mapdoublecategory}.
\end{proof}

\begin{proposition}[{\cite[Proposition 4.2.14]{hoshinoDoubleCategoriesRelations2023}}]
	\label{prop:univpropmap-complete}
    Suppose $\dbl{D}$ is a map-discrete cartesian equipment.
    Then, the canonical double functor $C\colon\dbl{D}\to\ol{\dbl{D}}$ in \Cref{lem:map-completion} exhibits $\ol{\dbl{D}}$ as the free unit-pure map-complete equipment on $\dbl{D}$, meaning the following universal property holds:
    \begin{enumerate}
        \item For any map-complete unit-pure equipment $\dbl{E}$ and a double functor $F\colon\dbl{D}\to\dbl{E}$, there exists a unique double functor $\wt F\colon\ol{\dbl{D}}\to\dbl{E}$ such that $F=\wt F\circ C$.
        \item For any 2-cell $\Psi\colon F\Rightarrow G$ in $\dbl{D}$, there exists a unique 2-cell $\wt\Psi\colon\wt F\Rightarrow\wt G$ in $\ol{\dbl{D}}$ such that $\Psi=\wt\Psi\circ C$.
    \qedhere
    \end{enumerate}
\end{proposition}

\begin{proof}
	In a unit-pure map-complete equipment, a tight arrow is uniquely determined by its representative adjoint pair.
	Therefore, for any pseudo double functor $F\colon\dbl{D}\to\dbl{E}$,
	assignment of the image of a tight arrow in $\ol{\dbl{D}}$ is uniquely determined by the image of its representative adjoint pair.
    We can also reduce general cells to a combination of the pair of tight arrows 
    and the corresponding globular cells in the loose bicategory,
	which implies that $\wt F$ is uniquely determined and also well-defined.

    Since $C$ is identity on the loose bicategory, the data of tightwise transformations 
    $\Psi\colon F\Rightarrow G$ and $\wt\Psi\colon\wt F\Rightarrow\wt G$ are the same.
    We prove that the naturality for the tight arrows in $\ol{\dbl{D}}$ automatically follows from the naturality of $F$.
    Given a tight arrow $f\colon I\to J$ in $\ol{\dbl{D}}$,
    we have the cell on the left below in $\dbl{E}$,
    which leads to the cell on the right below.
    \[
    \begin{tikzcd}
        FI
            \sar[r,"Ff_*"]
            \ar[d,"\Psi_I"']
            \ar[dr, phantom, "\Psi_{f_*}"]
        &
        FJ
            \ar[d,"\Psi_J"]
        \\
        GI
            \sar[r,"Gf_*"']
        &
        GJ
    \end{tikzcd}
    \hspace{2ex}
    \vline\!\vline
    \hspace{2ex}
    \begin{tikzcd}[row sep=tiny]
        &
        FI 
        \ar[dr, "Ff"]
        \ar[dl, "\Psi_I"']
        \\
        GI
        \ar[rr, phantom, "\wt\Psi_{f_*}"]
        \ar[dr, "Gf"']
        &&
        FJ
        \ar[dl, "\Psi_J"]
        \\
        &
        GJ
    \end{tikzcd}
    \]
    By the unit-pureness of $\ol{\dbl{D}}$, we have $\Psi_J\circ Ff = Gf\circ\Psi_I$.
    The naturality for the additional cells in $\ol{\dbl{D}}$ is also shown in a similar way.
\end{proof}

\begin{remark}
    \label{rem:Cauchisation}
    The proof only depends on the conclusion of \Cref{lem:map-completion}. The cartesian equipment with such a property is called the \textit{Cauchisation} of $\dbl{D}$ in \cite{hoshinoDoubleCategoriesRelations2023}.
    The functor should also posess a weak universal property for general map-complete equipments, but we do not discuss it here.
\end{remark}

\begin{example}
    \label{ex:mapdiscrete} 
    There are at least three appearances of map-completion in the literature.
    \begin{enumerate}
        \item \cite{Kel91} discusses the process of turning a stable (proper) factorization system into a regular factorization system, meaning it consists of the classes of regular epimorphisms and monomorphisms in a regular category.
        Since a stable factorization system $(\zero{E},\zero{M})$ on a category $\one{C}$ with finite limits is regular if and only if the double category of $\zero{M}$-relations $\Rel[\one{C}][(\zero{E},\zero{M})]$ is map-complete, this process can be understood as the map-completion; see \cite[Section 4.2]{hoshinoDoubleCategoriesRelations2023} for more details.
        Note that the cartesian equipment $\Rel[\one{C}][(\zero{E},\zero{M})]$ is map-discrete because it has Beck-Chevalley pullbacks, and hence is separable and Frobenius.
        \item In \cite[Section 8]{Pav96}, the author discusses the \textit{function comprehension completion} of regular fibrations.
        This is equivalent to map-completion up to the equivalence of regular fibrations and cartesian equipments with Beck-Chevalley pullbacks, \text{i.e.,} the completion is the same as $\uni[\Cau{\Bil[(-)]}]$ in our notation.
        For a regular fibration $\mf{p}\colon\one{E}\to\one{B}$, the equipment $\Bil[\mf{p}]$ is map-discrete since it is separable, but this was also proven in \cite[Proposition 4.2, Theorem 5.3]{Pav96}.
        \item In \cite{BSS21}, the authors construct a fully faithful functor from the category of locally posetal Frobenius cartesian bicategories to the category of elementary existential doctrines and a left adjoint to it.\footnote{
    Note that in this paper, the term \textit{cartesian bicategory} is used only for locally posetal ones, following the terminology in \cite{CW87}.
    However, it seems that they also assume the Frobenius law in the definition of cartesian bicategories, which is not compatible with the definition in \cite{CW87} nor \cite{CKWW07,WW08} (see \Cref{def:cartbicat}).
    The corresponding notion is rather called \textit{`bicategories of relations'} in \cite{CW87}.
    }
        Up to the equivalence between unit-pure map-complete cartesian equipments and map-discrete cartesian bicategories, the left adjoint gives the map-completion of the cartesian equipment $\Bil[\mf{p}]$ for an elementary existential doctrine $\mf{p}$.
        See \Cref{rem:map-completecompletion} for more details.
        \qedhere
    \end{enumerate}
\end{example}

\begin{remark}
    \label{rem:map-completecompletion}
    The map-completion of a map-discrete cartesian equipment $\dbl{D}$ together with the equivalence in \Cref{prop:mapdiscreteandmap-complete} gives the following left 2-adjoint, where we define 2-cells of cartesian bicategories to make the equivalence a 2-equivalence.
    \[
    \begin{tikzcd}
    {\Eqp\carttwo}_{,\textbf{m-disc}}
    \ar[r, shift left=1ex, "\ol{(-)}"]
    \ar[r, phantom, "\rotatebox{-90}{$\dashv$}"]
    \ar[r, shift right=1ex, hookleftarrow]    &
    {\Eqp\carttwo}_{,\textbf{m-comp},\textbf{unit-pure}}
    \ar[l, shift left=1ex, "\LBi{-}"]
    \ar[r, shift left=1ex, "\LBi{-}"]
    \ar[r, phantom, "\simeq"]
    &
    {\mathbf{CartBi}}_{\textbf{m-disc}}
    \ar[l, shift left=1ex, "\Map{-}"]
    \end{tikzcd}
    \]
    The composite left adjoint is $\LBi{-}$, whereas the composite right adjoint is $\Map{-}$.
    This shows the universal property of $\Map{\bi{B}}$ for a map-discrete bicategory $\bi{B}$ as the cofree map-discrete equipment.

    When we restrict to the Frobenius case, we can further connect this with the equivalence between the category of elementary existential fibrations and that of Frobenius cartesian equipments.
    \[
    \begin{tikzcd}
        {\Fib_\eef}_{,\textbf{m-disc}}
        \ar[r, shift left=1ex, "\Bil"]
        \ar[r, phantom, "\simeq"]
        &
        {\Eqp\Frob}_{,\textbf{m-disc}}
        \ar[l, shift left=1ex, "\uni"]
        \ar[r, shift left=1ex, "\LBi{-}"]
        \ar[r, phantom, "\rotatebox{-90}{$\dashv$}"]
        \ar[r, shift right=1ex, hookleftarrow]
        &
        {\mathbf{CartBi}\Frob}_{,\textbf{m-disc}}
        \ar[l, shift left=1ex, "\Map{-}"]
    \end{tikzcd}.
    \]
    Here, ${\Fib_\eef}_{,\textbf{m-disc}}$ is the 2-category of elementary existential fibrations that correspond to map-discrete cartesian equipments.
    Although we have not sorted out completely what this condition amounts to for fibrations, we know that  extensional elementary existential fibrations satisfy this from \Cref{prop:separable-extensional} and \Cref{lem:separableFrobeniusmapdiscrete}.
    In particular, elementary existential doctrines and regular fibrations are in ${\Fib_\eef}_{,\textbf{m-disc}}$.
    The image of the right adjoint is characterized as elementary existential fibrations with the three conditions, (1) having the comprehensive diagonals, (2) the unique choice principle, and (3) the condition corresponding to map-discreteness.
    In the locally posetal case, the above adjunction recovers the result in \cite{BSS21}.
\end{remark}

\printbibliography

\appendix

\section{A note on cartesian objects}\label{sec:cartesianobj}
\begin{definition}[{\cite[\S 5.1]{CKW91}}]
    \label{def:cartesianobj}
    A \emph{cartesian object} in a 2-category $\bi{K}$ with strict (2-dimensional) finite products $\monunit,\otimes$ is an object $x$ of $\bi{K}$ such that
    the canonical 1-cells $!\colon x\to \monunit$ and $\Delta\colon x\to x\otimes x$ have right adjoints $1\colon\monunit\to x$ and $\times\colon x\otimes x\to x$, respectively.
    A \emph{cartesian 1-cell} (or \emph{cartesian arrow}) in $\bi{K}$ is a 1-cell $f\colon x\to y$ between cartesian objects $x$ and $y$ of $\bi{K}$ such that
    the canonical 2-cells obtained by the mate construction $\times\circ(f\otimes f)\Rightarrow f\circ \times$ and $f\circ 1\Rightarrow 1$ are invertible. 

    For a 2-category $\bi{K}$ with strict finite products, 
    we write $\bi{K}_{\bi{cart}}$ for the 2-category of cartesian objects, cartesian 1-cells, and arbitrary 2-cells in $\bi{K}$.
\end{definition}

\begin{remark}
    \label{rem:cartesianobj}
    By strict finite products, we mean the most strict notion of finite products,
    that is, 
    the terminal object $\monunit$ and the binary product $\otimes$ 
    come with the isomorphisms in the 2-category of categories:
    \[
        \begin{aligned}
            \bi{K}(x,\monunit)&\cong \bi{1}\\
            \bi{K}(x,y\otimes z)&\cong \bi{K}(x,y)\times \bi{K}(x,z)
        \end{aligned}
    \]
    2-naturally in $x,y,z$.
    We call this kind of limits \emph{strict (2-)limits}.
\end{remark}

\begin{example}
    \label{ex:cartesianobj}
    In the 2-category $\Cat$ of categories, functors, and natural transformations, 
    the cartesian objects are the categories with finite products,
    where the right adjoints $1$ and $\times$ are the functors of
    the terminal object and the binary product, respectively.
\end{example}
\begin{lemma}
    \label{lemma:Cartesianequiv}
    Let $\bi{K}$ be a 2-category with strict finite products. 
    A 1-cell $f\colon x\to y$ in $\bi{K}_{\bi{cart}}$ is an equivalence in $\bi{K}_{\bi{cart}}$ 
    if and only if the underlying 1-cell of $f$ is an equivalence in $\bi{K}$.
\end{lemma}
\begin{proof}
    The only if part is clear since we have the forgetful 2-functor $\bi{K}_{\bi{cart}}\to\bi{K}$.
    For the if part, 
    take the right adjoint $g$ of the underlying 1-cell of $f$ as its inverse.
    Taking the right adjoint of both sides of the isomorphism 2-cells $!\circ f\cong\,!$ and $(f\otimes f)\circ\Delta\cong\Delta\circ f$,
    we obtain the isomorphism 2-cells $g\circ 1\cong 1$ and $\times\circ (g\otimes g)\cong g\circ\times$.
    This shows that $g$ gives a cartesian morphism from $y$ to $x$,
    and $g$ is indeed the inverse of $f$ in $\bi{K}_{\bi{cart}}$.
\end{proof}

\begin{lemma}
    \label{lem:cartobjinlffsub}
    Let $\bi{K},\bi{K}'$ be 2-categories with strict finite products $(\monunit,\otimes)$,
    and $\abs{-}\colon\bi{K}'\to\bi{K}$ be a 2-functor
    preserving strict finite products and locally full-inclusion, \textit{i.e.}, injective on 1-cells and bijective on 2-cells.
    For an object $x$ of $\bi{K}'$ to be cartesian,
    it is necessary and sufficient that $\abs{x}$ is cartesian in $\bi{K}$ and that
    the 1-cells $1\colon\monunit\to \abs{x}$ and $\times\colon \abs{x}\otimes \abs{x}\to \abs{x}$ 
    right adjoint to the canonical 1-cells are essentially in the image of $\abs{-}$.

    Moreover, for a 1-cell $f\colon x\to y$ of $\bi{K}'$ where $x$ and $y$ are cartesian in $\bi{K}'$,
    $f$ is cartesian in $\bi{K}'$ if and only if $\abs{f}$ is cartesian in $\bi{K}$.
\end{lemma}
\begin{proof}
    The necessity of the first condition follows from the fact that any 2-functor preserves adjunctions,
    that right adjoints are unique up to isomorphism, and that $\abs{-}$ preserves finite products.
    Since $\abs{-}$ is locally fully faithful, it also reflects units, counits, and the triangle identities with respect to the adjunctions,
    and hence the sufficiency of the first condition follows.

    The necessity of the second condition is again immediate from the fact that $\abs{-}$ preserves finite products.
    The sufficiency follws from the fact that $\abs{-}$ is locally fully faithful, in particular, reflects isomorphisms.
\end{proof}

\begin{lemma}
    \label{lem:cartobjinpb}
    Let $\bi{K}$, $\bi{L}$, and $\bi{M}$ 
    be a 2-category with strict finite products $(\monunit,\otimes)$,
    and $T\colon\bi{K}\to\bi{M}$ and $S\colon\bi{L}\to\bi{M}$ be 2-functors 
    preserving the finite products strictly, and locally isofibrations.
    Then, the canonical 2-functor 
    \[
        \left(\bi{K}\times_{\bi{M}}\bi{L}\right)\carttwo
        \to
        \bi{K}_{\bi{cart}}\times_{\bi{M}_{\bi{cart}}}\bi{L}_{\bi{cart}}
    \]
    is a 2-equivalence,
    where $-\times_{\bi{M}}-$ denotes the strict pullback of 2-categories,
    that is, the 2-category of pairs $(k,l)$ of 0-cells $k\in\bi{K}$ and $l\in\bi{L}$
    with $T(k)=S(l)$ in $\bi{M}$.
\end{lemma}
\begin{proof}
    finite products in $\bi{K}\times_{\bi{M}}\bi{L}$ are given by 
    pointwise finite products in $\bi{K}$ and $\bi{L}$,
    namely, $(k,l)\otimes (k',l')\coloneqq(k\otimes k',l\otimes l')$,
    and $\monunit\coloneqq(\monunit,\monunit)$.
    In addition, a 1-cell $(f,g)\colon (k,l)\to (k',l')$ has a right adjoint 
    if and only if $f$ and $g$ have right adjoints in $\bi{K}$ and $\bi{L}$, respectively.
    Here, we use the assumption that $T$ and $S$ are locally isofibrations.
    From this, we see that the canonical 2-functor is essentially surjective, 
    and locally fully faithful by the fact that natural isomorphisms in the pullback are pointwise.
\end{proof}

\section{Internal language of fibrations}\label{rem:intlang}
It is sometimes useful to have a formal language to describe what is going on in a categorical structure.
The internal language of fibrations with preordered fibers is given in \cite[Section 4.3]{jacobsCategoricalLogicType1999a} as first-order predicate logic, and that of general fibrations is a proof-relevant version of the former, which appears in  \cite[\S 2.1.2]{Her93} and \cite{Pav96}.
We present here a language in the style we will use in the following sections, but our use is restricted to regular logic.
This language is compatible with the type theory introduced in \cite{nasu2024internallogicvirtualdouble}.

The language, or the type theory, has algebraic type theory as its base.
It also has another kind called the \emph{proposition} depending on a context $\syn{\Gamma}= 
\syn{x}_1\colon\syn{I}_1,\dots,\syn{x}_n\colon\syn{I}_n$.
Given propositions $\syn{\alpha}_1,\dots,\syn{\alpha}_n$ and $\syn{\beta}$
in the same context $\syn{\Gamma}$, there is another syntactic entity 
called \emph{proof} of a Horn clause
$\syn{\alpha}_1,\dots,\syn{\alpha}_n\vdash\syn{\beta}$.
Here, we would prefer the following judgment declaration:
\[
\begin{aligned}
    & \vdash \syn{I}\ \textsf{type} \\
    \syn{\Gamma} & \vdash \syn{t}:\syn{I} \\
    \syn{\Gamma} & \vdash \syn{\alpha}\ \textsf{prop} \\
    \syn{\Gamma} \mid \syn{a}_1:\syn{\alpha}_1,\dots,\syn{a}_n:\syn{\alpha}_n & 
    \vdash \syn{\mu}:\syn{\beta}.
\end{aligned}
\]
The variables $\syn{a}_1,\dots,\syn{a}_n$ will serve as \emph{proof variables}.
If we make the variable dependency explicit, 
terms, propositions, and proofs are given by the following grammar:
\[
\syn{t}(\syn{x}_1,\dots,\syn{x}_n),\quad
\syn{\alpha}(\syn{x}_1,\dots,\syn{x}_n),\quad
\syn{\mu}(\syn{x}_1,\dots,\syn{x}_n)\{\syn{a}_1,\dots,\syn{a}_n\}.
\]
Here $\{-\}$ denotes the proof variable dependency,
but we will omit it,
or even drop the proof variables from the notation and write a proof as 
\[
    \syn{\Gamma}\vdash\syn{\alpha}_1,\dots,\syn{\alpha}_n\vdash^{\syn{\mu}}\syn{\beta}.
\]

\section[A diagrammatic proof]{A diagrammatic proof of \Cref{prop:CartdoubleCartbicat}}\label{sec:diagrammatic-proof}

\begin{proof}
We refer to the conditions presented in \Cref{def:cartbicat}.
(i) is clear from the definition of a cartesian equipment.
We prove (a) in (ii).
Since we have an obvious object $!_*$ in $\MapBi{\dbl{D}}(K,1)$ for any object $K$ in $\dbl{D}$,
we can define a functor $\bi{1}\to\MapBi{\dbl{D}}(K,1)$ by
sending the unique object of $\bi{1}$ to $!_*$.
There is a unique cell from $!_*$ to $!_*$ because of the universal property of $1$,
which implies that the functor is fully faithful.
We prove that it is essentially surjective.
Given a map $\gamma\colon K\sto 1$ in $\MapBi{\dbl{D}}(K,1)$,
then we have the unique cells $!$ from $\gamma$ and $\gamma^*$ into $\delta_1$,
as shown on the left below.
Define $\lambda$ and $\rho$ as follows.
\[
\begin{tangle}{(2,2)}
	\tgBlank{(0,0)}{\tgColour5}
	\tgBorderA{(1,0)}{\tgColour5}{\tgColour3}{\tgColour3}{\tgColour5}
	\tgBorderA{(0,1)}{\tgColour5}{\tgColour5}{\tgColour3}{\tgColour3}
	\tgBorderA{(1,1)}{\tgColour5}{\tgColour3}{\tgColour3}{\tgColour3}
	\tgCell{(1,1)}{!}
	\tgAxisLabel{(1.5,0)}{south}{\gamma}
	\tgAxisLabel{(0,1.5)}{east}{!
    }
    \node at (0.5,0.5) {$K$};
\end{tangle}
\begin{tangle}{(2,2)}
	\tgBorderA{(0,0)}{\tgColour3}{\tgColour5}{\tgColour5}{\tgColour3}
	\tgBlank{(1,0)}{\tgColour5}
	\tgBorderA{(0,1)}{\tgColour3}{\tgColour5}{\tgColour3}{\tgColour3}
	\tgBorderA{(1,1)}{\tgColour5}{\tgColour5}{\tgColour3}{\tgColour3}
	\tgCell{(0,1)}{!}
	\tgAxisLabel{(0.5,0)}{south}{\gamma}
	\tgAxisLabel{(2,1.5)}{west}{!}
    \node at (1.5,0.5) {$K$};
\end{tangle}
		\hspace{2em}
		\lambda\coloneqq
		\hspace{1em}
\begin{tangle}{(3,2)}
	\tgBorderA{(0,0)}{\tgColour5}{\tgColour5}{\tgColour3}{\tgColour5}
	\tgBorderA{(1,0)}{\tgColour5}{\tgColour5}{\tgColour5}{\tgColour3}
	\tgBorderA{(2,0)}{\tgColour5}{\tgColour3}{\tgColour3}{\tgColour5}
	\tgBorderA{(0,1)}{\tgColour5}{\tgColour3}{\tgColour3}{\tgColour5}
	\tgBorderA{(1,1)}{\tgColour3}{\tgColour5}{\tgColour3}{\tgColour3}
	\tgBorderA{(2,1)}{\tgColour5}{\tgColour3}{\tgColour3}{\tgColour3}
	\tgCell[(1,0)]{(0.5,0)}{\gamma}
	\tgCell{(1,1)}{!}
	\tgArrow{(1.5,1)}{0}
	\tgArrow{(2,0.5)}{1}
	\tgAxisLabel{(2.5,0)}{south}{!_*}
	\tgAxisLabel{(0.5,2)}{north}{\gamma}
\end{tangle},
		\hspace{1em}
		\rho\coloneqq
\begin{tangle}{(2,2)}
	\tgBorderA{(0,0)}{\tgColour5}{\tgColour5}{\tgColour3}{\tgColour5}
	\tgBorderA{(1,0)}{\tgColour5}{\tgColour3}{\tgColour3}{\tgColour3}
	\tgBorderA{(0,1)}{\tgColour5}{\tgColour3}{\tgColour3}{\tgColour5}
	\tgBlank{(1,1)}{\tgColour3}
	\tgCell{(1,0)}{!}
	\tgArrow{(0.5,0)}{0}
	\tgArrow{(0,0.5)}{1}
\end{tangle}
\]
Then we can deduce that they are mutually inverse by the following calculation. 
\[
\begin{tangle}{(4,2)}
	\tgBorderA{(0,0)}{\tgColour5}{\tgColour5}{\tgColour3}{\tgColour5}
	\tgBorderA{(1,0)}{\tgColour5}{\tgColour5}{\tgColour5}{\tgColour3}
	\tgBorderA{(2,0)}{\tgColour5}{\tgColour5}{\tgColour3}{\tgColour5}
	\tgBorderA{(3,0)}{\tgColour5}{\tgColour3}{\tgColour3}{\tgColour3}
	\tgBorderA{(0,1)}{\tgColour5}{\tgColour3}{\tgColour3}{\tgColour5}
	\tgBorderA{(1,1)}{\tgColour3}{\tgColour5}{\tgColour3}{\tgColour3}
	\tgBorderA{(2,1)}{\tgColour5}{\tgColour3}{\tgColour3}{\tgColour3}
	\tgBlank{(3,1)}{\tgColour3}
	\tgCell[(1,0)]{(0.5,0)}{\gamma}
	\tgCell{(1,1)}{!}
	\tgCell{(3,0)}{!}
	\tgArrow{(1.5,1)}{0}
	\tgArrow{(2,0.5)}{1}
	\tgArrow{(2.5,0)}{0}
\end{tangle}
=
\begin{tangle}{(3,2)}
	\tgBorderA{(0,0)}{\tgColour5}{\tgColour5}{\tgColour3}{\tgColour5}
	\tgBorderA{(1,0)}{\tgColour5}{\tgColour5}{\tgColour5}{\tgColour3}
	\tgBorderA{(2,0)}{\tgColour5}{\tgColour3}{\tgColour3}{\tgColour5}
	\tgBorderA{(0,1)}{\tgColour5}{\tgColour3}{\tgColour3}{\tgColour5}
	\tgBorderA{(1,1)}{\tgColour3}{\tgColour5}{\tgColour3}{\tgColour3}
	\tgBorderA{(2,1)}{\tgColour5}{\tgColour3}{\tgColour3}{\tgColour3}
	\tgCell[(1,0)]{(0.5,0)}{\gamma}
	\tgCell{(1,1)}{!}
	\tgCell{(2,1)}{!}
	\tgArrow{(1.5,1)}{0}
\end{tangle}
=
\begin{tangle}{(3,2)}
	\tgBorderA{(0,0)}{\tgColour5}{\tgColour5}{\tgColour3}{\tgColour5}
	\tgBorderA{(1,0)}{\tgColour5}{\tgColour5}{\tgColour5}{\tgColour3}
	\tgBorderA{(2,0)}{\tgColour5}{\tgColour3}{\tgColour3}{\tgColour5}
	\tgBorderA{(0,1)}{\tgColour5}{\tgColour3}{\tgColour3}{\tgColour5}
	\tgBorderA{(1,1)}{\tgColour3}{\tgColour5}{\tgColour3}{\tgColour3}
	\tgBorderA{(2,1)}{\tgColour5}{\tgColour3}{\tgColour3}{\tgColour3}
	\tgCell[(1,0)]{(0.5,0)}{\gamma}
	\tgCell[(1,0)]{(1.5,1)}{\gamma}
\end{tangle}
=
\begin{tangle}{(1,2)}
	\tgBorderA{(0,0)}{\tgColour5}{\tgColour3}{\tgColour3}{\tgColour5}
	\tgBorderA{(0,1)}{\tgColour5}{\tgColour3}{\tgColour3}{\tgColour5}
\end{tangle}
\]
\[
\begin{tangle}{(4,2)}
	\tgBlank{(0,0)}{\tgColour5}
	\tgBorderA{(1,0)}{\tgColour5}{\tgColour5}{\tgColour3}{\tgColour5}
	\tgBorderA{(2,0)}{\tgColour5}{\tgColour5}{\tgColour5}{\tgColour3}
	\tgBorderA{(3,0)}{\tgColour5}{\tgColour3}{\tgColour3}{\tgColour5}
	\tgBorderA{(0,1)}{\tgColour5}{\tgColour5}{\tgColour3}{\tgColour5}
	\tgBorderA{(1,1)}{\tgColour5}{\tgColour3}{\tgColour3}{\tgColour3}
	\tgBorderA{(2,1)}{\tgColour3}{\tgColour5}{\tgColour3}{\tgColour3}
	\tgBorderA{(3,1)}{\tgColour5}{\tgColour3}{\tgColour3}{\tgColour3}
	\tgCell[(1,0)]{(1.5,0)}{\gamma}
	\tgCell{(1,1)}{!}
	\tgCell{(2,1)}{!}
	\tgArrow{(0.5,1)}{0}
	\tgArrow{(2.5,1)}{0}
	\tgArrow{(3,0.5)}{1}
\end{tangle}
=
\begin{tangle}{(5,3)}
	\tgBlank{(0,0)}{\tgColour5}
	\tgBorderA{(1,0)}{\tgColour5}{\tgColour5}{\tgColour3}{\tgColour5}
	\tgBorderA{(2,0)}{\tgColour5}{\tgColour5}{\tgColour3}{\tgColour3}
	\tgBorderA{(3,0)}{\tgColour5}{\tgColour5}{\tgColour5}{\tgColour3}
	\tgBorderA{(4,0)}{\tgColour5}{\tgColour3}{\tgColour3}{\tgColour5}
	\tgBlank{(0,1)}{\tgColour5}
	\tgBorderA{(1,1)}{\tgColour5}{\tgColour3}{\tgColour5}{\tgColour5}
	\tgBorderA{(2,1)}{\tgColour3}{\tgColour3}{\tgColour5}{\tgColour5}
	\tgBorder{(2,1)}{0}{0}{1}{0}
	\tgBorderA{(3,1)}{\tgColour3}{\tgColour5}{\tgColour5}{\tgColour5}
	\tgBorderA{(4,1)}{\tgColour5}{\tgColour3}{\tgColour3}{\tgColour5}
	\tgBorderA{(0,2)}{\tgColour5}{\tgColour5}{\tgColour3}{\tgColour5}
	\tgBorderA{(1,2)}{\tgColour5}{\tgColour5}{\tgColour3}{\tgColour3}
	\tgBorderA{(2,2)}{\tgColour5}{\tgColour5}{\tgColour3}{\tgColour3}
	\tgBorder{(2,2)}{1}{0}{0}{0}
	\tgBorderA{(3,2)}{\tgColour5}{\tgColour5}{\tgColour3}{\tgColour3}
	\tgBorderA{(4,2)}{\tgColour5}{\tgColour3}{\tgColour3}{\tgColour3}
	\tgCell[(2,0)]{(2,0)}{\gamma}
	\tgCell{(2,2)}{!}
	\tgArrow{(0.5,2)}{0}
	\tgArrow{(3.5,2)}{0}
	\tgArrow{(4,0.5)}{1}
	\node[fill=black, inner xsep=3em, inner ysep=1.5pt] at (2.5,1.5) {}; 
\end{tangle}
=
\begin{tangle}{(2,2)}
	\tgBlank{(0,0)}{\tgColour5}
	\tgBorderA{(1,0)}{\tgColour5}{\tgColour3}{\tgColour3}{\tgColour5}
	\tgBorderA{(0,1)}{\tgColour5}{\tgColour5}{\tgColour3}{\tgColour5}
	\tgBorderA{(1,1)}{\tgColour5}{\tgColour3}{\tgColour3}{\tgColour3}
	\tgArrow{(0.5,1)}{0}
	\tgArrow{(1,0.5)}{1}
\end{tangle}
=
\begin{tangle}{(1,2)}
	\tgBorderA{(0,0)}{\tgColour5}{\tgColour3}{\tgColour3}{\tgColour5}
	\tgBorderA{(0,1)}{\tgColour5}{\tgColour3}{\tgColour3}{\tgColour5}
	\tgArrow{(0,0.5)}{1}
\end{tangle}
\]
This implies that $\gamma$ is isomorphic to $!_*$.
Here, nodes only with two vertical strings on one side,
like the ones labeled by $\gamma$ above,
represent the units and counits of some adjunctions.

Next, we prove (b) in (ii).
For a triple of objects $I,J,K$ in $\dbl{D}$,
we define two functors $\Phi$ and $\Psi$ as follows.
\[
	\begin{tikzcd}
		\MapBi{\dbl{D}}(K,I)\times\MapBi{\dbl{D}}(K,J)
		\ar[r,"\Phi",shift left=1ex]
		&
		\MapBi{\dbl{D}}(K,I\times J)
		\ar[l,"\Psi",shift left=1ex]
	\end{tikzcd},
	\hspace{1em}
	\begin{aligned}
		(\alpha,\beta)
		&\to<mapsto, "\Phi", > \tpl{0,0}_*(\alpha\times\beta),
		\\
		\gamma
		&\to<mapsto, "\Psi", > (\tpl{0}_*(\gamma),\tpl{1}_*(\gamma)).
	\end{aligned}
\]
Note that composites of left adjoints are also left adjoints.
We prove that $\Phi$ and $\Psi$ give an equivalence of categories.
For one direction, we prove that $\Psi\circ\Phi$ is naturally isomorphic to the identity. 
We only show this only on the first component of the product.

Suppose we are given a pair of maps
$\alpha\colon K\sto I$ and $\beta\colon K\sto J$ in $\dbl{D}$.
Let $\lambda$ and $\rho$ be defined as follows\footnote{
Since we run out of letters, we reuse $\lambda$ and $\rho$ for different cells.
}.
\[
\lambda\coloneqq
\hspace{1em}
\begin{tangle}{(6,3)}
	\tgBlank{(0,0)}{\tgColour5}
	\tgBorderA{(1,0)}{\tgColour5}{\tgColour4}{\tgColour4}{\tgColour5}
	\tgBorderA{(2,0)}{\tgColour4}{\tgColour2}{\tgColour2}{\tgColour4}
	\tgBlank{(3,0)}{\tgColour2}
	\tgBorderA{(4,0)}{\tgColour2}{\tgColour0}{\tgColour0}{\tgColour2}
	\tgBlank{(5,0)}{\tgColour0}
	\tgBorderA{(0,1)}{\tgColour5}{\tgColour5}{\tgColour4}{\tgColour5}
	\tgBorderA{(1,1)}{\tgColour5}{\tgColour4}{\tgColour4}{\tgColour4}
	\tgBorderA{(2,1)}{\tgColour4}{\tgColour2}{\tgColour2}{\tgColour4}
	\tgBlank{(3,1)}{\tgColour2}
	\tgBorderA{(4,1)}{\tgColour2}{\tgColour0}{\tgColour0}{\tgColour2}
	\tgBlank{(5,1)}{\tgColour0}
	\tgBorderA{(0,2)}{\tgColour5}{\tgColour4}{\tgColour5}{\tgColour5}
	\tgBorderA{(1,2)}{\tgColour4}{\tgColour4}{\tgColour5}{\tgColour5}
	\tgBorderA{(2,2)}{\tgColour4}{\tgColour2}{\tgColour0}{\tgColour5}
	\tgBorderA{(3,2)}{\tgColour2}{\tgColour2}{\tgColour0}{\tgColour0}
	\tgBorderA{(4,2)}{\tgColour2}{\tgColour0}{\tgColour0}{\tgColour0}
	\tgBlank{(5,2)}{\tgColour0}
	\tgCell{(2,2)}{\scriptstyle\tpl{0}}
	\tgAxisLabel{(1.5,0)}{south}{\tpl{0,0}_*}
	\tgAxisLabel{(2.5,0)}{south}{\alpha\times\beta}
	\tgAxisLabel{(4.5,0)}{south}{\tpl{0}_*}
	\tgAxisLabel{(2.5,3)}{north}{\alpha}
    \node at (0.5,0.5) {$K$};
    \node at (1.5,2.0) {$K\times K$};
    \node at (3.5,2.0) {$I\times J$};
    \node at (5.0,2.0) {$I$};
    \node[fill=black, inner xsep=1.5pt, inner ysep=1.5em] at (0.5,2.0) {}; 
\end{tangle}
\hspace{1em}
,
\hspace{1em}
\rho\coloneqq
\hspace{1em}
\begin{tangle}{(6,4)}
	\tgBorderA{(0,0)}{\tgColour5}{\tgColour5}{\tgColour4}{\tgColour5}
	\tgBorderA{(1,0)}{\tgColour5}{\tgColour5}{\tgColour4}{\tgColour4}
	\tgBorderA{(2,0)}{\tgColour5}{\tgColour5}{\tgColour4}{\tgColour4}
	\tgBorderA{(3,0)}{\tgColour5}{\tgColour5}{\tgColour4}{\tgColour4}
	\tgBorderA{(4,0)}{\tgColour5}{\tgColour5}{\tgColour5}{\tgColour4}
	\tgBorderA{(5,0)}{\tgColour5}{\tgColour0}{\tgColour0}{\tgColour5}
	\tgBorderA{(0,1)}{\tgColour5}{\tgColour4}{\tgColour4}{\tgColour5}
	\tgBorderA{(1,1)}{\tgColour4}{\tgColour4}{\tgColour2}{\tgColour4}
	\tgBorderA{(2,1)}{\tgColour4}{\tgColour4}{\tgColour2}{\tgColour2}
	\tgBorderA{(3,1)}{\tgColour4}{\tgColour4}{\tgColour4}{\tgColour2}
	\tgBorderA{(4,1)}{\tgColour4}{\tgColour5}{\tgColour5}{\tgColour4}
	\tgBorderA{(5,1)}{\tgColour5}{\tgColour0}{\tgColour0}{\tgColour5}
	\tgBorderA{(0,2)}{\tgColour5}{\tgColour4}{\tgColour4}{\tgColour5}
	\tgBorderA{(1,2)}{\tgColour4}{\tgColour2}{\tgColour2}{\tgColour4}
	\tgBorderA{(2,2)}{\tgColour2}{\tgColour2}{\tgColour0}{\tgColour2}
	\tgBorderA{(3,2)}{\tgColour2}{\tgColour4}{\tgColour5}{\tgColour0}
	\tgBorderA{(4,2)}{\tgColour4}{\tgColour5}{\tgColour5}{\tgColour5}
	\tgBorderA{(5,2)}{\tgColour5}{\tgColour0}{\tgColour0}{\tgColour5}
	\tgBorderA{(0,3)}{\tgColour5}{\tgColour4}{\tgColour4}{\tgColour5}
	\tgBorderA{(1,3)}{\tgColour4}{\tgColour2}{\tgColour2}{\tgColour4}
	\tgBorderA{(2,3)}{\tgColour2}{\tgColour0}{\tgColour0}{\tgColour2}
	\tgBorderA{(3,3)}{\tgColour0}{\tgColour5}{\tgColour0}{\tgColour0}
	\tgBorderA{(4,3)}{\tgColour5}{\tgColour5}{\tgColour0}{\tgColour0}
	\tgBorderA{(5,3)}{\tgColour5}{\tgColour0}{\tgColour0}{\tgColour0}
	\tgCell[(2,0)]{(2,1)}{\alpha\times\beta}
	\tgCell[(2,0)]{(4,3)}{\alpha}
	\tgCell{(3,2)}{\scriptstyle\tpl{0}}
    \node[fill=black, inner xsep=1.5pt, inner ysep=3em] at (4.5,1.5) {};
\end{tangle}
\]
Then, $\lambda$ and $\rho$ are mutually inverse
because of the following calculation.
\[
\begin{tangle}{(9,4)}
	\tgBorderA{(0,0)}{\tgColour5}{\tgColour5}{\tgColour4}{\tgColour5}
	\tgBorderA{(1,0)}{\tgColour5}{\tgColour5}{\tgColour4}{\tgColour4}
	\tgBorderA{(2,0)}{\tgColour5}{\tgColour5}{\tgColour4}{\tgColour4}
	\tgBorderA{(3,0)}{\tgColour5}{\tgColour5}{\tgColour4}{\tgColour4}
	\tgBorderA{(4,0)}{\tgColour5}{\tgColour5}{\tgColour5}{\tgColour4}
	\tgBorderA{(5,0)}{\tgColour5}{\tgColour5}{\tgColour4}{\tgColour5}
	\tgBorderA{(6,0)}{\tgColour5}{\tgColour4}{\tgColour4}{\tgColour4}
	\tgBorderA{(7,0)}{\tgColour4}{\tgColour2}{\tgColour2}{\tgColour4}
	\tgBorderA{(8,0)}{\tgColour2}{\tgColour0}{\tgColour0}{\tgColour2}
	\tgBorderA{(0,1)}{\tgColour5}{\tgColour4}{\tgColour4}{\tgColour5}
	\tgBorderA{(1,1)}{\tgColour4}{\tgColour4}{\tgColour2}{\tgColour4}
	\tgBorderA{(2,1)}{\tgColour4}{\tgColour4}{\tgColour2}{\tgColour2}
	\tgBorderA{(3,1)}{\tgColour4}{\tgColour4}{\tgColour4}{\tgColour2}
	\tgBorderA{(4,1)}{\tgColour4}{\tgColour5}{\tgColour5}{\tgColour4}
	\tgBorderA{(5,1)}{\tgColour5}{\tgColour4}{\tgColour4}{\tgColour5}
	\tgBlank{(6,1)}{\tgColour4}
	\tgBorderA{(7,1)}{\tgColour4}{\tgColour2}{\tgColour2}{\tgColour4}
	\tgBorderA{(8,1)}{\tgColour2}{\tgColour0}{\tgColour0}{\tgColour2}
	\tgBorderA{(0,2)}{\tgColour5}{\tgColour4}{\tgColour4}{\tgColour5}
	\tgBorderA{(1,2)}{\tgColour4}{\tgColour2}{\tgColour2}{\tgColour4}
	\tgBorderA{(2,2)}{\tgColour2}{\tgColour2}{\tgColour0}{\tgColour2}
	\tgBorderA{(3,2)}{\tgColour2}{\tgColour4}{\tgColour5}{\tgColour0}
	\tgBorderA{(4,2)}{\tgColour4}{\tgColour5}{\tgColour5}{\tgColour5}
	\tgBorderA{(5,2)}{\tgColour5}{\tgColour4}{\tgColour5}{\tgColour5}
	\tgBorderA{(6,2)}{\tgColour4}{\tgColour4}{\tgColour5}{\tgColour5}
	\tgBorderA{(7,2)}{\tgColour4}{\tgColour2}{\tgColour0}{\tgColour5}
	\tgBorderA{(8,2)}{\tgColour2}{\tgColour0}{\tgColour0}{\tgColour0}
	\tgBorderA{(0,3)}{\tgColour5}{\tgColour4}{\tgColour4}{\tgColour5}
	\tgBorderA{(1,3)}{\tgColour4}{\tgColour2}{\tgColour2}{\tgColour4}
	\tgBorderA{(2,3)}{\tgColour2}{\tgColour0}{\tgColour0}{\tgColour2}
	\tgBorderA{(3,3)}{\tgColour0}{\tgColour5}{\tgColour0}{\tgColour0}
	\tgBorderA{(4,3)}{\tgColour5}{\tgColour5}{\tgColour0}{\tgColour0}
	\tgBorderA{(5,3)}{\tgColour5}{\tgColour5}{\tgColour0}{\tgColour0}
	\tgBorderA{(6,3)}{\tgColour5}{\tgColour5}{\tgColour0}{\tgColour0}
	\tgBorderA{(7,3)}{\tgColour5}{\tgColour0}{\tgColour0}{\tgColour0}
	\tgBlank{(8,3)}{\tgColour0}
	\tgCell[(2,0)]{(2,1)}{\alpha\times\beta}
	\tgCell{(3,2)}{\scriptstyle\tpl{0}}
	\tgCell[(4,0)]{(5,3)}{\alpha}
	\tgCell{(7,2)}{\scriptstyle\tpl{0}}
    \node[fill=black, inner xsep=1.5pt, inner ysep=3em] at (4.5,1.5) {};
    \node[fill=black, inner xsep=1.5pt, inner ysep=3em] at (5.5,1.5) {};
\end{tangle}=
\begin{tangle}{(7,4)}
	\tgBorderA{(0,0)}{\tgColour5}{\tgColour5}{\tgColour4}{\tgColour5}
	\tgBorderA{(1,0)}{\tgColour5}{\tgColour5}{\tgColour4}{\tgColour4}
	\tgBorderA{(2,0)}{\tgColour5}{\tgColour5}{\tgColour4}{\tgColour4}
	\tgBorderA{(3,0)}{\tgColour5}{\tgColour5}{\tgColour4}{\tgColour4}
	\tgBorderA{(4,0)}{\tgColour5}{\tgColour4}{\tgColour4}{\tgColour4}
	\tgBorderA{(5,0)}{\tgColour4}{\tgColour2}{\tgColour2}{\tgColour4}
	\tgBorderA{(6,0)}{\tgColour2}{\tgColour0}{\tgColour0}{\tgColour2}
	\tgBorderA{(0,1)}{\tgColour5}{\tgColour4}{\tgColour4}{\tgColour5}
	\tgBorderA{(1,1)}{\tgColour4}{\tgColour4}{\tgColour2}{\tgColour4}
	\tgBorderA{(2,1)}{\tgColour4}{\tgColour4}{\tgColour2}{\tgColour2}
	\tgBorderA{(3,1)}{\tgColour4}{\tgColour4}{\tgColour4}{\tgColour2}
	\tgBlank{(4,1)}{\tgColour4}
	\tgBorderA{(5,1)}{\tgColour4}{\tgColour2}{\tgColour2}{\tgColour4}
	\tgBorderA{(6,1)}{\tgColour2}{\tgColour0}{\tgColour0}{\tgColour2}
	\tgBorderA{(0,2)}{\tgColour5}{\tgColour4}{\tgColour4}{\tgColour5}
	\tgBorderA{(1,2)}{\tgColour4}{\tgColour2}{\tgColour2}{\tgColour4}
	\tgBorderA{(2,2)}{\tgColour2}{\tgColour2}{\tgColour0}{\tgColour2}
	\tgBorderA{(3,2)}{\tgColour2}{\tgColour4}{\tgColour5}{\tgColour0}
	\tgBorderA{(4,2)}{\tgColour4}{\tgColour4}{\tgColour5}{\tgColour5}
	\tgBorderA{(5,2)}{\tgColour4}{\tgColour2}{\tgColour0}{\tgColour5}
	\tgBorderA{(6,2)}{\tgColour2}{\tgColour0}{\tgColour0}{\tgColour0}
	\tgBorderA{(0,3)}{\tgColour5}{\tgColour4}{\tgColour4}{\tgColour5}
	\tgBorderA{(1,3)}{\tgColour4}{\tgColour2}{\tgColour2}{\tgColour4}
	\tgBorderA{(2,3)}{\tgColour2}{\tgColour0}{\tgColour0}{\tgColour2}
	\tgBorderA{(3,3)}{\tgColour0}{\tgColour5}{\tgColour0}{\tgColour0}
	\tgBorderA{(4,3)}{\tgColour5}{\tgColour5}{\tgColour0}{\tgColour0}
	\tgBorderA{(5,3)}{\tgColour5}{\tgColour0}{\tgColour0}{\tgColour0}
	\tgBlank{(6,3)}{\tgColour0}
	\tgCell[(2,0)]{(2,1)}{\alpha\times\beta}
	\tgCell{(3,2)}{\scriptstyle\tpl{0}}
	\tgCell{(5,2)}{\scriptstyle\tpl{0}}
	\tgCell[(2,0)]{(4,3)}{\alpha}
\end{tangle}
\]
\[
=
\begin{tangle}{(7,4)}
	\tgBorderA{(0,0)}{\tgColour5}{\tgColour5}{\tgColour4}{\tgColour5}
	\tgBorderA{(1,0)}{\tgColour5}{\tgColour5}{\tgColour4}{\tgColour4}
	\tgBorderA{(2,0)}{\tgColour5}{\tgColour5}{\tgColour4}{\tgColour4}
	\tgBorderA{(3,0)}{\tgColour5}{\tgColour5}{\tgColour4}{\tgColour4}
	\tgBorderA{(4,0)}{\tgColour5}{\tgColour4}{\tgColour4}{\tgColour4}
	\tgBorderA{(5,0)}{\tgColour4}{\tgColour2}{\tgColour2}{\tgColour4}
	\tgBorderA{(6,0)}{\tgColour2}{\tgColour0}{\tgColour0}{\tgColour2}
	\tgBorderA{(0,1)}{\tgColour5}{\tgColour4}{\tgColour4}{\tgColour5}
	\tgBorderA{(1,1)}{\tgColour4}{\tgColour4}{\tgColour2}{\tgColour4}
	\tgBorderA{(2,1)}{\tgColour4}{\tgColour4}{\tgColour2}{\tgColour2}
	\tgBorderA{(3,1)}{\tgColour4}{\tgColour4}{\tgColour4}{\tgColour2}
	\tgBlank{(4,1)}{\tgColour4}
	\tgBorderA{(5,1)}{\tgColour4}{\tgColour2}{\tgColour2}{\tgColour4}
	\tgBorderA{(6,1)}{\tgColour2}{\tgColour0}{\tgColour0}{\tgColour2}
	\tgBorderA{(0,2)}{\tgColour5}{\tgColour4}{\tgColour4}{\tgColour5}
	\tgBorderA{(1,2)}{\tgColour4}{\tgColour2}{\tgColour2}{\tgColour4}
	\tgBlank{(2,2)}{\tgColour2}
	\tgBorderA{(3,2)}{\tgColour2}{\tgColour4}{\tgColour2}{\tgColour2}
	\tgBorderA{(4,2)}{\tgColour4}{\tgColour4}{\tgColour2}{\tgColour2}
	\tgBorderA{(5,2)}{\tgColour4}{\tgColour2}{\tgColour2}{\tgColour2}
	\tgBorderA{(6,2)}{\tgColour2}{\tgColour0}{\tgColour0}{\tgColour2}
	\tgBorderA{(0,3)}{\tgColour5}{\tgColour4}{\tgColour4}{\tgColour5}
	\tgBorderA{(1,3)}{\tgColour4}{\tgColour2}{\tgColour2}{\tgColour4}
	\tgBorderA{(2,3)}{\tgColour2}{\tgColour2}{\tgColour0}{\tgColour2}
	\tgBorderA{(3,3)}{\tgColour2}{\tgColour2}{\tgColour0}{\tgColour0}
	\tgBorderA{(4,3)}{\tgColour2}{\tgColour2}{\tgColour0}{\tgColour0}
	\tgBorderA{(5,3)}{\tgColour2}{\tgColour2}{\tgColour0}{\tgColour0}
	\tgBorderA{(6,3)}{\tgColour2}{\tgColour0}{\tgColour0}{\tgColour0}
	\tgCell[(2,0)]{(2,1)}{\alpha\times\beta}
	\tgCell[(2,0)]{(4,2)}{\alpha\times\beta}
\end{tangle}
=
\begin{tangle}{(3,4)}
	\tgBorderA{(0,0)}{\tgColour5}{\tgColour4}{\tgColour4}{\tgColour5}
	\tgBorderA{(1,0)}{\tgColour4}{\tgColour2}{\tgColour2}{\tgColour4}
	\tgBorderA{(2,0)}{\tgColour2}{\tgColour0}{\tgColour0}{\tgColour2}
	\tgBorderA{(0,1)}{\tgColour5}{\tgColour4}{\tgColour4}{\tgColour5}
	\tgBorderA{(1,1)}{\tgColour4}{\tgColour2}{\tgColour2}{\tgColour4}
	\tgBorderA{(2,1)}{\tgColour2}{\tgColour0}{\tgColour0}{\tgColour2}
	\tgBorderA{(0,2)}{\tgColour5}{\tgColour4}{\tgColour4}{\tgColour5}
	\tgBorderA{(1,2)}{\tgColour4}{\tgColour2}{\tgColour2}{\tgColour4}
	\tgBorderA{(2,2)}{\tgColour2}{\tgColour0}{\tgColour0}{\tgColour2}
	\tgBorderA{(0,3)}{\tgColour5}{\tgColour4}{\tgColour4}{\tgColour5}
	\tgBorderA{(1,3)}{\tgColour4}{\tgColour2}{\tgColour2}{\tgColour4}
	\tgBorderA{(2,3)}{\tgColour2}{\tgColour0}{\tgColour0}{\tgColour2}
\end{tangle},
\]
\[
\begin{tangle}{(7,4)}
	\tgBlank{(0,0)}{\tgColour5}
	\tgBorderA{(1,0)}{\tgColour5}{\tgColour5}{\tgColour4}{\tgColour5}
	\tgBorderA{(2,0)}{\tgColour5}{\tgColour5}{\tgColour4}{\tgColour4}
	\tgBorderA{(3,0)}{\tgColour5}{\tgColour5}{\tgColour4}{\tgColour4}
	\tgBorderA{(4,0)}{\tgColour5}{\tgColour5}{\tgColour4}{\tgColour4}
	\tgBorderA{(5,0)}{\tgColour5}{\tgColour5}{\tgColour5}{\tgColour4}
	\tgBorderA{(6,0)}{\tgColour5}{\tgColour0}{\tgColour0}{\tgColour5}
	\tgBorderA{(0,1)}{\tgColour5}{\tgColour5}{\tgColour4}{\tgColour5}
	\tgBorderA{(1,1)}{\tgColour5}{\tgColour4}{\tgColour4}{\tgColour4}
	\tgBorderA{(2,1)}{\tgColour4}{\tgColour4}{\tgColour2}{\tgColour4}
	\tgBorderA{(3,1)}{\tgColour4}{\tgColour4}{\tgColour2}{\tgColour2}
	\tgBorderA{(4,1)}{\tgColour4}{\tgColour4}{\tgColour4}{\tgColour2}
	\tgBorderA{(5,1)}{\tgColour4}{\tgColour5}{\tgColour5}{\tgColour4}
	\tgBorderA{(6,1)}{\tgColour5}{\tgColour0}{\tgColour0}{\tgColour5}
	\tgBorderA{(0,2)}{\tgColour5}{\tgColour4}{\tgColour4}{\tgColour5}
	\tgBlank{(1,2)}{\tgColour4}
	\tgBorderA{(2,2)}{\tgColour4}{\tgColour2}{\tgColour2}{\tgColour4}
	\tgBorderA{(3,2)}{\tgColour2}{\tgColour2}{\tgColour0}{\tgColour2}
	\tgBorderA{(4,2)}{\tgColour2}{\tgColour4}{\tgColour5}{\tgColour0}
	\tgBorderA{(5,2)}{\tgColour4}{\tgColour5}{\tgColour5}{\tgColour5}
	\tgBorderA{(6,2)}{\tgColour5}{\tgColour0}{\tgColour0}{\tgColour5}
	\tgBorderA{(0,3)}{\tgColour5}{\tgColour4}{\tgColour5}{\tgColour5}
	\tgBorderA{(1,3)}{\tgColour4}{\tgColour4}{\tgColour5}{\tgColour5}
	\tgBorderA{(2,3)}{\tgColour4}{\tgColour2}{\tgColour0}{\tgColour5}
	\tgBorderA{(3,3)}{\tgColour2}{\tgColour0}{\tgColour0}{\tgColour0}
	\tgBorderA{(4,3)}{\tgColour0}{\tgColour5}{\tgColour0}{\tgColour0}
	\tgBorderA{(5,3)}{\tgColour5}{\tgColour5}{\tgColour0}{\tgColour0}
	\tgBorderA{(6,3)}{\tgColour5}{\tgColour0}{\tgColour0}{\tgColour0}
	\tgCell[(2,0)]{(3,1)}{\alpha\times\beta}
	\tgCell[(2,0)]{(5,3)}{\alpha}
	\tgCell{(4,2)}{\scriptstyle\tpl{0}}
	\tgCell{(2,3)}{\scriptstyle\tpl{0}}
    \node [fill=black, inner xsep=1.5pt, inner ysep=3em] at (0.5,2.5) {};
    \node [fill=black, inner xsep=1.5pt, inner ysep=3em] at (5.5,1.5) {}; 
\end{tangle}
=
\begin{tangle}{(7,4)}
	\tgBorderA{(0,0)}{\tgColour5}{\tgColour5}{\tgColour4}{\tgColour5}
	\tgBorderA{(1,0)}{\tgColour5}{\tgColour5}{\tgColour4}{\tgColour4}
	\tgBorderA{(2,0)}{\tgColour5}{\tgColour5}{\tgColour4}{\tgColour4}
	\tgBorderA{(3,0)}{\tgColour5}{\tgColour5}{\tgColour4}{\tgColour4}
	\tgBorderA{(4,0)}{\tgColour5}{\tgColour5}{\tgColour4}{\tgColour4}
	\tgBorderA{(5,0)}{\tgColour5}{\tgColour5}{\tgColour5}{\tgColour4}
	\tgBorderA{(6,0)}{\tgColour5}{\tgColour0}{\tgColour0}{\tgColour5}
	\tgBorderA{(0,1)}{\tgColour5}{\tgColour4}{\tgColour4}{\tgColour5}
	\tgBlank{(1,1)}{\tgColour4}
	\tgBorderA{(2,1)}{\tgColour4}{\tgColour4}{\tgColour2}{\tgColour4}
	\tgBorderA{(3,1)}{\tgColour4}{\tgColour4}{\tgColour2}{\tgColour2}
	\tgBorderA{(4,1)}{\tgColour4}{\tgColour4}{\tgColour4}{\tgColour2}
	\tgBorderA{(5,1)}{\tgColour4}{\tgColour5}{\tgColour5}{\tgColour4}
	\tgBorderA{(6,1)}{\tgColour5}{\tgColour0}{\tgColour0}{\tgColour5}
	\tgBorderA{(0,2)}{\tgColour5}{\tgColour4}{\tgColour5}{\tgColour5}
	\tgBorderA{(1,2)}{\tgColour4}{\tgColour4}{\tgColour5}{\tgColour5}
	\tgBorderA{(2,2)}{\tgColour4}{\tgColour2}{\tgColour0}{\tgColour5}
	\tgBorderA{(3,2)}{\tgColour2}{\tgColour2}{\tgColour0}{\tgColour0}
	\tgBorderA{(4,2)}{\tgColour2}{\tgColour4}{\tgColour5}{\tgColour0}
	\tgBorderA{(5,2)}{\tgColour4}{\tgColour5}{\tgColour5}{\tgColour5}
	\tgBorderA{(6,2)}{\tgColour5}{\tgColour0}{\tgColour0}{\tgColour5}
	\tgBlank{(0,3)}{\tgColour5}
	\tgBlank{(1,3)}{\tgColour5}
	\tgBorderA{(2,3)}{\tgColour5}{\tgColour0}{\tgColour0}{\tgColour5}
	\tgBlank{(3,3)}{\tgColour0}
	\tgBorderA{(4,3)}{\tgColour0}{\tgColour5}{\tgColour0}{\tgColour0}
	\tgBorderA{(5,3)}{\tgColour5}{\tgColour5}{\tgColour0}{\tgColour0}
	\tgBorderA{(6,3)}{\tgColour5}{\tgColour0}{\tgColour0}{\tgColour0}
	\tgCell[(2,0)]{(3,1)}{\alpha\times\beta}
	\tgCell[(2,0)]{(5,3)}{\alpha}
	\tgCell{(4,2)}{\scriptstyle\tpl{0}}
	\tgCell{(2,2)}{\scriptstyle\tpl{0}}
    \node [fill=black, inner xsep=1.5pt, inner ysep=3em] at (0.5,1.5) {};
    \node [fill=black, inner xsep=1.5pt, inner ysep=3em] at (5.5,1.5) {};
\end{tangle}
\]
\[
=
\begin{tangle}{(5,4)}
	\tgBorderA{(0,0)}{\tgColour5}{\tgColour5}{\tgColour4}{\tgColour5}
	\tgBorderA{(1,0)}{\tgColour5}{\tgColour5}{\tgColour4}{\tgColour4}
	\tgBorderA{(2,0)}{\tgColour5}{\tgColour5}{\tgColour4}{\tgColour4}
	\tgBorderA{(3,0)}{\tgColour5}{\tgColour5}{\tgColour5}{\tgColour4}
	\tgBorderA{(4,0)}{\tgColour5}{\tgColour0}{\tgColour0}{\tgColour5}
	\tgBorderA{(0,1)}{\tgColour5}{\tgColour4}{\tgColour5}{\tgColour5}
	\tgBorderA{(1,1)}{\tgColour4}{\tgColour4}{\tgColour5}{\tgColour5}
	\tgBorderA{(2,1)}{\tgColour4}{\tgColour4}{\tgColour5}{\tgColour5}
	\tgBorderA{(3,1)}{\tgColour4}{\tgColour5}{\tgColour5}{\tgColour5}
	\tgBorderA{(4,1)}{\tgColour5}{\tgColour0}{\tgColour0}{\tgColour5}
	\tgBorderA{(0,2)}{\tgColour5}{\tgColour5}{\tgColour0}{\tgColour5}
	\tgBorderA{(1,2)}{\tgColour5}{\tgColour5}{\tgColour0}{\tgColour0}
	\tgBorderA{(2,2)}{\tgColour5}{\tgColour5}{\tgColour5}{\tgColour0}
	\tgBlank{(3,2)}{\tgColour5}
	\tgBorderA{(4,2)}{\tgColour5}{\tgColour0}{\tgColour0}{\tgColour5}
	\tgBorderA{(0,3)}{\tgColour5}{\tgColour0}{\tgColour0}{\tgColour5}
	\tgBlank{(1,3)}{\tgColour0}
	\tgBorderA{(2,3)}{\tgColour0}{\tgColour5}{\tgColour0}{\tgColour0}
	\tgBorderA{(3,3)}{\tgColour5}{\tgColour5}{\tgColour0}{\tgColour0}
	\tgBorderA{(4,3)}{\tgColour5}{\tgColour0}{\tgColour0}{\tgColour0}
	\tgCell[(2,0)]{(3,3)}{\alpha}
	\tgCell[(2,0)]{(1,2)}{\alpha}
    \node [fill=black, inner xsep=1.5pt, inner ysep=1.5em] at (0.5,1.0) {};
    \node [fill=black, inner xsep=1.5pt, inner ysep=1.5em] at (3.5,1.0) {};
\end{tangle}
=
\begin{tangle}{(1,4)}
	\tgBorderA{(0,0)}{\tgColour5}{\tgColour0}{\tgColour0}{\tgColour5}
	\tgBorderA{(0,1)}{\tgColour5}{\tgColour0}{\tgColour0}{\tgColour5}
	\tgBorderA{(0,2)}{\tgColour5}{\tgColour0}{\tgColour0}{\tgColour5}
	\tgBorderA{(0,3)}{\tgColour5}{\tgColour0}{\tgColour0}{\tgColour5}
\end{tangle}
\]

By this and the dual argument, we have that
$\Psi(\Phi(\alpha,\beta))$ is isomorphic to $(\alpha,\beta)$.

For the other direction, we prove that $\Phi\circ\Psi$ is naturally isomorphic to the identity. 
Suppose we are given a map $\gamma\colon K\sto I\times J$ in $\dbl{D}$.
Let $\lambda_i$ and $\rho_i$ be defined as follows for $i=0,1$,
where $\gamma_0\coloneqq\gamma\tpl{0}_*$ and $\gamma_1\coloneqq\gamma\tpl{1}_*$,
and thus, $\gamma_i^*$ can be taken as $\tpl{i}^*\gamma^*$ for $i=0,1$.
\[
\lambda_0\coloneqq
\begin{tangle}{(3,2)}
	\tgBorderA{(0,0)}{\tgColour5}{\tgColour2}{\tgColour2}{\tgColour5}
	\tgBlank{(1,0)}{\tgColour2}
	\tgBorderA{(2,0)}{\tgColour2}{\tgColour2}{\tgColour0}{\tgColour2}
	\tgBorderA{(0,1)}{\tgColour5}{\tgColour2}{\tgColour5}{\tgColour5}
	\tgBorderA{(1,1)}{\tgColour2}{\tgColour2}{\tgColour0}{\tgColour5}
	\tgBorderA{(2,1)}{\tgColour2}{\tgColour0}{\tgColour0}{\tgColour0}
	\tgArrow{(2,0.5)}{1}
	\tgAxisLabel{(0.5,0)}{south}{\gamma}
	\tgAxisLabel{(3,0.5)}{west}{\tpl{0}}
	\tgAxisLabel{(1.5,2)}{north}{\gamma_0}
    \node [fill=black, inner ysep=1.5pt, inner xsep=3em] at (1.5,1.5) {};
\end{tangle}
,
\hspace{0.5em}
\lambda_1\coloneqq
\begin{tangle}{(3,2)}
	\tgBorderA{(0,0)}{\tgColour5}{\tgColour2}{\tgColour2}{\tgColour5}
	\tgBlank{(1,0)}{\tgColour2}
	\tgBorderA{(2,0)}{\tgColour2}{\tgColour2}{\tgColour1}{\tgColour2}
	\tgBorderA{(0,1)}{\tgColour5}{\tgColour2}{\tgColour5}{\tgColour5}
	\tgBorderA{(1,1)}{\tgColour2}{\tgColour2}{\tgColour1}{\tgColour5}
	\tgBorderA{(2,1)}{\tgColour2}{\tgColour1}{\tgColour1}{\tgColour1}
	\tgArrow{(2,0.5)}{1}
	\tgAxisLabel{(0.5,0)}{south}{\gamma}
	\tgAxisLabel{(3,0.5)}{west}{\tpl{1}}
	\tgAxisLabel{(1.5,2)}{north}{\gamma_1}
    \node [fill=black, inner ysep=1.5pt, inner xsep=3em] at (1.5,1.5) {};
\end{tangle}
\]
\[
\rho_0\coloneqq
\begin{tangle}{(3,2)}
	\tgBorderA{(0,0)}{\tgColour2}{\tgColour2}{\tgColour2}{\tgColour0}
	\tgBlank{(1,0)}{\tgColour2}
	\tgBorderA{(2,0)}{\tgColour2}{\tgColour5}{\tgColour5}{\tgColour2}
	\tgBorderA{(0,1)}{\tgColour0}{\tgColour2}{\tgColour0}{\tgColour0}
	\tgBorderA{(1,1)}{\tgColour2}{\tgColour2}{\tgColour5}{\tgColour0}
	\tgBorderA{(2,1)}{\tgColour2}{\tgColour5}{\tgColour5}{\tgColour5}
	\tgArrow{(0,0.5)}{3}
	\tgArrow{(5.5,4)}{0}
	\tgAxisLabel{(2.5,0)}{south}{\gamma^*}
	\tgAxisLabel{(0,0.5)}{east}{\tpl{1}}
	\tgAxisLabel{(1.5,2)}{north}{\gamma_0^*}
    \node [fill=black, inner ysep=1.5pt, inner xsep=3em] at (1.5,1.5) {};
\end{tangle}
,
\hspace{0.5em}
\rho_0\coloneqq
\begin{tangle}{(3,2)}
	\tgBorderA{(0,0)}{\tgColour2}{\tgColour2}{\tgColour2}{\tgColour1}
	\tgBlank{(1,0)}{\tgColour2}
	\tgBorderA{(2,0)}{\tgColour2}{\tgColour5}{\tgColour5}{\tgColour2}
	\tgBorderA{(0,1)}{\tgColour1}{\tgColour2}{\tgColour1}{\tgColour1}
	\tgBorderA{(1,1)}{\tgColour2}{\tgColour2}{\tgColour5}{\tgColour1}
	\tgBorderA{(2,1)}{\tgColour2}{\tgColour5}{\tgColour5}{\tgColour5}
	\tgArrow{(0,0.5)}{3}
	\tgAxisLabel{(2.5,0)}{south}{\gamma^*}
	\tgAxisLabel{(0,0.5)}{east}{\tpl{1}}
	\tgAxisLabel{(1.5,2)}{north}{\gamma_1^*}
    \node [fill=black, inner ysep=1.5pt, inner xsep=3em] at (1.5,1.5) {};
\end{tangle}
\]
Using the universal property of the binary product in $\dbl{D}$,
we have the pairings $\tpl{\lambda_0,\lambda_1}$ and $\tpl{\rho_0,\rho_1}$.
Then, we can define $\lambda$ and $\rho$ as follows.
\[
\lambda\coloneqq
\hspace{1em}
\begin{tangle}{(3,3)}
	\tgBlank{(0,0)}{\tgColour5}
	\tgBorderA{(1,0)}{\tgColour5}{\tgColour2}{\tgColour2}{\tgColour5}
	\tgBlank{(2,0)}{\tgColour2}
	\tgBorderA{(0,1)}{\tgColour5}{\tgColour5}{\tgColour4}{\tgColour5}
	\tgBorderA{(1,1)}{\tgColour5}{\tgColour2}{\tgColour2}{\tgColour4}
	\tgBlank{(2,1)}{\tgColour2}
	\tgBorderA{(0,2)}{\tgColour5}{\tgColour4}{\tgColour4}{\tgColour5}
	\tgBorderA{(1,2)}{\tgColour4}{\tgColour2}{\tgColour2}{\tgColour4}
	\tgBlank{(2,2)}{\tgColour2}
	\tgCell[(1,0)]{(1,1)}{\scriptstyle\tpl{\lambda_0,\lambda_1}}
	\tgAxisLabel{(1.5,0)}{south}{\gamma}
	\tgAxisLabel{(0.5,3)}{north}{\tpl{0,0}_*}
	\tgAxisLabel{(1.5,3)}{north}{\tpl{\gamma_0,\gamma_1}}
\end{tangle},
\hspace{1em}
\rho\coloneqq
\hspace{1em}
\begin{tangle}{(5,3)}
	\tgBorderA{(0,0)}{\tgColour5}{\tgColour5}{\tgColour2}{\tgColour5}
	\tgBorderA{(1,0)}{\tgColour5}{\tgColour5}{\tgColour2}{\tgColour2}
	\tgBorderA{(2,0)}{\tgColour5}{\tgColour5}{\tgColour5}{\tgColour2}
	\tgBorderA{(3,0)}{\tgColour5}{\tgColour4}{\tgColour4}{\tgColour5}
	\tgBorderA{(4,0)}{\tgColour4}{\tgColour2}{\tgColour2}{\tgColour4}
	\tgBorderA{(0,1)}{\tgColour5}{\tgColour2}{\tgColour2}{\tgColour5}
	\tgBlank{(1,1)}{\tgColour2}
	\tgBorderA{(2,1)}{\tgColour2}{\tgColour5}{\tgColour4}{\tgColour2}
	\tgBorderA{(3,1)}{\tgColour5}{\tgColour4}{\tgColour4}{\tgColour4}
	\tgBorderA{(4,1)}{\tgColour4}{\tgColour2}{\tgColour2}{\tgColour4}
	\tgBorderA{(0,2)}{\tgColour5}{\tgColour2}{\tgColour2}{\tgColour5}
	\tgBlank{(1,2)}{\tgColour2}
	\tgBorderA{(2,2)}{\tgColour2}{\tgColour4}{\tgColour2}{\tgColour2}
	\tgBorderA{(3,2)}{\tgColour4}{\tgColour4}{\tgColour2}{\tgColour2}
	\tgBorderA{(4,2)}{\tgColour4}{\tgColour2}{\tgColour2}{\tgColour2}
	\tgCell[(2,0)]{(1,0)}{\gamma}
	\tgCell[(2,0)]{(3,2)}{\gamma_0\times\gamma_1}
	\tgCell[(1,0)]{(2,1)}{\scriptstyle\tpl{\rho_0,\rho_1}}
	\tgAxisLabel{(3.5,0)}{south}{\tpl{0,0}_*}
	\tgAxisLabel{(4.5,0)}{south}{\tpl{\gamma_0,\gamma_1}}
	\tgAxisLabel{(0.5,3)}{north}{\gamma}
\end{tangle}
\]
Before showing that $\lambda$ and $\rho$ are mutually inverse,
we observe the following equality:
\[
\begin{tangle}{(5,2)}
	\tgBlank{(0,0)}{\tgColour2}
	\tgBorderA{(1,0)}{\tgColour2}{\tgColour5}{\tgColour4}{\tgColour2}
	\tgBorderA{(2,0)}{\tgColour5}{\tgColour5}{\tgColour4}{\tgColour4}
	\tgBorderA{(3,0)}{\tgColour5}{\tgColour2}{\tgColour2}{\tgColour4}
	\tgBlank{(4,0)}{\tgColour2}
	\tgBlank{(0,1)}{\tgColour2}
	\tgBorderA{(1,1)}{\tgColour2}{\tgColour4}{\tgColour2}{\tgColour2}
	\tgBorderA{(2,1)}{\tgColour4}{\tgColour4}{\tgColour2}{\tgColour2}
	\tgBorderA{(3,1)}{\tgColour4}{\tgColour2}{\tgColour2}{\tgColour2}
	\tgBlank{(4,1)}{\tgColour2}
	\tgCell[(2,0)]{(2,1)}{\gamma_0\times\gamma_1}
	\tgCell[(1,0)]{(1,0)}{\scriptstyle\tpl{\rho_0,\rho_1}}
	\tgCell[(1,0)]{(3,0)}{\scriptstyle\tpl{\lambda_0,\lambda_1}}
	\tgAxisLabel{(1.5,0)}{south}{\gamma^*}
	\tgAxisLabel{(3.5,0)}{south}{\gamma}
	\node at (2.5,0.25) {$\scriptstyle\tpl{0,0}$};
\end{tangle}
=
\begin{tangle}{(2,2)}
	\tgBorderA{(0,0)}{\tgColour2}{\tgColour5}{\tgColour5}{\tgColour2}
	\tgBorderA{(1,0)}{\tgColour5}{\tgColour2}{\tgColour2}{\tgColour5}
	\tgBorderA{(0,1)}{\tgColour2}{\tgColour5}{\tgColour2}{\tgColour2}
	\tgBorderA{(1,1)}{\tgColour5}{\tgColour2}{\tgColour2}{\tgColour2}
	\tgCell[(1,0)]{(0.5,1)}{\gamma}
	\tgAxisLabel{(0.5,0)}{south}{\gamma^*}
	\tgAxisLabel{(1.5,0)}{south}{\gamma}
\end{tangle}
\]
This follows from that the both sides postcomposed by the projections
give the same result as shown below.
\[
\begin{tangle}{(5,3)}
	\tgBlank{(0,0)}{\tgColour2}
	\tgBorderA{(1,0)}{\tgColour2}{\tgColour5}{\tgColour4}{\tgColour2}
	\tgBorderA{(2,0)}{\tgColour5}{\tgColour5}{\tgColour4}{\tgColour4}
	\tgBorderA{(3,0)}{\tgColour5}{\tgColour2}{\tgColour2}{\tgColour4}
	\tgBlank{(4,0)}{\tgColour2}
	\tgBlank{(0,1)}{\tgColour2}
	\tgBorderA{(1,1)}{\tgColour2}{\tgColour4}{\tgColour2}{\tgColour2}
	\tgBorderA{(2,1)}{\tgColour4}{\tgColour4}{\tgColour2}{\tgColour2}
	\tgBorderA{(3,1)}{\tgColour4}{\tgColour2}{\tgColour2}{\tgColour2}
	\tgBlank{(4,1)}{\tgColour2}
	\tgBorderA{(0,2)}{\tgColour2}{\tgColour2}{\tgColour0}{\tgColour0}
	\tgBorderA{(1,2)}{\tgColour2}{\tgColour2}{\tgColour0}{\tgColour0}
	\tgBorderA{(2,2)}{\tgColour2}{\tgColour2}{\tgColour0}{\tgColour0}
	\tgBorderA{(3,2)}{\tgColour2}{\tgColour2}{\tgColour0}{\tgColour0}
	\tgBorderA{(4,2)}{\tgColour2}{\tgColour2}{\tgColour0}{\tgColour0}
	\tgCell[(2,0)]{(2,1)}{\gamma_0\times\gamma_1}
	\tgCell[(1,0)]{(1,0)}{\scriptstyle\tpl{\rho_0,\rho_1}}
	\tgCell[(1,0)]{(3,0)}{\scriptstyle\tpl{\lambda_0,\lambda_1}}
	\tgAxisLabel{(1.5,0)}{south}{\gamma^*}
	\tgAxisLabel{(3.5,0)}{south}{\gamma}
	\tgAxisLabel{(5,2.5)}{west}{\tpl{0}}
\end{tangle}
=
\begin{tangle}{(5,3)}
	\tgBlank{(0,0)}{\tgColour2}
	\tgBorderA{(1,0)}{\tgColour2}{\tgColour5}{\tgColour4}{\tgColour2}
	\tgBorderA{(2,0)}{\tgColour5}{\tgColour5}{\tgColour4}{\tgColour4}
	\tgBorderA{(3,0)}{\tgColour5}{\tgColour2}{\tgColour2}{\tgColour4}
	\tgBlank{(4,0)}{\tgColour2}
	\tgBorderA{(0,1)}{\tgColour2}{\tgColour2}{\tgColour0}{\tgColour0}
	\tgBorderA{(1,1)}{\tgColour2}{\tgColour4}{\tgColour5}{\tgColour0}
	\tgBorderA{(2,1)}{\tgColour4}{\tgColour4}{\tgColour5}{\tgColour5}
	\tgBorderA{(3,1)}{\tgColour4}{\tgColour2}{\tgColour0}{\tgColour5}
	\tgBorderA{(4,1)}{\tgColour2}{\tgColour2}{\tgColour0}{\tgColour0}
	\tgBlank{(0,2)}{\tgColour0}
	\tgBorderA{(1,2)}{\tgColour0}{\tgColour5}{\tgColour0}{\tgColour0}
	\tgBorderA{(2,2)}{\tgColour5}{\tgColour5}{\tgColour0}{\tgColour0}
	\tgBorderA{(3,2)}{\tgColour5}{\tgColour0}{\tgColour0}{\tgColour0}
	\tgBlank{(4,2)}{\tgColour0}
	\tgCell[(1,0)]{(1,0)}{\scriptstyle\tpl{\rho_0,\rho_1}}
	\tgCell[(1,0)]{(3,0)}{\scriptstyle\tpl{\lambda_0,\lambda_1}}
	\tgCell{(1,1)}{\scriptstyle\tpl{0}}
	\tgCell{(3,1)}{\scriptstyle\tpl{0}}
	\tgCell[(2,0)]{(2,2)}{\gamma_0}
\end{tangle}
\]
\[
=
\begin{tangle}{(5,3)}
	\tgBlank{(0,0)}{\tgColour2}
	\tgBorderA{(1,0)}{\tgColour2}{\tgColour5}{\tgColour5}{\tgColour2}
	\tgBlank{(2,0)}{\tgColour5}
	\tgBorderA{(3,0)}{\tgColour5}{\tgColour2}{\tgColour2}{\tgColour5}
	\tgBlank{(4,0)}{\tgColour2}
	\tgBorderA{(0,1)}{\tgColour2}{\tgColour2}{\tgColour0}{\tgColour0}
	\tgBorderA{(1,1)}{\tgColour2}{\tgColour5}{\tgColour5}{\tgColour0}
	\tgBlank{(2,1)}{\tgColour5}
	\tgBorderA{(3,1)}{\tgColour5}{\tgColour2}{\tgColour0}{\tgColour5}
	\tgBorderA{(4,1)}{\tgColour2}{\tgColour2}{\tgColour0}{\tgColour0}
	\tgBlank{(0,2)}{\tgColour0}
	\tgBorderA{(1,2)}{\tgColour0}{\tgColour5}{\tgColour0}{\tgColour0}
	\tgBorderA{(2,2)}{\tgColour5}{\tgColour5}{\tgColour0}{\tgColour0}
	\tgBorderA{(3,2)}{\tgColour5}{\tgColour0}{\tgColour0}{\tgColour0}
	\tgBlank{(4,2)}{\tgColour0}
	\tgCell[(2,0)]{(2,2)}{\gamma_0}
	\tgCell{(1,1)}{\rho_0}
	\tgCell{(3,1)}{\lambda_0}
\end{tangle}
=
\begin{tangle}{(5,3)}
	\tgBlank{(0,0)}{\tgColour2}
	\tgBorderA{(1,0)}{\tgColour2}{\tgColour5}{\tgColour5}{\tgColour2}
	\tgBlank{(2,0)}{\tgColour5}
	\tgBorderA{(3,0)}{\tgColour5}{\tgColour2}{\tgColour2}{\tgColour5}
	\tgBlank{(4,0)}{\tgColour2}
	\tgBlank{(0,1)}{\tgColour2}
	\tgBorderA{(1,1)}{\tgColour2}{\tgColour5}{\tgColour2}{\tgColour2}
	\tgBorderA{(2,1)}{\tgColour5}{\tgColour5}{\tgColour2}{\tgColour2}
	\tgBorderA{(3,1)}{\tgColour5}{\tgColour2}{\tgColour2}{\tgColour2}
	\tgBlank{(4,1)}{\tgColour2}
	\tgBorderA{(0,2)}{\tgColour2}{\tgColour2}{\tgColour0}{\tgColour0}
	\tgBorderA{(1,2)}{\tgColour2}{\tgColour2}{\tgColour0}{\tgColour0}
	\tgBorderA{(2,2)}{\tgColour2}{\tgColour2}{\tgColour0}{\tgColour0}
	\tgBorderA{(3,2)}{\tgColour2}{\tgColour2}{\tgColour0}{\tgColour0}
	\tgBorderA{(4,2)}{\tgColour2}{\tgColour2}{\tgColour0}{\tgColour0}
	\tgCell[(2,0)]{(2,1)}{\gamma}
\end{tangle}
\]

Then, with this equality, we can prove one direction of the mutual invertibility of $\lambda$ and $\rho$. 
\[
\begin{tangle}{(5,3)}
	\tgBorderA{(0,0)}{\tgColour5}{\tgColour5}{\tgColour2}{\tgColour5}
	\tgBorderA{(1,0)}{\tgColour5}{\tgColour5}{\tgColour5}{\tgColour2}
	\tgBorderA{(2,0)}{\tgColour5}{\tgColour5}{\tgColour4}{\tgColour5}
	\tgBorderA{(3,0)}{\tgColour5}{\tgColour2}{\tgColour2}{\tgColour4}
	\tgBlank{(4,0)}{\tgColour2}
	\tgBorderA{(0,1)}{\tgColour5}{\tgColour2}{\tgColour2}{\tgColour5}
	\tgBorderA{(1,1)}{\tgColour2}{\tgColour5}{\tgColour4}{\tgColour2}
	\tgBorderA{(2,1)}{\tgColour5}{\tgColour4}{\tgColour4}{\tgColour4}
	\tgBorderA{(3,1)}{\tgColour4}{\tgColour2}{\tgColour2}{\tgColour4}
	\tgBlank{(4,1)}{\tgColour2}
	\tgBorderA{(0,2)}{\tgColour5}{\tgColour2}{\tgColour2}{\tgColour5}
	\tgBorderA{(1,2)}{\tgColour2}{\tgColour4}{\tgColour2}{\tgColour2}
	\tgBorderA{(2,2)}{\tgColour4}{\tgColour4}{\tgColour2}{\tgColour2}
	\tgBorderA{(3,2)}{\tgColour4}{\tgColour2}{\tgColour2}{\tgColour2}
	\tgBlank{(4,2)}{\tgColour2}
	\tgCell[(1,0)]{(0.5,0)}{\gamma}
	\tgCell[(1,0)]{(1,1)}{\scriptstyle\tpl{\rho_0,\rho_1}}
	\tgCell[(1,0)]{(3,0)}{\scriptstyle\tpl{\lambda_0,\lambda_1}}
	\tgCell[(2,0)]{(2,2)}{\gamma_0\times\gamma_1}
	\tgArrow{(2,0.5)}{1}
\end{tangle}
=
\begin{tangle}{(5,3)}
	\tgBorderA{(0,0)}{\tgColour5}{\tgColour5}{\tgColour2}{\tgColour5}
	\tgBorderA{(1,0)}{\tgColour5}{\tgColour5}{\tgColour5}{\tgColour2}
	\tgBlank{(2,0)}{\tgColour5}
	\tgBorderA{(3,0)}{\tgColour5}{\tgColour2}{\tgColour2}{\tgColour5}
	\tgBlank{(4,0)}{\tgColour2}
	\tgBorderA{(0,1)}{\tgColour5}{\tgColour2}{\tgColour2}{\tgColour5}
	\tgBorderA{(1,1)}{\tgColour2}{\tgColour5}{\tgColour4}{\tgColour2}
	\tgBorderA{(2,1)}{\tgColour5}{\tgColour5}{\tgColour4}{\tgColour4}
	\tgBorderA{(3,1)}{\tgColour5}{\tgColour2}{\tgColour2}{\tgColour4}
	\tgBlank{(4,1)}{\tgColour2}
	\tgBorderA{(0,2)}{\tgColour5}{\tgColour2}{\tgColour2}{\tgColour5}
	\tgBorderA{(1,2)}{\tgColour2}{\tgColour4}{\tgColour2}{\tgColour2}
	\tgBorderA{(2,2)}{\tgColour4}{\tgColour4}{\tgColour2}{\tgColour2}
	\tgBorderA{(3,2)}{\tgColour4}{\tgColour2}{\tgColour2}{\tgColour2}
	\tgBlank{(4,2)}{\tgColour2}
	\tgCell[(1,0)]{(0.5,0)}{\gamma}
	\tgCell[(1,0)]{(1,1)}{\scriptstyle\tpl{\rho_0,\rho_1}}
	\tgCell[(2,0)]{(2,2)}{\gamma_0\times\gamma_1}
	\tgCell[(1,0)]{(3,1)}{\scriptstyle\tpl{\lambda_0,\lambda_1}}
\end{tangle}
\]
\[
=
\begin{tangle}{(3,3)}
	\tgBorderA{(0,0)}{\tgColour5}{\tgColour5}{\tgColour2}{\tgColour5}
	\tgBorderA{(1,0)}{\tgColour5}{\tgColour5}{\tgColour5}{\tgColour2}
	\tgBorderA{(2,0)}{\tgColour5}{\tgColour2}{\tgColour2}{\tgColour5}
	\tgBorderA{(0,1)}{\tgColour5}{\tgColour2}{\tgColour2}{\tgColour5}
	\tgBorderA{(1,1)}{\tgColour2}{\tgColour5}{\tgColour2}{\tgColour2}
	\tgBorderA{(2,1)}{\tgColour5}{\tgColour2}{\tgColour2}{\tgColour2}
	\tgBorderA{(0,2)}{\tgColour5}{\tgColour2}{\tgColour2}{\tgColour5}
	\tgBlank{(1,2)}{\tgColour2}
	\tgBlank{(2,2)}{\tgColour2}
	\tgCell[(1,0)]{(0.5,0)}{\gamma}
	\tgCell[(1,0)]{(1.5,1)}{\gamma}
\end{tangle}
=
\begin{tangle}{(1,3)}
	\tgBorderA{(0,0)}{\tgColour5}{\tgColour2}{\tgColour2}{\tgColour5}
	\tgBorderA{(0,1)}{\tgColour5}{\tgColour2}{\tgColour2}{\tgColour5}
	\tgBorderA{(0,2)}{\tgColour5}{\tgColour2}{\tgColour2}{\tgColour5}
\end{tangle}
\]

Using the universal property of the binary product in $\dbl{D}$ again,
the other direction follows from the following equality and its dual.
\[
\begin{tangle}{(3,4)}
	\tgBorderA{(0,0)}{\tgColour5}{\tgColour4}{\tgColour5}{\tgColour5}
	\tgBorderA{(1,0)}{\tgColour4}{\tgColour4}{\tgColour2}{\tgColour5}
	\tgBorderA{(2,0)}{\tgColour4}{\tgColour2}{\tgColour2}{\tgColour2}
	\tgBorderA{(0,1)}{\tgColour5}{\tgColour5}{\tgColour4}{\tgColour5}
	\tgBorderA{(1,1)}{\tgColour5}{\tgColour2}{\tgColour4}{\tgColour4}
	\tgBorderA{(2,1)}{\tgColour2}{\tgColour2}{\tgColour2}{\tgColour4}
	\tgBorderA{(0,2)}{\tgColour5}{\tgColour4}{\tgColour4}{\tgColour4}
	\tgBlank{(1,2)}{\tgColour4}
	\tgBorderA{(2,2)}{\tgColour4}{\tgColour2}{\tgColour2}{\tgColour4}
	\tgBorderA{(0,3)}{\tgColour4}{\tgColour4}{\tgColour5}{\tgColour5}
	\tgBorderA{(1,3)}{\tgColour4}{\tgColour4}{\tgColour5}{\tgColour5}
	\tgBorderA{(2,3)}{\tgColour4}{\tgColour2}{\tgColour0}{\tgColour5}
	\tgCell[(2,0)]{(1,0)}{\rho}
	\tgCell{(2,3)}{\scriptstyle\tpl{0}}
	\tgCell[(2,0)]{(1,1)}{\lambda}
	\tgAxisLabel{(0.5,0)}{south}{\tpl{0,0}_*}
	\tgAxisLabel{(2.5,0)}{south}{\tpl{\gamma_0,\gamma_1}}
	\tgAxisLabel{(0,2.5)}{east}{\tpl{0,0}}
	\tgAxisLabel{(0,3.5)}{east}{\tpl{0}}
	\tgAxisLabel{(3,3.5)}{west}{\tpl{0}}
	\tgAxisLabel{(2.5,4)}{north}{\gamma_0}
\end{tangle}
=
\begin{tangle}{(6,4)}
	\tgBlank{(0,0)}{\tgColour5}
	\tgBorderA{(1,0)}{\tgColour5}{\tgColour5}{\tgColour2}{\tgColour5}
	\tgBorderA{(2,0)}{\tgColour5}{\tgColour5}{\tgColour2}{\tgColour2}
	\tgBorderA{(3,0)}{\tgColour5}{\tgColour5}{\tgColour5}{\tgColour2}
	\tgBorderA{(4,0)}{\tgColour5}{\tgColour4}{\tgColour4}{\tgColour5}
	\tgBorderA{(5,0)}{\tgColour4}{\tgColour2}{\tgColour2}{\tgColour4}
	\tgBorderA{(0,1)}{\tgColour5}{\tgColour5}{\tgColour4}{\tgColour5}
	\tgBorderA{(1,1)}{\tgColour5}{\tgColour2}{\tgColour2}{\tgColour4}
	\tgBlank{(2,1)}{\tgColour2}
	\tgBorderA{(3,1)}{\tgColour2}{\tgColour5}{\tgColour4}{\tgColour2}
	\tgBorderA{(4,1)}{\tgColour5}{\tgColour4}{\tgColour4}{\tgColour4}
	\tgBorderA{(5,1)}{\tgColour4}{\tgColour2}{\tgColour2}{\tgColour4}
	\tgBorderA{(0,2)}{\tgColour5}{\tgColour4}{\tgColour4}{\tgColour4}
	\tgBorderA{(1,2)}{\tgColour4}{\tgColour2}{\tgColour2}{\tgColour4}
	\tgBlank{(2,2)}{\tgColour2}
	\tgBorderA{(3,2)}{\tgColour2}{\tgColour4}{\tgColour2}{\tgColour2}
	\tgBorderA{(4,2)}{\tgColour4}{\tgColour4}{\tgColour2}{\tgColour2}
	\tgBorderA{(5,2)}{\tgColour4}{\tgColour2}{\tgColour2}{\tgColour2}
	\tgBorderA{(0,3)}{\tgColour4}{\tgColour4}{\tgColour5}{\tgColour5}
	\tgBorderA{(1,3)}{\tgColour4}{\tgColour2}{\tgColour0}{\tgColour5}
	\tgBorderA{(2,3)}{\tgColour2}{\tgColour2}{\tgColour0}{\tgColour0}
	\tgBorderA{(3,3)}{\tgColour2}{\tgColour2}{\tgColour0}{\tgColour0}
	\tgBorderA{(4,3)}{\tgColour2}{\tgColour2}{\tgColour0}{\tgColour0}
	\tgBorderA{(5,3)}{\tgColour2}{\tgColour2}{\tgColour0}{\tgColour0}
	\tgCell[(1,0)]{(1,1)}{\scriptstyle\tpl{\lambda_0,\lambda_1}}
	\tgCell[(2,0)]{(2,0)}{\gamma}
	\tgCell[(1,0)]{(3,1)}{\scriptstyle\tpl{\rho_0,\rho_1}}
	\tgCell{(1,3)}{\scriptstyle\tpl{0}}
	\tgCell[(2,0)]{(4,2)}{\gamma_0\times\gamma_1}
\end{tangle}
\]
\[
\begin{tangle}{(6,4)}
	\tgBlank{(0,0)}{\tgColour5}
	\tgBorderA{(1,0)}{\tgColour5}{\tgColour5}{\tgColour2}{\tgColour5}
	\tgBorderA{(2,0)}{\tgColour5}{\tgColour5}{\tgColour2}{\tgColour2}
	\tgBorderA{(3,0)}{\tgColour5}{\tgColour5}{\tgColour5}{\tgColour2}
	\tgBorderA{(4,0)}{\tgColour5}{\tgColour4}{\tgColour4}{\tgColour5}
	\tgBorderA{(5,0)}{\tgColour4}{\tgColour2}{\tgColour2}{\tgColour4}
	\tgBorderA{(0,1)}{\tgColour5}{\tgColour5}{\tgColour4}{\tgColour4}
	\tgBorderA{(1,1)}{\tgColour5}{\tgColour2}{\tgColour2}{\tgColour4}
	\tgBlank{(2,1)}{\tgColour2}
	\tgBorderA{(3,1)}{\tgColour2}{\tgColour5}{\tgColour4}{\tgColour2}
	\tgBorderA{(4,1)}{\tgColour5}{\tgColour4}{\tgColour4}{\tgColour4}
	\tgBorderA{(5,1)}{\tgColour4}{\tgColour2}{\tgColour2}{\tgColour4}
	\tgBorderA{(0,2)}{\tgColour4}{\tgColour4}{\tgColour5}{\tgColour5}
	\tgBorderA{(1,2)}{\tgColour4}{\tgColour2}{\tgColour0}{\tgColour5}
	\tgBorderA{(2,2)}{\tgColour2}{\tgColour2}{\tgColour0}{\tgColour0}
	\tgBorderA{(3,2)}{\tgColour2}{\tgColour4}{\tgColour5}{\tgColour0}
	\tgBorderA{(4,2)}{\tgColour4}{\tgColour4}{\tgColour5}{\tgColour5}
	\tgBorderA{(5,2)}{\tgColour4}{\tgColour2}{\tgColour0}{\tgColour5}
	\tgBlank{(0,3)}{\tgColour5}
	\tgBorderA{(1,3)}{\tgColour5}{\tgColour0}{\tgColour0}{\tgColour5}
	\tgBlank{(2,3)}{\tgColour0}
	\tgBorderA{(3,3)}{\tgColour0}{\tgColour5}{\tgColour0}{\tgColour0}
	\tgBorderA{(4,3)}{\tgColour5}{\tgColour5}{\tgColour0}{\tgColour0}
	\tgBorderA{(5,3)}{\tgColour5}{\tgColour0}{\tgColour0}{\tgColour0}
	\tgCell[(1,0)]{(1,1)}{\scriptstyle\tpl{\lambda_0,\lambda_1}}
	\tgCell[(2,0)]{(2,0)}{\gamma}
	\tgCell[(1,0)]{(3,1)}{\scriptstyle\tpl{\rho_0,\rho_1}}
	\tgCell{(1,2)}{\scriptstyle\tpl{0}}
	\tgCell{(3,2)}{\scriptstyle\tpl{0}}
	\tgCell{(5,2)}{\scriptstyle\tpl{0}}
	\tgCell[(2,0)]{(4,3)}{\gamma_0}
\end{tangle}
=
\begin{tangle}{(7,4)}
	\tgBlank{(0,0)}{\tgColour5}
	\tgBorderA{(1,0)}{\tgColour5}{\tgColour5}{\tgColour2}{\tgColour5}
	\tgBorderA{(2,0)}{\tgColour5}{\tgColour5}{\tgColour2}{\tgColour2}
	\tgBorderA{(3,0)}{\tgColour5}{\tgColour5}{\tgColour5}{\tgColour2}
	\tgBlank{(4,0)}{\tgColour5}
	\tgBorderA{(5,0)}{\tgColour5}{\tgColour4}{\tgColour4}{\tgColour5}
	\tgBorderA{(6,0)}{\tgColour4}{\tgColour2}{\tgColour2}{\tgColour4}
	\tgBorderA{(0,1)}{\tgColour5}{\tgColour5}{\tgColour5}{\tgColour4}
	\tgBorderA{(1,1)}{\tgColour5}{\tgColour2}{\tgColour0}{\tgColour5}
	\tgBorderA{(2,1)}{\tgColour2}{\tgColour2}{\tgColour0}{\tgColour0}
	\tgBorderA{(3,1)}{\tgColour2}{\tgColour5}{\tgColour5}{\tgColour0}
	\tgBorderA{(4,1)}{\tgColour5}{\tgColour5}{\tgColour4}{\tgColour5}
	\tgBorderA{(5,1)}{\tgColour5}{\tgColour4}{\tgColour4}{\tgColour4}
	\tgBorderA{(6,1)}{\tgColour4}{\tgColour2}{\tgColour2}{\tgColour4}
	\tgBorderA{(0,2)}{\tgColour4}{\tgColour5}{\tgColour5}{\tgColour5}
	\tgBorderA{(1,2)}{\tgColour5}{\tgColour0}{\tgColour0}{\tgColour5}
	\tgBlank{(2,2)}{\tgColour0}
	\tgBorderA{(3,2)}{\tgColour0}{\tgColour5}{\tgColour5}{\tgColour0}
	\tgBorderA{(4,2)}{\tgColour5}{\tgColour4}{\tgColour5}{\tgColour5}
	\tgBorderA{(5,2)}{\tgColour4}{\tgColour4}{\tgColour5}{\tgColour5}
	\tgBorderA{(6,2)}{\tgColour4}{\tgColour2}{\tgColour0}{\tgColour5}
	\tgBlank{(0,3)}{\tgColour5}
	\tgBorderA{(1,3)}{\tgColour5}{\tgColour0}{\tgColour0}{\tgColour5}
	\tgBlank{(2,3)}{\tgColour0}
	\tgBorderA{(3,3)}{\tgColour0}{\tgColour5}{\tgColour0}{\tgColour0}
	\tgBorderA{(4,3)}{\tgColour5}{\tgColour5}{\tgColour0}{\tgColour0}
	\tgBorderA{(5,3)}{\tgColour5}{\tgColour5}{\tgColour0}{\tgColour0}
	\tgBorderA{(6,3)}{\tgColour5}{\tgColour0}{\tgColour0}{\tgColour0}
	\tgCell[(2,0)]{(2,0)}{\gamma}
	\tgCell{(6,2)}{\scriptstyle\tpl{0}}
	\tgCell{(1,1)}{\lambda_0}
	\tgCell{(3,1)}{\rho_0}
	\tgCell[(3,0)]{(4.5,3)}{\gamma_0}
	\node [fill=black, inner xsep=1.5pt, inner ysep=1.5em] at (0.5,2) {}; 
	\node [fill=black, inner xsep=1.5pt, inner ysep=1.5em] at (4.5,2) {};
\end{tangle}
\]
\[
\begin{tangle}{(6,4)}
	\tgBorderA{(0,0)}{\tgColour5}{\tgColour5}{\tgColour4}{\tgColour4}
	\tgBorderA{(1,0)}{\tgColour5}{\tgColour5}{\tgColour4}{\tgColour4}
	\tgBorderA{(2,0)}{\tgColour5}{\tgColour5}{\tgColour5}{\tgColour4}
	\tgBorderA{(3,0)}{\tgColour5}{\tgColour5}{\tgColour4}{\tgColour5}
	\tgBorderA{(4,0)}{\tgColour5}{\tgColour4}{\tgColour4}{\tgColour4}
	\tgBorderA{(5,0)}{\tgColour4}{\tgColour2}{\tgColour2}{\tgColour4}
	\tgBorderA{(0,1)}{\tgColour4}{\tgColour4}{\tgColour5}{\tgColour5}
	\tgBorderA{(1,1)}{\tgColour4}{\tgColour4}{\tgColour5}{\tgColour5}
	\tgBorderA{(2,1)}{\tgColour4}{\tgColour5}{\tgColour5}{\tgColour5}
	\tgBorderA{(3,1)}{\tgColour5}{\tgColour4}{\tgColour5}{\tgColour5}
	\tgBorderA{(4,1)}{\tgColour4}{\tgColour4}{\tgColour5}{\tgColour5}
	\tgBorderA{(5,1)}{\tgColour4}{\tgColour2}{\tgColour0}{\tgColour5}
	\tgBlank{(0,2)}{\tgColour5}
	\tgBorderA{(1,2)}{\tgColour5}{\tgColour5}{\tgColour0}{\tgColour5}
	\tgBorderA{(2,2)}{\tgColour5}{\tgColour5}{\tgColour0}{\tgColour0}
	\tgBorderA{(3,2)}{\tgColour5}{\tgColour5}{\tgColour5}{\tgColour0}
	\tgBlank{(4,2)}{\tgColour5}
	\tgBorderA{(5,2)}{\tgColour5}{\tgColour0}{\tgColour0}{\tgColour5}
	\tgBlank{(0,3)}{\tgColour5}
	\tgBorderA{(1,3)}{\tgColour5}{\tgColour0}{\tgColour0}{\tgColour5}
	\tgBlank{(2,3)}{\tgColour0}
	\tgBorderA{(3,3)}{\tgColour0}{\tgColour5}{\tgColour0}{\tgColour0}
	\tgBorderA{(4,3)}{\tgColour5}{\tgColour5}{\tgColour0}{\tgColour0}
	\tgBorderA{(5,3)}{\tgColour5}{\tgColour0}{\tgColour0}{\tgColour0}
	\tgCell[(2,0)]{(2,2)}{\gamma_0}
	\tgCell[(2,0)]{(4,3)}{\gamma_0}
	\tgCell{(5,1)}{\scriptstyle\tpl{0}}
	\node [fill=black, inner xsep=1.5pt, inner ysep=1.5em] at (2.5,1) {};
	\node [fill=black, inner xsep=1.5pt, inner ysep=1.5em] at (3.5,1) {};
\end{tangle}
=
\begin{tangle}{(3,3)}
	\tgBorderA{(0,0)}{\tgColour5}{\tgColour5}{\tgColour4}{\tgColour4}
	\tgBorderA{(1,0)}{\tgColour5}{\tgColour4}{\tgColour4}{\tgColour4}
	\tgBorderA{(2,0)}{\tgColour4}{\tgColour2}{\tgColour2}{\tgColour4}
	\tgBorderA{(0,1)}{\tgColour4}{\tgColour4}{\tgColour5}{\tgColour5}
	\tgBorderA{(1,1)}{\tgColour4}{\tgColour4}{\tgColour5}{\tgColour5}
	\tgBorderA{(2,1)}{\tgColour4}{\tgColour2}{\tgColour0}{\tgColour5}
	\tgBlank{(0,2)}{\tgColour5}
	\tgBlank{(1,2)}{\tgColour5}
	\tgBorderA{(2,2)}{\tgColour5}{\tgColour0}{\tgColour0}{\tgColour5}
	\tgCell{(2,1)}{\scriptstyle\tpl{0}}
\end{tangle}
\]
This implies that $\Phi(\Psi(\gamma))$ is isomorphic to $\gamma$.

Finally, we prove (iii).
The lax functors $1\colon\bi{1}\to\LBi{\dbl{D}}$ and $\times\colon\LBi{\dbl{D}}\times\LBi{\dbl{D}}\to\LBi{\dbl{D}}$
induced by the finite products in $\MapBi{\dbl{D}}$ are
in fact the same as the loose parts of the double functors $1\colon\dbl{1}\to\dbl{D}$ and $\times\colon\dbl{D}\times\dbl{D}\to\dbl{D}$.
This is because $\times\colon\LBi{\dbl{D}}\times\LBi{\dbl{D}}\to\LBi{\dbl{D}}$ 
induced by the finite products in $\MapBi{\dbl{D}}$
sends a pair of loose arrows $(\alpha,\beta)$ to $(\tpl{0}_*\alpha\tpl{0}^*)\land(\tpl{1}_*\beta\tpl{1}^*)$,
which in a certain equipment is isomorphic to $\alpha\times\beta$.
The laxity cells are also confirmed to be the same as the cells 
derived from the universal properties of the binary product in the cartesian equipment $\dbl{D}$. 
Therefore, the lax functors $1$ and $\times$ are pseudo since $\dbl{D}$ is a cartesian equipment.  
\end{proof}

\end{document}